\documentclass[10pt,a4paper]{these}
\usepackage{t1enc}
\usepackage[latin1]{inputenc}
\usepackage[french]{babel}
\usepackage{indentfirst}
\usepackage{amsthm}
\usepackage{amsmath}
\usepackage{latexsym}
\usepackage{amssymb}
\usepackage[all]{xy}
\usepackage{makeidx}

\makeindex

\newcommand{\surj}{-\!\!\!-\!\!\!\twoheadrightarrow}
\newcommand{\ie}{{\em i.e. }}
\newcommand{\sens}{\mathrm{SENS}}

\newcommand{\ens}{\mathrm{ENS}}

\newcommand{\uEEns}{\underline{\mathcal{ENS}}}
\newcommand{\catcat}{\mathrm{CAT}}
\newcommand{\gpd}{\mathrm{GPD}}
\newcommand{\sgpd}{\mathrm{SGPD}}
\newcommand{\gp}{\mathrm{GP}}

\newcommand{\cmg}{\mathrm{CM}^g}
\newcommand{\cmd}{\mathrm{CM}^d}
\newcommand{\cmsg}{\mathrm{CMS}^g}
\newcommand{\cmsd}{\mathrm{CMS}^d}
\newcommand{\scmg}{\mathrm{SCM}^g}
\newcommand{\scmd}{\mathrm{SCM}^d}
\newcommand{\scmsg}{\mathrm{SCMS}^g}
\newcommand{\scmsd}{\mathrm{SCMS}^d}
\newcommand{\scms}{\mathrm{SCMS}}
\newcommand{\catcateq}{\mathrm{CATEQ}}
\newcommand{\scatcat}{\mathrm{SCAT}}

\newcommand{\Hom}{{\mathrm{Hom}}}
\newcommand{\HHom}{{\mathcal{H}\mathrm{om}}}
\newcommand{\End}{{\mathrm{End}}}
\newcommand{\uHom}{{\underline{\mathrm{Hom}}}}
\newcommand{\uHHom}{{\underline{\mathcal{H}\mathrm{om}}}}
\newcommand{\RHom}{{\mathbb{R}\mathrm{Hom}}}
\newcommand{\RHHom}{{\mathbb{R}\mathcal{H}\mathrm{om}}}

\newcommand{\uEq}{{\underline{\mathrm{Eq}}}}
\newcommand{\Eq}{{\mathrm{Eq}}}

\newcommand{\N}{\mathbb{N}}
\newcommand{\B}{\mathcal{B}}
\newcommand{\G}{\mathcal{G}}

\newcommand{\kk}{\Bbbk}
\newcommand{\Sn}{\mathfrak{S}_n}
\newcommand{\SE}{\mathfrak{S}_E}

\newcommand{\cat}{\mathcal{C}at}

\newcommand{\kcat}{\mathcal{C}at}
\newcommand{\ukcat}{\underline{\mathcal{C}at}}
\newcommand{\catiso}{\mathcal{C}at^{Iso}}
\newcommand{\kcatiso}{\mathcal{C}at^{Iso}}
\newcommand{\cateq}{\mathcal{C}at^{Eq}}
\newcommand{\kcateq}{\mathcal{C}at^{Eq}}
\newcommand{\tkcateq}{\widetilde{\mathcal{C}at}^{Eq}}
\newcommand{\tkcatmor}{\widetilde{\mathcal{C}at}^{Mor}}

\newcommand{\kcatmor}{\mathcal{C}at^{Mor}}
\newcommand{\ukcatmor}{\underline{\mathcal{C}at}^{Mor}}

\newcommand{\CAT}{\mathcal{CAT}}

\newcommand{\ukCAT}{\underline{\mathcal{CAT}}}
\newcommand{\CATiso}{\mathcal{CAT}^{Iso}}
\newcommand{\kCATiso}{\mathcal{CAT}^{Iso}}
\newcommand{\CATeq}{\mathcal{CAT}^{Eq}}
\newcommand{\kCATeq}{\mathcal{CAT}^{Eq}}
\newcommand{\CATmor}{\mathcal{CAT}^{Mor}}

\newcommand{\kCATmor}{\mathcal{CAT}^{Mor}}
\newcommand{\ukCATmor}{\underline{\mathcal{CAT}}^{Mor}}
\newcommand{\tkCATmor}{\widetilde{\mathcal{CAT}}^{Mor}}

\newcommand{\ukcateq}{\underline{\mathcal{C}at}^{Eq}}
\newcommand{\ukcateqg}{\underline{\mathcal{C}at}^{Eq,g}}
\newcommand{\ukcateqd}{\underline{\mathcal{C}at}^{Eq,d}}
\newcommand{\kcateqg}{\mathcal{C}at^{Eq,g}}

\newcommand{\ukcatiso}{\underline{\mathcal{C}at}^{Iso}}

\newcommand{\ukcatison}{\underline{\mathcal{C}at}^{Iso,n}}

\newcommand{\kcatisod}{\mathcal{C}at^{Iso,(d)}}
\newcommand{\ukcatisod}{\underline{\mathcal{C}at}^{Iso,d}}

\newcommand{\ukcatisons}{\underline{\mathcal{C}at}^{Iso,n,strict}}
\newcommand{\ukCATisons}{\underline{\mathcal{CAT}}^{Iso,n,strict}}

\newcommand{\ukCATeqg}{\underline{\mathcal{CAT}}^{Eq,g}}
\newcommand{\ukCATeqd}{\underline{\mathcal{CAT}}^{Eq,d}}
\newcommand{\tkCATeqg}{\widetilde{\mathcal{CAT}}^{Eq,g}}
\newcommand{\tkCATeqd}{\widetilde{\mathcal{CAT}}^{Eq,d}}
\newcommand{\tkCATeq}{\widetilde{\mathcal{CAT}}^{Eq}}
\newcommand{\kCATeqg}{\mathcal{CAT}^{Eq,g}}
\newcommand{\kCATeqd}{\mathcal{CAT}^{Eq,d}}

\newcommand{\tkCATiso}{\widetilde{\mathcal{CAT}}^{Iso}}
\newcommand{\tkcatiso}{\widetilde{\mathcal{C}at}^{Iso}}

\newcommand{\ukcateqinfn}{\underline{\mathcal{C}at}^{Eq,\,\leq n}}
\newcommand{\ukCATeqinfn}{\underline{\mathcal{CAT}}^{Eq,\,\leq n}}

\newcommand{\ukcateqinfnplusun}{\underline{\mathcal{C}at}^{Eq,\,\leq n+1}}

\newcommand{\ukCATeq}{\underline{\mathcal{CAT}}^{Eq}}
\newcommand{\ukCATiso}{\underline{\mathcal{CAT}}^{Iso}}

\newcommand{\com}{\mathcal{C}om}
\newcommand{\COM}{\mathcal{COM}}
\newcommand{\ukCOM}{\underline{\mathcal{COM}}}
\newcommand{\kcomm}{\mathcal{C}om}
\newcommand{\kcom}{\mathcal{COM}}

\newcommand{\ukcom}{\underline{\mathcal{C}om}}
\newcommand{\aff}{\mathcal{AFF}}
\newcommand{\kaff}{\mathcal{AFF}}
\newcommand{\ukaff}{\underline{\mathcal{AFF}}}
\newcommand{\kaffw}{\mathcal{AFF}^\wedge}

\newcommand{\chaffk}{\mathcal{C}h(\mathcal{AFF})}

\newcommand{\ass}{\mathcal{A}ss}

\newcommand{\ukass}{\underline{\mathcal{A}ss}}
\newcommand{\ASS}{\mathcal{ASS}}
\newcommand{\kASS}{\mathcal{ASS}}
\newcommand{\ukASS}{\underline{\mathcal{ASS}}}

\newcommand{\uvect}{{\underline{\mathcal{V}ect}}}
\newcommand{\uvectd}{{\underline{\mathcal{V}ect}^{(d)}}}

\newcommand{\gln}{G\ell_n}
\newcommand{\bgln}{\mathcal{B}G\ell_n}
\newcommand{\ubgln}{\underline{\mathcal{B}G\ell}_n}
\newcommand{\ab}{\mathcal{A}b}
\newcommand{\kab}{\mathcal{A}b}
\newcommand{\tkab}{\widetilde{\mathcal{A}b}}
\newcommand{\ukab}{\underline{\mathcal{A}b}}
\newcommand{\AB}{\mathcal{AB}}
\newcommand{\kAB}{\mathcal{AB}}
\newcommand{\ukAB}{\underline{\mathcal{AB}}}
\newcommand{\tkAB}{\widetilde{\mathcal{AB}}}

\newcommand{\KAR}{\mathcal{KAR}}

\newcommand{\z}{\textrm{-}}
\newcommand{\Mod}{\textrm{-}Mod}
\newcommand{\zMMod}{\textrm{-}\mathcal{M}od}
\newcommand{\MMod}{\mathcal{M}od}
\newcommand{\proj}{\textrm{-}\textrm{proj}}

\newcommand{\UU}{\mathbb{U}}
\newcommand{\VV}{\mathbb{V}}
\newcommand{\WW}{\mathbb{W}}

\newenvironment{nota}%
{\refstepcounter{tot}\par\medskip \noindent{\sffamily \bfseries Notation~\thetot~}}%
{\par\medskip}
\newenvironment{cor}%
{\refstepcounter{tot}\par\medskip \noindent{\sffamily \bfseries Corollaire~\thetot~} \it}%
{\par\medskip\normalsize}
\newenvironment{corx}%
{\refstepcounter{totx}\par\medskip \noindent{\sffamily \bfseries Corollaire~\thetotx~} \it}%
{\par\medskip\normalsize}
\newenvironment{thm}%
{\refstepcounter{tot}\par\medskip \noindent{\sffamily \bfseries Théorème~\thetot~} \it}%
{\par\medskip\normalsize}
{\refstepcounter{totx}\par\medskip \noindent{\sffamily \bfseries Théorème~\thetotx~} \it}%
{\par\medskip\normalsize}
\newenvironment{conj}%
{\refstepcounter{tot}\par\medskip \noindent{\sffamily \bfseries Conjecture~\thetot~} \it}%
{\par\medskip\normalsize}

{\refstepcounter{tot}\par\medskip \noindent{\sffamily \bfseries Hypothèse~\thetot~} \it}%
{\par\medskip\normalsize}
{\refstepcounter{tot}\par\medskip \noindent{\sffamily \bfseries Affirmation~\thetot~} \it}%
{\par\medskip\normalsize}
\newenvironment{lemme}%
{\refstepcounter{tot}\par\medskip \noindent{\sffamily \bfseries Lemme~\thetot~} \it}%
{\par\medskip\normalsize}
\newenvironment{lemmex}%
{\refstepcounter{totx}\par\medskip \noindent{\sffamily \bfseries Lemme~\thetotx~} \it}%
{\par\medskip\normalsize}
{\refstepcounter{tot}\par\medskip \noindent{\sffamily \bfseries Sous-lemme~\thetot~} \it}%
{\par\medskip\normalsize}
\newenvironment{prop}%
{\refstepcounter{tot}\par\medskip \noindent{\sffamily \bfseries Proposition~\thetot~} \it}%
{\par\medskip\normalsize}
\newenvironment{propx}%
{\refstepcounter{totx}\par\medskip \noindent{\sffamily \bfseries Proposition~\thetotx~} \it}%
{\par\medskip\normalsize}
\newenvironment{rem}%
{\refstepcounter{tot}\par\medskip \noindent{\sffamily \bfseries Remarque~\thetot~}}%
{\par\medskip\normalsize}
{\refstepcounter{tot}\par\medskip \noindent{\sffamily \bfseries Fait~\thetot~} \it}%
{\par\medskip\normalsize}
\newenvironment{defi}%
{\refstepcounter{tot}\par\medskip \noindent{\sffamily \bfseries Définition~\thetot~}}%
{\par\medskip}
\newenvironment{defix}%
{\refstepcounter{totx}\par\medskip \noindent{\sffamily \bfseries Définition~\thetotx~}}%
{\par\medskip}

\newenvironment{rthm}[1][]%
{\refstepcounter{tot}\par\medskip \noindent{\sffamily \bfseries Théorème~\thetot~}{\sffamily (#1)} \it}%
{\par\medskip\normalsize}
\newenvironment{rlemme}[1][]%
{\refstepcounter{tot}\par\medskip \noindent{\sffamily \bfseries Lemme~\thetot~}{\sffamily (#1)} \it}%
{\par\medskip\normalsize}
\newenvironment{rprop}[1][]%
{\refstepcounter{tot}\par\medskip \noindent{\sffamily \bfseries Proposition~\thetot~}{\sffamily (#1)} \it}%
{\par\medskip\normalsize}
{\refstepcounter{totx}\par\medskip \noindent{\sffamily \bfseries Théorème~\thetotx~}{\sffamily (#1)} \it}%
{\par\medskip\normalsize}
{\refstepcounter{totx}\par\medskip \noindent{\sffamily \bfseries Lemme~\thetotx~}{\sffamily (#1)} \it}%
{\par\medskip\normalsize}
\newenvironment{rpropx}[1][]%
{\refstepcounter{totx}\par\medskip \noindent{\sffamily \bfseries Proposition~\thetotx~}{\sffamily (#1)} \it}%
{\par\medskip\normalsize}
{\refstepcounter{tot}\par\medskip \noindent{\sffamily \bfseries Définition~\thetot~}{\sffamily (#1)}}%
{\par\medskip}

\newenvironment{pr}%
{\noindent{\sffamily \slshape Preuve~}}%
{\hfill $\square$ \medskip}
\newenvironment{prpr}%
{\noindent{\sffamily \slshape Preuve~de~la~proposition~}}%
{\hfill $\square$ \medskip}

\newenvironment{preuve}%
{\noindent{\sffamily \slshape Preuve~}}%
{\hfill $\square$ \medskip}

\renewcommand{\subparagraph}[1]{\medskip\par\noindent{\bfseries\sffamily #1}~}

\begin{document}


\thispagestyle{empty}
\begin{center}
{\Large \sffamily \bfseries THÈSE \\
présentée pour le diplôme de \\
docteur de l'université Toulouse III -- Paul Sabatier\\
\vspace{5pt} en MATH\'EMATIQUES PURES\\}
\vspace{10pt}
{\large \sffamily \bfseries par\\}
\vspace{20pt}
{\huge \sffamily \bfseries Mathieu Anel \\}
\vspace{20pt}
{\large \sffamily \bfseries intitulée\\}
\vspace{50pt}
{\Huge \sffamily \bfseries CHAMPS DE MODULES\\ DES CAT\'EGORIES LIN\'EAIRES\\ ET AB\'ELIENNES\\}
\vspace{50pt}

{\large \sffamily \bfseries soutenue le 23 juin 2006\\}
{\large \sffamily \bfseries devant le jury composé de\\}

\vspace{20pt}
{\sffamily \Large
\begin{tabular}{llr}
	{\bfseries Lawrence Breen} & Professeur, Univ. Paris XIII & Président\\
	{\bfseries Bernhard Keller} & Professeur, Univ. Paris VII & Rapporteur\\
	{\bfseries Georges Maltsiniotis} & Directeur de Recherche, Univ. Paris VII & Examinateur\\
	{\bfseries Carlos Simpson} & Directeur de Recherche, Univ. Nice & Rapporteur\\
	{\bfseries Joseph Tapia} & Professeur, Univ. Paul Sabatier & Examinateur\\
	{\bfseries Bertrand Toën} & Chargé de Recherche, Univ. Paul Sabatier & Directeur\\
\end{tabular}
}
\vfill
{\sffamily \normalsize
Institut de Mathématiques de Toulouse, UMR 5580, UFR MIG\\
Laboratoire \'Emile Picard,\\
Université Paul Sabatier 31062 TOULOUSE Cédex 9
}
\end{center}
\newpage
\thispagestyle{empty}
~~
\newpage
\thispagestyle{empty}


\title{CHAMPS DE MODULES\\ DES CAT\'EGORIES LIN\'EAIRES\\ ET AB\'ELIENNES\\}
\author{{\Large Mathieu Anel}\\{\sf (anel@math.ups-tlse.fr)}}
\date{}
\maketitle
\thispagestyle{empty}


\newpage
\thispagestyle{empty}
~~
\newpage
\thispagestyle{empty}
~~
\vspace{60pt}
\begin{quote}
{\sffamily
\footnotesize
<<~J'ai trouvé mon île au trésor. Je l'ai trouvé dans mon monde intérieur, dans mes rencontres, dans mon travail. Passer ma vie avec un monde imaginaire a été mon île au trésor. Bien sûr, c'est vrai que les mondes que je visite au hasard de mes recherches peuvent parfois être jugés puérils ou inutiles, tant ils sont éloignés des préoccupations quotidiennes, mais quand aujourd'hui je repense à ceux qui m'accusaient d'être inutile, alors, vis-à-vis d'eux, je n'ai pas seulement le plaisir d'être inutile, mais aussi le désir d'être inutile.~>>

\hfill Hugo Pratt, {\em Le désir d'être inutile -- entretients avec Dominique Petitfaux}, \'Editions Robert Laffont, 1991.}
\end{quote}

\vspace{60pt}
\begin{quote}
{\sffamily
\footnotesize
<<~Mais je ne suis pas sûr que dans un univers où tous les phénomènes seraient régis par un schéma mathématiquement cohérent, mais dépourvu de contenu imagé, l'esprit humain serait peinement satisfait. Ne serait-on pas alors en pleine magie ? Dépourvu de toute possibilité d'intellection, c'est-à-dire d'interpréter géométriquement le schéma donné, ou l'homme cherchera à se créer malgré tout par des images appropriées une justification intuitive au schéma donné, ou sombrera dans une incompréhension résignée que l'habitude transformera en indifférence.~>>

\hfill René Thom, {\em Stabilité et morphogenèse}, \'Editions Dunod, 1984.}
\end{quote}

\vspace{60pt}
\begin{quote}
{\sffamily
\footnotesize
<<~L'une des conditions essentielles au travail mathématique est de savoir cerner ce à quoi on est exactement en train de penser.~>>

\hfill Michel Broué, {\em Une école de liberté} dans {\em Le goût de la science (textes rassemblés par Julie Clarini)}, \'Editions Alvik, 2005.} 
\end{quote}

\newpage
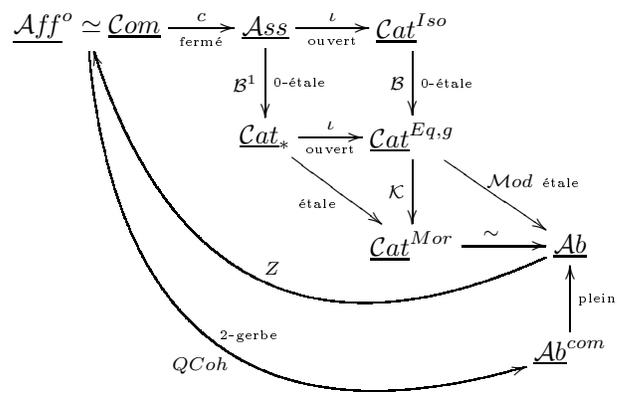
\thispagestyle{empty}
~
\newpage
\thispagestyle{empty}
~


\vfill
$$\xymatrix{
\underline{\mathcal{A}f\!f}^o\simeq \ukcom \ar[r]^-c_-{\textrm{{\tiny fermé}}}\ar@/_6pc/[rrrddd]_{QCoh}^{\textrm{{\tiny 2-gerbe}}} & \ukass \ar[d]_{\B^1}^-{\textrm{{\tiny 0-étale}}} \ar[r]^-\iota_-{\textrm{{\tiny ouvert}}} & \ukcatiso \ar[d]_\B^-{\textrm{{\tiny 0-étale}}}\\
& \ukcat_* \ar[rd]_-{\textrm{{\tiny étale}}} \ar[r]^-\iota_-{\textrm{{\tiny ouvert}}} & \ukcateqg \ar[d]_-{\mathcal{K}} \ar[rd]^-{\MMod\ \textrm{{\tiny étale}}}  \\
&  & \ukcatmor  \ar[r]^{\sim} & \ukab\ar@/^5pc/[llluu]_Z\\
&&&\ukab^{com}\ar[u]_-{\textrm{{\tiny plein}}}
}$$
\vfill
~
\vfill
~

\newpage
\thispagestyle{empty}


\chapter*{Remerciements}
\thispagestyle{chheadings}

\begin{quote}
{\sffamily
\footnotesize
\hfill <<~Without a doubt~>>

\hfill {\em Magic eight ball}, septembre 2003.}
\end{quote}
\vspace{30pt}

Ma thèse aura été un investissement psychologique, physique, éventuellement idéologique, largement ponctué de moments formidables, d'amitiés et de doutes que j'ai envie de retracer un peu ici. 
Je vais être long, mais, hormis la joie de faire des phrases de français, il faut bien donner de quoi lire à tous ceux qui ne comprendront rien à ce que je raconte dans ce travail ; je pense particulièrement, et affectueusement, à ma famille et à nombre d'amis\footnote{dont certains sont mathématiciens !} pour qui ce travail restera à jamais obscur. Tous se sont montrés curieux de mes travaux et c'était toujours à grande peine que je leur répondais qu'il m'était impossible de leur expliquer ce que je faisais ("but it is a heavy burden"), que ce texte soit pour eux le prétexte à découvrir au moins les circonstances de ma thèse.

En outre, ces remerciements sont aussi le lieu où se reposer un peu de la froideur académique des publications scientifiques ; la connaissance est joyeuse, sociale, elle est faite pour les humains par les humains, si elle possède ses indispensables protocoles, il est appréciable de temps en temps de les oublier.

\medskip
Ainsi, mes premiers remerciements seront pour ma famille. J'ai eu la chance que mes parents me laissent libre dans mon choix de faire des études et dans le choix d'icelles. Ils m'ont soutenu financièrement pendant mes années d'université, jusqu'à mon allocation de thèse ; mais plus que cela, je leur dois une bonne part de ma vocation scientifique. J'ai par exemple le souvenir d'un petit livre de cosmologie trouvé dans la bibliothèque, dont les pages, à force que je le feuillette, ne tenaient plus ; ou des nombreux jeux de constructions et autres microscopes qui ont éveillé mon esprit à une certaine logique ; ou encore des colonies de vacances passées à construire des fusées.

L'origine consciente de ma vocation aura été mon instituteur de {\sc cm2}, M. Barba, qui nous faisait faire de nombreuses expériences/observations scientifiques (électriques ou biologiques). Après lui, je suis rentré au collège attendant impatiemment les cours de physique. À l'époque, les mathématiques, où je ne rencontrais pas de difficultés particulières, m'indifféraient assez. Il a fallu attendre la fin du lycée pour que j'y rentre par l'intermédiaire de la théorie du chaos et des fractales, notamment grâce à Grégoire Winterstein. L'exposé fait sur le sujet dans le cours de philosophie de M. Bories ainsi que les encouragements de Mme Doucet et Mme Fouilleron, mes enseignantes de mathématiques et de physique en terminale, m'ont aussi largement assuré dans mon choix d'orientation.

C'est en classes préparatoires où ma relation avec la physique a été malmenée : ma rencontre avec les équations de Maxwell et particulièrement le laplacien a été fatale à ma capacité de compréhension en physique. D'une décision que je n'ai jamais regrettée, je me suis alors tourné entièrement vers les mathématiques, persuadé que la compréhension de la physique passait d'abord par celle des maths. Sans vouloir être catégorique dans cette opinion, je n'ai jamais trouvé à la repenser, et c'est la même volonté de bien comprendre les choses qui m'aura poussé jusqu'à finir une thèse dans des mathématiques (soi-disant) très pures.

Je voudrais que mes parents trouvent dans ce petit texte et dans ce travail, la récompense de leur efforts d'éducation et l'expression de ma gratitude à leur égard. Enfin, parce que je n'ai peut-être pas été très facile à vivre pendant ces années d'études, particulièrement celles de thèse, j'exprime à ma famille mes regrets et excuses.

\bigskip
C'est une grande joie pour moi que de pouvoir remercier mon directeur de thèse Bertrand Toën.
Je me souviens de moi lorsque je suis allé le voir en septembre 2003 pour lui demander un sujet de thèse, tout intimidé par la réputation qu'il avait laissé parmi les thésards du laboratoire ; il venait à peine de revenir sur l'université et a accepté du jour au lendemain de travailler avec moi.
Très vite, il s'est avéré que j'étais plein de questions et qu'il était plein de réponses ; Bertrand a eu durant nos années de collaboration cette compétence remarquable de comprendre mes points de vue et de devancer mes propres interrogations. Cela m'aura été d'une double économie : je ne perdais pas trop de temps à chercher les réponses à mes questions, et surtout à chercher les bonnes questions ! Il ne me restait plus qu'à ruminer quelques jours la remarque avant de me rendre compte qu'il avait touché le bon point.

Je le remercie de son incroyable disponibilité : à moins qu'il ne travaillât déjà avec quelqu'un d'autre, il ne m'a jamais demandé de repasser. Il s'est également montré particulièrement compréhensif et encourageant lors de mes moments de doutes.
\'Ecrivant ces lignes, il me revient une phrase de Bertrand issue de sa propre thèse, où il remerciait ses directeurs pour l'avoir encadré avec <<~beaucoup d'humanité~>> ; sept ans après, c'est une nouvelle génération qui lui retourne le compliment.
Je suis heureux d'avoir trouvé autour de lui toute une troupe de personnes très sympathiques et abordables (je pense particulièrement à Gabriele Vezzosi et Michel Vaquié), tranchant agréablement l'image de la sociologie officielle d'un milieu scientifique mené par les luttes d'influences des egos.

\medskip
Je me suis tourné vers mon domaine mathématique, la géométrie homotopique, à la suite d'une volonté de géométriser et de lutter pour le sens de toute cette algèbre abstraite (faisceaux, cohomologie, champs...) qui prétend qu'elle parle de géométrie mais qui n'explique jamais comment ! Au fond, je m'interrogeais maladroitement sur la notion d'espace et j'avais bien cru trouver dans les champs un début de réponse à cette question encore mal formulée. 
Travailler avec Bertrand ne m'a pas trompé, j'ai appris maintenant que la réponse croise étonnamment la topologie avec les catégories supérieures. Plus largement, j'ai réussi grâce à lui à répondre à de nombreuses questions que je me posais, et dont l'absence de réponses m'aurait condamné à une certaine angoisse sinon physique du moins métaphysique.

Travailler dans ce domaine m'aura donné un point de vue sur les mathématiques assez particulier où l'enjeu n'est pas mis dans les démonstrations -- ce que je conçois comme une vision logiciste très réductrice -- mais dans les définitions des objets. Pour philosopher le temps d'une phrase, les mathématiques ont cela de particulier que les objets qu'elles étudient sont exactement leurs définitions ; en conséquence, énormément de problèmes (irrégularité de propriétés, représentabilité d'objets, ...) peuvent en fait se résoudre "simplement" en trouvant une nouvelle définition, déformant légèrement l'initiale, mieux, {\sl élucidant} les limitations de l'initiale et les problèmes encourus. Un bel exemple d'une telle construction est le traitement homotopique des quotients, dont la version géométrique forme la théorie des champs. 

Dans ce domaine, les nombreuses conversations avec Joseph Tapia -- dont la compétence s'exerce particulièrement par la recherche d'un discours précis et juste sur les choses -- sur les topos, la logique, les maths quantiques et la physique, l'épistémologie et autres sujets auront toujours été enrichissantes. Avec cela, ses blagues inattendues m'auront toujours pris au dépourvu. Je suis particulièrement content qu'il participe à mon jury.

\bigskip

Carlos Simpson et Bernhard Keller ont accepté de rapporter cette thèse et je les remercie de leurs rapports.
Je les remercie aussi avec Georges Maltsiniotis de leur invitation à Nice et à Paris où j'ai pu exposer mes résultats, ce qui est toujours une stimulation importante. Je suis très heureux que tous les trois aient acceptés de faire parti du jury.
Enfin, je remercie beaucoup Lawrence Breen d'avoir accepté de présider mon jury.

\bigskip

Deux personnes méritent une place spéciale dans ces remerciements. 
Tout d'abord, je pense à David Gary, croisé au détour d'un cercle mathématique, qui se sera ouvert en un long chemin en ligne droite.
David est la personne qui a le plus partagé mes joies et mes angoisses mathématiques ; ses connaissances philosophiques et sa passion épistémologique m'ont, à travers nos nombreuses discussions, aidé à appréhender certaines questions que peu de mathématiciens se posent (et dont encore moins en cherchent les réponses). En m'en rendant compte assez tard, je me suis, parmi mes interrogations, largement lancé à la poursuite de cette chimère qu'est le rapport des mathématiques au réel. Je pense maintenant qu'il n'y a pas tant de mathématiques que de {\sl mathématisation} et réussir à (commencer de) dégager ce point de vue sur les mathématiques aura été une de mes plus grandes satisfactions de mon travail durant ces années de thèse. De surcroît, je pense, grâce à David, ne pas avoir oublié les raisons et questions qui m'avaient poussé à m'engager dans des études scientifiques, et de mathématiques en particulier. Je lui dédie la citation de Michel Broué.

Je pense ensuite à Damien Calaque, croisé sous ce soleil de Luminy si favorables aux mathématiques, à qui je dois peut-être un peu plus qu'il ne le croit. Damien a toujours montré un intérêt enthousiaste (pondéré par son flegme naturel) pour mon type de mathématiques. Là où la plupart des mathématiciens de son domaine, à propos de la théorie des champs, se contentent d'ânonner un intérêt formel, il a, lui, été plus sensible aux vertus d'une "vraie" approche géométrique. Son intérêt a souvent été stimulant pour moi et nos conversations, mathématiques ou autres, ont toujours été un plaisir. Je ne peux que nous souhaiter une prospère collaboration ! (si on se décide à s'y mettre...)

\bigskip

Il est temps de passer enfin aux remerciements de mes compagnons d'études au laboratoire \'Emile Picard.
Alors que j'étais en {\sc dea}, en 2002, il s'est créé une ambiance particulièrement chaleureuse et exceptionnelle entre les doctorants.
Sans regrouper l'ensemble des thésards du laboratoire, nous étions néanmoins un grand groupe uni, d'abord autour de notre groupe de travail (le {\sc gdte}, prononcé "g'deuteu"), puis après par amitié. Je leur dois à tous beaucoup sur les plans mathématiques, affectifs et humain.

Mes premières salutations seront pour mes camarades de promotions (et leur compagne) : 
Guillaume Rond (et la douce Opaline Chesneau) "chef" des thésards aux attributs supérieurs et grand parfumeur de rhum dont la compétence à mentir m'aura toujours impressionné ; Nicolas Puignau (et... on ne va citer qu'Aurélie Cavaillé et son saxo) qui fût le premier à m'intégrer dans le groupe du {\sc dea} et à qui je dois un certain retour à la vie sociale ; et Olivier Drévillon, mon cobureau pendant trois ans, tout aussi avide que moi de comprendre la physique ou le calcul stochastique.

Il me faut souligner le rôle important qu'auront eu Guy Casale-Zawichost, Emmanuel Opshtein et Julien Keller (feu le bureau 31) durant leurs années au laboratoire. 
Guy a été l'animateur charismatique de notre groupe ; bon vivant, il sait générer autour de lui une ambiance extraordinairement enthousiaste, propre à la réalisation de tous les projets ! Je lui dois une certaine décomplexion et je le remercie en particulier pour tous les conseils pertinents qu'il m'a donné.\footnote{Je remercie aussi Belinda Zawichost-Casale, sa simplicité et sa bonne humeur auront toujours été apaisantes et appréciées.}
Manu a lui été l'animateur des premières années du "g'deuteu", sa viabilité ainsi que la diversité des thèmes abordés sont dû en bonne part à sa soif de mathématiques. Je lui dois quelques bonnes luttes dans la neige et  randos à pieds ou en raquettes. Je tiens aussi à mentionner son formidable style d'exposé.
Enfin, Julien (avec qui on m'a si souvent confondu...) aura eu cette vertu admirable, en acceptant la présidence du Collectif de Doctorants Toulousain, de politiser notre bande et de la sensibiliser aux problèmes de l'université et de la recherche en général. Je le remercie de son engagement, qui a été courageux et fructueux, il m'a enseigné beaucoup de choses et je lui présente mes excuses quant à mon militantisme approximatif.

Leur passion des maths à tous les trois a été la source de la création et de la richesse du "g'deuteu" ainsi que du groupe de travail sur le programme de Mori, où nous avons tous appris, et fait apprendre, beaucoup de belles mathématiques.

\medskip

Je remercie également Laurent Mazet, coconstructeur du site des doctorants picardiens et partenaire de trop peu de balades en montagne ; Grégoire Montcouquiol, pour \verb+gy!be+, ses soirées jeux, son pistolet à fléchettes et pas pour ce cocktail (gris!) qu'il me fit un jour ; Vanessa Vitse, pour son caractère, ses soirées et le pousse-rapière ! Yohann Genzmer pour son humeur chantante et jamais renfrognée (et Johanna Chappuis à qui je dédie un sourire et un clin d'oeil) ; ainsi que Mathieu Fructus (dont je suis très content de ne pas avoir suivi le conseil qu'il me donna un jour) pour les rigolades sur les sites web débiles et pour tout le travail effectué dans ses fonctions de représentant des thésards auprès du laboratoire.
Je termine sur une pensée pour d'autres thésards et post/pré-docs, moins fréquentés mais que je me reprocherais d'oublier : Magali Bouffet, Arnaud Hilion, Sonia Rafi, Eva Miranda et Sergi Simon.

Je garde d'eux tous d'inoubliables souvenirs de randonnées à la Llagonne, au Canigou ou en Corse, de descentes en luge, de cafés et discussions sur la terrasse de l'upsidum, de squattage de bureaux, de déconnades sur le web, de bières du café pop', de manifs, de goûters, de barbecues et autres soirées.

\medskip
Les circonstances ont voulu que je fréquente également la nouvelle génération de thésards : Cécile Poirier, dont l'énergie et la gentillesse ont parfumé notre couloir d'une nouvelle fraîcheur, Julien Roques, à qui je lègue le miroir de mon bureau, Anne Granier, Landry Salle, Benjamin Audoux, Michaël Ayoul, Maxime Rebout et les autres.
Je leur lègue ces conseils : investissez-vous dans le labo, ne laissez pas mourir le "g'deuteu" et boycottez le {\sc cies} (oui, c'est dur quand on a tous l'agreg mais fallait pas la passer !)

\bigskip
J'en arrive enfin à remercier mes amis proches, qui ont suivit mes travaux sans pouvoir les comprendre. 
La première en liste est Anne-Lise Marchand que je remercie vivement et affectueusement pour sa profonde amitié et son soutien constant. Elle fait partie des gens qui ont le mieux suivi ma vocation et qui l'ont le mieux comprise. Je lui dédie la citation de Hugo Pratt.
Le second est Grégoire Winterstein, que je dois remercier pour trop de raisons (du livre de Gleick aux cours de l'{\sc esm2} en passant par les Monty Python et le Pendule) qui font notre complicité depuis toutes ces années ("it's truly a real honorable experience, a wonderful and warm and emotional moment for all of us") et promis ! la prochaine fois j'enlève mes chaussures. Je lui dédie le diagramme kabbalistique de la page {\sf 5} : lu dans un sens il guérit le cancer et dans l'autre la bêtise.
D'une époque révolue, je remercie Paul Anel, Gwenhaël et Mikaël Finkel, vieux camarades de jeux de rôles, pour les nombreuses conversations sur tout et n'importe quoi qui ont exercé mon esprit à l'interrogation.
Je remercie Mathieu Richardoz, Mathieu \'Eludut, autres vieux camarades rôlistes, et plus récemment Yann Sauvage, à qui je dois de nombreux hébergements et des conversations toujours stimulantes. J'ai aussi une pensée pour \'Eliette Ferré qu'un curieux destin aura voulu lier au début et à la fin des mes études mathématiques.

Je remercie Ivane Pairaud de son investissement dans le {\sc cdt} et de son amitié, notamment ces derniers mois.

Enfin je voudrais mentionner mes compagnons d'études universitaires, dont certains que je n'ai pas vu depuis longtemps, avec qui j'ai largement travaillé mes {\sc td}s, exams et mémoires : Julien Orus, Caroline Begaud, Renaud Marty et Cédric Pecqueur.

\bigskip

J'estime avoir été très chanceux par rapport aux choix ou aux rencontres que j'ai faits, aux idées que j'ai apprises ou dégagées ; globalement, tout cela me satisfait et je le dois à toutes les personnes mentionnées.

\bigskip
\hfill --- \hfill~
\bigskip

<<~Mais à quoi sert tout ça~?~>> Question embarrassante, qu'on m'a beaucoup posé, et à laquelle je crois important d'apporter quelques éléments de réponses. 

La première réponse que je donne régulièrement aux non mathématiciens est que ça ne sert à rien. D'abord, leur expliquer les intérêts disciplinaires serait vain et puis, au fond, pas très juste, ceux-ci ressemblant souvent à des prétextes, déguisant des intérêts personnels du scientifique, qui ne travaille pas tant pour la science que pour sa curiosité et son plaisir. On peut aussi arguer d'éventuelles applications en physique ou dans d'autres domaines, mais, à moins de travailler déjà dans cette perspective, cela ressemble plus à un souhait qu'à une motivation.
Le travail d'un scientifique, au fond, lui sert d'abord à lui : cherchant ses propres réponses, écoutant celles de ses collègues, inventant des unifications, il construit sa vision des choses, et c'est bien là ce qu'on doit attendre de lui : qu'il ait un point de vue sur les choses, qu'il soit {\sl savant}. La question n'est plus alors <<~à quoi ça sert~?~>> mais <<~à qui ça sert~?~>>.

Dans sa forme enthousiaste, la question est souvent une forme maladroite de la question <<~de quoi ça parle~?~>> à laquelle il est encore plus dur de répondre ! Les mathématiques n'étant pas une science de la nature, elles n'ont pas de "problématique" facile à résumer et ne forment un tout uni que via le partage d'une méthode et non d'un substrat d'étude.

Plus insidieusement, la question peut-être aussi le reflet d'un certain discours économico-utilitariste, parfois posée pour humilier le chercheur, qui, sauf s'il travaille sur le cancer ou une autre généralité médiatique, aura généralement du mal à justifier l'utilité de ses recherches. Ce type de discours oublie insidieusement le mot {\sl savant} pour le terme asservissant de {\sl chercheur}.
Hugo Pratt, à qui on reprochait de n'être qu'un artiste, répondait lui par un <<~désir d'être inutile~>>. Je le suivrais assez dans cette provocation, si la question est mal intentionnée, elle ne mérite qu'une réponse narguante, et afficher une satisfaction à son inutilité me semble tout à fait adéquat.

\medskip

Enfin, la même question dans la bouche d'un mathématicien est encore plus embarrassante, parce qu'il faut avec les collègues encore plus sauver les apparences ! Mais ceci est une autre histoire.

\newpage
\thispagestyle{empty}
~~
\newpage
\thispagestyle{empty}

\vspace*{4cm}
\begin{flushright}
{\sffamily \slshape Je dédie cette thèse à toutes les bonnes idées que j'ai eues pendant son élaboration\\
et qui n'en furent pas.}
\end{flushright}

\newpage
\thispagestyle{empty}

\newpage
~~
\thispagestyle{empty}

\tableofcontents


\setcounter{tot}{0}
\newpage
\newenvironment{irthm}[1][]%
{\refstepcounter{tot}\par\medskip \noindent{\sffamily \bfseries TH\'EOR\`EME~\thetot~}{\sffamily (#1)} \it}%
{\par\medskip\normalsize}

\newpage
\thispagestyle{empty}
\chapter*{Introduction}
\addcontentsline{toc}{chapter}{Introduction}
\thispagestyle{chheadings}

\begin{quote}
{\sffamily
\footnotesize
\hfill <<~Tu verras, quand on sait ce qu'est un champ, on en voit partout !~>>

\hfill Bertrand Toën, automne 2003.}
\end{quote}
\vspace{30pt}

Ce travail étudie plusieurs problèmes de modules pour les catégories linéaires et abéliennes, nos {\sl résultats principaux} sont des théorèmes d'existence pour les espaces de modules concernés (théorème \ref{thmintrogeom}, {\em infra}) avec le calcul de leurs tangents (théorème \ref{thmintrotangent}), ainsi qu'un théorème de comparaison entre les deux plus importants de ces espaces de modules (théorème \ref{thmintroetale}).

\bigskip

La première partie de cette introduction replace le contexte des champs supérieurs puis celui des déformations des structures associatives et des catégories linéaires ; la seconde présente les résultats obtenus dans cette étude ; et la dernière explique le plan des chapitres de la thèse.

\section{Contexte}

On retrace d'abord le développement des champs dans le cadre des problèmes de modules.
Puis, on définit le problème des modules de catégories linéaires, pour lequel on rappelle les résultats de W.~Lowen et M.~Van den Bergh.


\subsection{Problèmes de modules et champs supérieurs}

\paragraph{Problèmes de modules et faisceaux}
Comme les ensembles, les espaces géométriques peuvent se définir en extension ou en compréhension. La définition en extension correspond, par exemple, à la donnée explicite d'un atlas pour l'espace, la définition en compréhension consiste, elle, en la donnée d'un sens aux points de cet espace, \ie d'un problème duquel les points de l'espace forment les solutions. Un tel problème, où les solutions forment un continuum, s'appelle un {\sl problème de modules} et sa solution un {\sl espace de modules}.


Les espaces de modules sont légions (espaces projectifs, grassmaniennes, schémas de Hilbert, modules de fibrés vectoriels, de courbes, espaces de Teichmüller, plan de Mandelbrot, déploiements de singularités, etc.) et depuis Riemann et ses surfaces, qui sont à l'origine de la notion, la construction-description de ces espaces est devenu un enjeu.

Nombre de problèmes de modules apparaissent particulièrement en géométrie algébrique où la "grossièreté" de la topologie de Zariski complique l'étude en extension ; et depuis les techniques développées dans les années 60 autour de A.~Grothendieck, l'enjeu de la construction s'est déplacé. Ces techniques, utilisant largement le langage des catégories à travers celui des {\sl faisceaux}, permettent de toujours définir l'objet répondant à un problème de modules : le point de vue du {\sl foncteur des points} caractérise un objet par ses ensembles de points, et c'est presque tautologiquement qu'un problème de modules donne un foncteur des points. L'enjeu n'est plus alors de construire l'espace des modules, qui existe toujours en tant que faisceau, mais de savoir si ce faisceau est équivalent à celui d'un espace géométrique.

\paragraph{Champs}
Malgré ses succès, l'approche faisceautique a néanmoins vite trouvé ses limites, et s'est rapidement ouverte sur la {\sl théorie des champs}. D'une manière sous-estimée, nombre de problèmes de modules ont des solutions qui forment naturellement des catégories et non plus des ensembles, et ne retenir de ces catégories que les ensembles de classes d'isomorphie d'objets s'est avéré créer des irrégularités dans les espaces de modules déduits, leur interdisant d'être des espaces classiques (schémas, variétés).

L'idée s'est alors imposée que, pour régulariser les propriétés des espaces de modules, il fallait conserver la trace de la structure catégorielle des points de l'espace des modules. Les {\sl champs} sont ainsi ce que deviennent la notion d'espace topologique (en fait celle, plus subtile, de faisceau) lorsque leurs points ne forment plus des ensembles {\sl mais des catégories}. Ceci dit, comme précédemment avec les faisceaux, un problème de modules définit toujours canoniquement un champ, et, avant même de se demander s'il possédait de bonnes propriétés,  il a fallu inventer, parce que les champs étaient de nouveaux objets, la classe des {\sl champs géométriques}, \ie ceux parmi les champs suffisamment proches des espaces géométriques habituels pour que les constructions connus sur ces derniers s'y généralisent. Ce sont, par exemple, les définitions de P.~Deligne et D.~Mumford \cite{dm} ou de M.~Artin \cite{ar} (voir aussi \cite{lmb}).

Ainsi certains problèmes de modules trouvent leur solution non plus sous la forme d'un espace de modules mais d'un {\sl champ de modules}. 

\paragraph{Champs supérieurs}
D'une manière peut-être un peu décevante, les champs ont, comme les faisceaux, trouvé leur limites. Inventés, entre autres raisons, pour représenter certains quotients de groupoïdes, il se trouve que le quotient d'un groupoïde de la catégorie des champs n'est pas en général un champ, l'obstruction étant la possibilité d'obtenir pour les points du quotient des 2-catégories et non plus seulement des (1-)catégories \cite{breen}. \'Egalement, certains problèmes de modules comme celui des complexes traité dans \cite{hsimpson}, où celui des catégories traité dans ce travail, qui sont des objets de catégories supérieures ne trouvent des réponses qu'hors des champs.

Le dernier ingrédient de la théorie des champs a donc été la considération de catégories supérieures : les champs supérieurs sont les objets topologiques dont les points forment des catégories supérieures. Bien que les définitions de ces objets aient stagné un moment faute d'une bonne axiomatisation des catégories supérieures, une notion de {\sl champ supérieur géométrique} a finalement été proposé par C.~Simpson \cite{simpson}, reprenant la définition d'A.~Joyal des champs comme certains préfaisceaux simpliciaux.

\medskip

Il est notable que les exemples de catégories supérieures en tant que telles ne sont pas fréquents, et qu'on les obtient beaucoup plus souvent par la considération de 1-catégories dans lesquelles on a une notion d'équivalence plus faible que l'isomorphie. L'inversion formelle des flèches de la catégorie qui sont des équivalences donne en effet lieu à un phénomène des plus remarquables : l'ambiguïté à décrire explicitement les flèches de la catégorie localisée enrichit celle-ci en une catégorie simpliciale, \ie un modèle de catégorie supérieure pour lequel toute les flèches de degré supérieur à 2 sont inversibles.

Afin de travailler avec ces localisations, il a été inventé par Quillen le formalisme des catégories de modèles \cite{quillen, hovey, hirschhorn} qui sont devenues des objets de pratique quotidienne pour les manipulateurs de catégories supérieures.
Peut-être encore plus pertinent que la notion de champ supérieur semble alors être la notion de {\sl champ de Quillen} introduite par A.~Hirschowitz et C.~Simpson dans \cite{hsimpson}, qui permet de manipuler les champs supérieurs à partir de simples préfaisceaux en catégories. Le cadre d'apparition des problèmes de modules (axiomatisé un peu plus bas) fait de la notion de préfaisceau de Quillen une notion merveilleusement adaptée\footnote{Cette notion a tout de même ses limites et s'avère inadaptée pour l'étude de certains points, par exemple les auto-équivalences d'un champ.}.



\medskip

La notion de champ supérieur est cette fois stable par quotient de groupoïdes, précisément, la catégorie des champs supérieurs possède une caractérisation par des axiomes copiant ceux de Giraud pour les topos : comme les catégories de faisceaux sont universelles (essentiellement) pour l'effectivité et l'universalité des relations d'équivalence, les catégories de champs sont universelles pour l'effectivité et l'universalité des relations d'équivalence supérieures, \ie des groupoïdes. C'est la notion de topos supérieur de \cite{lurie1, hag1}.

\paragraph{Après les champs}
Malgré que cette caractérisation en termes topossiques semble clore le développement des objets nécessaires à la géométrie et aux problèmes de modules, la nécessité d'autres élargissements se fait déjà sentir.
Par exemple, la notion de champs supérieur dont nous parlons ici est celle de champs en infini-groupoïdes, et il semblerait utile de développer une notion de {\sl champs en catégories supérieures} qui ne soient pas forcément des groupoïdes (par exemple en catégories simpliciales ou de Segal)\footnote{Dans les problèmes de modules des catégories étudié dans cette thèse, il n'est, au fond, pas très naturel de se limiter aux seuls morphismes qui soient des équivalences.}.

Dans d'autres directions, la philosophie de la {\sl géométrie algébrique dérivée} à la Kontsevich-Deligne-Kapranov, qui appelle une notion de {\sl champs dérivés} formalisée récemment dans \cite{lurie2, hag2}), ou la géométrie algébrique sous $Spec(\mathbb{Z})$ de \cite{tvaq}) attendent leurs développements.\footnote{Aura-t-on besoin un jour d'autres objets que des {\sl champs-en-catégories-supérieures-dérivés-sous-$Spec(\mathbb{Z})$} ?}

\paragraph{Complexe tangent}
Les champs (supérieurs ou non) sont des objets très différents des espaces géométriques habituels et attraper leur intuition est difficile\footnote{Peut-être, d'ailleurs, un élément retenu contre l'utilisation des champs dans certains domaines mathématiques est le fait que, contrairement aux espaces, ils ne se "dessinent" pas ; si on pense aux champs avec une intuition géométrique, c'est sans autre support dessiné que celui des diagrammes commutatifs.}.
Outre la nature de leurs points, la particularité la plus facile à attraper des champs est sans doute la structure de leur tangent en un point. Comme ensemble de certains points, cet objet doit être une catégorie et elle doit de plus posséder la structure linéaire propre aux tangents. De tels objets sont modélisables par des complexes de modules (en degrés négatifs) et on parlera donc du {\sl complexe tangent} d'un champ en un point.

De tels objets sont en fait fréquents en mathématiques : tous les complexes de cohomologie obtenus pour résoudre des problèmes de déformations infinitésimales sont des complexes tangents de certains champs, cette thèse l'illustre avec le complexe de Hochschild.


\paragraph{Champs de modules}\label{introquillen}
L'analyse des exemples amène à voir qu'un problème de modules consiste essentiellement en les données suivantes :
\begin{itemize}
\item un {\sl site} $S$ d'objets géométriques de référence (\ie connus) ;
\item une notion de {\sl famille} de solutions au problème paramétrée par un objet du site, laquelle va toujours de pair avec une notion de morphisme entre telle familles ;
\item une notion d'{\sl équivalence} sur ces familles, traduisant la relation pour laquelle on veut confondre deux familles (souvent enrichissable en une structure de modèles) ;
\item et des foncteurs de {\sl changement de base} le long de chaque morphisme du site, compatible avec les équivalences (éventuellement quitte à les dériver), pour lesquels la pratique montre qu'ils ne vérifient en général que des fonctorialités faibles par rapport à la composition dans $S$.
\end{itemize}
Cette axiomatisation basique associe à un problème de modules un préfaisceau sur le site de base dont les valeurs sont des catégories avec équivalences (voire des catégories de modèles), à partir duquel on peut construire un champ de modules\footnote{L'expression {\sl champs de modules} prend le mot {\sl champ} au sens des champs supérieurs.} répondant au problème des modules (cette construction est détaillée au {\sf\S \ref{chpmodchpmod}}).

\subsection{Modules de catégories linéaires}

\paragraph{Des modules associatifs aux abéliens}

La formalisation des problèmes de modules en termes de champs a eu pour origine des problèmes de modules géométriques  dans les années 60 (espaces de fibrés, de courbes, etc.), et il est assez normal que ce soient des géomètres qui aient dégagé la notion de champ. D'un autre côté, les problèmes de modules se posaient aussi pour les structures algébriques, et, en l'absence des champs, ils ont été particulièrement étudié à un niveau infinitésimal où la linéarisation des phénomènes permet la résolution du problème par des complexes de {\sl cohomologie de déformations}. Les exemples sont nombreux : groupes Ext pour les modules, cohomologie des groupes, cohomologie de Hochschild pour les algèbres associatives, de Harrison pour les algèbres commutatives, de Chevalley-Eilenberg pour les algèbres de lie, d'André-Quillen pour les points dans un schéma, etc.

Dans ce cadre, le problème des modules d'algèbres associatives est posé implicitement depuis très longtemps (l'article de Hochschild date de 1945) mais il faut attendre les travaux de M.~Gerstenhaber et S.~Schack \cite{gerst} pour élargir la cohomologie de Hochschild aux diagrammes d'anneaux et aux catégories linéaires.
Depuis la cohomologie de Hochschild a connu des expansions à d'autres notions essentiellement associatives ($A_\infty$-algèbres, dg-catégories, etc.) dont la plus récente est peut-être la définition de W.~Lowen et M.~Van den Bergh de la cohomologie de Hochschild d'une catégorie abélienne dans \cite{lvdb2}.

\medskip

Il est curieux de noter que les degrés des groupes de cohomologie classifiant les déformations infinitésimales varient avec le problème : c'est le 0-ième pour la cohomologie d'André-Quillen, le premier pour celle des Ext, et le second pour les autres exemples cités.

\paragraph{Modules associatifs et champs supérieurs}
Une des motivations primordiales pour ce travail a été la construction d'exemples de problèmes de modules dont les solutions soient données par des champs supérieurs.

Contrairement aux champs en groupoïdes, ou 1-champs, relativement bien adoptés par la communauté mathématique, les champs supérieurs restent encore des objets peu vulgarisés et dont on connaît peu d'exemples. À cet égard, celui des modules associatifs est peut-être le plus frappant. 

Il est en effet bien connu que les déformations au premier ordre de structures associative sont classifiées par le {\sl deuxième} module de cohomologie de Hochschild, et que cela s'explique par le fait que la catégorie des algèbres associatives peut se voir comme une 2-catégorie, en distinguant isomorphismes intérieurs et extérieurs. Si, dans le simple langage des algèbres, la construction parait artificielle, elle trouve un cadre naturel dans un élargissement aux catégories linéaires et leurs équivalences (et c'est bien dans ce cadre que se placent M.~Gerstenhaber et S.~Schack dans leur papier fondateur \cite{gerst}). Les catégories linéaires formant des 2-catégories, leurs modules doivent donc former un 2-champ.

Généralisant la situation de cet exemple, on peut se risquer à établir le {\sl slogan} suivant pour les champs supérieurs : les objets dont les modules infinitésimaux sont classifiés par un $n$-ième module de cohomologie vivent naturellement dans des $n$-catégories et leurs modules forment des $n$-champs. Ceci élucide la remarque sur les degrés fait au paragraphe précédent.

Plus précisément, le complexe de cohomologie classifiant les déformations au premier ordre d'une structure $A$ s'interprète comme le {\sl complexe tangent} du champ des modules au point définit par $A$ (cf. {\S\ref{cplxtangent}}).


\paragraph{Déformations infinitésimales}
La théorie des déformations infinitésimales des catégories linéaires et abéliennes a été développée, dans les deux papiers de W.~Lowen et M.~Van den Bergh \cite{lvdb1,lvdb2}.

Un cadre général pour les déformations est celui des catégories fibrées (ou des préfaisceaux en catégories) :
si $C\to S$ est une catégorie fibrée et qu'on a un diagramme 
$$\xymatrix{
x\ar@{-->}[d]\ar[r]& y\ar@{-->}[d]& \in C\\
s\ar[r]^u& t &\in S
}$$
on dit que $y$ est une {\sl déformation de $x$ le long de $u$} si $x$ est la {\sl fibre de $y$ en $u$}, \ie si $x\to y$ induit un isomorphisme $x\to u^*y$ (\footnote{On peut aussi considérer un isomorphisme $u^*y\to x$, suivant les contextes il apparaît les deux sens.}) où $u^*$ est un foncteur de tiré-en-arrière le long de $u$ pour $C$.
Si $S$ est la catégorie des schémas affines, une déformation infinitésimale est une déformation le long d'un épaississement infinitésimal $s\to t$ de $s$, \ie d'un morphisme d'anneaux de noyau nilpotent.

W.~Lowen et M.~Van den Bergh définissent deux notions de déformations : la première pour les catégories linéaires et la seconde pour les catégories abéliennes. Dans un langage qu'ils n'adoptent pas, et conformément à notre présentation des problèmes de modules, ces deux définitions reviennent à considérer deux préfaisceaux en catégories différents au-dessus du site des schémas affines : le premier est celui dont les valeurs sont les catégories de catégories linéaires et où le changement de base le long d'un morphisme d'anneaux $A\to B$ est donné par le foncteur $-\otimes_AB$ (cf. {\sf\S\ref{modcateq}}) ; le second est celui dont les catégories de valeurs sont celles des catégories abéliennes et où le changement de base le long de $A\to B$ est donné par le foncteur $\uHom_A(B,-)$ (catégorie des $B$-modules) (cf. {\sf\S\ref{modcatab}}).

Leurs résultats principaux dans \cite{lvdb1,lvdb2} sont de montrer que, sous certaines restrictions indispensables de platitude, les déformations infinitésimales d'une catégorie linéaire ou abélienne sont contrôlées par un complexe de cohomologie de Hochschild \cite[thm. 3.1]{lvdb2}, et que les problèmes des déformations d'une catégorie linéaire $C$ et de la catégorie abélienne de ses modules $C\Mod$ sont équivalents \cite[thm. 8.16]{lvdb1}.

\paragraph{De la théorie infinitésimale à la globale}
Voulant construire un objet global classifiant les catégories linéaires et abéliennes, la contrainte naturelle sur celui-ci est que sa théorie infinitésimale (l'étude de son tangent en un point) redonne la théorie de W.~Lowen et M.~Van-den-Bergh. Or, une fois transposées dans le langage des préfaisceaux en catégories, leurs définitions offrent des candidats naturels pour les foncteurs des points des champs classifiant les catégories linéaires et abéliennes ; ce sont ces définitions qu'on reprend au chapitre {\ref{chapmodcat}}.

\section{Résultats}

Dans cette section le mot {\sl champ} est employé par défaut au sens de {\sl champ en infini-groupoïdes}.

\paragraph{Les différents modules des catégories}
Il a paru intéressant, tant qu'à construire les modules de catégories linéaires, de remarquer que celles-ci possédaient naturellement trois notions d'équivalences :
\begin{itemize}
\item l'isomorphie,
\item l'équivalence de catégorie,
\item et l'équivalence de Morita.
\end{itemize}
Nous avons donc construit aux {\sf \S\S \ref{modcatiso}, \ref{modcateq}} et {\ref{modcatmor}} les champs de modules pour ces trois relations d'identification, notés respectivement $\ukCATiso$, $\ukCATeq$ et $\ukCATmor$. En outre, chacune de ces relations élargissant la suivante, on dispose de morphismes naturellement surjectifs entre ces champs :
$$
\ukCATiso \surj \ukCATeq \surj \ukCATmor.
$$
Il est notable que, dans l'approche utilisée pour les construire, les foncteurs des points de ces objets ont les mêmes catégories de valeurs (celles des catégories linéaires) et les mêmes changements de base ; leur différence se fait seulement par les {\sl structures de modèles} mises sur les catégories de valeurs (cf. {\sf \S\ref{compstmod}}).

\medskip

Nous nous sommes également intéressés aux modules de catégories abéliennes (à équivalence de catégorie près) pour lesquels le champ $\ukAB$ correspondant est construit au {\sf \S\ref{modcatab}}. 
La notion de catégorie abélienne à laquelle on se limite est telle qu'elles possèdent toutes des petits générateurs et le morphisme naturel, issu de l'association à une catégorie de la catégorie de ses modules, est donc surjectif
$$
\MMod : \ukCATeq\surj \ukAB.
$$
Par définition de l'équivalence de Morita, se factorise $\ukCATmor$, et offre donc un diagramme de morphismes surjectifs
$$
\ukCATiso \surj \ukCATeq \surj \ukCATmor \surj \ukAB.
$$

\medskip

Il est remarquable qu'on ait l'équivalence suivante, qui est une traduction géométrique de la version à la Freyd de l'équivalence de Morita.

\begin{irthm}[thm. {\ref{thmequivmodulemoritaabel}}]
Le morphisme $\ukCATmor \surj \ukAB$ est une équivalence de champs.
\end{irthm}

L'équivalence passe par la construction d'une structure de modèles sur la catégorie des catégories linéaires pour laquelle les équivalences sont les équivalences de Morita, \ie les foncteurs $f:C\to D$ qui prolongés aux catégories de modules $f_!:C\Mod\to D\Mod$ sont des équivalences de catégories.

\begin{irthm}[thm. {\ref{thmstmodelemorita}}]
Pour un anneau commutatif $A$, il existe sur la catégorie $A\z\CAT$ des catégories $A$-linéaires une structure de modèles au sens de  engendrée par cofibrations (cf. \cite{hovey}) pour laquelle les équivalences sont les équivalences de Morita et les objets fibrants les catégories karoubiennes (cf. définition \ref{defkar}).
\end{irthm}

\medskip

Cette équivalence entre $\ukCATmor$ et $\ukAB$ doit être mise en perspective avec la démonstration par W.~Lowen et M.~Van den Bergh de l'équivalence des déformations d'une catégorie linéaire $C$ et sa catégorie des modules $C\Mod$ \cite[thm. 8.16]{lvdb1}. Leur raisonnement passe par la considération de la sous-catégorie $inj(C)$ des injectifs de $C\Mod$ qui, lorsqu'on déforme $C\Mod$ comme catégorie abélienne, se déforme, elle, linéairement.
L'équivalence $\ukCATmor\simeq \ukAB$ dit essentiellement la même chose mais avec la sous-catégorie $\widehat{C}$ des $C$-modules projectifs de type fini à la place de $inj(C)$ : les déformations linéaires (et pas seulement infinitésimales) de $\widehat{C}$ sont équivalentes à celles abéliennes de $C\Mod$.

\paragraph{Géométricité}

Afin d'obtenir des résultats de géométricité sur les champs précédents, il est nécessaire d'imposer des conditions de finitude sur les objets classifiés, celles choisies sont les suivantes : 
\begin{itemize}
\item pour $\ukCATiso$ on se limite au sous-champ $\ukcatiso$ classifiant les catégories dont les Hom sont projectifs de type fini et qui ont un nombre fini d'objets ;
\item pour $\ukCATeq$ on se limite au sous-champ $\ukcateq$ classifiant les catégories équivalentes à celles de $\ukcatiso$ ;
\item pour $\ukCATmor$ on se limite au sous-champ $\ukcatmor$ classifiant les catégories Morita-équivalentes à celles de $\ukcateq$ ;
\item et pour $\ukAB$ on se limite au sous-champ $\ukab$ image du morphisme $\MMod : \ukcateq\longrightarrow \ukAB$.
\end{itemize}

\bigskip

Le résultat principal de ce travail est le suivant.
\begin{irthm}[thm. {\ref{geomiso}}, {\ref{thmgeommod}}, {\ref{thmgeomeq}}]\label{thmintrogeom}
Les champs $\ukcatiso$, $\ukcateq$, $\ukcatmor$ et $\ukab$ sont géométriques.
\end{irthm}

\paragraph{Tangents}

Comme on s'y attend, nous démontrons que les complexes tangents de nos champs sont des tronqués du complexe de Hochschild. Si $C$ est une catégorie linéaire, notons $HC^0(C)\to HC^1(C)\to HC^2(C)\to \dots$ son complexe de Hochschild de $C$ (définit en Annexe {\ref{hoch}}) et $HZ^2(C)$ le sous-module de $HC^2(C)$ formé des 2-cocycles.

\begin{irthm}[thm. {\ref{tangentiso}}, {\ref{tangentequiv}}, {\ref{tangentabel}}]\label{thmintrotangent}
\begin{itemize}
\item Le complexe tangent en un point de $\ukcatiso$ correspondant à une catégorie linéaire $C$ est 
$$
Der^{\leq1}:= HC^1(C)\longrightarrow  HZ^2(C)
$$
où $HZ^2(C)$ est en degré 0.
\item Le complexe tangent en un point de $\ukcateq$ correspondant à une catégorie linéaire $C$ est 
$$
Hoch^{\leq2}:= HC^0(C)\longrightarrow HC^1(C)\longrightarrow  HZ^2(C)
$$
où $HZ^2(C)$ est en degré 0.
\item Le complexe tangent en un point de $\ukab$ correspondant à la catégorie $C\Mod$ des modules sur une catégorie linéaire $C$ est 
$$
Hoch^{\leq2}:= HC^0(C)\longrightarrow HC^1(C)\longrightarrow  HZ^2(C)
$$
où $HZ^2(C)$ est en degré 0.
\end{itemize}
\end{irthm}

En particulier nous démontrons le résultat suivant, élucidant la coïncidence des complexes tangents de $\ukcateq$ et $\ukab$. Ce résultat est à comparer avec \cite[thm. 8.16]{lvdb1} qu'il précise un peu puisque les auteurs n'y considère que les foncteurs de déformation en groupoïdes et non en 2-groupoïdes.

\begin{irthm}[thm. {\ref{thmetale}}]\label{thmintroetale}
$\MMod: \ukcateq\surj \ukab$ est étale.
\end{irthm}

Il est tout à fait remarquable que la cohomologie de Hochschild possède ainsi deux champs très différents l'admettant comme tangent.

\paragraph{Diagramme}
Essentiellement les résultats de la thèse se résument par les diagrammes suivant entre champs géométriques et leurs tangents (voir {\ref{diagrammefinal}} pour un diagramme plus complet) :
$$
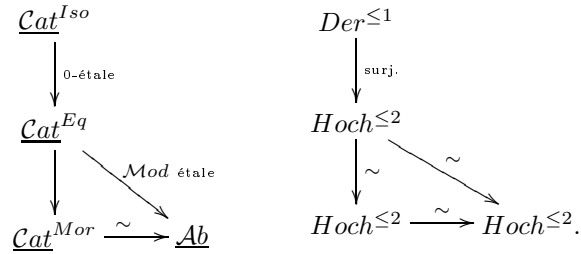
\xymatrix{
\ukcatiso \ar[d]^-{\textrm{{\tiny 0-étale}}}  \\
\ukcateq \ar[d] \ar[rd]^-{\MMod\ \textrm{{\tiny étale}}}  \\
\ukcatmor  \ar[r]^{\sim} & \ukab 
}\quad\quad\quad
\xymatrix{
Der^{\leq1} \ar[d]^{\textrm{{\tiny surj.}}} \\
Hoch^{\leq2}\ar[d]^-{\sim} \ar[rd]^-{\sim}  \\
Hoch^{\leq2}  \ar[r]^-{\sim} & Hoch^{\leq2}.
}$$

\section{Développements}

\paragraph{Conséquences}
En conséquence directe de la géométricité des champs $\ukass$, $\ukcat_*$ et $\ukab$, on peut espérer montrer que les champs $\uHHom(X,\ukass)$, $\uHHom(X,\ukcat_*)$ et $\uHHom(X,\ukab)$, classifiant les morphismes depuis un type d'homotopie $X\in\sens$ vérifiant certaines conditions de finitude ou depuis un certain schéma projectif $X$, sont également géométriques.
Sans mentionner les conditions de finitudes sur les objets, on peut dire que $\uHHom(X,\ukass)$ classifie les faisceaux d'algèbres associative sur $X$, $\uHHom(X,\ukcat_*)$ les gerbes linéaires sur $X$ et $\uHHom(X,\ukass)$ les formes tordues de catégories abéliennes sur $X$.

Les tangents de $\uHHom(X,\ukcat_*)$ et $\uHHom(X,\ukab)$ en un point devraient être donnés par la cohomologie de $X$ à coefficients dans une cohomologie de Hochschild, intégrant d'autres cadres de définition de la cohomologie de Hochschild de W.~Lowen et M.~Van den Bergh (cf. \cite[ch. 7]{lvdb2}) et donnant un cadre à la décomposition de Hodge de \cite[ch. ?]{gerst}. 

Les morphismes surjectifs $\ukass \longrightarrow \ukcat_* \longrightarrow \ukab$ induisent des morphismes
$$
\uHHom(X,\ukass) \longrightarrow \uHHom(X,\ukcat_*) \longrightarrow \uHHom(X,\ukab)
$$
dont seul le second est surjectif (il doit même rester étale car les tangents sont quasi-isomorphes), le défaut de surjectivité du premier traduisant le fait que toutes les gerbes sur $X$ ne sont pas neutres.

%
%
%

\paragraph{Champification}
Il se trouve que les objets classifiant les catégories linéaires (pour toutes les équivalences considérés) et abéliennes ne sont pas naturellement des champs et qu'il faut les champifier, opération correspondant au remplacement fibrant dans la catégorie de modèles des champs.

Pour chaque problème de modules abordé dans ce travail, il est posé la question de la description du champ associé (c'est à dire de trouver un autre problème de modules, élargissant le premier, dont le classifiant soit naturellement un champ, qui soit équivalent à celui qu'on étudie) et, dans les paragraphes qui en parlent, il est fait plusieurs fois référence au travail en cours sur la champification \cite{anel1}.

L'étude des problèmes de champification de cette thèse a conduit à trouver quelques énoncés généraux (\ie adaptables pour les champs en groupoïdes, en catégories, en 2-catégories, simpliciaux, etc.) de champification dont le point clé est le fait que les modules de champs forment eux-même des champs. 
Ce travail ayant été entamé trop récemment, il n'a pas pu être achevé pour figurer dans cette thèse et, n'étant encore totalement démontrés, nous avons choisi de ne pas donner les énoncés explicites de description des champification des problèmes de modules des catégories linéaires et abéliennes. Nous donnons seulement quelques indications rapides de ce que cela peut être (être plus locace prendrait trop de place dans ce travail).

\paragraph{Le reste du complexe}
L'un de nos résultats est une interprétation géométrique du début du complexe de Hochschild comme complexe tangent aux champs $\ukcateq$ et $\ukab$ et il est logique de se demander si le reste du complexe possède ou pas une interprétation géométrique.

Il est déjà prouvé dans \cite{penkava} que le $(n+2)$-ième module de cohomologie de Hochschild classifie les déformations au premier ordre de la structure associative d'une algèbre (discrète) $A$ en une structure de $A_\infty$-algèbre sur une déformation de l'ensemble $A$ en un type d'homotopie n'ayant de non-trivial que son $\pi_n$ (son $\pi_0$ restant, bien sûr, $A$).

Depuis quelques années, se mettent en place des techniques de {\sl géométrie algébrique dérivée} qui peuvent donner un sens à cette assertion \cite{lurie2, hag2}, et on peut espérer construire un {\sl champ dérivé} des catégories linéaires pour lequel le complexe tangent serait tout le complexe de Hochschild.
Un tel champ, définit, par exemple, sur le site de modèles des anneaux simpliciaux \cite{hag2}, devrait classifier certaines catégories simpliciales, conformément à l'analyse infinitésimale précédente.

\paragraph{dg-catégories}
Des notions de cohomologie de Hochschild sont définies dans le cadre très général des dg-catégories et il est normal d'en espérer une interprétation en terme du tangent à un champ des dg-catégories.

Un champ de modules de dg-catégories serait, par exemple, utile pour étudier l'association de sa catégorie dérivée à un schéma et donnerait un contexte pour étudier la conjecture de Kawamata sur la finitude de l'ensemble des schémas projectifs ayant des catégories dérivées équivalentes \cite{rouquier}.





\paragraph{En dehors de la géométrie algébrique}
Les résultats de cette thèse concernent la géométricité de certains champs définis au-dessus du site des schémas affines, et cela la rattache au domaine de la géométrie algébrique ; néanmoins, tant les définitions des champs que de la géométricité transcendent ce domaine et on peut souhaiter des résultats similaires à ceux de cette thèse dans le domaine des déformations d'algèbres normées, qui infinitésimalement sont données par la cohomologie du sous-complexe de Hochschild formé des cochaînes qui sont des opérateurs multidifférentiels.

\paragraph{Algèbres de Lie et autres opérades}
Comme les algèbres associatives, les algèbres de Lie possèdent une cohomologie de déformation (celle de Chevalley-Eilenberg) dont il est assez évident que le champ classifiant ces algèbres à isomorphisme près admet un double-tronqué comme complexe tangent. Peut-être un peu moins évident est qu'il doit exister un 2-champ $\underline{\mathcal{L}ie}$ qui serait l'équivalent pour les algèbres de lie du champ $\ukcat$ pour les algèbres associatives, et qui, en particulier récupérerait tout le début de la cohomologie de Chevalley-Eilenberg comme complexe tangent.

Le champ $\ukcat$ se définit en termes de catégories et de leurs équivalences et la notion correspondante semble être des {\sl algèbres de lie à plusieurs objets}, qu'il faut sans doute mieux penser comme les objets infinitésimaux correspondant à des {\sl groupoïdes de Lie}, ces objets doivent former naturellement une 2-catégorie, d'où le fait que le champ de modules soit un 2-champ.

La condition naturelle de finitude naturelle pour obtenir des champs géométriques serait de se limiter aux algèbres de lie dont le module sous-jacent est projectif de type fini.

\medskip

La notion correspondante à celle de catégorie abélienne pour les algèbres de lie serait les catégories de représentations des algèbres de Lie et, compte tenu du contexte, on peut conjecturer que le champ classifiant ces catégories est géométrique et que le morphisme depuis le champ $\underline{\mathcal{L}ie}$ est étale.

\bigskip

Plus généralement, on devine que la notion nécessaire pour copier la situation des algèbres associative et de Lie est une certaine opérade $\mathcal{O}$ dont les algèbres possèdent une notion de morphisme et de dérivation intérieures. Une telle structure permet de définir des 2-flèches et des équivalences (1-flèches inversibles à des 2-flèches près) dans la catégorie des $\mathcal{O}$-algèbres, autorisant la double classification à isomorphisme près et équivalence près. De telles algèbres doivent avoir leurs déformations infinitésimales à équivalence près contrôlées par une cohomologie généralisant celle de Hochschild, dont on peut s'attendre, comme dans le cas associatif, à ce qu'elle coïncide avec celle de Quillen sauf pour les premiers degrés.

\section{Plan de la thèse}

\noindent Le chapitre {\ref{chapchamp}} pose les définitions et énonce les principales propriétés des champs.
La section {\ref{secchpsplx}} étudie les champs simpliciaux ou supérieurs, définit leur géométricité et établit quelques proposition pour le calcule de leurs complexes tangents.
La section {\ref{secchpquillen}} étudie l'objet technique principal de ce travail : les préfaisceaux de Quillen.

\noindent Le chapitre {\ref{chapcatlin}} définit les catégories linéaires et abéliennes et construit les principaux schémas affines qui serviront au chapitre {\ref{chapgeom}} dans les preuves de géométricité. La fin du chapitre est dédiée à l'étude des structures de modèles sur les catégories de catégories.

\noindent Le chapitre {\ref{chapmodcat}} définit les champs associé aux différents problèmes de modules étudié dans ce travail.

\noindent Le chapitre {\ref{chapgeom}} contient les preuves de géométricité des champs précédents ainsi que des études de leur tangent. À la fin du chapitre se trouve un diagramme commutatif récapitulant les principaux objets et la nature des morphismes entre eux.

\noindent L'annexe {\ref{hoch}} définit la cohomologie de Hochschild d'une catégorie linéaire et interprète les premiers modules de cocycles et de cobord en les termes habituels.

\noindent L'annexe {\ref{cmf}} regroupe quelques résultats sur les nerfs d'équivalences de catégories de modèles et la localisation simpliciale de Dwyer-Kan.

\noindent Enfin, un index contient le vocabulaire principal utilisé tout au long de l'étude.

\newpage
\thispagestyle{empty}
\chapter*{Notations}
\thispagestyle{chheadings}
\addcontentsline{toc}{chapter}{Notations}

Dans tout ce travail on se fixe trois univers $\UU\in\VV\in\WW$, et on fait la convention que, si le contraire n'est pas précisé, tout les objets considérés seront toujours $\UU$-petits. Essentiellement, les différents univers n'apparaîtrons que lors de la considérations de catégories de catégories, et on oublie de signaler la dépendance en l'univers des autres catégories utilisées.

On utilise les notions suivantes pour les principales catégories utilisées :
\begin{itemize}
\item $\ens$\index{$\ens$} est la catégorie des ensembles $\UU$-petits ;
\item $\sens$\index{$\sens$} est la catégorie des ensembles simpliciaux $\UU$-petits ;
\item $\Delta$\index{$\Delta$} est la catégorie simpliciale ;
\item $\gp$\index{$\gp$} est la catégorie des groupes $\UU$-petits ;
\item $\gpd$\index{$\gpd$} est la catégorie des groupoïdes $\UU$-petits ;
\item $\catcat_\UU$\index{$\catcat$}\index{$\catcat_\UU$} est la catégorie des catégories $\UU$-petites.
\end{itemize}

Dans $\Delta$ les objets sont notés $[n]$, ou simplement $n$, et leurs foncteurs des points dans $\sens$ sont notés $\Delta^n$. Si $X$ est un ensemble simplicial, on note $X_n$ l'ensemble des $n$-simplexes de $X$.

\medskip

Dans tout ce qui suit $\Bbbk$ sera toujours un anneau commutatif fixé dans $\UU$, servant de référent absolu pour toutes les algèbres.
\begin{itemize}
\item $\kcom$\index{$\kcom$} désigne la catégorie des algèbres associatives et commutatives dans les $\kk$-modules $\UU$-petits ; 
\item $\kaff$\index{$\kaff$} désigne la catégorie des schémas affines sur $\kk$ (catégorie équivalente à $\kcom^o$) et, si $A\in\kcom$, $\aff_A$\index{$\aff_A$} désigne la catégorie des schémas sur $Spec(A)$ ;
\item pour $A\in\kcom$, $A\Mod$\index{$A\Mod$} est la catégorie des $A$-modules $\UU$-petits ;
\item $A\z\ASS$\index{$A\z\ASS$} désigne la catégorie des algèbres associatives unitaires dans $A\Mod$ ;
\item $A\z\ass$\index{$A\z\ass$} désigne la catégorie des algèbres associatives unitaire dans la sous-catégorie de $A\Mod$ des modules projectifs de type fini ;
\item pour $B\in A\z\ASS$, $B\Mod$ est la catégorie des $B$-modules $\UU$-petits ;
\end{itemize}

Par défaut, une algèbre commutative sera toujours supposée associative et unitaire.

\medskip

$pt$\index{$pt$} désigne toujours l'objet final de la catégorie à laquelle il appartient ; en particulier dans $\kaff$, $pt$ désigne $spec(\kk)$, dans $\ens$ il désigne un ensemble à un élément et dans $\sens$ l'ensemble simplicial constant à un élément.

Si $C$ est une catégorie $C^o$\index{$C^o$} désigne la catégorie opposée et $C^{int}$\index{$C^{int}$} le sous-groupoïde maximal formé des isomorphismes.

Dans une catégorie de modèles "holim"\index{holim} et "hocolim"\index{hocolim} désignent les limites et colimites homotopiques et $-\times^h_--$ le produit fibré homotopique.

Si $G$ est un groupe ou $A$ une algèbre associative unitaire, $\B G$\index{$\B G$} et $\B A$\index{$\B A$} désignent les catégories à un seul objet ayant $G$ ou $A$ comme endomorphismes.

\setcounter{tot}{0}
\newpage
\thispagestyle{empty}
\chapter{Champs}\label{chapchamp}
\thispagestyle{chheadings}

Il existe plusieurs notions de champs qui se distinguent par la structure de leurs points : classiquement, ceux-ci forment des groupoïdes \cite{lmb}, mais ils peuvent plus généralement être des catégories, et des groupoïdes ou des catégories supérieurs\footnote{Et même, au fond, tout objet d'une catégorie supérieure.}.

Ce chapitre est constitué de deux parties. La première définit une notion de champ en infini-groupoïdes, nommé champs simpliciaux car ils utilisent le formalisme des préfaisceaux simpliciaux ; puis, suivant \cite{simpson} et \cite{hag2} définit une notion de géométricité pour ces champs (cf. {\sf\S \ref{pqgeom}} pour une discussion sur les motivations d'une telle définition). La fin de cette première partie définit le complexe tangent d'un champs et établit quelques résultats pour les calculer.

La deuxième partie du chapitre définit et étudie les préfaisceaux de Quillen, qui sont l'outil technique fondamental de ce travail. On explique notamment comment un tel préfaisceau définit toujours un champ simplicial. L'essentiel de l'étude consiste à établir des propositions permettant de calculer les champs de morphismes ou les fibres des champs associés aux préfaisceaux de Quillen.

\medskip

On rappelle que, par défaut, l'usage du mot {\sl champ} a toujours le sens de {\sl champ en infini-groupoïdes} ; les champs au sens de \cite{lmb} seront qualifiés de {\sl 1-champs} ou de {\sl champs en groupoïdes}.

\section{Champs simpliciaux}\label{secchpsplx}

Comme les ensembles simpliciaux modélisent les infini-groupoïdes, les champs simpliciaux modélisent les champs en infini-groupoïdes. 
Les références principales pour les définitions et les propositions de la section sont les sections 3 ou 4 de \cite{hag1} et les sections 1.3, 1.4.1 et 2.1 de \cite{hag2}.
Les références pour les catégories de modèles sont \cite{hovey} et \cite{hirschhorn}, les notations sont celles de l'annexe {\ref{cmf}}.
On considère sur $\sens$ la structure de modèles classique décrite dans \cite{hovey}.

\subsection{La catégorie des champs}

\subsubsection{Préchamps}

\begin{defi}
La catégorie des {\em préchamps simpliciaux}\index{préchamp simplicial} sur $S$ est la catégorie $\kaffw$\index{$\kaffw$} des {\em préfaisceaux simpliciaux}\index{préfaisceau simplicial} sur la catégorie $\kaff$ des schémas affines. 
\end{defi}
Comme la structure de modèles de $\sens$ est engendrée par cofibrations \cite{hovey}, elle en induit une sur $\kaffw$ telle que les équivalences et les fibrations soient définies termes à termes \cite{gj, dhi}. En particulier, un préchamp simplicial est fibrant s'il est fibrant terme à terme.

\paragraph{Structure simpliciale} $\kaffw$ est une catégorie monoïdale pour le produit terme à terme dans $\sens$ et $\sens$ s'y plonge naturellement en les préfaisceaux constants. Ceci permet de définir un enrichissement de $\kaffw$ sur $\sens$ : pour $F$ et $G$ deux préfaisceaux simpliciaux, si $\Hom_{\kaffw}(F,G)$ désigne l'ensemble des morphismes de $F$ vers $G$ dans $\kaffw$, on définit l'espace simplicial des morphismes entre $F$ et $G$ par :
\begin{eqnarray*}
\Hom^{\Delta}_{\kaffw}(C,D)& = & n\mapsto \Hom_{\kaffw}(C\times\Delta^n,D) \quad\in SEns.
\end{eqnarray*}
Le foncteur dérivé de $\Hom^{\Delta}_{\kaffw}(-,-)$ est noté $\RHom^{\Delta}_{\kaffw}(-,-)$.

\paragraph{Yoneda} Via le plongement $\ens\to \sens$, la catégorie des préfaisceaux sur $\kaff$ se plonge dans celle des préfaisceaux simpliciaux, en particulier les objets de $\kaff$ peuvent être vu comme tels et on dispose d'un lemme de Yoneda \cite{hag1} : si $X\in\kaff$ et $F\in\kaffw$ on a une équivalence naturelle dans $\sens$ :
$$
F(X)\simeq \Hom^\Delta_{\kaffw}(X,F) \simeq \RHom^\Delta_{\kaffw}(X,F).
$$

\paragraph{Hom interne}
Soient $C,D\in\kaffw$ le foncteur $C\mapsto C\times D$ admet un adjoint à droite noté $D^C$ ou $\HHom_{\kaffw}(C,D)$ définit par :
\begin{eqnarray*}
D^C := \HHom_{\kaffw}(C,D) : \kaff & \longrightarrow & \sens \\
X & \longmapsto & \Hom^{\Delta}_{(\kaff_{|X})^\wedge}(C_{|X},D_{|X}).
\end{eqnarray*}
En particulier, si $K\in\sens$ est vu comme préfaisceau constant sur $\kaff$, on a 
$$
C^K(X) = \HHom_{\kaffw}(K,C)(X)=\Hom^{\Delta}_{\sens}(K,C(X)) = C(X)^K
$$
Le foncteur dérivé de $\HHom^{\Delta}_{\kaffw}(-,-)$ est noté $\RHHom^{\Delta}_{\kaffw}(-,-)$.

\paragraph{Produits fibrés homotopiques}
Soient $C_1\overset{a}{\to} C_3\overset{b}{\leftarrow} C_2\in \kaffw$ et $C_1\to C'_1\to C_3$ une factorisation de $a$ en une cofibration triviale puis une fibration. Comme la structure de modèles de $\kaffw$ est propre à droite, le produit fibré homotopique $C_1\times^h_{C_3}C_2$ peut se calculer comme le produit fibré $C'_1\times_{C_3}C_2$ \cite[prop. 13.3.7]{hirschhorn}.

\begin{lemme}\label{prodfib}
Soient $C_1\overset{a}{\to} C_3\overset{b}{\leftarrow} C_2\in \kaffw$. Le produit fibré homotopique $C_1\times^h_{C_3}C_2$ admet comme modèle
$$
(C_1\times C_2)\times_{C_3\times C_3}C_3^{\Delta^1}.
$$
\end{lemme}
\begin{pr}
Le produit fibré homotopique $C_1\times^h_{C_3}C_2$ est équivalent au produit fibré homotopique $(C_1\times C_2)\times^h_{C_3\times C_3}C_3$ où $d:C_3\to C_3\times C_3$ est le morphisme diagonal. et 
$C_3\to C_3^{\Delta^1}\to C_3\times C_3$, où le premier morphisme est issu de $\Delta^1\to \Delta^0$ et le second est tiré de $\Delta^0\times \Delta^0\to \Delta^1$, est une factorisation de $d$ en une cofibration triviale puis une fibration.
\end{pr}

\subsubsection{Faisceaux d'homotopie}\label{faischot}

À un préfaisceau simplicial $C$ on associe ses {\em préfaisceaux d'homotopies} :
\begin{eqnarray*}
\widehat{\pi}_0(C) : (\kaff)^o &\longrightarrow & \ens \\
X &\longmapsto & \pi_0(C(X))
\end{eqnarray*}
($\pi_0(C)$ est dit l'{\em espace grossier}\index{espace grossier} de $C$) et si $x\in C(X)$ est un point de $C$ :
\begin{eqnarray*}
\widehat{\pi}_n(C,x) : (\aff_{/X})^o &\longrightarrow & \gp \\
u:Y\to X &\longmapsto & \pi_n(C(Y),u^*x)
\end{eqnarray*}
où les foncteurs $\pi_n$ sont ceux définis sur $\sens$.

\begin{defi}
Ces préfaisceaux ne sont, en général, ni des faisceaux, ni même séparés (même si $C$ est un champ) et on définit les {\em faisceaux d'homotopie} de $C$, qu'on note $\pi_n(C)$ ($n\geq 0$), comme les faisceaux associés aux $\widehat{\pi}_n(C)$\footnote{Lorsqu'on parle des préfaisceaux ou des faisceaux d'homotopie en général on oublie, pour simplifier les notations, de noter la dépendance en un éventuel point base.}.
\end{defi}

\paragraph{Interprétation}
Si on pense à un ensemble simplicial comme à un infini-groupoïde, le préfaisceau $\widehat{\pi}_n(C,x)$ correspond exactement aux $n$-endomorphismes de l'objet $x$. L'interprétation des faisceaux d'homotopie est, elle, plus subtile : les points des $\pi_n(C)$ ne correspondent pas aux simples recollements de $n$-flèches, la faute à l'opération de séparation.
Considérons, par exemple, le faisceau $\pi_0(C)$ ; $\widehat{\pi}_0(C)$ classifie les objets de $C$ est sa séparation oblige à confondre certains objets, d'abord les isomorphes, mais également ceux qui sont seulement localement isomorphes : un objet et toutes ses formes tordues seront confondus dans $\pi_0(C)$. Il en est de même pour les flèches : la séparation de $\pi_n(C)$ oblige à confondre une flèche et ses formes tordues dans une classe qu'on qualifie de {\em forme locale} de la flèche ou de l'objet.
Après séparation, on sait qu'il suffit pour faisceautiser de saturer le préfaisceau pour les recollements de ses objets et les faisceaux $\pi_n(C)$ classifient donc les {\em recollements de formes locales}, lesquels ne correspondent pas forcément à des recollements d'objets (puisque cela revient à recoller sans la condition de cocycle).
C'est, notamment, le fait responsable de ce que le $\pi_0$ d'une gerbe non neutre n'est pas vide.

\paragraph{Suite exacte longue d'homotopie}
Soient $C\to D$ et $pt\to D$ dans $\sens$, il est associé à un produit fibré homotopique $F := pt\times^h_DC$ et $x:pt\to F$ une suite exacte longue d'homotopie
$$
\dots \longrightarrow \pi_1(F,x)\longrightarrow \pi_1(C,x)\longrightarrow \pi_1(D,x)
\longrightarrow \pi_0(F)\longrightarrow \pi_0(C)\longrightarrow \pi_0(D).
$$
En utilisant le fait que les produits homotopiques de $\kaffw$ se calculent terme à terme et l'exactitude du foncteur de faisceautisation, on associe à $C\to D\leftarrow pt\in \kaffw$ et $x:pt\to F := pt\times^h_DC$ la suite exacte 
$$
\dots \longrightarrow \pi_1(F,x)\longrightarrow \pi_1(C,x)\longrightarrow \pi_1(D,x)
\longrightarrow \pi_0(F)\longrightarrow \pi_0(C)\longrightarrow \pi_0(D).
$$

\subsubsection{\'Equivalences locales et hyper-recouvrements}

Les équivalences de $\kaffw$, sont les morphismes $f$ qui induisent des équivalences faibles de $\sens$ terme à terme, \ie tels que tous les $\widehat{\pi}_n(f)$ soient des isomorphismes de préfaisceaux.

\begin{defi}
On dit qu'un morphisme $f:C \to D\in \kaffw$ est une {\em équivalence locale}\index{equivalence locale@{équivalence locale}} si tous les $\pi_n(f):\pi_n(C)\to\pi_n(D)$ sont des isomorphismes de faisceaux, \ie des isomorphismes locaux de préfaisceaux. 
\end{defi}

\medskip

L'adjectif {\sl local} est utilisé pour rappeler la dépendance via-à-vis de la topologie ; par opposition, les équivalences de $\kaffw$ sont qualifiés d'{\em équivalences globales}\index{equivalence globale@{équivalence globale}}.

\begin{defi}
Un {\em hyper-recouvrement}\index{hyper-recouvrement} d'un objet $X\in\kaff$ est un morphisme $U\to X\in\kaffw$, qui soit une équivalence locale et où, pour tout $n$, le préfaisceau $U_n$ est une réunion quelconque de préfaisceaux représentables. 
\end{defi}

\begin{defi}
Un préchamp $C$ est dit {\em de descente}\index{préchamp de descente}\index{descente} si pour tout hyper-recouvrement $U\to X$, la flèche naturelle 
$$
C(X)\simeq \RHom^{\Delta}_{\kaffw}(X,C) \longrightarrow \RHom^{\Delta}_{\kaffw}(U,C)
$$
est une équivalence dans $\sens$.
\end{defi}

\subsubsection{Champs}\index{champ}

\begin{rthm}[{\cite[thm. 6.2]{dhi}}]
Il existe une localisation de Bousfield à gauche de $\kaffw$ pour laquelle les équivalences sont les équivalences locales. Cette localisation se définit en inversant la classe $H$ des hyper-recouvrements.

Les objets fibrants de cette localisation sont exactement les préchamps qui sont fibrant dans $\kaffw$ et de descente.
\end{rthm}

\begin{defi}
On note $\chaffk$\index{$\chaffk$} la catégorie $\kaffw$ munie de la structure de modèles issue de la localisation précédente et on s'y réfère comme à la {\em catégorie des champs simpliciaux}. Les objets fibrants $\chaffk$ sont appelés les {\em champs simpliciaux}\index{champ simplicial} sur $\kaff$, on conserve pour les autres le terme de préchamp simplicial.
\end{defi}

Le foncteur naturel $a(=id) :\kaffw \longrightarrow \chaffk$ est un adjoint de Quillen à gauche, commutant avec les limites homotopiques finies, dont l'adjoint à droite induit un foncteur pleinement fidèle $\mathbb{R}j:Ho(\chaffk)\to Ho(\kaffw)$. $a$ commute avec les colimites homotopiques finies et les produits fibrés homotopiques de champs peuvent donc se calculer dans la catégorie des préchamps.


\begin{defi}
Les préchamps localement équivalents à un champ $C$ seront qualifiés de {\em modèles locaux}\index{modèle local} pour $C$. Réciproquement, si $C$ est un préchamp, un remplacement fibrant $D$ sera qualifié de {\em champification} de $C$.

Si $C\to D$ est une équivalence locale, les points de $C$ seront qualifiés de {\em points locaux}\index{point local} de $D$.
\end{defi}

\begin{defi}
Un {\em $n$-champ}\index{n-champ@$n$-champ} (resp. un {\em $n$-préchamp}\index{n-préchamp@$n$-préchamp}) est un champ $C$ tel que ses faisceaux d'homotopies $\pi_k(C)$ (resp. ses préfaisceaux $\widehat{\pi}_k(C)$) soient triviaux pour $k>n$.
\end{defi}

Un champ peut posséder des modèles locaux qui ne sont pas des $n$-préchamps, mais tout champ équivalent à un $n$-préchamp est un $n$-champ.

\paragraph{Descente}\index{descente}

Pour $U\in\kaffw$, on a toujours $U\simeq \textrm{hocolim}_n\ U_n$ dans $\kaffw$.
Soit $U\to X$ un hyper-recouvrement et $C$ un préchamp, on a
$$
\RHom^{\Delta}_{\chaffk}(U,C) \simeq \textrm{holim}_n\ \RHom^{\Delta}_{\chaffk}(U_n,C)\simeq C(U_n).
$$
La condition de descente pour un préchamp $C$ se récrit sous la forme habituelle en demandant que, pour tout hyper-recouvrement $U\to X$, le morphisme canonique
\begin{eqnarray*}\label{cdtdescente}
C(X) \longrightarrow \textrm{holim}_n\ C(U_n)
\end{eqnarray*}
soit une équivalence dans $\sens$.

\subsubsection{Champ diagonal}\label{chpdiago}
\index{champ diagonal}\index{diagonale d'un champ}\index{champ des morphismes}


Comme le foncteur $a:\kaffw\longrightarrow \chaffk$ commute avec les limites homotopiques finies, on peut calculer les produits fibrés homotopiques de $\chaffk$ par la formule du lemme {\ref{prodfib}}.

\medskip
\begin{defi}
Soient $C\in\kaffw$ et $X\in\kaff$. À tout couple de points $x,y:\in C(X)$ on associe le {\em préchamp diagonal} de $C$ en $(x,y)$ ou {\em préchamp des morphismes} (des {\em endomorphismes} si les extrémités coïncident) de $x$ vers $y$ :
$$
\Omega_{x,y}C := X\times^h_{C}X = (X\times X)\times_{C\times C}C^{\Delta^1}.
$$
\end{defi}
Dans le cas où $x=y$ on note simplement $\Omega_xC$.

\smallskip

Le terme de "préchamp diagonal" vient de ce que $\Omega_{x,y}C$ est la fibre en $x\times y: X\times X\longrightarrow C\times C$ du morphisme diagonal $C\longrightarrow C\times C$.

\medskip

On déduit de la suite exacte longue d'homotopie associée à ce produit fibré que, pour tout $n\in\mathbb{N}$ :
$$
\pi_n(\Omega_xC,e_x)\simeq \pi_{n+1}(C,x)
$$
(où $e_x$ désigne l'identité de $x$).

\medskip

\begin{lemme}
Si $C$ est un champ, alors, pour tout $X\in\kaff$ et tout $x\in C(X)$, $\Omega_xC$ est un champ.
\end{lemme}
\begin{pr}
On utilise le fait que les champs soient les objets locaux par rapport aux équivalences locales. 
Soit $u:A\to B$ une équivalence locale entre deux préchamps, on veut savoir si 
$\RHom^{\Delta}_{Ch(\kaff)}(B,\Omega_xC)\to \RHom^{\Delta}_{Ch(\kaff)}(A,\Omega_xC)$ est une équivalence dans $\sens$.
Or $\RHom^{\Delta}_{Ch(\kaff)}(B,\Omega_xC)\simeq\Omega_x\RHom^{\Delta}_{Ch(\kaff)}(B,C)\simeq \Omega_xC(B)$ et $u$ induit un diagramme commutatif
$$\xymatrix{
pt\ar[r]\ar[d] &C(B)\ar[d]&\ar[l]\ar[d]pt\\
pt\ar[r] &C(A)&\ar[l]pt
}$$
où les flèches verticales sont des équivalences : c'est clair pour les deux extrêmes et cela résulte du fait que $C$ soit local pour celle du milieu. La proposition 13.3.4 de \cite{hirschhorn} assure alors que $\Omega_xC(B)\to \Omega_xC(A)$ est une équivalence.
\end{pr}

\begin{lemme}\label{commdiagochpif}
Si $C$ est un champ et $D$ un de ses modèles locaux, les espaces de chemins de $D$ sont des modèles locaux pour ceux de $C$.
\end{lemme}
\begin{pr}
Soit $D\to C$ une équivalence locale où $C$ est un champ, on veut montrer que, pour tout point $x$ de $D$, $\Omega_xD\to \Omega_xC$ est une équivalence locale.
C'est encore une application du lemme 13.3.4 de \cite{hirschhorn}.
\end{pr}


\subsubsection{Morphismes}

On définit quelques qualificatifs des morphismes de champs.

\begin{defi}
Un morphisme $f: C\to D\in\chaffk$ est dit
\begin{enumerate}
\item un {\em épimorphisme}\index{epimorphisme de champs@{épimorphisme de champs}} ou {\em essentiellement surjectif}\index{essentielle surjectivité} si $\pi_0(f):\pi_0(C)\to \pi_0(D)$ est un épimorphisme de faisceaux ;
\item une {\em gerbe}\index{gerbe} si $\pi_0(f):\pi_0(C)\to \pi_0(D)$ est un isomorphisme de faisceaux ;
\item {\em pleinement fidèle}\index{morphisme pleinement fidèle} si pour tout $n>0$ et tout $x\in C$, $\pi_n(f):\pi_n(C,x)\to \pi_n(D,f(x))$ est un isomorphisme de faisceaux ;
\item {\em ouvert}\index{morphisme ouvert} (resp. {\em fermé}\index{morphisme fermé}) si, pour tout $X\in\kaff$ et tout $X\to C\in\chaffk$, $X\times_CD\to X$ est un morphisme ouvert (resp. {\em fermé}) de schémas. Dans ce cas on dit que $C$ est un {\em sous-champ ouvert}\index{sous-champ ouvert} (resp. {\em fermé}\index{sous-champ fermé} de $D$.
\end{enumerate}
\end{defi}


\medskip
\begin{defi}
\begin{itemize}
\item Un morphisme de champs $f:C\to D$ est dit une {\em gerbe}\index{gerbe} ou {\em connexe}\index{morphisme connexe} si $\pi_0(f):\pi_0(C)\to\pi_0(D)$ est un isomorphisme.
\item Une gerbe $f:C\to D$ est dite {\em neutre}\index{gerbe neutre} si elle possède une section, \ie un morphisme $s:D\to C$ tel que $fs=id_D$.
\item Un champ $C$ est dit une {\em gerbe} ou un {\em connexe}\index{champ connexe} si le morphisme structural $C\to pt$ est connexe, \ie si $\pi_0(C)\simeq pt$. Il est {\em neutre} s'il possède un point global.
\end{itemize}
\end{defi}

\subsection{Champs géométriques}
\index{géométricité}
\index{champ géométrique}

\subsubsection{Pourquoi la géométricité ?}\label{pqgeom}

Ce paragraphe détaille les raisons de la définition des champs géométriques ; afin d'en faire ressortir l'essentiel,
l'analyse est faite dans un contexte plus général que celui de la géométrie algébrique.

\medskip

L'idée essentielle de l'invention des champs est de rendre effectifs et universels les quotients par des groupoïdes (cf. la caractérisation à la Giraud de la catégorie des champs \cite[\S 4.9]{hag1}), mais la complétion du site de base par les champs est trop forte en un certain sens.

Comme dans le cas des relations d'équivalences, pour lesquelles la catégorie universelle d'effectivité-universalité est celle des faisceaux, on va distinguer dans la catégorie des champs plusieurs sous-catégories correspondant à des notions plus ou moins proches de celle d'objets du site $S$. On rappelle que si $S$ est un site s'enrichissant en un contexte géométrique (cf. \cite{courstoen} et \cite{hag2} pour une version plus générale), on distingue dans la catégorie des faisceaux sur $S$, les {\em variétés} (ou {\em schémas} dans le cas algébrique) correspondant intuitivement au recollement d'objets de $S$ le long de la topologie ; et les {\em espaces géométriques} (ou {\em espaces algébriques} dans le cas algébrique), définis comme les faisceaux ayant {\sl localement} la structure des objets de $S$, \ie possédant un morphisme {\sl étale} (isomorphisme local) depuis une variété ; un tel morphisme est dit une {\em carte}. Les variétés sont caractérisables comme la clôture des objets de $S$ par les relations d'équivalences ouvertes et les espaces géométriques comme la clôture pour les relations d'équivalence étales (cf. \cite[prop. 1.3]{lmb} pour le cas algébrique). Ceci donne les inclusions pleines de catégories :
$$
S\subset \textrm{ variétés } \subset \textrm{ espaces géométriques }\subset \textrm{ faisceaux}
$$
Les espaces géométriques, forment une catégorie d'objets où se représentent convenablement, par exemple, les quotients de $\mathbb{R}/\mathbb{Q}$, ou du tore par une action de $\mathbb{R}$ donné par une pente irrationnelle, et c'est précisément la raison de leur définition\footnote{Ces exemples ne rentrent pas dans le cadre de la géométrie algébrique ; pour des exemples dans ce contexte, cf. \cite{knutson}.}.

Dans le cas des champs, la motivation et l'idée de la définition sont similaires ; pour représenter effectivement certains quotients de groupoïdes, on cherche à définir une classe de champs qui localement ressemblent à des objets de $S$ et il apparaît essentiellement deux notions différentes pour généraliser les espaces géométriques.
Tout d'abord, la notion d'isomorphisme local garde un sens dans la catégorie des champs et on peut vouloir définir un champ ayant localement la structure de $S$ comme un champ possédant une carte au sens précédent, \ie un morphisme étale depuis une variété, on obtient ainsi la notion de {\em champ de Deligne-Mumford}. 

Mais ces champs ne permettent pas de représenter tous les quotients de groupoïdes, à commencer par  $\underline{\B\mathbb{G}_a}$ et $\underline{\B\gln}$, l'obstruction étant que les cartes qui sont des isomorphisme locaux ne permettent de récupérer que des groupes discrets comme stabilisateurs des points :
si $C$ et un champ et $X\to C$ un morphisme surjectif, $X\times_CX\rightrightarrows X$ est naturellement muni d'une structure de groupoïde dans la catégorie des champs et son quotient est équivalent à $C$ ;
si $X\to C$ est supposé étale, les morphismes $X\times_CX\rightrightarrows X$ sont également étales et cela force les points du quotient à avoir des stabilisateurs qui sont des groupes étales.

Il faut autoriser des cartes ayant plus de "fibre" pour récupérer des stabilisateurs "continus". On définit alors une classe de champs ayant une carte dans $S$ tel que le morphisme soit une {\sl submersion} (morphisme {\sl lisses} dans le cas algébrique), ce sont les {\em champs géométriques} (ou {\em champs d'Artin} dans le cadre algébrique, cf. p.ex. \cite{simpson}). En résumé, sur le site $\aff$, on a les foncteurs pleinement fidèles suivants :
$$
\kaff \subset\textrm{schémas } \subset\textrm{espaces algébriques } \subset\textrm{ champs D-M} \subset\textrm{ champs géométriques } \subset\textrm{ champs}.
$$
Comme les faisceaux généraux, les champs non géométriques, sont rejetés comme des objets à la géométrie trop pathologique.

Il est notable que dans le cas des groupoïdes algébriques agissant librement, les deux notions de champs redonnent la même notions d'espace géométrique (cf. \cite[cor. 8.1.1]{lmb} pour le cas des 1-champs).

\subsubsection{Géométricité -- lissité}

La définition choisie pour la géométricité d'un champ est celle de Simpson mais on utilisera le vocabulaire de \cite{hag2} pour parler de $n$-géométricité, où $n$ est un entier traduisant la complexité à décrire l'objet à partir des objets de base que sont les schémas affines : un champ 0-géométrique est obtenu comme quotient d'un groupoïde affine, un champ 1-géométrique est obtenu comme quotient d'un groupoïde 0-géométrique, etc.

Les champs auront donc deux indices, le premier traduisant la nature de leurs groupoïdes de points et le second leur complexité géométrique. Un champ $n$-géométrique est toujours un $(n+1)$-champ.

\medskip

On se contente de définir les (-1)-, 0- et 1-géométricités qui suffiront pour les applications, la définition générale ainsi que les preuves des propriétés utilisées se trouvent dans \cite[\S 1.3.3]{hag2}.
La définition utilise la notion de morphisme lisse, qui est définie après\footnote{Les définitions des géométricités et lissités se font par récurrence et sont interdépendantes, d'où une présentation, toujours un peu maladroite, où l'on définit l'un avant l'autre tout en y faisant référence.}.

Les définitions sont données au cas par cas, n'utilisant pas la notion de $n$-atlas de \cite[\S 1.3.3]{hag2}.

\begin{defi}\index{champ géométrique}
\begin{enumerate}
\item Un champ est dit {\em (-1)-géométrique} si c'est un schéma affine.
\item Un champ est dit {\em 0-géométrique}
\begin{itemize}
\item si sa diagonale est affine, \ie si tous les $\Omega_{x,y}C$ sont des schémas affines ;
\item et s'il possède un morphisme surjectif et lisse depuis une petite réunion de schémas affines.
\end{itemize}
\item Un champ est dit {\em 1-géométrique}
\begin{itemize}
\item si sa diagonale est 0-géométrique, \ie si tous les $\Omega_{x,y}C$ sont 0-géométriques ;
\item et s'il possède un morphisme surjectif et lisse depuis une petite réunion de schémas affines.
\end{itemize}
\item Un champ est dit {\em géométrique} s'il est $n$-géométrique pour un certain $n$.
\end{enumerate}
\end{defi}

\bigskip

\begin{defi}
\'Etant donné un morphisme de champ $f:C\to D$ et un morphisme $x:X\to D$ depuis un schéma affine, la {\em fibre affine}\index{fibre affine} de $f$ en $x$ est le produit fibré homotopique $X\times^h_DC$.
\end{defi}

\begin{defi}
Si $n=-1,0,1$ un morphisme est dit {\em $n$-géométrique}\index{morphisme géométrique} (ou {\em $n$-représentable}\index{morphisme représentable}) si ses fibres affines sont $n$-géométriques. Un morphisme {\em représentable} ou {\em géométrique} est un morphisme $n$-géométrique ou $n$-géométrique pour un certain $n$.
\end{defi}

\begin{rlemme}[{\cite[prop. 1.3.3.4]{hag2}}]\label{basefibngeom}
Soit $C\to D$ un morphisme $n$-géométrique où $D$ est $n$-géométrique, alors $C$ est $n$-géométrique.
\end{rlemme}

\begin{rlemme}[{\cite[prop. 1.3.4.5]{hag2}}]\label{topfibngeom}
Soit $C\to D$ un morphisme $(n-1)$-géométrique où $C$ est $n$-géométrique, alors $D$ est $n$-géométrique.
\end{rlemme}

On déduit de leur définition les propriétés suivantes des morphismes ouverts et fermés.
\begin{cor}\label{ouvfergeom}
Soit $C$ un sous-champ ouvert ou fermé d'un champ géométrique $D$, alors $C$ est géométrique
\end{cor}

\begin{lemme}\label{fibnnungeom}
Soit $f:C\to D\in Ch(\aff)$ où $C$ est $n$-géométrique et $D$ est $(n+1)$-géométrique, alors les fibres affines de $f$ sont $n$-géométriques.
\end{lemme}
\begin{pr}
Soit $X\in\kaff$ et $X\to D$ on considère $X\times_CD\simeq (X\times D)\times_{C\times C}C$.
$C\to C\times C$ étant $n$-géométrique $X\times_CD\to X\times D$ est $n$-géométrique et comme $X\times D$ est $n$-géométrique, le lemme {\ref{basefibngeom}} assure que $X\times_CD$ est $n$-géométrique.
\end{pr}

\begin{defi}
Soit $C$ un champ, un morphisme représentable, lisse et surjectif $X\to C$ depuis une petite réunion disjointe de schémas (affines) est dit une {\em carte (affine)}\index{carte} de $C$.
\end{defi}

\begin{defi}
Un morphisme $f:C\to D$ est {\em quasi-compact}\index{morphisme quasi-compact} si toute ses fibres affines admettent des cartes qui sont des réunions finies de schémas affines.
\end{defi}

\medskip

On définit maintenant la lissité.

\begin{defi}\index{morphisme lisse}\label{deflisse}
Essentiellement, un morphisme entre champs est dit {\em lisse} si toutes ses fibres affines sont lisses. 
On se fixe un morphisme $f:C\to D$ de champs et on précise la définition dans plusieurs cas.
\begin{enumerate}
\item Si $C$ et $D$ sont (-1)-géométriques, on utilise la notion classique de lissité.
\item Si $C$ est (-1)-géométrique et $D$ est 0-géométrique, les fibres affines de $f$ sont affines et on utilise la notion précédente.
\item Si $C$ est 0-géométrique et $D$ est (-1)-géométrique, on dit que $f$ est lisse s'il existe $Y\to C$, une carte affine, telle que le composé $Y\to C\to X$ soit lisse.
\item Si $C$ et $D$ sont 0-géométriques, les fibres affines de $f$ sont 0-géométriques par le lemme {\ref{fibnnungeom}} et on utilise la notion précédente.
\item Si $C$ est 0-géométrique et $D$ est 1-géométrique, les fibres affines de $f$ sont 0-géométriques par le lemme {\ref{fibnnungeom}} et on se sert des notions précédentes.
\item Si $C$ est 1-géométrique et $D$ est (-1)- ou 0-géométrique, on dit que $f$ est lisse s'il existe $Y\to C$, une carte affine, telle que le composé $Y\to C\to X$ soit lisse.
\item Si $C$ et $D$ sont 1-géométriques, les fibres affines de $f$ sont 1-géométriques par le lemme {\ref{fibnnungeom}} et on utilise la définition précédente.
\end{enumerate}
\end{defi}

\bigskip

On définit maintenant la notion de {\em morphisme localement de présentation finie}.
\begin{defi}\index{morphisme localement de présentation finie}\label{deflocpresfini}
Essentiellement un morphisme est {\em localement de présentation finie} s'il est géométrique et si toute ses fibres affines le sont. Soit un morphisme $f:C\to D$ de champs, on précise la définition en deux temps :
\begin{enumerate}
\item si $D$ est affine $f$ est localement de présentation finie s'il existe une carte $C_0$ de $C$ telle que le morphisme composé $C_0\to D$ soit localement de présentation fini comme morphisme de schémas.
\item si $D$ n'est pas affine, on, demande que toutes ses fibres affines soient localement de présentation finie. 
\end{enumerate}
\end{defi}

\begin{defi}
Un morphisme $f:C\to D$ est dit de {\em présentation finie}\index{morphisme de présentation fini} s'il est localement de présentation finie et quasi-compact.
\end{defi}

\medskip

\begin{defi}\index{champ lisse}\index{champ localement de présentation finie}
Un champ $C$ est dit {\em lisse} ou {\em (localement) de présentation finie} si son morphisme structural $C\to pt$ l'est.
\end{defi}

\subsubsection{Présentation}\label{presentation}

\begin{defi} Un {\em groupoïde géométrique lisse}\index{groupoïde géométrique lisse} est un groupoïde $s,b:G_1\rightrightarrows G_0$ dans $\chaffk$ où $s$ ou $b$ est lisses et où $G_0$ et $G_1$ sont des champs géométriques.
\end{defi}
Le morphisme d'inversion des flèches d'un groupoïde étant une involution, il est toujours lisse, on en déduit que la lissité de $s$ entraîne celle de $b$ et réciproquement ; $s$ et $b$ étant lisses, on en déduit que le morphisme de multiplication est également lisse.

\medskip

Un tel groupoïde possède une colimite homotopique (\ie un quotient) dans $\chaffk$ noté $[G_0/G_1]$ pour lequel $G_0\to [G_0/G_1]$ est une carte.

\begin{defi} Soit $C$ un champ, une {\em présentation de $C$ par un groupoïde lisse}\index{présentation} est la donnée d'un groupoïde lisse $G_1\rightrightarrows G_0$ de $\chaffk$ et d'une équivalence $[G_0/G_1]\to C$.
\end{defi}

\begin{prop}
Un champ $(n+1)$-géométrique $C$ admet toujours une présentation par un groupoïde $n$-géométrique lisse.
\end{prop}
\begin{pr}
Par hypothèse, il existe une carte 0-géométrique $x:X\to C$ et $\Omega_xC$ est aussi $n$-géométrique.
$\Omega_xC\rightrightarrows X$ s'enrichit naturellement en une structure de groupoïde $n$-géométrique, et les morphismes source et but sont lisses car tirés en arrière de $X\to C$ qui est lisse.
\end{pr}

Si on itère la proposition à $\Omega_xC$ avec une carte $x_1:X_1\to \Omega_xC$ on obtient un $(n-1)$-groupoïde $\Omega_{x_1}\Omega_xC\rightrightarrows X_1$ et un 2-graphe 
$$
\Omega_{x_1}\Omega_xC\rightrightarrows X_1\rightrightarrows X
$$
dont on peut se demander s'il s'enrichit en un 2-groupoïde.
La réponse est en général négative, car la multiplication $\Omega_xC\times_X\Omega_xC\to \Omega_xC$
de $\Omega_xC\rightrightarrows X$ est un morphisme de champs, et ne se remonte que localement aux cartes $X_1\times_XX_1$ et $X_1$ de $\Omega_xC\times_X\Omega_xC$ et $\Omega_xC$.

Toutefois dans les exemples traités dans ce travail, un tel relevé existe et on pose les définitions suivantes.
(On renvoie à \cite[ch. 6]{maclane} pour les définitions de 2-groupoïde.)

\begin{defi}
Un 2-groupoïde $G_2\rightrightarrows G_1\rightrightarrows G_0$ est dit {\em affine}, (resp. {\em schématique, $n$-géométrique}) si les $G_i$ le sont.
Il est dit {\em lisse} si les quatre morphismes sources et buts sont lisses.
\end{defi}

\begin{defi}
Un champ $C$ admet une {\em présentation par un 2-groupoïde}\index{présentation} affine (resp. $n$-géométrique) lisse $G_2\rightrightarrows G_1\rightrightarrows G_0$ si 
$$
[G_1/G_1] \rightrightarrows G_0
$$
est un groupoïde lisse qui est une présentation de $C$ par un groupoïde.
\end{defi}

\subsection{Complexe tangent et lissité}\label{cplxtangent}

On définit le complexe tangent d'un champ qu'on compare avec sa version dérivée. Puis on compare la notion de lissité avec son pendant dérivé.

\bigskip

Soit $F$ un champ et $x:Spec(A)\to F$ un point affine.
Pour un $A$-module $M$ on note $A\oplus M$ l'extension infinitésimale au premier ordre de $A$ associée (cf. \cite[\S 1.2.1]{hag2}).
On définit $\mathbb{D}er_x(F,M)$ comme la fibre en $x$ du morphisme d'ensemble simpliciaux $Map(Spec(A\oplus M),F)\to Map(Spec(A),F)$.

\begin{defi}
Le {\em complexe cotangent}\index{complexe cotangent} de $F$ au point $x$, noté $\mathbb{L}_{F,x}$\index{$\mathbb{L}_{F,x}$}, est le complexe de $A$-modules en degrés positifs, représentant le foncteur 
\begin{eqnarray*}
\mathbb{D}er_x(F,-) : A\Mod &\longrightarrow & \sens\\
M & \longmapsto & \mathbb{D}er_x(F,M)
\end{eqnarray*}
au sens où le foncteur
\begin{eqnarray*}
h^{\mathbb{L}_{F,x}}:A\Mod &\longrightarrow & \sens\\
M & \longmapsto & DP\left(\RHom_A(\mathbb{L}_{F,x},M)^{\leq0}\right)
\end{eqnarray*}
où 
\begin{itemize}
\item $\RHom_A(-,-)$ est le hom interne dérivé de la catégorie des complexes de $A$-modules ;
\item $C^{\leq0}$ désigne la troncation du complexe $C$ ne conservant que les degrés négatifs\footnote{Pour un complexe $C = \dots C^{-1}\overset{d^0}{\to} C^0\overset{d^1}{\to} C^1\dots$, on définit $C^{\leq0} := \dots C^{-1}\overset{d^0}{\to} \ker(d^1)\to 0 \dots $\index{$C^{\leq0}$} et $C^{\geq0} := \dots 0\to \textrm{im}(d^0)\to C^1\to \dots $\index{$C^{\geq0}$}.} ;
\item et $DP$ désigne la transformation de Dold-Puppe associant un ensemble simplicial à un complexe en degré négatif \cite{illusie}
\end{itemize}
est équivalent à $\mathbb{D}er_x(F,M)$ dans la catégorie des préfaisceaux simpliciaux sur $A\Mod^o$.
\end{defi}

\bigskip

Cette notion de complexe cotangent, malgré la similitude de son nom, n'a rien à voir avec le complexe cotangent de Quillen-Illusie \cite{illusie}, qu'on qualifie de {\sl dérivé} afin de le distinguer du précédent. 
Une différence notable est que le complexe cotangent dérivé déborde a priori en degrés négatifs alors que le complexe cotangent d'un champ est complètement en degrés positifs.

On rappelle maintenant les définitions du complexe tangent dérivé ; puis, le lemme {\ref{comparcotan}}  compare les deux notions, rendant claire la remarque précédente.

\medskip

Soit $A\Mod_s$ la catégorie des $A$-modules simpliciaux, si $M\in A\Mod_s$, comme dans le cas non simplicial, on a une notion d'extension au premier ordre de $A$ par $M$ (cf. \cite[\S 1.2.1]{hag2}) et on définit $\mathbb{D}er_x(F,M)$ de même que précédemment.

\begin{defi}
Le {\em complexe cotangent dérivé}\index{complexe cotangent dérivé} de $F$ en $x$, est définit comme le complexe de $A$-modules (en degrés quelconques cette fois) $\mathbb{L}^{der}_{F,x}$ représentant le foncteur
\begin{eqnarray*}
\mathbb{D}er_x(F,-) : A\Mod_s &\longrightarrow & \sens\\
M & \longmapsto & \mathbb{D}er_x(F,M)
\end{eqnarray*}
au même sens que précédemment (cf. \cite[\S 1.4.1]{hag2}).
\end{defi}

On a le résultat suivant de comparaison entre les deux complexes tangents.
\begin{lemme}\label{comparcotan}
Soit $F$ un champ et $x:Spec(A)\to F$ un point affine, on a 
$$
\mathbb{L}_{F,x}\simeq (\mathbb{L}^{der}_{F,x})^{\geq0}.
$$
\end{lemme}
\begin{pr}
Il est clair que restreint à $A\Mod$, $(\mathbb{L}^{der}_{F,x})^{\geq 0}$ représente $\mathbb{D}er_x(F,-)$ dans les complexes en degrés quelconques, l'équivalence voulue vient du fait que, parce que $M$ est concentré en degré 0, on a $\RHom_A(\mathbb{L}^{der}_{F,x},M)\simeq \RHom_A((\mathbb{L}^{der}_{F,x})^{\geq0},M)$.
\end{pr}

\medskip

\begin{defi}
Le {\em complexe tangent dérivé}\index{complexe tangent dérivé} de $F$ au point $x$ est le complexe de $A$-modules $\mathbb{T}^{der}_{F,x} := \mathbb{R}\Hom_A(\mathbb{L}^{der}_{F,x},A)$.

Le {\em complexe tangent}\index{complexe tangent} de $F$ au point $x$ est le complexe de $A$-modules en degrés négatifs
$$
\mathbb{T}_{F,x}\index{$\mathbb{T}_{F,x}$}:=
\mathbb{R}\Hom_A(\mathbb{L}_{F,x},A)
\simeq \mathbb{R}\Hom_A((\mathbb{L}^{der}_{F,x})^{\geq0},A)
\simeq \mathbb{R}\Hom_A(\mathbb{L}^{der}_{F,x},A)^{\leq0}
=(\mathbb{T}^{der}_{F,x})^{\leq0}.
$$
\end{defi}
Par définition de $\mathbb{L}_{F,x}$ on a $DP(\mathbb{T}_{F,x})\simeq \mathbb{D}er_x(F,A)$.

\bigskip

La suite de cette section est dévoué à comparer les notions de lissité dans les cadres classiques et dérivés ; on cite abondamment \cite[ch. 2.2]{hag2}.

\medskip

On rappelle qu'on s'est fixé un anneau commutatif $\kk$ de référence et que $\kcom$ est la catégorie des $\kk$-algèbres commutatives.
\begin{defi}
On note $D^-\z\chaffk$ la catégorie de modèles des champs sur le site de modèles des $\kk$-algèbres commutatifs simpliciales (cf. \cite[ch. 2.2]{hag2}). Les objets de $D^-\z\chaffk$ sont qualifiés de {\em préchamps dérivés} et les fibrants de {\em champs dérivés}.
\end{defi}

On dispose d'une adjonction de Quillen comparant champs et champs dérivés (\cite[\S 2.2.4]{hag2}) :
$$
i_! : D^-\z\chaffk \leftrightarrows \chaffk : i^*
$$
dont il nous suffira de savoir que $\mathbb{L}i_!$ est pleinement fidèle (\cite[lem. 2.2.4.1]{hag2}), \ie que la catégorie des champs se plonge dans celle des champs dérivés.

\medskip

On dispose dans $D^-\z\chaffk$ d'une notion de lissité, définit essentiellement comme à la définition {\ref{deflisse}} (cf. \cite[\S 2.2.2]{hag2}), dont on va prouver que, appliqué aux objets non dérivés, elle redonne la définition {\ref{deflisse}}.

\begin{lemme}\label{lissvslissder}
Un morphisme de $\chaffk$ est lisse ssi il est lisse dans $D^-\z\chaffk$.
\end{lemme}
\begin{pr}
Compte tenu des définitions, il est clair qu'un morphisme lisse dans $\chaffk$ le sera dans $D^-\z\chaffk$. Il faut alors montrer qu'un morphisme entre champs non dérivé lisse dans $D^-\z\chaffk$, l'est en fait dans $\chaffk$. Soit donc $f:C\to D\in \chaffk$ lisse dans $D^-\z\chaffk$, on y va au cas par cas.
Si $C$ et $D$ sont des schémas, la notion de lissité de $D^-\z\chaffk$ est la classique les deux notions coïncident bien. Si $C$ est affine et $D$ un champ 0-géométrique, les fibres de $f$ sont des schémas affines et là encore les deux notions coïncident.
Si $C$ est un champ 0-géométrique et $D$ est affine, $f$ est lisse s'il existe un schéma représentable et $D_1$ et un morphisme lisse surjectif $D_1\to C$ dans $D^-\z\chaffk$ telle que la composition $D_1\to D$ soit lisse. Les raisonnement précédents montrent que $D_1\to D$ et $D_1\to C$ sont en fait lisses dans $\chaffk$, on en déduit, par définition d'être lisse dans $\chaffk$, que $f$ aussi.
Si $C$ et $D$ sont 0-géométriques, l'étude des fibres affines ramène le problème au cas précédent.

Le raisonnement se poursuit identiquement pour les cas où $C$ et $D$ sont 1-géométriques.
\end{pr}

\bigskip

On définit deux notions plus forte de lissité ; la première dit que le champs des relèvements d'une extension infinitésimale est connexe, \ie que le relèvement est unique à isomorphisme près ; et la seconde que le champs des relèvements est contractile, \ie que le relèvement est unique à une isomorphisme près, lui-même unique à un 2-isomorphisme près, etc.
\begin{defi}\label{defetale}
\begin{itemize}
\item Un morphisme $f:C\to D$ est dit {\em 0-étale}\index{morphisme 0-étale} s'il est lisse et si pour tout $x:Spec(A)\to C$, l'application tangente $T^0f: H^0(\mathbb{T}_{C,x})\longrightarrow H^0(\mathbb{T}_{D,f(x)})$ est un isomorphisme.
\item Un morphisme $f:C\to D$ est dit {\em étale}\index{morphisme étale} s'il est lisse et si pour tout $x:Spec(A)\to C$, l'application tangente $\mathbb{T}f: \mathbb{T}_{C,x})\longrightarrow \mathbb{T}_{D,f(x)})$ est un quasi-isomorphisme.
\end{itemize}
\end{defi}

\subsubsection{Calcul de complexes tangents}

Le but de cette section est d'établir le corollaire {\ref{tangpd}} permettant le calcul des complexes tangents de champs définit comme quotients de groupoïdes.

\bigskip

On rappelle que $\kcom$ désigne la catégorie des $\kk$-algèbres commutative pour un anneaux commutatif $\kk$ fixé.

Soit $A\in\kcom$ le {\em module des différentielles de Kähler}, noté $\Omega_A$, est le $A$-module obtenu par le quotient du $A$-module libre engendré par les symboles $da$ où $a\in A$ par les relations $d(aa')=a(da')+(da)a'$ et $d(ka+a')=kda+da'$ où $a,a'\in A$ et$k\in\kk$.

Comme tout module $\Omega_A$ définit un faisceau quasi-cohérent sur $Spec(A)$ dont la valeur en $A\to B\in \kcom$ est $\Omega_A\otimes_AB$.

\begin{prop}
Si $A\in\kcom$, le complexe cotangent de $Spec(A)$ en un point affine $x:Spec(B)\to Spec(A)$ s'identifie au $B$-module des différentielles de Kähler $\Omega_A\otimes_AB$, en particulier son dual s'identifie au $B$-module dérivations de $A$ (cf. \cite[ch. 16]{eisenbud}).
\end{prop}
\begin{pr}
Il suffit de remarquer que, dans le cas affine, $\mathbb{D}er_x(Spec(A),-)$ est à valeurs discrètes. On peut donc appliquer le résultat classique de représentabilité par le module des différentielles de Kähler (cf. p.ex. le début de \cite[ch. 16]{eisenbud}).
\end{pr}

\begin{nota}
Dans le cas d'un schéma affine $X$, on notera $\mathrm{T}_{X,x}$\index{$\mathrm{T}_{X,x}$} le complexe tangent à $X$ en $x$ afin de rappeler que c'est un simple module.
\end{nota}

\bigskip

La proposition suivante est le fait essentiel qui permet le calcul des complexes tangents.

\begin{prop}\label{calcultangent}
Soient $X=Spec(A)\in\kaff$ et un carré homotopiquement cartésien entre champs (munie d'une section de $b$) : 
$$\xymatrix{
F\ar[r]^a\ar[d]^b &G\ar[d]^c \\
X\ar@/^1pc/[u]^s \ar[r]_x &H
}$$
tel que $c$ soit lisse, alors on a la relation
$$
\mathbb{T}_{H,x} \simeq \textrm{hocolim}\ (\mathbb{T}_{F,s}\longrightarrow \mathbb{T}_{G,as}).
$$
\end{prop}
\begin{pr}
On utilise le (1) du lemme {\sf 1.4.1.16} de \cite{hag2}. Avec ses notations, ce lemme donne deux triangles de complexes
$$\xymatrix{
\mathbb{L}^{der}_{H,x}\ar[d]_{x^*}\ar[r]^{c^*} & \mathbb{L}^{der}_{G,bs}\ar[d]_{a^*}\ar[r] & \mathbb{L}^{der}_{G/H,bs}\ar[d]_{e}\\
\mathbb{L}^{der}_{X,id}\ar[r]^{b^*}& \mathbb{L}^{der}_{F,s}\ar[r]^f &\mathbb{L}^{der}_{F/X,s}
}$$
or, $\mathbb{L}^{der}_{X,id}=0$ donc $f$ est un quasi-isomorphisme.
Par cartésianité du carré, le (2) du lemme {\sf 1.4.1.16} assure que $e$ est un isomorphisme.
On a donc un triangle 
$$
\mathbb{L}^{der}_{H,x} \longrightarrow \mathbb{L}^{der}_{G,bs} \longrightarrow \mathbb{L}^{der}_{F,s}
$$
dont on va montrer que a suite des troncations
$$
\mathbb{L}_{H,x} \longrightarrow \mathbb{L}_{G,bs} \longrightarrow \mathbb{L}_{F,s}
$$
est encore un triangle. Il suffit pour cela de voir que la suite
$$
H^0(\mathbb{L}^{der}_{G,bs}) \longrightarrow H^0(\mathbb{L}^{der}_{F,s}) \longrightarrow 0
$$
est exacte ; or, parce que $b$ est lisse (comme tiré en arrière de $c$), on en déduit ({\cite[cor. 2.2.5.3]{hag2}}), que $H^1(\mathbb{L}^{der}_{F,s})=0$ et la suite est bien exacte.
On en déduit le triangle dual
$$
\mathbb{T}_{F,s} \longrightarrow \mathbb{T}_{G,bs} \longrightarrow \mathbb{T}_{H,x}.
$$
\end{pr}

\begin{cor}\label{tangpd}
Soient $s,b:G_1\rightrightarrows G_0$ est un groupoïde schématique lisse, et $x:X=Spec(A)\to G_0$, alors, $\mathbb{T}_{[G_0/G_1],x}$ est quasi-isomorphe au complexe de $A$-modules
$$
\textrm{T}_{G_1\times_{G_0}X,x}\overset{b^*}{\longrightarrow} \textrm{T}_{G_0,x}
$$
où, les $\textrm{T}_{S,s}$ désigne le module tangent du schéma $S$ au point $s$, et où $b^*$ est la différentielle de l'application but $b$ restreinte à $G_1\times_{G_0}X$.
\end{cor}
\begin{pr}
On applique la proposition {\ref{calcultangent}} au carré 
$$\xymatrix{
G_1\times_{G_0}X\ar[r]\ar[d] &G_0\ar[d]^c \\
X\ar@/^1pc/[u] \ar[r] & {[}G_1/G_0{]}
}$$
où la section est donnée par les flèches identités et où $c$ est lisse parce que le groupoïde est lisse.
On conclut en remarquant que, comme $G_0$, $G_1$, $X$ et donc $G_1\times_{g_0}X$ sont affines, leur complexe tangent s'identifie à leur module tangent.
\end{pr}

\subsubsection{Critère de lissité d'un morphisme}

Le but de cette section est d'établir la proposition {\ref{lissite}} donnant un critère de lissité des morphismes de champs par relèvement d'extensions infinitésimales.

\bigskip

En conséquence du lemme {\ref{lissvslissder}} on peut utiliser la caractérisation de la lissité dans $D^-\z\chaffk$ de \cite[\S 2.2.5]{hag2}. 

Pour un $A$-module simplicial $M$, $\Omega M$ désigne le produit fibré homotopique $0\times^h_M0$ dans $A\Mod_s$ (qui se calcule dans $\sens$) et $SM$, noté aussi $M[1]$, désigne la somme amalgamée homotopique $0\coprod^{\mathbb{L}}_M0$ dans $A\Mod_s$. On dit que $M$ est {\em connexe} si $\pi_0(M)=*$, en particulier on a alors $M \simeq S\Omega M$ dans $A\Mod_s$ (cf. \cite[\S 2.2.1]{hag2}).

On rappelle le foncteur $i^*:D^-\z\chaffk\longrightarrow \chaffk$.

\begin{rprop}[{\cite[prop. 2.2.5.1]{hag2}}]\label{criterelissder}
Soit $f:C\to D$ un morphisme géométrique entre champs dérivés. $f$ est lisse si et seulement si \begin{itemize}
\item $i^*f$ est localement de présentation finie dans $\chaffk$ ;
\item pour tout $A$-module simplicial connexe $M$ et toute dérivation $d:A\to A\oplus M$ pour lesquels on note $A\oplus_d\Omega M$ l'extension au premier ordre associée, tout diagramme
$$\xymatrix{
Spec(A)\ar[r]\ar[d] & C\ar[d]\\
Spec(A\oplus_d\Omega M)\ar[r]\ar@{-->}[ru] & D.
}$$
admet un relèvement, \ie le morphisme naturel $A\oplus_d\Omega M\to A$ induit une surjection
$$
\pi_0(C(A\oplus_d \Omega M)) \longrightarrow \pi_0\left(C(A)\times^h_{D(A)}D(A\oplus_d \Omega M)\right).
$$
\end{itemize}
\end{rprop}

En conséquence du lemme {\ref{lissvslissder}} on peut utiliser la caractérisation ci-dessus pour les morphismes lisses de $\chaffk$ ; le lemme suivant permet de voir comment la situation se simplifie lorsque le morphisme $f$ est dans $\chaffk$.

\begin{lemme}\label{pizeroext}
Soit $A\in\kcom$ et soient $M$ un $A$-module simplicial connexe ($\pi_0(M)=*$), $d:A\to A\oplus M$ une dérivation et $A\oplus_d\Omega M$ l'extension au premier ordre associée. 
Alors, on a $\pi_0(A\oplus_d\Omega M) \simeq A\oplus_{d_0}\pi_0(\Omega M)$ (où $d_0$ est définit dans la preuve).
\end{lemme}
\begin{pr}
Soit $M_{\leq1}$ le premier étage de la tour de Postnikov de $M$ (\cite[\S 2.2.1]{hag2}), on a $\Omega M_{\leq1}\simeq \pi_1(M)\simeq \pi_0(\Omega M)$.
Le morphisme naturel $M\to M_{\leq1}$ induit un morphisme d'extensions au premier ordre $A\oplus M\to A\oplus M_{\leq1}$ et donc, en composant par $d$, une dérivation $d_0:A\to A\oplus M\to A\oplus M_{\leq1}$. On en déduit un morphisme $\tau:A\oplus_d\Omega M\to A\oplus_{d_0}\Omega M_{\leq1}$.
L'algèbre $A\oplus_{d_0}\Omega M_{\leq1}$ est une extension de $A$, qui est discret, par $\Omega M_{\leq1}\simeq \pi_1(M)$, discret également, elle est donc discrète.
Par définition du $\pi_0$, on déduit alors de $\tau$ un morphisme d'algèbres $\tau_0:\pi_0(A\oplus_d\Omega M)\to A\oplus_{d_0}\Omega M_{\leq1}$.
Comme $A$ est discret, on tire de la longue suite d'homotopie de la fibration $\Omega M\to A\oplus_d\Omega M\to A$ que $\pi_0(A\oplus_d\Omega M)$ est une extension de $\pi_0(A)$ par $\pi_0(\Omega M)$. $\tau_0$ induit donc un morphisme d'extensions, qui est l'identité sur $\pi_0(A)$ et $\pi_0(\Omega M)$ :
c'est un isomorphisme.
\end{pr}

\begin{cor}\label{ptextchp}
Si $C\in\chaffk$ alors, avec les notations du lemme précédent, $C(A\oplus_d\Omega M) = C(\pi_0(A\oplus_d\Omega M)) = C(A\oplus_d\pi_0(\Omega M))$
\end{cor}
\begin{pr}
$C(A\oplus_d\Omega M) = C(\pi_0(A\oplus_d\Omega M))$ est dû à la définition du foncteur $i^*:\chaffk\longrightarrow D^\z\chaffk$, le reste est une application du lemme {\ref{pizeroext}}.
\end{pr}

On reformule alors la proposition {\ref{criterelissder}} de la manière suivante.
\begin{cor}\label{lissiterelvt}
Un morphisme $f:C\longrightarrow D$ de $\chaffk$ entre deux champs est lisse si et seulement s'il est localement de présentation finie et si, pour tout $I\in A\Mod$ et toute dérivation $d:A\to A\oplus I[1]$ dont note $A\oplus_d I \ (\in\kcom)$ l'extension au premier ordre associée, tout diagramme
$$\xymatrix{
Spec(A)\ar[r]\ar[d] & C\ar[d]\\
Spec(A\oplus_dI)\ar[r]\ar@{-->}[ru] & D.
}$$
admet un relèvement.
\end{cor}
\begin{pr}
D'après le corollaire {\ref{ptextchp}} il suffit de tester le relèvement sur tout extension du type $A\oplus_d\Omega M$ où $\Omega M$ est discret, comme $M$ est supposé connexe un tel module s'écrit toujours sous la forme $SI = I[1]$ où $I\in A\Mod$.
\end{pr}

Afin d'avoir le critère voulu de relèvement, il reste à montrer le lemme suivant.
\begin{lemme}\label{extclassico}
Toute extension infinitésimale au premier ordre $a:A'\to A\in \kcom$ est de la forme $A\oplus_d I$, où $I=\ker(A'\to A)$ et $d:A\to A\oplus I[1]$ est une dérivation.
\end{lemme}
\begin{pr}
On considère un diagramme
$$\xymatrix{
I\ar[r]\ar[d] \ar@{}[rd]|1& A'\ar[r]^a\ar[d]^a \ar@{}[rd]|2& A\ar[d]^s\\
0\ar[r] & A\ar[r]^d & A''
}$$
où le carré 2 est le pushout homotopique de $a,a:A'\rightrightarrows A$ dans la catégorie des algèbres commutatives simpliciales et où le carré 1 décrit $A$ comme quotient de $A'$ par $I$, c'est aussi un pushout dans cette catégorie car $I\to A'$ est un monomorphisme. On en déduit que le carré 
$$\xymatrix{
I\ar[r]^0\ar[d] & A\ar[d]^s\\
0\ar[r] & A''
}$$
est également un pushout, $A''$ est donc isomorphe à $A\oplus SI=A\oplus I[1]$ et $s:A\to A\oplus I[1]$ correspond par cet isomorphisme à la dérivation nulle.
\end{pr}

Le lemme {\ref{extclassico}} et le corollaire {\ref{lissiterelvt}} permettent enfin d'énoncer le critère voulu.
\begin{prop}\label{lissite}
Un morphisme $f:C\longrightarrow D$ de $\chaffk$ entre deux champs est lisse si et seulement s'il est localement de présentation fini et si, pour toute l'extension au premier ordre $A'\to A\in\kcom$, tout diagramme
$$\xymatrix{
Spec(A)\ar[r]\ar[d] & C\ar[d]\\
Spec(A')\ar[r]\ar@{-->}[ru] & D.
}$$
admet un relèvement, \ie si le morphisme suivant est surjectif
$$
\pi_0(C(A')) \longrightarrow \pi_0\left(C(A)\times^h_{D(A)}D(A')\right).
$$
\end{prop}

\subsection{Exemples de champs}

\subsubsection{Schémas}

L'inclusion $\ens\to\sens$ permet de voir tout préfaisceau en ensembles comme un préfaisceau simplicial, en particulier un tel préfaisceau, dit {\em discret}\index{préfaisceau simplicial discret}, est un champ si et seulement c'est un faisceau. Ce provient de ce que dans la condition de descente {\ref{cdtdescente}} 
$\textrm{holim}_n\ F(U_n)$ se réduit à $\textrm{lim}(F(U_0)\rightrightarrows F(U_1)$ lorsque les valeurs de $F$ sont des ensembles.

En particulier tout schéma définit un champ et, avec la terminologie choisie, les schémas à diagonale affine (\ie séparés) sont 0-géométriques, et les schémas à diagonale schématique séparée sont 1-géométriques.

\subsubsection{Champ associé à un groupoïde}\label{groupoide}

Soit $G := s,b:G_1\rightrightarrows G_0$ un groupoïde schématique tel que les deux flèches $s$ et $b$ soient des morphismes lisses de schémas, alors $G$ définit un 1-champ géométrique $[G_0/G_1]$\index{$[G_0/G_1]$}, tel que le morphisme $G_0\to [G_0/G_1]$ soit une carte. En particulier, $[G_0/G_1]$ est localement de présentation finie si $G_0$ et $G_1$ le sont.

$[G_0/G_1]$ est 0-géométrique si $s$ ou $b$ sont des morphismes affines, 1-géométrique si leurs fibres sont  des schémas séparés.

\bigskip

Si $G_0=pt$ et $G_1=G\ell_n$, on note $\underline{\B G\ell_n}$\index{$\underline{\B G\ell_n}$} le champ $[pt/G\ell_n]$. Il classifie les $G\ell_n$-torseurs et, comme $G\ell_n$ a une représentation fidèle et transitive dans le groupe des automorphismes d'un module libre de rang $n$, les $G\ell_n$-torseurs sont en équivalence avec le groupoïde des fibrés vectoriels de rang $n$. On note $\uvect^n$\index{$\uvect^n$} le champ classifiant ces derniers, on a donc $\underline{\B G\ell_n}\simeq \uvect^n$ et si $\uvect$ désigne le champ classifiant les fibrés vectoriels de tout rang on a $\uvect\simeq \coprod_n\underline{\B G\ell_n}$.


\section{Champs de modèles}\label{secchpquillen}

\subsection{Champs de modules -- Champs de modèles}\label{chpmodchpmod}

On formalise l'analyse commencée en introduction sur l'association d'un champ à un problème de modules.

\medskip

L'axiomatisation primitive d'un problème de modules donne un préfaisceau faible\footnote{Comme tout est strictifiable en de vrais préfaisceaux (cf. {\sf\S\ref{strictification}}) on ne parlera plus que de préfaisceaux stricts.} en catégories $M$ (ou une catégorie fibrée clivée) sur un site $S$ (correspondant au type de famille choisie : algébrique, analytique, différentiable, etc.), chacune des catégories-points $M(x)$ étant munie d'une sous-catégorie d'équivalences $W(x)$. 

Terme à terme, on peut localiser $M(x)$ par $W(x)$ et obtenir une catégorie simpliciale (cf. \cite{dk}), toutefois ces catégories simpliciales ne forment un préfaisceau que si les équivalences sont conservés par les changements de base, \ie si les $W(x)$ forment un sous-préfaisceau en catégories. Dans l'affirmative, on note $L(M,W)$ le préfaisceau en catégories simpliciales obtenu. Les catégories simpliciales sont des modèles pour les catégories supérieures (cf. {\sf\S\ref{catsplx}}) dont les $n$-flèches pour $n\geq2$ sont toutes inversibles ; or, comme on se limite à travailler avec des champs en infini-groupoïdes, on ne va retenir des catégories simpliciales $L(M(x),W(x))$ que leurs sous-infini-groupoïdes maximaux $L(M(x),W(x))^{int}$ (cf. {\sf\S\ref{catsplx}}). Les restrictions conservant l'inversibilité des flèches, les $L(M(x),W(x))^{int}$ forment un sous-préfaisceau $L(M,W)^{int}$ de $L(M,W)$, on obtient un préfaisceau simplicial, noté simplement $|M|$, en composant par le foncteur réalisation géométrique du nerf (cf. {\sf\S\ref{catsplx}}).

\medskip

Cette construction étant donnée, plusieurs problèmes se posent a priori, à commencer par le fait que des exemples de problèmes de modules (complexes, dg-catégories) montrent que les équivalences peuvent ne pas être préservées par les restrictions. Ceci amène à devoir dériver les foncteurs de restriction. Ensuite, il faut pouvoir calculer explicitement les morphismes dans les catégories localisés, notamment pour établir la géométricité des champs associés. Ces deux problèmes amènent à la considération préfaisceaux en catégories de modèles ou préfaisceaux de Quillen, qui est un bon cadre tant pour le calcul des foncteurs dérivés que pour celui des espaces de morphismes. On redéveloppe ici leur théorie en suivant \cite{hsimpson,hag2}.

Dernier problème enfin, les préfaisceaux simpliciaux obtenus par ce procédé ne sont en général pas des champs (pour des raisons souvent multiples). Ceci est illustré par les exemples de cette thèse (cf. {\sf\S \ref{champifiso}, \S \ref{champifeq}, \S \ref{champifabel}}). Si tel est le cas, l'alternative est la suivante : soit on renonce à ce que l'objet classifiant soit un champ, \ie on renonce à l'opération topologique de recollement des points et morphismes (descente) ; soit le problème de modules était mal posé, puisqu'il n'est pas stable par l'opération de recollement des points et morphismes. Comme il ne parait pas intéressant de perdre le recollement, on prendra le second point de vue. Par exemple, on considérera que le problème des modules des catégories linéaires est mal posé (cf. {\sf  \S\ref{champifiso}}) ; le bon problème des modules associé étant quelque chose comme celui des champs en catégories linéaires (qu'on traitera dans \cite{anel1}).

\subsection{Préfaisceaux de Quillen}

Les catégories de familles et les équivalences considérées dans les problèmes de modules peuvent souvent s'enrichir en catégories de modèles ; toutefois, les objets dont on espère que leur modules forment des champs géométriques doivent posséder des conditions de finitudes, et cela n'en fait presque jamais des catégories avec les limites ou colimites, et, a fortiori, des catégories de modèles. Heureusement ces catégories sont souvent des sous-catégories pleines de catégories de modèles et cela justifie la notion de sous-catégorie de modèles (définition {\ref{sscatmod}}) et sous-préfaisceau de Quillen (définition {\ref{ssprefq}}).

Ainsi, les champs de modules construit à partir d'un sous-préfaisceau de Quillen seront toujours des sous-champs pleins du champ associé au préfaisceau de Quillen.


\subsubsection{Définitions}\label{prefQ}

On renvoie à l'annexe {\ref{cmf}} pour les définitions concernant les catégories de modèles. On définit les catégories suivantes de catégories de modèles.
\begin{itemize}
\item $\cmg$\index{$\cmg$} est la catégorie des catégories de modèles avec comme morphismes les adjoints de Quillen à gauche.
\item $\cmsg$\index{$\cmsg$} est la catégorie des catégories de modèles simpliciales avec comme morphismes les adjoints de Quillen à gauche.
\item $\scmg$\index{$\scmg$} est la catégorie des sous-catégories de modèles avec comme morphismes les adjoints de Quillen à gauche.
\item $\scmsg$\index{$\scmsg$} est la catégorie des sous-catégories de modèles simpliciales avec comme morphismes les adjoints de Quillen à gauche.
\item $\scmsg_{cf}$\index{$\scmsg_{cf}$} est la catégorie des sous-catégories de modèles simpliciales avec comme morphismes les adjoints de Quillen à gauche préservant les objets fibrants et cofibrants.
\end{itemize}
Les mêmes catégories avec un exposant $d$ à la place du $g$ désignent les catégories formées des mêmes objets avec les adjoints de Quillen à droite comme morphismes.

\begin{defi}
Un {\em préfaisceau de Quillen à gauche}\index{préfaisceau de Quillen} (resp. {\em à droite}) $M$ sur une catégorie $S$ est un foncteur 
$$
M : S^o\longrightarrow \cmg \textrm{ (resp. $\cmd$)}.
$$
Un {\em préfaisceau de Quillen simplicial à gauche}\index{préfaisceau de Quillen simplicial} (resp. {\em à droite}) $M$ sur une catégorie $S$ est un foncteur 
$$
M : S^o\longrightarrow \cmsg \textrm{ (resp. $\cmsd$)}.
$$
Les morphismes de préfaisceaux de Quillen sont définis comme les transformations naturelles de foncteurs.
On note $Pr(S,\cmg)$ (resp. $Pr(S,\cmd)$) la catégorie des préfaisceaux de Quillen à gauche (resp. à droite) et 
$Pr(S,\cmsg)$ et $Pr(S,\cmsd)$ leurs équivalents simpliciaux.
\end{defi}

On a de même des notions de {\em préfaisceau faible de Quillen}\index{préfaisceau faible de Quillen} (simplicial ou pas) en relâchant les conditions de fonctorialité ; le strictifié d'un préfaisceau faible de Quillen est un préfaisceau de Quillen, de même dans le cas simplicial (cf. {\sffamily \S\ref{strictification}}).

\begin{defi}\label{ssprefq}
Un {\em sous-préfaisceau de Quillen à gauche}\index{sous-préfaisceau de Quillen} (resp. {\em à droite}) $C$ 
est un foncteur 
$$
C : S^o\longrightarrow \scmg \textrm{ (resp. $\scmd$)}.
$$
Un {\em sous-préfaisceau de Quillen simplicial à gauche}\index{sous-préfaisceau de Quillen simplicial} (resp. {\em à droite}) $C$ 
est un foncteur 
$$
C : S^o\longrightarrow \scmsg \textrm{ (resp. $\scmsd$)}.
$$
En particulier, un sous-préfaisceau de Quillen est toujours un sous-préfaisceau d'un préfaisceau de Quillen de même orientation.

Les morphismes de sous-préfaisceaux de Quillen sont définis comme les transformations naturelles de foncteurs.
On note $Pr(S,\scmg)$ (resp. $Pr(S,\scmd)$) la catégorie des préfaisceaux de Quillen à gauche (resp. à droite) et 
$Pr(S,\scmsg)$ et $Pr(S,\scmsd)$ leurs équivalents simpliciaux.
\end{defi}

\medskip
Soit $C$ un sous-préfaisceau de Quillen.
Pour $s\in S$ on note $C(s)$ la valeur de $C$ en $s$. On définit $W_C(s)$ comme la sous-catégorie des équivalences de $C(s)$ et $C(s)^c$, $C(s)^f$ et $C(s)^{cf}$ respectivement comme les sous-catégories des objets cofibrants, fibrants et fibrants-cofibrants de $C(s)$. Pour $*=c$, $f$ ou $cf$, $W_C^*(s)$ désigne $W_C(s)\cap C(s)^*$.

\subsubsection{Préfaisceaux de morphismes}

Soient $M:S^o\to \scmg$ un sous-préfaisceau de Quillen simplicial gauche et $x,y\in M(s)$ pour $s\in S$. 

D'après le {\sf\S\ref{fonctomorquillen}}, on peut associer à cette donnée un préfaisceau faible à valeurs dans $Ho(\sens)$, dit {\em préchamp des morphismes dans $M$ de $x$ vers $y$}\index{préchamp des morphismes} essentiellement donné par :
\begin{eqnarray*}
\mathcal{M}ap_M(x,y) : (S/s)^o &\longrightarrow & Ho(\sens)\\
u:t\to s &\longmapsto & Map^{eq}_{M_t}(\mathbb{L}u^*x,\mathbb{L}u^*y)\\
t'\overset{v}{\to}t \overset{u}{\to} s &\longmapsto & Map(v)
\end{eqnarray*}
où $Map(v)$ est la flèche dans $Ho(\sens)$ :
$$
Map^{eq}_{M_t}(\mathbb{L}u^*x,\mathbb{L}u^*y)\to Map^{eq}_{M_t}(\mathbb{L}v^*\mathbb{L}u^*x,\mathbb{L}v^*\mathbb{L}u^*y)\overset{\sim}{\to} Map^{eq}_{M_t}(\mathbb{L}(uv)^*x,\mathbb{L}(uv)^*y)
$$

Dans le cas où $M:S^o\to \scmg$ est en fait à valeurs dans $\scmsg_{cf}$, on a un modèle simplicial pour $\mathcal{M}ap_M(x,y)$.
\begin{eqnarray*}
\mathcal{E}q^\Delta_M(x,y) : (S/s)^o &\longrightarrow & \sens\\
u:t\to s &\longmapsto & \Eq^\Delta_{M_t}(u^*(x^{cf}),u^*(y^{cf})).
\end{eqnarray*}
où, pour $x\in M$, $x^{cf}$ est un remplacement fibrant cofibrant de $x$.

\subsubsection{Classifiants}\label{classifquillen}

Soit $C$ un sous-préfaisceau de Quillen, on pose plusieurs notations pour des objets associés qui reprennent celles du {\sf\S\ref{quillen}}.

\medskip

Si $C$ est de Quillen à gauche (resp. à droite), la composition par le foncteur $|W^c_-|:\cmg\to \sens$ (resp. $|W^f_-|:\cmg\to \sens$) permet d'associer à $C$ un préfaisceau simplicial noté $|C|$.
Dans le cas gauche, on a :
\begin{eqnarray*}
|C| : S^o & \longrightarrow & \sens\\
	s & \longmapsto & |W_C^c(s)|\\
	u:t\to s & \longmapsto & |u^*|:|W_C^c(s)|\to |W_C^c(t)|.
\end{eqnarray*}

\begin{defi}
Par la proposition \ref{classifmodele} chaque $|W_C^c(s)|$ est un classifiant pour $(C(s),W_C(s))$ et $|C|$ sera dit le {\em préfaisceau (ou préchamp) classifiant}\index{préfaisceau classifiant}\index{préchamp classifiant} canonique du sous-préfaisceau de Quillen $C$. On dira aussi que $|C|$ est un {\em sous-préchamp de Quillen}\index{sous-préchamp de Quillen}. Dans le cas où $C$ est en fait un préfaisceau de Quillen on parlera de $|C|$ comme d'un {\em préchamp de Quillen}\index{préchamp de Quillen}.
Si $S$ est un site, le champ associé à $|C|$ est noté $\underline{C}$ et est appelé le {\em champ associé au problème de module} $C$.
\end{defi}


\bigskip

Dans le cas d'un sous-préfaisceau de Quillen à valeurs dans $\scms_{cf}^*$ (ou *=$g$ ou $d$), on dispose d'un autre préfaisceau classifiant pour $C$, aidant au calcul des préchamps de morphismes. Si $C$ est un sous-préfaisceau de Quillen, en utilisant le foncteur $\mathcal{G}:\scms^*_{cf}\to \scatcat$ (cf. {\ref{G}}) et $|-|:\scatcat\to \sens$ (cf. {\ref{realgeom}}), on définit le préfaisceau simplicial $|\G(C)|$ par :
\begin{eqnarray*}
|\G(C)| : S^o & \longrightarrow & \sens\\
	s & \longmapsto & |\G(C(s))|.
\end{eqnarray*}

\begin{prop}\label{CGC}
Si $C\in Pr(S,\scms^*_{cf})$, où *=$g$ ou $d$, alors $|C|$ et $|\G(C)|$ sont équivalents (pour les équivalences globales de préfaisceaux simpliciaux).
\end{prop}
\begin{pr}
C'est un corollaire immédiat de la proposition \ref{comparclassifmodele} et du fait que les changements de base respectent les objets fibrants et cofibrants.
\end{pr}

\subsubsection{Espaces de chemins}

On renvoie au {\S\ref{chpdiago}} pour la définition du préchamp $\Omega_{x,y}C$ associé un préchamp simplicial $C$ pointé en deux objets $x$ et $y$.

\begin{prop}\label{cheminquillen}
Soit $C\in Pr(S,\scmsg_{cf}))$, $s\in S$ et $x,y\in C(s)$ 
$$
\Omega_{x,y}|C| \simeq \mathcal{E}q^\Delta_C(x,y).
$$
\end{prop}
\begin{pr}
Cela provient de la proposition \ref{comparclassifmodele} et du fait que les changements de base respectent les objets fibrants et cofibrants. On a les équivalences de préchamps :
$$
\Omega_{x,y}|C| \simeq \Omega_{x,y}|\mathcal{G}(C)| \simeq \mathcal{E}q^\Delta_C(x,y).
$$
\end{pr}

\subsection{Préfaisceaux faibles de Quillen}

Le cadre naturel d'apparition, ou de formalisation des problèmes de modules n'est pas tant celui des préfaisceaux que des préfaisceaux faibles. De ce fait, historiquement, les champs ont plutôt été développés en termes de catégories fibrées scindées clivées, langage équivalent à celui des préfaisceaux faibles (cf. \cite[ch. VI]{sga1}). Ce double langage correspond à modéliser un objet par sa catégorie des points ou son foncteur des points ; on a préféré dans ce travail le point de vue fonctoriel.


\medskip

Dans cette section on définit les préfaisceaux faibles de Quillen et on prouve que leur strictifié est un préfaisceau de Quillen. Puis on prouve que les espaces de morphismes d'un préfaisceau faible de Quillen simplicial et de sont strictifié sont canoniquement isomorphes.

\subsubsection{Strictification}\label{strictification}

On se contente de rappeler brièvement certaines définitions et de poser quelques notations, la référence prise est {\cite[ch. 5 \& App. B]{hollander}}.

\medskip
 \begin{defi}\label{prefaible}
Soit $S$ une catégorie, si $H$ est une 2-catégorie (faible ou stricte), on définit un {\em préfaisceau faible}\index{préfaisceau faible}  $\widetilde{C}$ sur $S$ à valeur dans $H$ par essentiellement la même donnée qu'un préfaisceau mais où la condition de fonctorialité est donnée par un 2-isomorphisme $\widetilde{c}_{v,u} : \widetilde{c}_u\circ \widetilde{c}_v \Rightarrow \widetilde{c}_{vu}$ rendant commutatifs tous les carrés
\begin{eqnarray}\label{cdtpreff}
\xymatrix{
\widetilde{c}_u\circ \widetilde{c}_v\circ \widetilde{c}_w \ar[r]^{1_{\widetilde{C}_w}*\widetilde{c}_{v,w}}\ar[d]_{\widetilde{c}_{v,u}*1_{\widetilde{c}_w}} & \widetilde{c}_u\circ \widetilde{c}_{wv} \ar[d]^{\widetilde{c}_{wv,u}}\\
\widetilde{c}_{vu}\circ \widetilde{c}_w \ar[r]_{\widetilde{c}_{w,vu}} & \widetilde{c}_{wvu}
}
\end{eqnarray}
\end{defi}

\medskip

Soient $\widetilde{C}$ et $\widetilde{D}$ deux préfaisceaux faibles sur $S$, une {\em transformation naturelle faible} $(a,\alpha)$ est la donnée
\begin{itemize}
\item pour tout objet $s$ de $S$ d'un foncteur $a_s : \widetilde{C}(s)\to \widetilde{D}(s)$ ;
\item et pour tout morphisme $u:t\to s\in S$ d'un isomorphisme naturel $\alpha_u$ faisant commuter le diagramme
$$\xymatrix{
\widetilde{C}(s)\ar[r]^{a_i} \ar[d]_{\widetilde{c}_u} & \widetilde{D}(s) \ar[d]^{\widetilde{d}_u} \\
\widetilde{C}(t)\ar[r]^{a_j} & \widetilde{D}(t)
}$$
(\ie $\alpha_u : a_t\circ \widetilde{c}_u \Rightarrow \widetilde{d}_u\circ a_s $)
\item tels que, pour toute paire de morphismes $t'\overset{v}{\to} t\overset{u}{\to} s$ de $S$ le diagramme suivant soit commutatif
$$\xymatrix{
& a_k \circ \widetilde{c}_v \circ \widetilde{c}_u \ar[rd]^{\alpha_v*\widetilde{c}_u} \ar[ld]_{a_k\circ \widetilde{c}_{v,u}} & \\
a_k \circ \widetilde{c}_{uv} \ar[d]_{\alpha_{uv}} & & \widetilde{d}_v\circ a_j \circ \widetilde{d}_u \ar[d]^{\widetilde{d}_v*\alpha_u} \\
\widetilde{d}_{uv}\circ a_i  & & \widetilde{d}_v\circ \widetilde{d}_u\circ a_i\ar[ll]^{ \widetilde{d}_{u,v}}\\
}$$
Ce diagramme correspond aux cinq façons d'aller de $\widetilde{C}(s)$ à $\widetilde{C}(t')$ dans le diagramme suivant :
$$\xymatrix{
\widetilde{C}(s)\ar[r]^{a_i}\ar[d]_{\widetilde{c}_u} \ar@/_2pc/[dd]_{\widetilde{c}_{vu}} & \widetilde{D}(s) \ar[d]^{\widetilde{d}_u} \ar@/^2pc/[dd]^{\widetilde{d}_{vu}}\\
\widetilde{C}t)\ar[r]^{a_j}\ar[d]_{\widetilde{c}_v} & \widetilde{D}(t) \ar[d]^{\widetilde{d}_v} \\
\widetilde{C}(t')\ar[r]^{a_k} & \widetilde{D}(t').\\
}$$
\end{itemize}

Une transformation naturelle faible est un {\em isomorphisme naturel faible} si tous les $\alpha_u$ sont des isomorphismes.

\medskip
On a également une notion de {\em 2-transformation naturelle} qu'on ne détaille pas. 

Toutes ces définitions définissent une 2-catégorie des préfaisceaux faibles sur $S$ à valeurs dans $H$, notée $Prf(S,H)$\index{$Prf(S,H)$}.
Par opposition au préfaisceaux faibles, les préfaisceaux seront qualifiés de {\em stricts}\index{préfaisceau strict}. La catégorie $Pr(S,H)$ des préfaisceaux stricts de $S$ à valeurs dans $H$ est munie d'un morphisme pleinement fidèle $\iota : Pr(S,H) \to Prf(S,H)$.

\medskip

Si $H=\catcat$, la catégories des catégorie, vue comme une 2-catégorie, un préfaisceau (strict ou faible) à valeurs dans $\catcat$ est dit un préfaisceau (strict ou faible) en catégories.
Dans ce cas une transformation naturelle faible $(a,\alpha) : \widetilde{C}\Rightarrow \widetilde{D}$ est la donnée pour tout objet $s\in S$ de foncteurs $a_s :\widetilde{C}(s)\Rightarrow \widetilde{D}(s)$ vérifiant les conditions de la définition. $(a,\alpha)$ est dite une {\em équivalence} si tous les $\alpha_u$ sont des équivalence de catégories.
Ces sous-catégories d'équivalences de $Pr(S,\catcat)$ et $Prf(S,\catcat)$ s'enrichissent en des structures de modèles (cf. \cite[ch. 7]{hollander}).

\medskip

Si $\widetilde{C}$ est un préfaisceau faible en catégories on rappelle (cf. \cite[ch. 5]{hollander}) que son {\em strictifié}\index{strictifié d'un préfaisceau faible} $C$ est le préfaisceau 
$$
C:= \Hom_{Prf(S,\catcat)}(h_{(-)},\widetilde{C})
$$
où $\hom_{Prf(S,\catcat)}(-,\widetilde{C})$ est le foncteur des points de $\widetilde{C}$ dans la catégorie $Prf(S,\catcat)$ et où $h_{(-)}$ est le foncteur $S\to Pr(S,\ens)\to Pr(S,\catcat)\to Prf(S,\catcat)$ (la première flèche est le plongement de Yoneda). Cette construction définit un foncteur $\mathcal{S} : Prf(S,\catcat) \to Pr(S,\catcat)$ adjoint à droite de $\iota$.

\begin{rprop}[{\cite[cor. 7.4]{hollander}}]\label{strictifholl}
Avec les notations précédentes, l'adjonction $\iota : Pr(S,\catcat) \leftrightarrow Prf(S,\catcat) : \mathcal{S}$ 
est une équivalence de Quillen.
\end{rprop}

\begin{defi}
Un {\em préfaisceau faible de Quillen à gauche} (resp. {\em à droite})\index{préfaisceau faible de Quillen} sur $S$ est un objet de $Prf(S,\cmg)$ (resp. $Prf(S,\cmd)$) où $Prf(S,\cmg)$ (resp. $Prf(S,\cmd)$) est pris avec sa structure naturelle de 2-catégorie donnée par les transformations naturelles de foncteurs.
\end{defi}

\begin{prop}\label{strictifquillen}
Le strictifié d'un préfaisceau faible de Quillen (resp. simplicial) est canoniquement un préfaisceau de Quillen (resp. simplicial).
\end{prop}
\begin{pr}
Soit $\widetilde{C}$ un préfaisceau faible en catégorie sur $S$ et $C$ sont strictifié. Pour chaque $s\in S$ on a une équivalence de catégories $\widetilde{C}(s)\overset{\sim}{\to}C(s)$ qu'on va expliciter et par laquelle on peut transporter la structure de modèles de $\widetilde{C}(x)$ sur $C(s)$.

Un objet $(a,\alpha)$ de $C(s)$ consiste en la donnée 
\begin{itemize}
\item pour tout $u:t\to s\in S$ d'un objet $a_u\in \widetilde{C}(t)$ ;
\item pour tout $t'\overset{v}{\to} t\overset{u}{\to} s$ d'un isomorphisme $\alpha_v :\widetilde{c}_v(a_u)\overset{\sim}{\to}a_{uv}\in\widetilde{C}(t')$
\item tels que pour tout $t''\overset{w}{\to} t'\overset{v}{\to} t\overset{u}{\to} s\in S$ le diagramme
$$\xymatrix{
\widetilde{c}_w\widetilde{c}_v(a_u) \ar[d]_{\widetilde{c}_{v,w}}\ar[r]^{\widetilde{c}_w(\alpha_v)} & \widetilde{c}_w(a_{uv})\ar[d]^{\alpha_w} \\
\widetilde{c}_{vw}(a_u) \ar[r]_{\alpha_{vw}} & a_{uvw}
}$$
soit commutatif.
\end{itemize}
Un morphisme  $m:(a,\alpha)\to (b,\beta)$ de $\widetilde{C}(s)$ est la donnée
\begin{itemize}
\item pour tout $u:t\to s\in S$ de morphismes $m_u:a_u\to b_u\in \widetilde{C}(y)$ ;
\item tels que pour tout $t'\overset{v}{\to} t\overset{u}{\to} s\in S$ le diagramme
$$\xymatrix{
\widetilde{c}_v(a_u) \ar[d]_{\alpha_v}\ar[r]^{\widetilde{c}_v(m_u)} &\widetilde{c}_v(b_u) \ar[d]^{\beta_v} \\
a^{vu} \ar[r]_{m_{vu}} & b^{vu}.
}$$
soit commutatif.
\end{itemize}

Les équivalences entre $\widetilde{C}(s)$ et $C(s)$ sont données par les unités et co-unités de l'adjonction de strictification sont les foncteurs {\em extension} et {\em troncation} :
\begin{eqnarray*}
\epsilon : \widetilde{C}(s) & \longrightarrow & C(s) \\
a & \longmapsto & (a,\alpha) := (\{\widetilde{c}_u(a)\}_u,\{\widetilde{c}_{u,v} : \widetilde{c}_v\widetilde{c}_u(a)\to \widetilde{c}_{uv}(a) \}_v)\\
m : a\to b & \longmapsto & m := \{\widetilde{c}_u(m) : \widetilde{c}_u(a)\to \widetilde{c}_u(b) \}_u\\
\ \\
\tau : C(s) & \longrightarrow & \widetilde{C}(s) \\
(a,\alpha) = (\{a_u\}_u,\{\alpha_v\}_v) & \longmapsto & a_{1_x}\\
m=\{m_u\}_u : (\{a_u\}_u,\{\alpha_v\}_v) \to (\{b_u\}_u,\{\beta_v\}_v) & \longmapsto & m_{1_x}:a_{1_x}\to b_{1_x}\\
\end{eqnarray*}

On se sert de ces équivalences pour transporter la structure de modèles de $\widetilde{C}(s)$ sur $C(s)$ en les transformant en équivalences de Quillen.
Si $u:t\to s$ et $v:t'\to s'$ sont deux flèches d'une catégorie $C$, on dit qu'elles sont {\em isomorphes} s'il existe $a:x\to z$ et $b:y\to t$ des isomorphismes tels que $va=bu$ ; les sous-catégories d'équivalences, de cofibrations et de fibrations de $C(s)$ sont alors données par les clôtures pour l'isomorphisme des flèches des sous-catégories image par $\epsilon$ des sous-catégories des équivalences, des cofibrations et des fibrations de $\widetilde{C}(s)$.

En particulier une flèche $m$ de $C(s)$ est une équivalence, une cofibration ou une fibration ssi $\tau(m)=m_{1_x}$ en est une. En revanche les autres $m_u$ peuvent ou pas être des équivalences, des cofibration ou des fibrations.
Toutefois, si $\widetilde{C}$ est de Quillen à gauche (resp. à droite) $m$ est une cofibration (resp. une fibration) ssi tous les $m_u$ en sont ; il en est de même avec les cofibrations et les fibrations triviales.


\medskip

À un morphisme $\chi:t\to s\in C$ est associé un foncteur de changement de base $c_{\chi}:C(s)\to C(t)$ qui se décrit au niveau des objets par 
$$
(a,\alpha) = (\{a_u\}_u,\{\alpha_v\}_v) \longmapsto  c_\chi(a,\alpha) = (\{a'_{u'}\}_{u'},\{\alpha'_{v'}\}_{v'}) = (\{a_{u'\chi}\}_{u'},\{\alpha_{v'}\}_{v'})
$$
et au niveau des flèches par
$$
m=\{m_u\}_u \longmapsto c_\chi(m)=\{m'_{u'}\}_{u'} = \{m_{u'\chi}\}_{u'}.
$$
Ainsi, si $\widetilde{C}$ est de Quillen à gauche et que $m$ est une cofibration (resp. une cofibration triviale), \ie si tous les $m_u$ sont des cofibrations (resp. des cofibration triviales), il en est de même de $c_\chi(m)$ : les foncteurs $c_{\chi}$ préservent les cofibrations et les cofibrations triviales.

Il reste à voir que les changements de base sont adjoints à gauche ; pour cela on considère le diagramme de foncteurs : 
$$\xymatrix{
\widetilde{C}(s) \ar[d]_{\widetilde{c}_\chi}\ar[r]^\epsilon & C(s) \ar[d]^{c_\chi} \\
\widetilde{C}(t) \ar[r]_\epsilon & C(t) \\
}$$
ce diagramme est commutatif à un isomorphisme naturel près donné essentiellement par les $\widetilde{c}_{v,u}:\widetilde{c}_u\circ \widetilde{c}_v\to \widetilde{c}_{vu}$.
Ainsi $c_\chi$ est adjoint à gauche (ou a droite) car équivalent au foncteur $\widetilde{c}_\chi$ qui est adjoint à gauche (ou à droite) si $\widetilde{C}$ est de Quillen à gauche (ou à droite).

\medskip

Pour la structure simpliciale, on la transporte des $\widetilde{C}(s)$ vers les $C(s)$ par $\epsilon$ en posant, pour tous $K\in \sens$ et $(a,\alpha)=(\{\widetilde{c}_u(a)\}_u,\{\widetilde{c}_{u,v} : \widetilde{c}_v\widetilde{c}_u(a)\to \widetilde{c}_{uv}(a) \}_v)\in C(s)$
$$
K\otimes (a,\alpha) := (\{K\otimes \widetilde{c}_u(a)\}_u,\{c^K_{u,v} : \widetilde{c}_v(K\otimes \widetilde{c}_u(a))\to K\otimes \widetilde{c}_{uv}(a) \}_v)
$$
où $c^K_{u,v}$ est la composition
$$
\widetilde{c}_v(K\otimes \widetilde{c}_u(a))\overset{k_v}{\longrightarrow}
K\otimes \widetilde{c}_v(\widetilde{c}_u(a))\overset{\widetilde{c}_{u,v}}{\longrightarrow}
K\otimes \widetilde{c}_{uv}(a)
$$
où $k_v$ est l'isomorphisme de la structure simpliciale de $\widetilde{c}_v$.
$\epsilon$ étant essentiellement surjectif, on prolonge cette définition aux objets qui ne sont pas dans son image module le choix d'un isomorphisme vers un objet de l'image.

Les axiomes $M6$ et $M7$ de la définition \cite[def. 9.1.6]{hirschhorn} des catégories de modèles simpliciales  sont vérifiées car $\epsilon$ est une équivalence de catégorie.
\end{pr}

\subsubsection{Préchamps de morphismes}

Le résultat principal de cette section est le lemme {\ref{strictifhomsplx}} établissant qu'un préfaisceau faible en catégories simpliciales (définition {\ref{prefsplxfaible}}) et son strictifié ont les mêmes préfaisceaux de morphismes entre deux objets. Il est l'outil technique principal de la proposition {\ref{cheminquillenfaible}} autorisant le calcul de la diagonale des champs associés aux préfaisceaux de Quillen.

\bigskip

Le lemme suivant établit que les ensembles de morphismes entre deux objets d'un préfaisceau faible en catégories forment naturellement des préfaisceaux.
\begin{lemme}
Soit $\widetilde{C}$ est un préfaisceau faible en catégories discrètes sur $S$, et soient, pour $s\in S$, deux objets $x$ et $y$ de $C(s)$, Alors 
\begin{eqnarray*}
\Hom_{\widetilde{C}}(x,y) : (S/s)^o & \longrightarrow & \ens\\
u:t\to s & \longmapsto & \Hom_{\widetilde{C}(t)}(u^*x,u^*y)\\
t'\overset{v}{\to} t\overset{u}{\to}s & \longmapsto & \psi\circ v^* : \Hom_{\widetilde{C}(t)}(u^*x,u^*y)\to \Hom_{\widetilde{C}(t')}((vu)^*x,(vu)^*y)
\end{eqnarray*}
où $\psi$ est issu de l'isomorphisme $v^*u^* \to (vu)^*$, est un préfaisceau.
\end{lemme}
\begin{pr}
On note $\widetilde{c}_u$ les transitions et $\widetilde{c}_{u,v}$ les isomorphismes de fonctorialité de $\widetilde{C}$. On considère un objet $s\in S$ et deux objets $x,y\in C(s)$, pour tout $t'\overset{v}{\to} t \overset{u}{\to}s\in S$ on a un isomorphisme 
$$
\gamma_{u,v}^x:f_vf_ux\longrightarrow f_{uv}x 
$$
duquel on tire un isomorphisme 
\begin{eqnarray*}
\psi_{u,v} : \Hom_{\widetilde{C}(t')}(\widetilde{c}_v\widetilde{c}_ux,\widetilde{c}_v\widetilde{c}_uy)&\longrightarrow &\Hom_{\widetilde{C}(t')}(\widetilde{c}_{uv}x,\widetilde{c}_{uv}y)\\
a & \longmapsto & \gamma_{u,v}^y\circ a\circ (\gamma_{u,v}^x)^{-1}
\end{eqnarray*}

$\Hom_{\widetilde{C}}(x,y)$ est un préfaisceau ssi pour tout triplet de morphismes $t''\overset{w}{\to}t'\overset{v}{\to} t \overset{u}{\to} s$ on a 
$$
(\psi_{uv,w}\circ c_{w})\circ (\psi_{u,v}\circ c_v) = \psi_{u,vw}\circ c_v
$$
\ie dans le diagramme suivant on veut montrer que $3=2\circ 1$
$$\xymatrix{
\Hom_{\widetilde{C}(t)}(c_ux,c_uy)\ar[r]^{c_v} \ar@/_1pc/[rdd]_{c_{vw}} \ar@/^2pc/@{-->}[rr]^1 \ar@/_6pc/@{-->}[rrdd]_3 \ar@{}[rd]|A & \Hom_{\widetilde{C}(t')}(c_vc_ux,c_vc_uy)\ar[r]^{\psi_{u,v}} \ar[d]^{c_w} \ar@{}[rd]|B & \Hom_{\widetilde{C}(t')}(c_{uv}x,c_{uv}y) \ar[d]^{f_w}\ar@/^5pc/@{-->}[dd]^2 \\
& \Hom_{\widetilde{C}(t'')}(c_wc_vc_ux,c_wc_vc_uy)\ar[r]^{c_w*\psi_{u,v}} \ar[d]^{\psi_{v,w}*c_u} \ar@{}[rd]|C & \Hom_{\widetilde{C}(t'')}(c_wc_{uv}x,c_wc_{uv}y)\ar[d]^{\psi_{uv,w}}\\
& \Hom_{\widetilde{C}(t'')}(c_{wv}c_ux,c_{wv}c_uy)\ar[r]^{\psi_{u,vw}} & \Hom_{\widetilde{C}(t'')}(c_{wvu}x,c_{wvu}y)
}$$
Il suffit pour cela de montrer que les trois carrés $A,\ B$ et $C$ sont commutatifs : pour $A$, c'est la définition de $\gamma_{v,w}$ ; pour $C$, c'est la relation de cohérence des $\gamma$ et pour $B$, c'est la fonctorialité de $c_w$.
\end{pr}

\begin{lemme}
Deux faibles préfaisceaux faiblement isomorphes donnent lieu à des préfaisceaux de morphismes isomorphes.
\end{lemme}
\begin{pr}
Soit $\widetilde{C}$ et $\widetilde{D}$ deux préfaisceaux sur $S$ et $(a,\alpha):\widetilde{C}\to \widetilde{D}$ un isomorphisme naturel faible entre eux.
Soient $s\in S$ et $x,y\in \widetilde{C}(s)$ dont on note $a_x$ et $a_y$ les objets correspondant dans $\widetilde{D}(s)$.
On a les deux faisceaux des morphismes :
\begin{eqnarray*}
\Hom_{\widetilde{C}}(x,y) : (S/s)^o & \longrightarrow & \ens \\
t\overset{u}{\to}  s &\longmapsto & \Hom_{\widetilde{C(t)}}(\widetilde{c}_ux,\widetilde{c}_uy)\\
\Hom_{\widetilde{D}}(a_x,a_y) : (S/s)^o & \longrightarrow & \ens \\
t\overset{u}{\to}  s &\longmapsto & \Hom_{\widetilde{D(t)}}(\widetilde{d}_ua_x,\widetilde{d}_ua_y).
\end{eqnarray*}
entre lesquels $(a,\alpha)$ induit un morphisme :
$$
\alpha_u\circ a : \widetilde{C}(t)(\widetilde{c}_ux,\widetilde{c}_uy) \overset{a}{\longrightarrow} \widetilde{D}(t)(a\widetilde{c}_ux,a\widetilde{c}_uy)\overset{\alpha_u}{\longrightarrow} \widetilde{D}(t)(\widetilde{d}_ua_x,\widetilde{d}_ua_y)
$$
Pour que ce morphisme soit un morphisme de préfaisceaux il faut vérifier que, pour tout $v:t'\to t$, le diagramme
$$\xymatrix{
\widetilde{C}(t)(\widetilde{c}_ux,\widetilde{c}_uy) \ar[r]^-{\alpha_u\circ a}\ar[d]_{\psi_{v,u} \circ f_v} & \widetilde{D}(t)(\widetilde{d}_ua_x,\widetilde{d}_ua_y) \ar[d]^{\psi_{v,u} \circ \widetilde{d}_v} \\
\widetilde{C}(t')(\widetilde{c}_{uv}x,\widetilde{c}_{uv}y) \ar[r]_-{\alpha_{uv}\circ a} & \widetilde{D}(t')(\widetilde{d}_{uv}a_x,\widetilde{d}_{uv}a_y)
}$$
est commutatif. Pour cela il suffit de montrer que les sous-diagrammes $A,\ B,\ C$ et $D$ du diagramme suivant sont commutatifs :
$$\xymatrix{
\widetilde{C}(t)(\widetilde{c}_ux,\widetilde{c}_uy) \ar@{}[rrd]|A \ar[rr]^-a\ar[d]_{\widetilde{c}_v} & & \widetilde{D}(t)(a\widetilde{c}_ux,a\widetilde{c}_uy) \ar@{}[rd]|B\ar[d]^{\widetilde{d}_v} \ar[r]^-{\alpha_u} & \widetilde{D}(t)(\widetilde{d}_ua_x,\widetilde{d}_ua_y) \ar[d]^{\widetilde{d}_v} \\
\widetilde{C}(t')(\widetilde{c}_v\widetilde{c}_ux,\widetilde{c}_v\widetilde{c}_uy) \ar@{}[rd]|C \ar[r]^-a \ar[d]_{\psi_{v,u}} & \widetilde{D}(t')(a\widetilde{c}_v\widetilde{c}_ux,a\widetilde{c}_v\widetilde{c}_uy) \ar@{}[rrd]|D \ar[d]^{a\circ\phi_{u,v}} \ar[r]^-{\alpha_v*\widetilde{c}_u} & \widetilde{D}(t')(\widetilde{d}_va\widetilde{c}_ux,\widetilde{d}_va\widetilde{c}_uy) \ar[r]^-{\widetilde{d}_v*\alpha_u} & \widetilde{D}(t')(\widetilde{d}_v\widetilde{d}_ua_x,\widetilde{d}_v\widetilde{d}_ua_y)\ar[d]^{\gamma_{u,v}}\\
\widetilde{C}(t')(\widetilde{c}_{uv}x,\widetilde{c}_{uv}y) \ar[r]_-a & \widetilde{D}(t')(a\widetilde{c}_{vu}x,a\widetilde{c}_{vu}y)\ar[rr]_-{\alpha_{vu}} & & \widetilde{D}(t')(\widetilde{d}_{uv}a_x,\widetilde{d}_{uv}a_y)
}$$
Or, la commutativité de $A$ est due à la définition des $\alpha$ ;
celle de $B$ est de à la fonctorialité de $g_v$ ;
celle de $C$ est de à la fonctorialité de $a$ ;
celle de $D$, enfin, est de à la condition de cohérence pentagonale de $(a,\alpha)$.
\end{pr}

Compte tenu qu'un préfaisceau faible et son strictifié sont toujours faiblement isomorphes, on a le corollaire suivant.
\begin{lemme}\label{strictifhom}
Soient $\widetilde{C}$ un préfaisceau faible sur $S$ en catégories et, pour $s\in S$ fixé, deux objets $x,y$ de $C(s)$. On note $C$ le strictifié de $\widetilde{C}$. On a un isomorphisme canonique de préfaisceaux en ensembles :
$$
\Hom_{\widetilde{C}}(x,y) \overset{\sim}{\longrightarrow} \Hom_C(x,y).
$$
\end{lemme}

\bigskip

$\sens$ est une catégorie monoïdale ; on note $\sens\z\catcat$ la 2-catégorie des $\sens$-modules dans $\catcat$. Une catégorie $C$ qui est un module sur $\sens$ s'enrichit naturellement en une catégorie  simpliciale, notée $C^{\Delta}$\index{$C^\Delta$} (cf. \cite[ch. IV]{hovey} et {\sf\S\ref{catsplx}}).

\begin{defi}\label{prefsplxfaible}
Un {\em préfaisceau faible en catégories simpliciales}\index{préfaisceau faible en catégories simpliciales} est un préfaisceau faible à valeur dans $\sens\z\catcat$ (cf. définition {\sf\S\ref{prefaible}}).
Une telle donnée définit toujours un préfaisceau faible à valeurs dans $\scatcat$ et c'est ce qui en justifie le nom.

Le {\em strictifié} d'un préfaisceau faible en catégorie simpliciale est son strictifié comme préfaisceau en catégorie, comme démontré dans la preuve de la proposition {\ref{strictifquillen}}, il est canoniquement muni d'une structure de $\sens$-module.
\end{defi}

\medskip

Si $C\in\scatcat$, on note $C_n$ la catégorie de ses n-simplexes.
La définition {\ref{prefsplxfaible}} est telle que, pour tout $n\in\Delta$ et tout préfaisceau faible en catégories simpliciales $\widetilde{C}$, les préfaisceaux en catégories $\widetilde{C}_n : s\mapsto C(s)^\Delta_n$ soient faibles au sens de la définition {\ref{prefaible}}.

\bigskip

Si $\widetilde{C}$ est un préfaisceau faible en catégories simpliciales sur $S$, et si, pour $s\in S$, $x$ et $y$ sont deux objets de $C(s)$, on définit le préfaisceau simplicial
\begin{eqnarray*}
\Hom_{\widetilde{C}}^\Delta(x,y) : (S/s)^o & \longrightarrow & \sens\\
	u:t\to s & \longmapsto & \Hom_{\widetilde{C}(t)}^\Delta(u^*x,u^*y)\\
	t'\overset{v}{\to} t\overset{u}{\to}s & \longmapsto & \alpha\circ v^* : \Hom_{\widetilde{C}(t)}^\Delta(u^*x,u^*y)\to \Hom_{\widetilde{C}(t')}^\Delta((vu)^*x,(vu)^*y)
\end{eqnarray*}
où, $\Hom_{\widetilde{C}(t)}^\Delta(-,-)$ désigne l'ensemble simplicial des morphismes dans $C(t)$ et où $\alpha$ désigne l'isomorphisme $v^*u^* \to (vu)^*$.

\begin{lemme}\label{strictifhomsplx}
Soient $\widetilde{C}$ un préfaisceau faible sur $S$ en catégories simpliciales et, pour $s\in S$ fixé, deux objets $x,y$ de $C(s)$. On note $C$ le strictifié de $\widetilde{C}$. On a un isomorphisme canonique de préfaisceaux simpliciaux :
$$
\Hom_{\widetilde{C}}^\Delta(x,y) \overset{\sim}{\longrightarrow} \Hom_C^\Delta(x,y).
$$
\end{lemme}
\begin{pr}
Cela découle du lemme {\ref{strictifhom}} en remarquant que, pour tout $n\in \Delta$, $\Hom_{\widetilde{C}}^\Delta(x,y)_n = \Hom_{\widetilde{C}}(x\otimes \Delta^n,y)\simeq \Hom_{C}(x\otimes \Delta^n,y) = \Hom_{C}^\Delta(x,y)_n$.
\end{pr}

\medskip

\begin{prop}\label{cheminquillenfaible}
Soient ${\widetilde{C}}$ un sous-préfaisceau faible de Quillen simplicial dont les restrictions conservent les objets fibrants et cofibrants et, pour $i\in I$ fixé, deux objets $x,y$ de $\widetilde{C}_i^{cf}$. On note $C$ le préfaisceau de Quillen strictifié de $\widetilde{C}$. On a une équivalence (globale) de préfaisceaux simpliciaux :
$$
\Omega_{xy}|C| \simeq \Eq_{\widetilde{C}}^\Delta(x,y).
$$
\end{prop}
\begin{pr}
C'est un corollaire immédiat de la proposition {\ref{cheminquillen}} et du lemme {\ref{strictifhomsplx}}.
\end{pr}

\begin{cor}\label{cheminchampquillen}
Si $\underline{C}$ est le champ associé à un sous-préfaisceau faible de Quillen simplicial $\widetilde{C}$ dont les restrictions conservent les objets fibrants et cofibrants et si $x,y$ deux $i$-points de $\widetilde{C}^{cf}$
et si  alors $\underline{\Eq_{\widetilde{C}}^\Delta(x,y)}$ est le champs associé au préchamp $\Eq_{\widetilde{C}}^\Delta(x,y)$, on a 
$$
\Omega_{xy}\underline{C} \simeq \underline{\Eq_{\widetilde{C}}^\Delta(x,y)}.
$$
\end{cor}
\begin{pr}
D'après le lemme {\ref{commdiagochpif}}, on a
$$
\Omega_{xy}\underline{C} \simeq \underline{\Omega_{xy}|C|}.
$$
\end{pr}

\medskip

Pour deux objets $x$ et $y$ non fibrants-cofibrants d'un sous-préfaisceau de Quillen simplicial $C$, on considère 
$Tx,$ et $Ty$ où $T$ est un foncteur de remplacement fibrant-cofibrant ; comme pour tout objet $z$, $z$ et $Tz$ sont fonctoriellement équivalents, on a des équivalences (canoniques) $\Omega_{x,y}\underline{C}\simeq \Omega_{Tx,Ty}\underline{C}\simeq \underline{\Eq_{\widetilde{C}}^\Delta(Tx,Ty)}$.

\subsection{Critère de lissité}

Soient $f:M\to N$ un morphisme entre deux sous-préfaisceaux de Quillen. Pour $A$ un anneau et $A'\to A$ une extension infinitésimale au premier ordre de $A$, on a le diagramme commutatif suivant :
$$\xymatrix{
M(A') \ar[d]_{\mu} \ar[r]^{f_{A'}} & N(A') \ar[d]^{\nu} \\
M(A) \ar[r]_{f_A} & N(A)
}$$
où les flèches horizontales sont celles de $f$ et les verticales les restrictions de $M$ et $N$.

D'après la proposition {\ref{lissite}} le morphisme de champs associé à $f$, $\underline{f}:\underline{M}\to \underline{N}$ est lisse ssi l'application naturelle $|W^c_{M(A')}|\to |W^c_{M(A)}|\times^h_{|W^c_{N(A)}|}|W^c_{N(A')}|$ est surjective à homotopie près.

On note $(x,y,z,u,v)$ les objets de $M(A)\times_{N(A)}N(A')$ (cf. {\S \ref{profibnerf}}).
Un objet $m$ de $M(A')$ donne dans $M(A)\times_{N(A)}N(A')$ l'objet $(\mu(m), f_A(\mu(m))=\nu(f_{A'}(m)), f_{A'}(m),id,id)$.

\begin{prop}\label{lissequillen}
(On utilise les notations précédentes et celles du {\S\ref{profibnerf}}.)
Soit $f:M\to N$ un morphisme de Quillen à gauche entre deux préfaisceaux de Quillen à gauche, le morphisme entre les champs associés est formellement lisse si et seulement si, pour tout anneau henselien $A$ et toute extension au premier ordre $A'$, tout objet $(x,y,z,u,v)$ de $(M(A)\times_{N(A)}N(A'))_{cart}$ est connecté par une chaîne d'équivalences à un objet du type $(\mu(m), f_A(\mu(m)), f_{A'}(m),id,id)$.
\end{prop}
\begin{pr}
On utilise la proposition {\ref{profibnerfprop}} pour caractériser $|W^c_{M(A)}|\times^h_{|W^c_{N(A)}|}|W^c_{N(A')}|$, avec ses notations $f$ est lisse ssi $|W^c_{M(A')}|\to |W^c_{cart}|$ est surjective. La proposition découle de la description des objets de $|W^c_{cart}|$.
La considération des anneaux henseliens, \ie locaux pour la topologie étale, est là pour éviter d'avoir à traiter les relèvements des points modulo un recouvrement.
\end{pr}

\setcounter{tot}{0}

\newpage
\thispagestyle{empty}
\chapter{Catégories linéaires}\label{chapcatlin}
\thispagestyle{chheadings}

Ce chapitre définit les catégories linéaires et les notions classiques les accompagnant. La partie originale consiste en la construction de schémas et champs classifiant les structures de catégories et d'autres objets, ainsi qu'en celle de la structure de modèles sur les catégories des catégories linéaires pour laquelle les équivalences sont les équivalences de Morita.
\footnote{La volonté de construire les classifiants des structures définies sitôt après leur définition, nous a obligé à quelques renvois dans l'avant du texte.}

Les équations des structures de catégorie, de foncteur, ou de modules sont exprimés en convention d'Einstein :
les vecteurs des bases ont leur indices en base, les indices contravariants sont notés en exposants et il y a toujours une somme implicite sur un même indice qui apparaît en haut et en bas.

\section{Graphes linéaires}
\label{graphes}\label{graphe}

\begin{defi}
Soient $E$ un ensemble et $A$ un anneau commutatif, on définit un {\em $A$-graphe sur $E$}, $C$, comme la donnée pour tout $(x,y)\in E^2$ d'un $A$-module $C(x,y)$ ; $E$ est l'ensemble des {\em objets du graphe}.
Un $A$-graphe est un graphe sur un ensemble non précisé.
Un {\em graphe linéaire}\index{graphe linéaire} est un $A$-graphe pour un anneau $A$ imprécisé.

On dit qu'un $A$-graphe est {\em libre}\index{graphe libre}, {\em projectif}\index{graphe projectif} ou de {\em type fini}\index{graphe de type fini} si les $A$-modules le constituant le sont tous (avec, pour le dernier cas, la condition supplémentaire que l'ensemble des objets soit fini).
Le {\em type}\index{type d'un graphe linéaire} d'un $A$-graphe projectif de type fini est le couple de son ensemble d'objets et de l'application $d:E^2\to\N$ qui à $(x,y)$ associe le rang de $C(x,y)$.

Un {\em morphisme strict}\index{morphisme strict de graphes} $f$ entre deux $A$-graphes $C$ et $D$ sur $E$ est la donnée pour tout $(x,y)\in E^2$ de morphismes de $A$-modules $f_{xy}:C(x,y)\to D(x,y)$ ; c'est un isomorphisme ssi tous les $f_{xy}:C(x,y)\to D(x,y)$ le sont.

Un {\em morphisme (non strict)}\index{morphisme de graphes} $f$ entre deux $A$-graphes $C$ et $D$ sur $E$ et $F$ respectivement est la donnée d'une application $f_0:E\to F$ (qui souvent sera noté $f$ dans la suite) et pour tout $(x,y)\in E^2$ de morphismes de $A$-modules $f_{xy}:C(x,y)\to D(f_0(x),f_0(y))$.
Les isomorphismes se caractérisent par le fait que $f$ soit une bijection et les $f_{xy}:C(x,y)\to D(f_0(x),f_0(y))$ des isomorphismes de $A$-modules.
\end{defi}

\medskip

\paragraph{Champs classifiants}
Soit $\uvect$\index{$\uvect$} le champ classifiant des fibrés vectoriels (\ie les modules projectifs) de rang fini. $\uvect$ est la réunion disjointe des champs classifiant les fibrés vectoriels de rang fixé, dont une présentation est donné par les groupoïdes associés aux $G\ell_n$, on a $\uvect \simeq \coprod_n \ubgln$ (cf. {\sf \S \ref{groupoide}}).
Pour $E$ un ensemble de cardinal $n$ et une fonction $d:E^2\to\N$ on définit les objets suivants :
\begin{eqnarray*}
G\ell_d \index{$G\ell_d$}&:=& \prod_{e\in E^2} G\ell_{d(e)},\\
\uvect^{(d)}\index{$\uvect^{(d)}$} := \underline{\B G\ell}_d\index{$\underline{\B G\ell}_d$} &:=& \prod_{e\in E^2} \underline{\B G\ell}_{d(e)},\\
\uvect^{(n)}\index{$\uvect^{(n)}$} := \uvect^{(E)} &:=& \prod_{E^2} \uvect = \coprod_{d:E^2\to\mathbb{N}} \uvect^{(d)}
\end{eqnarray*}
où la dernière réunion est prise sur toutes les applications $\delta:E^2\to\mathbb{N}$ avec $E$ fixé de cardinal $n$.
Les champs $\ubgln$ étant 1-géométriques, il en est de même de $\underline{\B G\ell}_d$ et $\uvect^{(E)}$.

Si $A$ est une algèbre commutative, les $A$-points de $\uvectd$ sont exactement les $A$-graphes projectifs de type $d$. Réciproquement tout $A$-graphe de type $d$ définit un point de $\uvectd$, mais non forcément unique (à cause des isomorphismes échangeant les objets). En fait, il est facile de voir qu'un isomorphisme entre deux points de $\uvectd$ fixe les objets des graphes.

\medskip
Soit $\SE$ le groupe des permutations de $E$, son action naturelle sur les applications $d:E^2\to\mathbb{N}$
le fait agir sur $\uvect^{(n)}$ par permutation des facteurs $\uvect^{(d)}$. On définit :
\begin{eqnarray*}
\underline{Gr}_n\index{$\underline{Gr}_n$} := \underline{Gr}_E := \uvect^{(E)}/\SE.
\end{eqnarray*}
L'action de $\SE$ est clairement lisse et le champ $\underline{Gr}_E$ est géométrique.
Les $A$-points de $\underline{Gr}_E$ sont aussi des $A$-graphes, mais cette fois les isomorphismes entre points permutent les objets.


\begin{prop}
Pour $(E,d)$ un type de graphe fixé, $\uvectd$ est le champ des modules des graphes linéaires sur $E$ de type $d$
à isomorphisme strict près ; $\uvect^{(E)}$ est le champ des modules des graphes linéaires sur $E$ à isomorphismes stricts près et $\underline{Gr}_E$ est le champ des modules des graphes linéaires sur $E$ à isomorphismes non stricts près. On a le produit fibré homotopique suivant :
$$\xymatrix{
\uvect^{(E)} \ar[d]\ar[r] & \underline{Gr}_E \ar[d] \\
pt \ar[r] & \underline{\B\SE}. \\
}$$
\end{prop}

\bigskip

Si $E$ et $F$ sont deux ensembles en bijection, les champs $\underline{Gr}_E$ et $\underline{Gr}_F$
sont en équivalence canonique, on désignera par $\underline{Gr}_n$ {\sl le} champ des graphes sur un ensemble à $n$ éléments.

\section{Catégories linéaires}\label{categolin}

\subsection{Définition et classifiant}\label{classifstcat}

Soient $C$ et $D$ deux graphes $A$-linéaires sur $E$, on définit leur produit\footnote{à ne pas confondre avec le produit tensoriel de catégories, définit au {\sffamily \S \ref{prodtens}}} $C\otimes_E D$ qui est un $A$-graphe sur $E$ en posant pour $(x,z)\in E^2$ : $$(C\otimes_ED)(x,z) := \bigoplus_{y\in E}C(x,y)\otimes_AD(y,z).$$

Ceci définit une structure monoïdale sur les graphes $A$-linéaires sur $E$ et permet de définir les {\em catégories $A$-linéaires d'objets $E$} (ou {\em $A$-algèbre sur $E$}) comme les monoïdes de cette catégorie. La loi monoïdale se décompose en une famille de morphismes $A$-bilinéaires :
$$
m_{xyz} : C(x,y)\otimes_AC(y,z) \longrightarrow C(x,z)
$$
et le neutre en éléments distingués $1_x$ dans chaque $C(x,x)$.
En particulier, chaque $C(x,x)$ est une $A$-algèbre associative unitaire et chaque $C(x,y)$ est un double module à gauche sous $C(x,x)$ et à droite sous $C(y,y)$, ces deux actions commutant.

Une {\em catégorie $A$-linéaire}\index{catégorie linéaire} (ou une {\em $A$-algèbre à plusieurs objets}\index{algèbre à plusieurs objets}) est une catégorie $A$-linéaire sur un ensemble d'objets imprécisé.
Une {\em catégorie linéaire} (ou une {\em algèbre à plusieurs objets}) est une catégorie $A$-linéaire pour un anneau commutatif $A$ et un ensemble d'objets imprécisés.
Une catégorie linéaire est {\em libre}\index{catégorie linéaire libre}\footnote{à ne pas confondre avec la liberté relativement à la loi monoïdale, qui, d'ailleurs, n'interviendra jamais dans ce travail}, {\em projective}\index{catégorie linéaire projective} ou de {\em type fini}\index{catégorie linéaire de type fini} si son graphe sous-jacent l'est.
Le {\em type}\index{type d'une catégorie linéaire} d'une catégorie linéaire est celui de son graphe sous-jacent.

\bigskip

Soit $(E,d)$ un type de graphe ; pour tout $A\in \kcom$, on définit $S_A$ comme l'ensemble des structures de catégories $A$-linéaire sur le $A$-graphe libre de type $(E,d)$. Si $A\to B\in \kcom$ on définit une application $S_A\to S_B$ en associant à $C\in\S_A$ la structure $C_B$\label{reststcat} de catégorie $B$-linéaire sur le $B$-graphe libre de type $(E,d)$ obtenue à partir de celle de $C$ en prolongeant par linéarité à $B$ le produit.

Le foncteur suivant classifie les structures de catégorie linéaire sur un graphe libre de type $(E,d)$.
\begin{eqnarray*}
S : \kaff^o &\longrightarrow & \ens\\
A&\longmapsto & S_A\\
A\to B &\longmapsto & S_A\to S_B
\end{eqnarray*}

\medskip

Si le graphe sous-jacent est libre, les structures de catégories s'expriment par des équations quadratiques explicites. Soit $C$ un $A$-graphe libre de type $(E,d)$, on note $e^{xy}_i$ les vecteurs des bases canoniques de chaque $C(x,y)$. Dans les bases canoniques, $m_{xyz}$ se caractérise par $d(x,y)d(y,z)d(x,z)$ éléments de $A$ :
$$
m_{xyz}(e^{xy}_i,e^{yz}_j) = \sum_{1\leq k \leq d(x,z)} (m_{xyz})^k_{ij}\ e^{xz}_k
$$
et les équations d'associativité s'écrivent :
$$
Ass_{xyzt}^{ijkp}:=\sum_{1\leq \ell \leq d(x,t)}
(m_{xzt})^{p}_{\ell k}(m_{xyz})^{\ell}_{ij} - (m_{xyt})^{p}_{i\ell}(m_{yzt})^{\ell}_{jk} = 0.
$$
En notant $1_x^i$ les coordonnées de l'identité dans la base $e_i^{xx}$ de $C(x,x)$, les équations le concernant sont :
\begin{eqnarray*}
Id_{xxy}^{jk}:=\sum_{1\leq i \leq d(x,x)} (m_{xxy})^{k}_{ij}1_x^i - \delta_{jk} = 0\\
Id_{yxx}^{jk}:=\sum_{1\leq i \leq d(x,x)} (m_{yxx})^{k}_{ji}1_x^i - \delta_{jk} = 0
\end{eqnarray*}
où $\delta_{ij}$ est symbole de Kronecker.
Ainsi les structures de catégorie $A$-linéaire sur $C$ sont en bijection avec les morphismes de $\kk$-algèbres :
$$
\kk[(m_{xyz})^{k}_{ij},1_x^i]/(Ass_{xyzt}^{ijkp}, Id_{xxy}^{jk},Id_{yxx}^{jk}) \longrightarrow A.
$$
Ce qui prouve le fait suivant.
\begin{prop}\label{structurecat}
Avec des abréviations évidentes, le schéma affine
$$
\kcat^{(d)}\index{$\kcat^{(d)}$} := Spec\left(\kk[m,1]/(A,I,I) \right)
$$
représente le foncteur $S$ classifiant les structures de catégorie linéaire sur les graphes libres de type $(E,d)$.
\end{prop}

\bigskip

Dans le cas où on se limite à des graphes ayant un seul objet, \ie à des algèbres associatives unitaires, on note $\kcat^{(n)}$ le schéma classifiant les structures d'algèbre associative unitaire sur un module libre de rang $n$.

Si on se limite en plus aux seules algèbres associatives unitaires qui soient commutative, il faut imposer en plus les équations :
$$
C^k_{ij} := m^k_{ij}-m^k_{ji} =0.
$$
On en déduit la proposition suivante. 

\begin{prop}\label{schstcom}
Avec les mêmes abréviations que précédemment, le schéma affine
$$
\com^{(n)}\index{$\com^{(n)}$} := Spec\left(\kk[m,1]/(A,I,I,C) \right)
$$
représente le foncteur classifiant les structures d'algèbre associative unitaire commutative sur un module libre de rang $n$. C'est un sous-schéma fermé de $\kcat^{(n)}$.
\end{prop}

\subsection{Foncteurs linéaires}

Un {\em morphisme strict de catégories $A$-linéaires}  (d'objets $E$) est un morphisme de monoïdes sur $E$ dont le morphisme sous-jacent de graphes soit strict.
Un {\em morphisme (non strict) de catégories $A$-linéaires} (ou {\em foncteur $A$-linéaire}\index{foncteur linéaire}) $f:C\to D$ est un morphisme de $A$-graphes compatible aux structures monoïdales au sens où tous les diagrammes suivants commutent :
$$\xymatrix{
C(x,y)\otimes_AC(y,z) \ar[d]_{f_{xy}\otimes f_{yz}} \ar[r]^-{m_{xyz}} & C(x,z) \ar[d]^{f_{xz}} \\
D(fx,fy)\otimes_AD(fy,fz) \ar[r]_-{m'_{fxfyfz}} & D(fx,fz)\\
}$$
($m$ est le produit de $C$, $m'$ celui de $D$, et $fx$ abrège $f_0(x)$).
Les {\em isomorphismes de catégories $A$-linéaires} sont les foncteurs qui sont des isomorphismes (non stricts) des graphes sous-jacents.

\medskip

Lorsque les catégories but et source sont libres de type fini, un foncteur s'exprime par des équations quadratiques. Reprenant les notations de ci-dessus, on note $m,1$ les constantes de structure de $C$ et $m',1'$ celles de $D$, un foncteur linéaire $f:C\to D$ est la donnée d'une application $f$ entre les objets de $C$ et $D$ et pour chaque paire $x,y$ d'objets de $C$ d'une application $f_{xy}:C(x,y)\to D(fx,fy)$ caractérisé dans les bases canoniques par des scalaires $(f_{xy})_j^{i}$. La fonctorialité de $f$ se traduit par les équations :
\begin{eqnarray*}
Fct_{xyz}^{ijk} & := & \sum_\ell (f_{xz})^i_\ell (m_{xyz})^\ell_{jk} - \sum_{\ell p} (m'_{fxfyfz})^i_{\ell p} (f_{xy})^\ell_j (f_{yz})^p_k = 0\\
FctId_x^{i} & := & \sum_\ell (f_{xx})^i_j1_x^j-1_{fx}^i=0
\end{eqnarray*}
qu'on abrège $F$ et $FI$. L'anneau des constantes de structures des foncteurs $A$-linéaires entre deux catégories de graphe libre est donc $\bigoplus_f A[f]/(F,FI)$ où la somme est prise sur les applications entre les objets.

\begin{prop}
\begin{enumerate}
\item Les foncteurs linéaires entre deux catégories $A$-linéaires libres de types finis $C$ et $D$ fixées sont classifiés par le schéma affine :
$$
Fct(C,D)_0\index{$Fct(C,D)_0$} := \coprod_f Spec \left( A[f]/(F,FI) \right).
$$
où la réunion est prise sur les applications $f$ entre les objets de $C$ et ceux de $D$. 
\item Les triplets formés de deux catégories linéaires libres de types finis et d'un foncteur entre elles sont classifiés par le schéma affine :
$$
CL_1\index{$CL_1$} := \coprod_{(E,d),(F,d')} \coprod_{f:E\to F} Spec(\kk[m,1,m',1',f]/(A,I,I,A',I',I',F,FI).
$$
où la première réunion est prise sur les couples de types de graphes finis et la seconde sur les applications $f$ entre les objets de $C$ et ceux de $D$. 
\end{enumerate}
\end{prop}

\begin{pr}
Le premier point est prouvé par le raisonnement précédent, le schéma est affine car la réunion se faisant sur des applications entre ensembles finis, elle est finie. Le second point s'en déduit en ne fixant pas les constantes de structures de $C$ et $D$.
\end{pr}

\bigskip

On note $A\z\mathcal{CAT}$\index{$A\z\CAT$} la catégorie des catégories $A$-linéaires et des foncteurs $A$-linéaires.
On y distingue la sous-catégorie des objets projectifs de type fini, qu'on note $A\z\cat^f$\index{$A\z\cat^f$}. 
Les catégories équivalentes à des catégories projectives de type fini forment la sous-catégorie $A\z\cat$\index{$A\z\cat^f$}.

\subsection{Isomorphismes d'objets}

Dans une catégorie linéaire $C$, deux objets $x,y$ sont dits {\em isomorphes}\index{objets isomorphes} s'il existe
deux éléments de $f\in C(x,y)$ et $g\in C(x,y)$ tels que $m_{xyx}(f,g)=1_x$ et $m_{yxy}(g,f)=1_y$.
Cela implique que $C(x,x)$, $C(y,y)$, $C(x,y)$ et $C(y,x)$ sont des $A$-modules isomorphes.
Dans le cas où $C$ est projective de type fini $d$, on en tire : $d(x,x)=d(x,y)=d(y,x)=d(y,y)$.
La relation d'isomorphie est une relation d'équivalence et on définit l'ensemble des {\em classes d'isomorphie d'objets de $C$} comme le quotient de ses objets par cette relation.
Si cet ensemble et de cardinal $n$, on dit que $C$ a {\em essentiellement $n$ objets}\index{nombre essentiel d'objets}.

\begin{lemme}\label{isoouvert}
Le lieu dans $Spec(A)$ d'inversibilité d'un morphisme d'une catégorie $A$-linéaire projective est un ouvert (éventuellement vide).
\end{lemme}
\begin{pr}
Soient $C$ une catégorie $A$-linéaire projective, $x$ et $y$ deux objets et $f\in C(x,y)$. $f$ est inversible s'il existe $g\in C(y,x)$ tel que $fg=1_y$ et $gf=1_x$. $f$ définit un morphisme de $A$-module projectifs $m_f:C(y,x)\to C(y,y) ; g\mapsto fg$ qui est un isomorphisme si et seulement si $f$ est inversible à droite.
Si ces modules n'ont pas le même rang, $m_f$ ne peut pas être un isomorphisme et l'ouvert cherché est vide. Si ces modules ont même rang, comme $m_f$ est de rang maximal exactement sur un ouvert $O_f$ de $Spec(A)$, c'est un isomorphisme sur cet ouvert. En particulier, sur $O_f$, l'image réciproque de $1_y$ donne $g$ un inverse à droite de $f$. Le raisonnement précédent appliqué à $g$ donne un ouvert $O'_f\subset O_f$ sur lequel $g$ est inversible à droite par un élément $h$. Par $gh=1_x$ et $fg=1_y$ on trouve $h=f$ et le lieu d'inversibilité de $f$ est $O'_f$. 
\end{pr}

\begin{cor}
Le lieu d'isomorphie de deux objets d'une catégorie $A$-linéaire projective est un ouvert (éventuellement vide) de $Spec(A)$.
\end{cor}
\begin{pr}
Si le lieu d'isomorphie n'est pas vide, il contient au moins un point de $Spec(A)$. 
Soit $a$ un tel point de corps résiduel $\kappa$, $x$ et $y$ y sont isomorphes, \ie il existe $f\in C(x,y)\otimes_A\kappa$ inversible. D'après le lemme précédent $f$ reste inversible sur un voisinage de $a$ et $x$ et $y$ sont donc isomorphes sur un ouvert.
\end{pr}

\begin{cor}\label{nbessobjets}
Le nombre essentiel d'objets d'une catégorie $A$-linéaire projective est une fonction semi-continue supérieurement sur $Spec(A)$. En particulier, le nombre essentiel d'objets est minimal sur un ouvert de $Spec(A)$.
\end{cor}
\begin{pr}
Soit $C$ une catégorie $A$-linéaire, on note $N$ le nombre essentiel d'objets vu comme fonction sur $Spec(A)$ à valeur dans $\mathbb{N}$. 
Soit $a$ un point de $Spec(A)$ d'après le corollaire précédent, il existe un voisinage de $a$ dans lequel deux objets isomorphes en $a$ restent isomorphes, ainsi $N$ ne varie pas au voisinage de $a$ et les ensembles $\{a : N(a)\leq n\}$ sont ouverts.
\end{pr}

\subsection{\'Equivalences de catégories}

Tout foncteur envoie deux objets isomorphes dans deux objets isomorphes et, en conséquence, définit une application entre les classes d'isomorphie d'objets des deux catégories. Un foncteur est dit {\em essentiellement surjectif}\index{foncteur essentiellement surjectif} si cette application est une surjection.

Soit $f:C\to D$ un foncteur $A$-linéaire, $f$ est dit {\em plein}\index{foncteur plein} (resp. {\em fidèle}\index{foncteur fidèle}) si tous les $f_{x,y}$ sont surjectifs (resp. injectifs). On dit que $f$ est {\em pleinement fidèle} si $f$ et plein et fidèle. 

\medskip

S'il n'existe aucun isomorphisme entre deux objets d'une catégorie, il n'en existera aucun entre leur image par un foncteur pleinement fidèle. Un tel foncteur induit donc une application injective entre les classes d'isomorphie d'objets. Dans le cas où les catégories sont projectives de type $d$ et $d'$, l'application $f_0$ entre les objets d'un foncteur pleinement fidèle $f$ induit l'égalité $d=d'\circ f_0$.
Un foncteur pleinement fidèle est aussi tel que tous les $f_{xy}$ soient des isomorphismes, si les catégories sont libres cela équivaut à ce que les $\det(f_{xy})$ soient tous non nuls. Ceci fournit la classification suivante.

\begin{prop}\label{schplfid}
\begin{enumerate}
\item Les foncteurs linéaires pleinement fidèles entre deux catégories $A$-linéaires $C$ et $D$ fixées, libres de types respectifs $d$ et $d'$, sont classifiés par le schéma affine :
$$
FPF(C,D)\index{$FPF(C,D)$} := \coprod_f Spec \left( A[f,\det(f)^{-1}]/(F,FI) \right).
$$
où la réunion est prise sur les applications $f$ entre les objets de $C$ et ceux de $D$ vérifiant $d=d'\circ f$.
En particulier ils forment un ouvert dans le schéma des foncteurs de $C$ vers $D$.
\item Les triplets formés de deux catégories linéaires et d'un foncteur pleinement fidèle entre elles sont classifiés par le schéma affine :
$$
CL_{3 \over 2}\index{$CL_{3 \over 2}$} := \coprod_{(E,d),(F,d')} \coprod_{f:E\to F,d=d'\circ f} Spec(\kk[m,1,m',1',f,\det(f)^{-1}]/(A,I,I,A',I',I',F,FI).
$$
où la première réunion est prise sur les couples de types de graphes finis (avec $E$ et $F$ fixés) et la seconde sur les applications $f:E\to F$ vérifiant $d=d'\circ f$. C'est un ouvert de $CL_1$.
\end{enumerate}
\end{prop}

\medskip

Une {\em équivalence de catégories}\index{equivalence de catégories@équivalence de catégories} entre deux catégories linéaires est un foncteur linéaire qui est pleinement fidèle et essentiellement surjectif. En particulier il induit une bijection entre les classes d'isomorphie d'objets. Tout isomorphisme est une équivalence.

\begin{prop}\label{schequiv}
\begin{enumerate}
\item Les équivalences entre deux catégories $A$-linéaires $C$ et $D$ fixées libres de types finis sont classifiés par un schéma $Eq(C,D)_0$\index{$Eq(C,D)_0$} qui est un ouvert de $FPF(C,D)$.
\item Les triplets formés de deux catégories linéaires et d'une équivalence entre elles sont classifiés par un schéma $EQ_1$\index{$EQ_1$} qui est un ouvert de $CL_{3 \over 2}$.
\end{enumerate}
\end{prop}
\begin{pr}
Pour $B$ une $A$-algèbre commutative, soit $f$ un $B$-point de $FPF$, \ie un foncteur pleinement fidèle $C_B\to D_B$ (cf. {\sf \S \ref{reststcat}} ou {\sf \S \ref{operations}} pour la définition de $C_B$), qui soit une équivalence une fois spécialisée en un point $b$ de $Spec(B)$ de corps résiduel $\kappa(b)$, on prouve qu'il est encore une équivalence sur un voisinage de $b$. 
Soit $f_0$ l'application entre les objets de $C_B$ et $D_B$ donnée par $f$, si $f_0$ est surjective alors $f$ est une équivalence sur tout $Spec(B)$. Sinon, soit $y$ un objet de $D_B$ non dans l'image de $f$, il existe $x$ un objet de $C_B$ et un élément $u$ de $C_B(f_0(x),y)$ tels qu'au point $b$, $u_b:f_0(x) \to y$ soit un isomorphisme. Par le lemme {\ref{isoouvert}} $x$ et $y$ sont en fait isomorphes sur un ouvert $U_y$ de $Spec(B)$. Le nombre d'objets de $D$ étant fini, l'intersection des $U_y$ est le lieu où $f$ est une équivalence.
\end{pr}

\medskip
Les applications surjectives entre ensembles ayant toujours des sections, les équivalences de catégories ont toujours des quasi-inverses. On en déduit le lemme suivant.

\begin{lemme}\label{quasiinv}
Si dans un triangle commutatif de foncteurs, deux sont des équivalences, le troisième aussi.
\end{lemme}

\subsection{Transformations naturelles}

On note $[1]_A$\index{$[1]_A$} la catégorie $A$-linéaire ayant deux objets $0$ et $1$ et tel que $[1]_A(0,0)=[1]_A(1,1)=A$, $[1]_A(1,0)=0$ et $[1]_A(0,1)=A$ vu comme $A$-bimodule, les compositions étant définies par l'identité de $A$.
On définit aussi $\Delta^1_A$\index{$\Delta^1_A$} le "groupoïde" $A$-linéaire par $\Delta^1_A(x,y)=[1]_A(x,y)$ sauf pour $\Delta^1_A(1,0)=A$ ; les produits sont tous donnés par l'identité de $A$. L'unité de $A$ de $\Delta^1_A(1,0)$ correspond à l'inverse de l'unité de $A$ dans $\Delta^1_A(0,1)$.

\medskip

Une {\em transformation naturelle}\index{transformation naturelle} $\alpha:f\to g$ entre deux foncteurs linéaires de $C$ dans $D$ est un foncteur
$$
\alpha : C\otimes_A [1]_A \longrightarrow D
$$
tel que ses restrictions aux deux objets 0 et 1 de $[1]_A$ soient $f$ et $g$ respectivement.
(Le produit tensoriel de catégories est redéfinit au {\S\ref{prodtens}}.)

Un {\em isomorphisme naturel}\index{isomorphisme naturel} $\alpha:f\to g$ entre deux foncteurs linéaires de $C$ dans $D$ est un foncteur
$$
\alpha : C\otimes_A \Delta^1_A \longrightarrow D
$$
tel que ses restrictions aux deux objets 0 et 1 de $\Delta_A$ soient $f$ et $g$.

Plus concrètement, une transformation naturelle est la donnée pour tout objet $x$ de $C$ de morphismes $\alpha_x\in D(fx,gx)$ tels que, pour tout $u\in C(x,y)$, les diagrammes suivants commutent :
$$\xymatrix{
f(x) \ar[r]^{\alpha_x}\ar[d]_{f(u)} & g(x)\ar[d]^{g(u)}\\
f(y) \ar[r]^{\alpha_y}& g(y).\\
}$$

\begin{lemme}\label{tnlin}
L'ensemble des transformations naturelles entre deux foncteurs $A$-linéaires est un $A$-module.
\end{lemme}
\begin{pr}
En utilisant les notations précédentes, il suffit de montrer que les $\{\alpha_x\}_x$ vérifiant la condition de commutation à $f$ et $g$ forment un sous-$A$-module de $\oplus_x D(fx,gx)$.
Soient deux transformations naturelles $\alpha,\beta : f \to g$ données par des morphismes $\alpha_x$ et $\beta_x$ de $D(fx,gx)$, et soit $a\in A$ ; la condition de commutativité s'écrit $g(u)\circ\alpha_x=\alpha_y\circ f(u)$ :
la $A$-linéarité des compositions de $C$ et $D$ donne $g(u)\circ(a.\alpha+\beta)=a.(g(u)\circ\alpha)+g(u)\circ\beta=a.(\alpha\circ f(u))+\beta\circ f(u)=(a.\alpha+\beta)\circ f(u)$.
\end{pr}

Les isomorphismes naturels se caractérisent, par le fait que tous les $\alpha_x$ soient des isomorphismes.
On a alors : 
$$
g(u)=\alpha_y \circ f(u) \circ \alpha_x^{-1},
$$
$\alpha$ et $f$ déterminent $g$.

\bigskip

Si les catégories buts et sources sont libres, les transformations naturelles entre foncteurs s'expriment par des équations linéaires. Soient $f$ et $g$ deux foncteurs entre $C$ et $D$ deux catégories libres. En coordonnées, une transformation naturelle $\alpha:f \to g$ est la donné de scalaires $\alpha_x^i$ soumis aux conditions :
$$
TN^i_{xy} := \sum_{j,k} (m_{fxgxgy})^i_{jk} (g_{xy})^j_\ell\alpha_x^k - (m_{fxfygy})^i_{jk}\alpha_y^j (f_{xy})^k_\ell = 0.
$$
L'anneau des constantes de structures des transformations naturelles entre deux foncteurs $A$-linéaires $f$ et $g$ fixés est abrégé $A[\alpha]/(TN)$. Pour obtenir celui des isomorphismes naturels, il est nécessaire que, pour tout $x$, $D(fx,fx)$, $D(fx,gx)$, $D(gx,fx)$ et $D(gx,gx)$ aient le même rang ; dans ce cas, l'opérateur $m_{\alpha_x}$ de multiplication par $\alpha_x$ est inversible si sa matrice dans les bases canoniques a un déterminant inversible, cela donne une localisation de l'anneau précédent.

\begin{prop}\label{schisoequiv}
Soient $C$ et $D$ deux catégories $A$-linéaires libres de types finis et $f$ et $g$ deux foncteurs linéaires de $C$ vers $D$. 
\begin{enumerate}
\item Les transformations naturelles de $f$ vers $g$ sont classifiés par le schéma affine :
$$
TN(f,g,C,D) := Spec(A[\alpha]/(TN).
$$
\item Les isomorphismes naturels de $f$ vers $g$ sont classifiés par le schéma affine :
$$
Iso(f,g,C,D) := Spec(A[\alpha,\det(m_\alpha)^{-1}]/(TN).
$$
\item Les triplets formés de deux foncteurs linéaires de $C$ vers $D$ et d'un isomorphisme entre ces deux foncteurs sont classifiés par le schéma affine :
$$
Fct(C,D)_1\index{$Fct(C,D)_1$} := \coprod_{f} Spec(A[f,\alpha,\det(m_\alpha)^{-1}]/(F,FI,TN).
$$
\item Les triplets formés de deux équivalences de $C$ vers $D$ et d'un isomorphisme entre ces deux équivalences forment un ouvert $Eq(C,D)_1$\index{$Eq(C,D)_1$} du schéma $Fct(C,D)_1$.
\item Les quintuplets formés de deux catégories linéaires $C$ et $D$ libres de types finis, de deux foncteurs de $C$ vers $D$ et d'une transformation naturelle entre ces deux foncteurs sont classifiés par le schéma affine :
$$
CL_2\index{$CL_2$} := \coprod_{(E,d),(F,d')} \coprod_{f:E\to F} Spec(\kk[m,1,m',1',f,\alpha] / (A,I,I,A',I',I',F,FI,TN).
$$
où $E$ et $F$ sont deux ensembles fixés.
\item Les quintuplets formés de deux catégories linéaires $C$ et $D$, de deux équivalences de $C$ vers $D$ et d'un isomorphisme entre ces deux équivalences sont classifiés par un schéma affine $EQ_2$\index{$CL_2$} qui est un ouvert de $CL_2$.
\end{enumerate}
\end{prop}
\begin{pr}
Les points {\em 1.} et {\em 2.} sont prouvés par le raisonnement précédent la proposition.
Le schéma classifiant les triplets de {\em 3.} est naturellement :
$$
\coprod_{f,g} Spec(A[f,\det(f)^{-1},g,\det(g)^{-1},\alpha,\det(m_\alpha)^{-1}]/(F_f,FI_f,F_g,FI_g,TN)\ ,
$$
mais, comme $\alpha$ et $f$ caractérisent entièrement $g$ par conjugaison, ce schéma se réduit à
$$
\coprod_{f,\alpha} Spec(A[f,\det(f)^{-1},\alpha,\det(m_\alpha)^{-1}]/(F_f,FI_f,TN).
$$
Le point {\em 4.} se prouve en remarquant que $Fct(C,D)_1$ se surjecte sur $Fct(C,D)_0$. La fibre de cette projection au-dessus de $Eq(C,D)_0$, qui est un ouvert comme image réciproque d'un tel, est formé des triplets d'une équivalence, d'un foncteur et d'un isomorphisme entre les deux ; or, tout foncteur isomorphe à une équivalence est une équivalence.
Cet ouvert est le schéma $Eq(C,D)_1$ cherché.
Le point {\em 5.} se déduit en ajoutant les équations des constantes de structures de $C$ et $D$.
Le point {\em 6.} se démontre avec des arguments similaires aux preuves des propositions {\ref{schplfid}} et {\ref{schequiv}}.
\end{pr}

On note $s$ la projection définie dans la preuve précédente, c'est l'application qui à un isomorphisme naturel associe sa source ; on note $b$ l'application but.

On pose $CL_0:=EQ_0:=\coprod_d\kcat^{(d)}$\index{$CL_0$}\index{$EQ_0$}

\medskip

On renvoie à \cite[ch. 6]{maclane} pour les définitions de 2-catégories et 2-groupoïdes.

\begin{prop}\label{gpdequiv}
\begin{enumerate}
\item Il existe une structure de groupoïde sur le graphe $s,b:Eq(C,D)_1 \rightrightarrows Eq(C,D)_0$.
\item Il existe une structure de 2-catégorie sur $CL_2 \rightrightarrows CL_1 \rightrightarrows CL_0$.
\item Il existe une structure de 2-groupoïde sur $EQ_2 \rightrightarrows EQ_1 \rightrightarrows EQ_0$.
\end{enumerate}
\end{prop}
\begin{pr}
La démonstration est évidente compte tenu de ce que classifient ces objets.
\end{pr}

\subsection{Catégories de foncteurs}
On pose les notations suivantes :%
\begin{itemize}
\item $\Hom_{A\z\CAT}(C,D)$\index{$\Hom_{A\z\CAT}(C,D)$} est l'ensemble des foncteurs $A$-linéaires de $C$ vers $D$ ;
\item $\Eq_{A\z\CAT}(C,D)$\index{$\uEq_{A\z\CAT}(C,D)$} est l'ensemble des équivalences de $C$ vers $D$ ;
\item $\uHom_{A\z\CAT}(C,D)$\index{$\uHom_{A\z\CAT}(C,D)$} est la catégorie des foncteurs de $C$ vers $D$ et de leur transformations naturelles, c'est une catégorie $A$-linéaire par le lemme \ref{tnlin} ;
\item $\uHom^{int}_{A\z\CAT}(C,D)$\index{$\uHom^{int}    _{A\z\CAT}(C,D)$} est le groupoïde des foncteurs de $C$ vers $D$ et de leurs isomorphismes naturels ;
\item $\uEq_{A\z\CAT}(C,D)$\index{$\uEq_{A\z\CAT}(C,D)$} est le groupoïde des équivalences de $C$ vers $D$ et de leurs isomorphismes naturels.
\end{itemize}

\section{Opérations sur les catégories}\label{operations}

\subsection{Catégories de catégories}

Soit $\UU$ un univers et $A\in \kcom_\UU$, on définit les catégories suivantes
\begin{itemize}
\item $A\z\CAT_\UU$\index{$A\z\CAT_\UU$} est la catégorie des catégories $A$-linéaires $\UU$-petites et des foncteurs $A$-linéaires ;
\item $A\z\cat^f_\UU$\index{$A\z\cat^f_\UU$} est la sous-catégorie pleine de $A\z\CAT_\UU$ des catégories dont le graphe est de type fini ;
\item $A\z\cat_\UU$\index{$A\z\cat_\UU$} est l'image essentielle de l'inclusion naturelle $A\z\cat^f_\UU\to A\z\CAT_\UU$ ;
\item $A\z\ASS_\UU$\index{$A\z\ASS_\UU$} est la catégorie des catégories $A$-linéaires ayant un unique objet, \ie des algèbres associatives unitaires, et des foncteurs $A$-linéaires ;
\item $A\z\ass_\UU$\index{$A\z\ass_\UU$} est la sous-catégorie pleine de $A\z\ASS_\UU$ dont le module sous-jacent à l'algèbre des endomorphismes est projectif de type fini sur $A$.
\end{itemize}

La plupart du temps, on oubliera de noter la référence à l'univers $\UU$, on ne la rappellera que lorsqu'on en aura besoin.

On dispose de plusieurs opérations sur les catégories $A\z\CAT$.

\subsection{Produit tensoriel et Hom interne}\label{prodtens}
Soient $C, D\in A\z\CAT$, on définit $C\otimes_A D$ comme la catégorie ayant pour objets les couples formés d'un objet de $C$ et d'un de $D$ et pour morphismes : 
$$
(C\otimes_AD)((c,d),(c',d')) =C(c,c')\otimes_AD(d,d').
$$
On a un isomorphisme canonique $\tau:C\otimes_AD\to D\otimes_AC$ tel que $\tau^2=id$ .

On définit également $\uHom_{A\z\CAT}(C,D)$\index{$\uHom_{A\z\CAT}(C,D)$} comme la catégorie ayant pour objets les foncteurs linéaires entre $C$ et $D$ et pour morphismes les transformations naturelles entre tels foncteurs, lesquelles forment naturellement des $A$-modules.

Ces deux foncteurs sont adjoints :
$$
\Hom_{A\z\CAT}(C\otimes_AD,E) \simeq \Hom_{A\z\CAT}(C,\uHom_{A\z\CAT}(D,E))
$$
et on a un isomorphisme canonique de catégories :
$$
\uHom_{A\z\CAT}(C\otimes_AD,E) \simeq \uHom_{A\z\CAT}(C,\uHom_{A\z\CAT}(D,E)).
$$

\begin{lemme}\label{homcat}
Si $C,C'$ et $D,D'$ sont des couples de catégories équivalentes alors les catégories $\uHom_{A\z\CAT}(C,D)$ et $\uHom_{A\z\CAT}(C',D')$ sont équivalentes et les équivalences sont envoyés sur des équivalences.
\end{lemme}
\begin{pr}
Immédiat, la démonstration repose sur l'existence de quasi-inverses pour les équivalences. La deuxième assertion résulte de ce que la composition des équivalences de catégorie en soit encore une.
\end{pr}

\begin{lemme}\label{oppo}
Soient $C$ et $D$ dans $A\z\CAT$, si $C^o$ désigne la catégorie opposée à $C$, on a les isomorphismes canoniques suivants :
\begin{eqnarray*}
(C\otimes_AD)^o & \simeq & C^o\otimes_AD^o \\
\uHom(C,D)^o & \simeq & \uHom(C^o,D^o)
\end{eqnarray*}
\end{lemme}
\begin{pr}
Les isomorphismes sont évidents au niveau des objets. Pour les morphismes on a d'une part $(C\otimes_AD)^o((c,d),(c',d')) = (C\otimes_AD)((c',d'),(c,d)) = C(c',c)\otimes_AD(d',d) = (C^o\otimes_AD^o)((c,d,),(c',d'))$ et d'autre que $\alpha\in\uHom(C,D)^o(f,g)$ est une famille $\alpha_x\in D(g(x),f(x))=D^o(f(x),g(x))$, \ie $\alpha$ définit canoniquement un élément de $\uHom(C^o,D^o)(f,g)$.
\end{pr}

\subsection{Changements de bases}\label{chgtbase}

Soient $u:A\to B$ un morphisme d'anneaux commutatifs, en considérant la catégorie $A$-linéaire $\B B$ ayant un objet et $B$ comme algèbre d'endomorphismes, on lui associe les trois foncteurs suivants.
\begin{itemize}
\item 'Oubli' $F: B\z\CAT \longrightarrow A\z\CAT$.
\item 'Catégorie $B$-linéaire libre' $(-)_u$\index{$(-)_u$} :
\begin{eqnarray*}
(-)_u : : A\z\CAT & \longrightarrow & B\z\CAT \\
C & \longmapsto & C_u := C\otimes_A\B B \\
f:C\to D & \longmapsto & f_u:C_u\to D_u
\end{eqnarray*}
où $f_u$ est l'extension $B$-linéaire de $f$.
\item 'Catégorie des $B$-modules' $(-)^u$\index{$(-)^u$} :
\begin{eqnarray*}
(-)^u : A\z\CAT & \longrightarrow & B\z\CAT \\
C & \longmapsto & C^u := \uHom_A(\B B,C) \\
f:C\to D & \longmapsto & f^u:C^u\to D^u
\end{eqnarray*}
où $f^u$ est défini sur les objets par les triangles commutatifs suivants en posant $f^u(x):=f\circ x$ :
$$\xymatrix{
\B B \ar[r]^x\ar[rd]_{f\circ x} & C\ar[d]^f \\
& D
}$$
et sur les morphismes $\alpha:x\to x'$ par $f^u(\alpha) := f(\alpha)$.
\end{itemize}

\begin{lemme}\label{adjcat}
On a les adjonctions suivantes entre les foncteurs précédents :
\begin{eqnarray*}
\Hom_{A\z\CAT}(C,FD) & \simeq & \Hom_{B\z\CAT}(C_u,D) \\
\Hom_{A\z\CAT}(FC,D) & \simeq & \Hom_{B\z\CAT}(C^u,D)
\end{eqnarray*}
qui donnent des isomorphismes canoniques de catégories :
\begin{eqnarray*}
\uHom_{A\z\CAT}(C,FD) & \simeq & \uHom_{B\z\CAT}(C_u,D) \\
\uHom_{A\z\CAT}(FC,D) & \simeq & \uHom_{B\z\CAT}(C^u,D)
\end{eqnarray*}
\end{lemme}
\begin{preuve}
\end{preuve}

\medskip

\begin{defi}
Les foncteurs $(-)_u$ et $(-)^u$ seront qualifiés de foncteurs de {\em restriction}\index{restriction} ou de {\em changement de base}\index{changement de base} car il seront utilisés dans la construction de certains préfaisceaux. On a les propriétés suivantes de ces foncteurs.
\end{defi}

\begin{lemme}\label{fonctofaible}
Soient $A\overset{u}{\to}B\overset{v}{\to}C$ deux morphismes d'anneaux. On a les isomorphismes canoniques suivants :
\begin{enumerate}
\item $((-)_u)_v \overset{\sim}{\longrightarrow} (-)_{vu}$\ ;
\item $((-)^u)^v \overset{\sim}{\longrightarrow} (-)^{vu}$.
\end{enumerate}
Pour $A\overset{u}{\to}B\overset{v}{\to}C\overset{w}{\to}D\in\kcom$, ces isomorphismes vérifient la condition de cohérence {\sf \ref{cdtpreff}} de la définition {\sf \ref{prefaible}}.
\end{lemme}
\begin{pr}
Soit $D$ une catégorie $A$-linéaire. L'isomorphisme $((-)_u)_v \overset{\sim}{\to} (-)_{vu}$ provient immédiatement des isomorphismes canoniques de modules $D(x,y)\otimes_AB\otimes_BC\overset{\sim}{\to}D(x,y)\otimes_AC$. La condition de cohérence résulte de celle du produit tensoriel d'anneau.

Pour $((-)^u)^v \overset{\sim}{\to} (-)^{vu}$, on explicite les objets et les morphismes de $D^u$.

Un objet de $D^u$ est la donnée d'un objet de $D$ et d'un morphisme $\chi:B\to D(x,x)$ de $A$-algèbres et un morphisme de $(x,\chi)$ vers $(y,\upsilon)$ est un élément de $D(x,y)$ commutant avec la double action de $B$ sur $D(x,y)$ donnée par $\chi$ et $\gamma$ ; on note $D^u(x,y)$ l'ensemble de ces éléments. Comme $B$ est commutatif, $D^u(x,x)$ est une $B$-algèbre pour le morphisme induit par $\chi$.

Un objet de $(D^u)^v$ est la donnée d'un objet de $D$ et d'un morphisme $\xi:C\to D^u(x,x)$ de $B$-algèbres ; un tel objet est caractérisé par le diagramme commutatif :
$$\xymatrix{
A\ar[r]^u\ar[rd] & B\ar[r]^v\ar[rd]\ar[d]^\chi & C\ar[d]^\xi\\
& D(x,x) & D^u(x,x)\ar[l]
}$$
qui est équivalent à la simple donnée du morphisme de $A$-algèbres $C\to D(x,x)$ induit par $\xi$ (la factorisation par $D^u(x,x)$ se déduisant de ce que $v$ envoie $B$ dans le centre de $C$). Ainsi les objets de $(D^u)^v$ sont en bijection avec ceux de $D^{vu}$.

Un morphisme de $(x,\chi,\xi)$ vers $(y,\upsilon,\gamma)$ est un élément de $D(x,y)$ commutant aux doubles actions de $B$ et $C$, mais l'action de $B$ se factorisant par celle de $C$, l'ensemble de ces morphismes est simplement $D^{vu}(x,y)$.
\end{pr}

\medskip

\begin{lemme}\label{monoidal}
Soit $u:A\to B$ un morphisme d'anneaux.
\begin{enumerate}
\item Le foncteur $(-)_u$ vérifie $(C\otimes_AD)_u \simeq C_u\otimes_BD_u$ (\ie sont monoïdaux).
\item Le foncteur $(-)^u$ vérifie $\uHom_{A\z\CAT}(C,D)^u \simeq \uHom_{B\z\CAT}(C_u,D^u)$.
\end{enumerate}
Ces deux équivalences sont des isomorphismes canoniques.
\end{lemme}
\begin{pr}
Le caractère monoïdal de $(-)_u$ résulte des isomorphisme canoniques $(M\otimes_AN)\otimes_A B \simeq (M\otimes_AB)\otimes_B(N\otimes_AB)$ pour deux $A$-modules $M$ et $N$.
Pour la deuxième assertion on a 
\begin{eqnarray*}
\uHom_{A\z\CAT}(C,D)^u &=& \uHom_{A\z\CAT}(\B B,\uHom_{A\z\CAT}(C,D)) \\
	&=& \uHom_{A\z\CAT}(\B B\otimes_AC,D)) \\
	&=& \uHom_{A\z\CAT}(C,D^u)) \\
	&=& \uHom_{B\z\CAT}(C_u,D^u))
\end{eqnarray*}
(où on utilise le symbole d'égalité pour des isomorphismes canoniques).
\end{pr}

\begin{lemme}\label{compatfibcofib}
Soit $u:A\to B$ un morphisme d'anneaux.
\begin{enumerate}
\item Les foncteurs $(-)_u$ et $(-)^u$ préservent les isomorphismes de catégories.
\item Le foncteur $(-)_u$ préserve les équivalences et les foncteurs induisant une inclusion sur les objets.
\item Le foncteur $(-)^u$ préserve les équivalences et les foncteurs surjectifs sur leur image essentielle.
\end{enumerate}
\end{lemme}
\begin{pr}
{\em 1.~}
Les ensembles d'objets ne changeant pas par restriction, il suffit de vérifier les isomorphismes des modules de morphismes.
Soit $f:C\to D$ un isomorphisme $A$-linéaire et $x$ et $y$ deux objets de $C$, $(f_u)_{xy}:C(x,y)\otimes_AB\to D(fx,fy)\otimes_AB$ est définit comme le prolongement $B$-linéaire de $f_{xy}:C(x,y)\to D(fx,fy)$ et reste donc un isomorphisme si $f_{xy}$ en était un.

{\em 2.~}
On garde les mêmes notations, mais cette fois $f$ est une équivalence de catégories.
Les identités étant conservées par restriction, il en est de même des inverse et un isomorphisme entre deux objets reste un isomorphisme entre les restrictions.
Si $f$ est essentiellement surjectif, soit $z\in D$ et $\alpha$ un isomorphisme $z\to f(x)$, $(f_u)_{xy}(\alpha)$ reste un isomorphisme et $z$ dans l'image essentielle. Comme les objets de $D$ et de $D_u$ sont les mêmes $f_u$ reste essentiellement surjectif.
Si $f$ est pleinement fidèle $f_{xy}$ est un isomorphisme de $A$-module et donc $(f_u)_{xy}$ est encore une bijection, et $f_u$ reste pleinement fidèle. 
Enfin, comme $C_u$ et $D_u$ ont respectivement les mêmes ensembles d'objets que $C$ et $D$ et comme $f_u$ est définit avec la même application que $f$ entre les objets, $f_u$ reste injectif sur les objets si $f$ l'était.

{\em 3.~}
$f^u$ préserve les équivalences d'après le lemme {\ref{homcat}}. La surjectivité sur l'image essentielle découle du fait suivant : soit $z\in D^u$ et $f^u(x)\in D^u$ pour $x\in C^u$, chacun de ces objets détermine un objet de $D$ (image de l'objet unique de $\B B$) et un isomorphisme $z\to f^u(x)$ est équivalent à un isomorphisme dans $D$ entre ces deux objets.
\end{pr}

\begin{lemme}\label{cdtfinitude}
Soit $u:A\to B$ un morphisme d'anneaux. Le foncteur $(-)_u$ conserve les caractères suivants des catégories :
\begin{itemize}
\item le nombre d'objets (en particulier sa finitude) ;
\item liberté des modules de morphismes ;
\item type fini ;
\item projectivité ;
\end{itemize}
et le nombre essentiel d'objets ne peut que diminuer par composition par $(-)_u$.
\end{lemme}
\begin{pr}
Le premier point vient de ce que, pour une catégorie $C$, l'ensemble d'objets de $C_u$ est pris le même que celui de $C$ ; les trois points suivants de ce que si $M$ est un $A$-module libre (resp. de type fini, resp. projectif), $M\otimes_AB$ est un $B$-module libre (resp. de type fini, resp. projectif). Le dernier point vient de ce que deux objets isomorphes restent isomorphes par restriction : si $a$ est un isomorphisme de $C$, d'inverse $b$, les restrictions des flèches $a$ et $b$ restent inverse l'une de l'autre.
\end{pr}

\section{Modules}

\subsection{Anneau matriciel d'une catégorie}\label{matrice}

À une catégorie $A$-linéaire $C$ ayant un nombre fini d'objets, on associe son {\em anneau matriciel} :
$$[C]\index{$[C]$} :=\bigoplus_{x,y} C(x,y)\ ;$$
c'est une $A$-algèbre pour le produit :
$$
m : [C]\otimes_A [C] = \bigoplus_{x,y,z,t} C(x,y)\otimes_A C(z,t)\longrightarrow \bigoplus_{x,y} C(x,y) =[C]
$$
définit par $m_{xyz}$ sur la composante indexée par $(x,y,y,z)$ et par 0 sur les composantes non de ce type.
L'unité est donnée par : $$1_{[C]}= \oplus_x 1_x\in \bigoplus_x C(x,x).$$

Notons $\{x_1,\dots,x_n\}$ les objets de $C$. Si on représente $[C]$ par le tableau carré indexé par les couples d'objets de $C$
$$[C]:=\left(\begin{array}{cccc}
C(x_1,x_1) & C(x_1,x_2) & \dots & C(x_1,x_n) \\
C(x_2,x_1) & C(x_2,x_2) & \dots & C(x_2,x_n) \\
\vdots & & \ddots & \vdots \\
C(x_n,x_1) & C(x_n,x_2) & \dots & C(x_n,x_n) \\
\end{array}\right),$$
le produit est exactement un produit matriciel (d'où le nom de $[C]$).

\medskip

Les objets de $C$ fournissent par les éléments $id_x$ une suite complète d'idempotents orthogonaux de $[C]$ (on a  $id_x^2=id_x$, $id_xid_y=0$ et $\sum_x i_dx = 1_{[C]}$) et réciproquement, un anneau $D$ muni d'une telle suite peut s'écrire comme l'anneau matriciel d'une catégorie. Si $p_1,\dots,p_n$ sont des idempotents orthogonaux on définit une catégorie $C$ ayant $n$ objets $x_i$ telle que $C(x_i,x_j)=p_iDp_j$ et où la composition est donnée par les morphismes naturels $p_iDp_j\otimes p_jDp_k\to p_iDp_k$.

\bigskip

Si les objets $C$ et $[C]$ se déterminent l'un l'autre, la correspondance entre les deux n'est pas fonctorielle. (Par exemple, un foncteur définit un morphisme entre les anneaux matriciels seulement s'il est injectif sur les objets.)

\bigskip

\begin{lemme}\label{propannmat}
On a des isomorphismes 
\begin{eqnarray*}
[C^o] & \simeq & [C]^o\\
{[}C\otimes_AD{]} & \simeq & [C]\otimes_A[D].
\end{eqnarray*}
\end{lemme}
\begin{pr}
Le premier point est évident compte tenu que $[C^o]$ se modélise par le tableau transposé de celui de $[C]$.
Le second se démontre en remarquant que la définition du produit tensoriel de deux catégories correspond au produit matriciel de leur représentation en tableau.
\end{pr}

\subsection{Modules}
Soit $C\in A\z\CAT$, un {\em module à gauche (resp. à droite) sur $C$} est un foncteur $C\to A\Mod$ (resp. $C^o\to A\Mod$).
On note
$$
\overline{C}\index{$\overline{C}$}:=C\Mod\index{$C\Mod$}:=\uHom_{A\z\CAT}(C,A\Mod)
$$ la catégorie des $C$-modules à gauche.
Si $C=\B B$ pour une $A$-algèbre $B$, on note $B\Mod$ la catégorie $\B B\Mod$.

\begin{lemme}\label{surjassab}
Il existe une équivalence canonique de catégories :
$$
C\Mod \simeq [C]\Mod.
$$
\end{lemme}
\begin{pr}
Explicitement un $C$-module à gauche est la donnée pour tout objet $x$ de $C$ d'un $A$-module $M(x)$ et de morphismes
$$
C(x,y)\otimes_A M(y) \overset{m_{xy}}{\longrightarrow} M(x)
$$
faisant commuter les diagrammes évidents.
On note $[M]=\bigoplus_x M(x)$, le $[C]$-module à gauche tiré de $M$ ; on peut représenter $[M]$ comme une matrice ligne $\left(M(x_1), \dots , M(x_n) \right)$ la structure de $[C]$-module est donné par un produit matriciel à gauche par $[C]$.
Soient $M$ et $N$ deux $C$-modules, un morphisme de $C$-modules $\alpha : M\to N$
est une transformation naturelle de foncteurs $C\to A\Mod$, \ie la donnée pour tout objet $x$ de $C$ d'$\alpha_x:M(x)\to N(x)$ commutant aux actions des $C(x,y)$. Alors
$[\alpha] := \bigoplus_x \alpha_x$
définit un morphisme de $[C]$-modules $[M]\to [N]$.
Ceci établit un foncteur (clairement fidèle) $C\Mod\longrightarrow [C]\Mod$.

On a également un foncteur en sens inverse. Pour $M$ un $[C]$-module, la décomposition $1_{[C]}=\oplus_{x\in C} 1_x$ permet de définir $M(x):=1_x.M$ tels que $M=\oplus_x M(x)$ et le fait que le produit soit matriciel assure que les produits de $C(x,y)\otimes_A M(y)$ soient dans $M(x)$.
Si $f: M\to N$ est un morphisme de $[C]$-modules, la linéarité de $f$ assure que
$$
f(1_x.m) = 1_x.f(m)
$$
et permet de définir $f_x:=1_xf1_x(=f1_x=1_xf) :M(x)\to N(x)$ et un morphisme de $C$-modules.
Tout ceci définit un foncteur $[C]\Mod\longrightarrow C$, fidèle car $f=\oplus f_x$, dont il est aisé de vérifier qu'il inverse le précédent.
\end{pr}

\medskip

\begin{lemme}\label{restricmodule}
Soient $u:A\to B\in\kcom$ et $C\in A\z\CAT$, on a une équivalence :
$$
\uHom(C,A\Mod)^u = \uHom(C_u,B\Mod).
$$
\end{lemme}
\begin{pr}
On remarque d'abord qu'il existe une équivalence de catégories $A\Mod^u\simeq B\Mod$.
Pour la voir on explicite deux foncteurs quasi-inverses l'un de l'autre. Les éléments de $A\Mod^u$ sont des couples $(M,\beta)$ où $M\in A\Mod$ et où $\beta:B\to End_A(M)$ est un morphisme de $A$-algèbres.
Un tel couple définit une structure de $B$-module sur $M$ et donc un certain objet de $B\Mod$ et les morphismes dans $A\Mod^u$ correspondent exactement aux morphismes de $B$-modules ; ceci fournit un premier foncteur. Réciproquement, soit $M\in B\Mod$ dont on note $u_*M\in A\Mod$ son image par le foncteur d'oubli ; la structure de $B$-module de $M$ fournit un morphisme $\beta:B\to End_A(u_*M)$ le second foncteur est définit en associant à $M$ le couple $(u_*M,\beta)$. Ces foncteurs sont quasi-inverses l'un de l'autre.

Pour avoir le lemme on utilise le lemme {\ref{monoidal}} avec $D=A\Mod$ et l'équivalence ci-dessus.
\end{pr}

\bigskip



Pour $C\in A\z\CAT$, le 1-préchamp $C\zMMod$ des $C$-modules, projectifs de type fini sur la base est définit par
\begin{eqnarray*}
C\zMMod^{proj} : (\aff_A)^o & \longrightarrow & \gpd\\
u:A\to B &\longmapsto & (C_u\Mod^{proj})^{int}
\end{eqnarray*}
où $(C_u\Mod^{proj})^{int}$ est le groupoïde intérieur de la catégorie des $C_u$-modules projectifs comme $B$-modules.

\begin{prop}
Le préchamp $C\zMMod^{proj}$ est un 1-champ 1-géométrique.
\end{prop}
\begin{pr}
Le fait que $C\zMMod^{proj}$ soit un champ résulte de la descente des modules et des structures d'algèbres \cite{sga1}.

On considère $A\proj$, le champ des modules projectifs sur $A$, \ie la restriction de $\bgln$ à $\aff_A$ ; d'après le lemme {\ref{projproj}} on a des morphismes
$$
s : C\zMMod^{proj} \longrightarrow  \uvect
$$
où $\uvect\simeq \coprod_n\bgln$ est le champ des modules projectifs de types finis, dont on va montrer la représentabilité.

Soit $u:A\to B$ et $M$ un $B$-module projectif, vu comme un morphisme $M:Spec(B)\to \uvect$, on considère la fibre $S(C_u,M)$ de $s$ le long de $M$
$$\xymatrix{
S(C_u,M)\ar[d]\ar[r] & C\zMMod^{proj}\ar[d]^s\\
Spec(B) \ar[r] & \uvect
}$$
c'est le faisceau classifiant les structures de $C_u$-module sur $M$.
Les morphismes suivants entre $B$-modules projectifs forment, une fois faisceautisés, un morphisme de fibrés vectoriels sur $Spec(B)$ dont l'égalisateur est $S(C_u,M)$ :
\begin{eqnarray*}
\Hom_B(C_B\otimes_BM,M) & \longrightarrow & \Hom_B(C_B\otimes_BC_B\otimes_BM,M)\\
f & \longmapsto & (f : a\otimes b\otimes m \mapsto f(a\otimes f(b\otimes m)) \\
f & \longmapsto & (f : a\otimes b\otimes m \mapsto f(ab\otimes m)
\end{eqnarray*}
où $\Hom_B$ désigne le Hom interne de la catégorie des $B$-modules.
Ceci prouve que $S(C_u,M)$ est représentable, ainsi que le morphisme $s$.
Comme $\uvect$ est 1-géométrique, on en déduit que $C\zMMod^{proj}$ est 1-géométrique.

\end{pr}

\begin{lemme}\label{defmodproj}
Soient $A\in\kcom$, $B\to A\in\kcom$ une extension infinitésimale au premier ordre, $C\in A\z\ass$, $M$ un $C$-module projectif de rang fini $r$, et $D\in B\z\ass$ tel que $D\otimes_BA\simeq C$ ; alors il existe un $D$-module projectif $N$ de rang $r$, tel que $N\otimes_DC$ soit isomorphe à $M$.
\end{lemme}
\begin{pr}
Soient les $B$-modules $I:=\ker (B\to A)$ et $J:=\ker(D\to D\otimes_BA\simeq C)$, on a $J=ID$ et, comme $I^2=0$, $J^2=0$. On déduit de $D\to C$ un morphisme entre les algèbres de matrices $n\times n$ : $M_n(D)\to M_n(C)$ dont le noyau $K=M_n(I)$ vérifie $K^2=0$.
$M$ est isomorphe à un facteur direct dans un certain $C^n$ et un tel facteur correspond à un projecteur $p\in M_n(C)$. Un tel $p$ se relève toujours en $p'\in M_n(D)$ qui n'est plus forcément un projecteur, mais 
on montre qu'il existe toujours $k\in K$, commutant avec $p'$, tel que $p'+k$ soit un projecteur.
$(p'+k)^2=p'+k \iff k(1-2p') = p'^2-p'$ et il faut montrer que $1-2p'$ est inversible, or $(1-2p')^2=1+4(p'^2-p')$ est inversible d'inverse $1-4(p'^2-p')$ car $p'^2-p'\in K$ ; l'inverse de $1-2p'$ est donc $(1-2p')(1-4(p'^2-p'))$ et on peut poser $k=(p'^2-p')(1-2p')(1-4(p'^2-p'))$.
\end{pr}

\subsection{Bimodules}
\label{modulesinversibles}

Soient $A\in\kcom$ et $C$ et $D$ des catégories $A$-linéaires, un {\em $C\z D$-bimodule}\index{bimodule} est un $(C^o\otimes_AD)$-module.

Pour des catégories ayant un nombre fini d'objets, les lemmes {\ref{propannmat}} et{\ref{surjassab}} caractérisent la catégorie des $C\z D$-bimodules comme étant $[C^o\otimes_AD]\Mod \simeq ([C]^o\otimes_A[D])\Mod$, et on peut se limiter au cas où $C$ et $D$ sont des $A$-algèbres associatives.

Les unités de $C$ et $D$ permettent de définir des morphismes de $A$-algèbres $C^o\to C^o\otimes_AD$ et $D\to C\otimes_AD$ et un $C\z D$-bimodule peut toujours se voir comme un $C^o$-module ou un $D$-module.
$C$ avec son action à droite et gauche par multiplication définit un $C\z C$-bimodule.

Un $C\z D$-bimodule est dit {\em biprojectif}\index{bimodule biprojectif} s'il est projectif comme $C^o$-module et comme $D$-module. Le $C\z C$-bimodule associé à $C$ est biprojectif.

\medskip

Si $C$ et $D$ sont dans $A\z\ass$ de modules sous-jacent libres de rangs $c$ et $d$ de bases $e_i$ et $e'_i$
et si $M=A^m$, de base $x_i$, une structure $\mu:C\otimes_AM\otimes_AD\to M$ de $C\z D$-bimodule sur $M$ est la donnée de $cdm^2$ éléments de $A$, $\mu_{ijk}^\ell$, tels que $\mu(e_i,x_j,e'_k)=\mu_{ijk}^\ell x_\ell$ ; en notant $m$ et $m'$ les produits respectifs de $C$ et $D$, ces éléments vérifient les équations suivantes :
\begin{eqnarray*}
Bim_{ijkpq}^\ell := \mu_{irq}^\ell \mu_{jkp}^r - \mu_{rks}^\ell m_{ij}^r(m')_{pq}^s = 0
\end{eqnarray*}

Si $M$ et $N$ sont deux structures de $C\z D$-bimodule sur $A^m$ et $A^n$ respectivement, dont on note les bases $x_i$ et $y_i$, un morphisme de $C\z D$-bimodules $f:M\to N$ est , en coordonnées, la donnée de $nm$ éléments de $A$, $f_i^j$ tels que $f(x_i) = f_i^jy_j$, ces nombres vérifient les équations suivantes
\begin{eqnarray*}
Morbim_{ijk}^\ell := \mu_{ijk}^pf_p^\ell -f_j^p\mu_{ipk}^\ell =0.
\end{eqnarray*}

\begin{prop}\label{presbim}
Soient $C,D\in A\z\ASS$ de $A$-modules sous-jacent libre de type fini. 
\begin{enumerate}
\item Avec des abréviations évidentes, les structures de $C\z D$-bimodules sur un $A$-module libre de rang $m$ sont classifiées par le schéma affine :
$$
Bim(C,D,m) := Spec(A[\mu]/(Bim)).
$$
\item Les triplets formés d'une structure $M$ de $C\z D$-bimodule sur $A^m$, d'une structure $M'$ de $C\z D$-bimodule sur $A^n$ et d'un morphisme $f:M\to M'$ de $C\z D$-bimodule sont classifiées par le schéma affine :
$$
Morbim(C,D,m,n) := Spec(A[\mu,\mu',f]/(Bim,Bim',Morbim)).
$$
\item Si on se limite aux seuls morphismes qui sont des isomorphismes on trouve un ouvert du schéma $Morbim(C,D,m,m)$ :
$$
Isobim(C,D,m) := Spec(A[\mu,\mu',f,\det(f)^{-1}]/(Bim,Bim',Morbim)).
$$
En outre, si $G\ell_m$ désigne le schéma automorphismes d'un module libre de rang $m$, on a :
$$
Isobim(C,D,m) = Spec(A[\mu,f,\det(f)^{-1}]/(Bim)) = Bim(C,D,m)\times G\ell_m.
$$
\end{enumerate}
\end{prop}
\begin{pr}
Les deux premiers résultats sont évidents. Pour le troisième il suffit de remarquer qu'un isomorphisme de $C\z D$-bimodule est en particulier un isomorphisme de $A$-modules. Enfin la dernière égalité résulte de ce que la structure $M$ et l'isomorphisme $f$ caractérisent la structure $M'$.
\end{pr}

\medskip

Si $M$ est un $C\z D$-bimodule et $N$ un $D\z C$-bimodule on leur associe $M\otimes_DN$ et $N\otimes_CM$ qui sont respectivement un $C\z C$-bimodule et un $D\z D$-bimodule.
En particulier, le $C\z C$-bimodule associé à $C$ agit comme un élément neutre à gauche sur les $C\z D$-bimodules ; de même le $D\z D$-bimodule associé à $D$ agit comme un élément neutre à droite.

Un $C\z D$-bimodule est dit {\em inversible}\index{bimodule inversible} s'il existe un $D\z C$-bimodule $N$ et des isomorphismes de bimodules $M\otimes_DN\simeq C$ et $N\otimes_CM\simeq D$.

%

\begin{lemme}\label{biminvouv}
Soient $C,D\in A\z\ass$ de $A$-modules sous-jacent libre de type fini. 
Le schéma $Inv(C,D,m)$ classifiant les structures de $C\z D$-bimodules inversibles sur un $A$-module libre de rang $m$, est un ouvert du schéma affine $Bim(C,D,m)$.
\end{lemme}
\begin{pr}
Le lemme {\ref{eqcateqmor}} établit que les bimodules inversibles sont en particulier projectifs de type fini sur chacun des anneaux de base (biprojectifs).
Dans une premier étape on montre donc que les structures de bimodules biprojectifs forment un ouvert de $Bim(C,D,m)$.
$C$ et $D$ jouant des rôles symétriques, il suffit de vérifier que la condition d'être un module projectif sur $C$ est ouverte dans le classifiant des structures de $C$-modules.

Soit donc $M$ un $C$-module et $p:C^n\to M$ une présentation libre (qui existe car $M$ et $C$ sont de type fini sur $A$), $M$ est projectif sur $C$ ssi $p$ admet une section. Soit $K$ le noyau de $p$, la section existe ssi $Ext^1_C(M,K)=0$. On va montrer que cette condition est vérifiée sur un ouvert de $Spec(A)$.
On considère le faisceau $E$ sur $Spec(A)$ définit pour $A\to B$ par $Ext^1_{C\otimes_AB}(M\otimes_AB,K\otimes_AB)$.

On complète $p$ en une résolution libre : $R^*:= \dots C^\ell\to C^m\to C^n\to M$. Pour tout morphisme $u:A\to B\in\kcom$, $R^*\otimes_AB$ reste une résolution libre de $M\otimes_AB$.

Les modules $Ext^1_{C\otimes_AB}(M\otimes_AB,K\otimes_AB)$ se calculent alors comme les $H^1$ des complexes $\uHom_{C\otimes_AB}(R^*\otimes_AB,K\otimes_AB) = (K^n\to K^m\to K\ell\to \dots)\otimes_AB$. Les différentielles de ce complexe étant $C\otimes_AB$- et donc $B$-linéaires, on a 
$H^1(R^*\otimes_AB) = H^1(R^*)\otimes_AB$ et le faisceau $E$ est donc cohérent. En particulier il est nul  exactement sur une localisation de $B$.

\medskip

On note $BimBip(C,D,m)$ l'ouvert de $Bim(C,D,m)$ ainsi obtenu.
On prouve maintenant que la condition d'inversibilité est ouverte dans $BimBip(C,D,m)$.

Soit $M$ un $A$-point de $BimBip(C,D,m)$. $\Hom_C(M,C)$ est naturellement un $D\z C$-bimodule, on étudie les flèches d'évaluations
\begin{eqnarray*}
\alpha: M\otimes_{D}\Hom_{C}(M,C)\longrightarrow C\quad \textrm{et}\\
\beta: \Hom_{C}(M,C)\otimes_{C}M\longrightarrow D
\end{eqnarray*}
La caractérisation des bimodules inversibles du lemme {\ref{eqcateqmor}} assure que $M$ est inversible si et seulement si ces flèches sont des isomorphismes.

\medskip

Sous l'hypothèse que $M$ est biprojectif on montre que $M\otimes_D\Hom_C(M,C)$ et $\Hom_C(M,C)\otimes_CM$ sont projectifs de type fini sur $A$.

Comme $M$ est $D$-projectif il existe un $D$-module $M'$ et un entier $k$ tels que $M\oplus M'\simeq D^k$ dont on tire que $(M\otimes_D\Hom_C(M,C))\oplus (M'\otimes_D\Hom_C(M,C))\simeq \Hom_C(M,C)^k$ comme $C$-modules. Ce dernier module est projectif sur $C$, en effet comme $M$ est projectif de type fini sur $C$ il existe un $C$-module $M'$ et un entier $n$ tels que $M\oplus M'\simeq C^n$ d'où  $\Hom_C(M,C)\oplus \Hom_C(M',C) \simeq C^n$. Ensuite, comme $C$ est projectif de type fini comme $A$-module, tout $C$-module projectif de type fini l'est aussi sur $A$ (si $N$ est un tel module $N\oplus N'\simeq C^\ell$ et $C\oplus C'\simeq A^p$ donnent $N\oplus N'\oplus C'^\ell\simeq A^{\ell p}$).
Finalement $M\otimes_D\Hom_C(M,C)$ est facteur direct d'un $A$-module libre, donc il est projectif sur $A$. Le raisonnement est analogue pour $\Hom_C(M,C)\otimes_CM$.

Le lieu de rang maximal d'une application linéaire entre module projectifs de type finis cohérents est un ouvert. 
On a donc l'alternative suivante : soit les modules $M\otimes_{D_u}\Hom_{C_u}(M,C_u)$ et $\Hom_{C_u}(M,C_u)\otimes_{C_u}M$ ont les mêmes rangs que $C$ et $D$ respectivement et l'ouvert de $Spec(A)$ cherché est l'intersection des ouverts de rang maximal de $\alpha$ et $\beta$ ; soit les modules n'ont pas les mêmes rangs et l'ouvert cherché est vide.
\end{pr}

On définit le champ $\underline{\mathcal{I}nv}(C,D)$ classifiant, à isomorphisme près, les $C\z D$-bimodules inversibles comme la champ associé au préchamp suivant.
\begin{eqnarray*}
\mathcal{I}nv(C,D) : \aff_A &\longrightarrow & \sens\\
u:A\to B &\longmapsto & \textrm{Inv}\left(C_u,D_u\right)\\
A\overset{u}{\to} B\overset{v}{\to} B' &\longmapsto & \textrm{Inv}_v : \textrm{Inv}\left(C_u,D_u\right)\to \textrm{Inv}\left(C_{vu},D_{vu}\right).
\end{eqnarray*}
où $\textrm{Inv}\left(C_u,D_u\right)$ est le nerf du groupoïde des isomorphismes de $C_u\z D_u$-bimodules inversibles.

\begin{prop}\label{geominv}
Soient $C,D\in A\z\ASS$ de $A$-modules sous-jacent libre de type fini. 
Le champ $\underline{\mathcal{I}nv}(C,D)$ des bimodules inversibles à isomorphisme près est un 1-champ géométrique et une présentation est donnée par le groupoïde affine lisse :
$$
\coprod_m Inv(C,D,m)\times G\ell_m \rightrightarrows \coprod_m Inv(C,D,m).
$$
\end{prop}
\begin{pr}
La lissité du groupoïde résulte de ce que les fibres de $s$ sont des $G\ell_n$-torseurs et que $G\ell_n$ est un groupe lisse.
$\coprod_m Inv(C,D,m)\rightarrow \underline{\mathcal{I}nv}(C,D)$ est surjectif car on s'est limité à des bimodule projectifs sur la base : si $A\to B\in\kcom$ tout $C\otimes B\z D\otimes B$-bimodule est, sur un recouvrement étale $B\to B'$, libre comme $B'$-module, \ie un point de $\coprod_m Inv(C,D,m)$.
Il reste à voir que $\left(\coprod_m Inv(C,D,m)\right)\times_{\underline{\mathcal{I}nv}(C,D)}\left(\coprod_m Inv(C,D,m)\right)\simeq \coprod_m Inv(C,D,m)\times G\ell_m$, or, parce que $\underline{\mathcal{I}nv}(C,D)$ est un 1-champ, le premier terme est le schéma classifiant les isomorphismes du "$C\z D$-bimodule libre sur la base universel" et est donc équivalent au second.
\end{pr}

\bigskip
La proposition suivante se déduit immédiatement de la proposition {\ref{presbim}}.

\begin{prop}\label{presabel}
\begin{enumerate}
\item Soient $c,d$ et $n$ des entiers. Les triplets formés d'une structure $C$ de $A$-algèbre associative sur $A^c$, d'une structure $C'$ d'algèbre associative sur $A^d$ et d'une structure de $C\z C'$-bimodule sur $A^n$ sont classifiés par le schéma affine :
$$
Bim(c,d,n) := Spec(A[m,m',\mu]/(A,I,I,A',I',I',Bim))
$$
(on renvoie au {\sf\S\ref{classifstcat}} pour les abréviations des structures associatives).
\item Si on se limite aux seuls bimodules inversibles, les triplets en question sont classifiés par un ouvert $Inv(c,d,n)$ du schéma $Bim(c,d,n,n)$.
\item Les quintuplets formés d'une structure $C$ de $A$-algèbre associative sur $A^c$, d'une structure $C'$ d'algèbre associative sur $A^d$, d'une structure $M$ de $C\z C'$-bimodule sur $A^m$, d'une structure $M'$ de $C\z C'$-bimodule sur $A^n$ et d'un morphisme de $C\z C'$-modules $f:M\to M'$ sont classifiées par le schéma affine :
$$
MorBim(c,d,m,n) := Spec(A[m,m',\mu,\mu',f]/(A,I,I,A',I',I',Bim,Bim',Morbim)).
$$
\item Si on se limite aux seuls isomorphismes de bimodules dans le problème précédent, le schéma affine classifiant, qui est un ouvert de $MorBim(c,d,m,m)$, est :
$$
IsoBim(c,d,m) := Spec(A[m,m',\mu,\mu',f,\det(f)^{-1}]/(A,I,I,A',I',I',Bim,Bim',Morbim)).
$$
Ce schéma se récrit :
$$
IsoBim(c,d,m) = Spec(A[m,m',\mu,f,\det(f)^{-1}]/(A,I,I,A',I',I',Bim)) = Bim(c,d,m)\times G\ell_m.
$$
\item Si on se limite en plus aux bimodules inversibles dans le problème précédent, le schéma classifiant est un ouvert $IsoInv(c,d,m)$, éventuellement vide, de $MorBim(c,d,m,m)$. On a la caractérisation :
$$
IsoInv(c,d,m) = Inv(c,d,m)\times G\ell_m.
$$
\item Il existe une structure de 2-groupoïde (cf. \cite[ch. 6]{maclane}) sur le 2-graphe
$$
\coprod_{c,d,m}Inv(c,d,m)\times G\ell_m\rightrightarrows \coprod_{c,d,m}Inv(c,d,m)
\rightrightarrows \coprod_n\ass^n.
$$
\end{enumerate}
\end{prop}

\section{Catégories karoubiennes et abéliennes}
\label{karoubi}\label{abel}

\begin{defi}\label{defkar}
Une catégorie $A$-linéaire est dite {\em additive}\index{catégorie additive} si elle possède un objet nul et si toutes les sommes finies existent.
Dans une telle catégorie, somme et produit coïncident et se caractérisent par les biproduits.

Un {\em projecteur}\index{projecteur} dans une catégorie $C$ est un endomorphisme idempotent d'un objet de $C$. Une catégorie $A$-linéaire est dite {\em karoubienne}\index{catégorie karoubienne} si elle est additive et si elle possède les noyaux et conoyaux de tout projecteur. Il revient au même de demander que tout projecteur scinde son objet en une somme deux facteurs directs (l'un étant le noyau-conoyau l'autre l'image).
On note $A\z\KAR_\UU$ la sous-catégorie pleine de $A\z\CAT$ formée par les catégories karoubiennes.

Une catégorie $A$-linéaire est dite {\em abélienne}\index{catégorie abelienne} si elle est additive, possède les noyaux et conoyaux de toutes flèches et si tout épimorphisme (resp. monomorphisme) est un conoyau (resp. un noyau). En particulier une catégorie abélienne est karoubienne.
\end{defi}

Si $C$ est une catégorie linéaire, la catégorie $\overline{C}=C\Mod$ est abélienne (et donc karoubienne et additive). On se limite dans cette étude aux seules catégories abéliennes qui sont équivalentes à des catégories de modules sur une catégorie linéaire $\UU$-petite.

\begin{defi}
La {\em catégorie des catégories abéliennes $A$-linéaires} est notée $A\z\AB_\UU$, ses objets sont les catégories équivalentes à une catégorie de modules sur une catégorie de $A\z\CAT_\UU$ et ses morphismes sont les foncteurs commutant aux colimites.
\end{defi}

Contrairement à ceux de $A\z\CAT_\UU$, les objets de $A\z\AB_\UU$ ne sont pas des catégories $\UU$-petites et on ne dispose d'un foncteur d'oubli de la structure abélienne qu'à valeurs dans
$A\z\CAT_\VV$.

\medskip

\begin{defi}
En utilisant le plongement de Yoneda $C\to C^o\Mod$, on définit un $C$-module {\em libre}\index{module libre} (resp. {\em libre de type fini}\index{module libre de type fini}) comme une somme directe (resp. d'un nombre fini) d'objets de $C$. Un $C$-module est dit {\em projectif}\index{module projectif} (resp. {\em projectif de type fini}\index{module projectif de type fini}) s'il est facteur direct d'un objet libre (resp. libre de type fini). On note $C$-proj la sous-catégorie pleine de $C^o\Mod$ formée des modules projectifs de type fini.
\end{defi}

Si $C$ est $\UU$-petite, $C$-proj ne l'est pas mais elle est toujours équivalence à une catégorie $\UU$-petite $\widehat{C}$.

\begin{defi}
\begin{itemize}
\item On définit l'{\em l'additivisation} de $C$ comme la catégorie $A$-linéaire $C^+$ dont les objets sont les familles $(x_1,\dots, x_n)$ d'objets de $C$ où $n\in\mathbb{N}$ est variable ; dont le module des morphismes de $(x_1,\dots,x_n)$ vers $(y_1,\dots,y_m)$ est $\oplus_{i,j} C(x_i,y_j)$ ; et dont la composition est donnée, pour trois objets $(x_1,\dots, x_n)$, $(y_1,\dots, y_m)$, $(z_1,\dots, z_\ell)$ et en notant $(f_i^j)_{i=1..n}^{j=1..m}$ les éléments de $\oplus_{i,j} C(x_i,y_j)$ par 
\begin{eqnarray*}
\oplus_{i,j} C(x_i,y_j) \otimes \oplus_{i,j} C(y_j,z_k) & \longrightarrow & \oplus_{i,k} C(x_i,z_k)\\
((f_i^j),(f_j^k)) & \longmapsto & (\sum_{j=1..m}f_j^kf_i^j).
\end{eqnarray*}
\item On définit la {\em karoubianisation} de $C$ comme la catégorie $A$-linéaire $\widehat{C}$ dont les objets sont les couples $(x,p)$ où $x\in C^+$ et $p$ est un projecteur de $x$ et dont le module de morphismes de $(x,p)$ vers $(y,q)$ ast donné par $\widehat{C}(x,y) := qC^+(x,y)p$) ; et dont la composition est donnée, si $p,q,r$ sont trois projecteurs, par le morphisme naturel $qC(x,y)p\otimes rC(y,z)q\to rC(x,y)p:(f,g)\mapsto gf$.
\end{itemize}
\end{defi}

Il est clair par construction que ces catégories restent $\UU$-petites.

\medskip

\begin{lemme}
$C^+$ est équivalente à la sous-catégorie pleine de $\overline{C}$ formée des modules libres de types finis.
\end{lemme}
\begin{pr}
On a un foncteur $C^+\to \overline{C}$ qui associe à $(x_1,\dots,x_n)$ le module $\oplus_i x_i$, 
compte tenu des isomorphismes de $A$-modules $\overline{C}(\oplus_i x_i,\oplus_j y_j) \simeq \oplus_{i,j}\overline{C}(x_i,y_j)\simeq \oplus_{i,j}C(x_i,y_j)$ s'étend canoniquement en un foncteur.
On dispose d'un foncteur réciproque qui à un module libre $L$ isomorphe à $\oplus_ix_i$ associe $\oplus_ix_i$, qui, afin d'être bien défini sur les flèches, demande de fixer un isomorphisme $L\simeq\oplus_ix_i$. Il est clair que ces deux foncteurs sont quasi-inverses.
\end{pr}

\begin{lemme}\label{hatc}
$\widehat{C}$ est équivalente à la sous-catégorie pleine de $\overline{C}$ formée des modules projectifs de types finis.
\end{lemme}
\begin{pr}
Si $P$ est un module projectif de $\overline{C}$ facteur d'un module libre de type fini $X$, il définit, modulo le choix d'un supplémentaire $P'$, un projecteur $p$ de $X$ dont il est l'image, $P'$ en étant le noyau ; le projecteur ayant $P$ comme noyau et $P'$ comme image est noté $p'$. Si $Q$ est un autre module projectif de type fini de supplémentaire $Q'$ et de projecteurs associés $q$ et $q'$, en posant $X:=P\oplus P'$ et $Y:=Q\oplus Q'$ on a
$$
\overline{C}(X,Y) = (q\oplus q')\overline{C}(X,Y)(p\oplus p') \simeq q\overline{C}(X,Y)p\oplus q\overline{C}(X,Y)p'\oplus q'\overline{C}(X,Y)p\oplus q'\overline{C}(X,Y)p'
$$
d'un autre côté, on a $\overline{C}(X,Y)\simeq \overline{C}(P,Q)\oplus \overline{C}(P',Q)\oplus \overline{C}(P,Q')\oplus \overline{C}(P',Q')$ et en comparant les noyaux et images des morphismes on trouve $\overline{C}(P,Q)\simeq q\overline{C}(X,Y)p$.

On construit un foncteur $\widehat{C}\to \overline{C}$ en associant à $((x_1,\dots,x_n),p)$ le noyau de $p$ vu comme projecteur de $x^n\in\overline{C}$, un tel module est bien projectif de type fini ; l'analyse précédente assure que cette association se prolonge au niveau des flèches en un foncteur (pleinement fidèle).

Réciproquement un module de type fini est toujours isomorphe à l'image d'un $((x_1,\dots,x_n),p)$ par le foncteur précédent, le choix d'un isomorphisme avec un $((x_1,\dots,x_n),p)$ pour tout module projectif de type fini définit un quasi-inverse au premier foncteur.
\end{pr}

\begin{lemme}
$\widehat{C}$ est une catégorie karoubienne.
\end{lemme}
\begin{pr}
En vertu du lemme {\ref{hatc}} il suffit de montrer que la catégorie $\widehat{C}'$des $C$-module projectifs de types finis est karoubienne.
Si $M$ et $N$ sont deux tels modules, il existe $M'$ et $N'$ tels que $M\oplus M'$ et $N\oplus N'$ soient libres.
d'où $M\oplus N\oplus M'\oplus N$ est libre et $M\oplus N$ est dans $\widehat{C}'$.
Soit $p$ un projecteur d'un objet $M$, $p$ induit un projecteur du module libre $M\oplus M'$ et ses noyaux et conoyaux correspondent à un facteur direct de $M\oplus M'$, \ie à un élément de $\widehat{C}'$.
\end{pr}

%
%

\begin{prop}\label{univcat}
\begin{enumerate}
\item On a une 2-adjonction :
$$
\widehat{(-)}\ : A\z\CAT_\UU \leftrightarrows A\z\KAR_\UU : F
$$
où $F$ est l'oubli de la structure abélienne. En particulier, si $C\in A\z\CAT_\UU$ et $D\in A\z\KAR_\UU$, on a : 
$$
\uHom_{A\z\CAT_\UU}(C,D) \simeq \uHom_{A\z\KAR_\UU}(\widehat{C},D).
$$
\item On a deux foncteurs
$$
A\z\CAT_\UU \overset{\overline{(-)}}{\longrightarrow} A\z\AB_\UU \overset{F}{\longrightarrow} A\z\CAT_\VV  
$$
vérifiant, pour $C\in A\z\CAT_\UU$ et $D\in A\z\AB_\UU$
$$
\uHom_{A\z\AB_\UU}(\overline{C},D) \simeq \uHom_{A\z\CAT_\VV}(C^o,D).
$$
\end{enumerate}
\end{prop}
\begin{pr}
{\em 1.} Soit $C\in A\z\CAT$ et $D\in A\z\KAR$, et $f:C\to D$, $f$ s'étend en un foncteur $\hat{f}:\widehat{C}\to D$. On l'étend d'abord aux modules libres : si $M=\oplus_{x\in C} x$ est libre on pose $\hat{f}(M) = \oplus_x f(x)$. Puis, si $M,M'\in\widehat{C}$ sont tels que $M\oplus M'$ soit libre, on considère le projecteur $p$ de $M\oplus M'$ associé à $M$, et on définit $f(M)=\ker\ f(p)$.
Réciproquement, étend donné un foncteur $f:\widehat{C}\to D$, sa composition à la source avec $C\to \widehat{C}$ fournit un foncteur $C\to D$.
Ces deux foncteurs sont clairement quasi-inverse l'un de l'autre et donnent l'équivalence $\Hom_{A\z\CAT_\UU}(C,D) \simeq \Hom_{A\z\KAR_\UU}(\widehat{C},D)$.

{\em 2.} Soit $C\in A\z\CAT$ et $D\in A\z\AB$, et $f:C^o\to D$, $f$ s'étend uniquement en un foncteur $f_!:\overline{C}\to D$\index{$f_{"!}$}. 
Pour $M=\oplus_{x\in I}x$ un $C$-module libre, on définit $f_!(M):=\oplus_x f(x)$ ; puis pour $M\in\overline{C}$ quelconque, dont $L_1\to L_2\to M\to 0$ est une présentation libre de $M$, on définit l'extension de Kan $f_!$ de $f$ par $f_!(M):= \textrm{colim}\ f(L_1)\to f(L_2)$.
Réciproquement, la composition d'un foncteur $f:\overline{C}\to D$ à la source avec le plongement de Yoneda $C^o\to\overline{C}$ donne un foncteur $C\to D\in A\z\CAT_\VV$. Ces deux foncteurs sont quasi-inverse l'un de l'autre et, croisé avec le fait que $A\z\CAT_\UU$ est une sous-catégorie pleine de $A\z\CAT_\VV$, donnent la bijection cherchée.
\end{pr}

\subsection{\'Equivalences de Morita}\label{morita}

\begin{lemme}\label{moritamorphisme}
Soient $C,D\in A\z\CAT$, on a l'équivalence de catégorie suivante 
$$
\underline{\mathrm{Hom}}_{A\z\AB}(\overline{C},\overline{D})\simeq \overline{C^o\otimes_AD}.
$$
\end{lemme}
\begin{pr}
\begin{eqnarray*}
\underline{\mathrm{Hom}}_{A\z\AB}(\overline{C},\overline{D}) & \simeq & \underline{\mathrm{Hom}}_{A\z\CAT}(C^o,\overline{D})\\
&\simeq & \underline{\mathrm{Hom}}_{A\z\CAT}(C^o,\uHom_{A\z\CAT}(D,\overline{A}))\\
&\simeq & \underline{\mathrm{Hom}}_{A\z\CAT}(C^o\otimes_AD,\overline{A})\\
&\simeq & \overline{C^o\otimes_AD}.
\end{eqnarray*}
\end{pr}

Explicitement, la correspondance est la suivante. Soit $f:\overline{C}\to \overline{D}$, pour $x\in C$, $f(x)$ est un $D$-module est se décompose en une famille $\{f(x)(y)\}_{y\in D}$ indexée par les objets de $D$.
Le bimodule $F$ associé à $f$ est donné, pour $x\in C$ et $y\in D$, donné par $F(x,y):= f(x)(y)$.

En particulier, le bimodule associé à l'extension de $f:C\to D$ est donné, pour $x\in C$ et $y\in D$, par $F(y,x)=D(y,f(x))$. L'action de $D(y',y)$ est simplement le produit à gauche et celle de $C(x,x')$ est le produit à droite par son image dans $D(f(x),f(x'))$. L'identité de $\overline{C}$ correspond donc à $C$ vu comme bimodule sur elle-même.

Réciproquement, si $F\in \overline{C^o\otimes_AD}$, on lui associe un foncteur $-\otimes_CF:\overline{C}\to \overline{D}$ qui à $M\in \overline{C}$ associe le $D$-module $y\mapsto \oplus_{x\in C}M(x)\otimes_AF(x,y)$.

\begin{lemme}\label{eqcateqmor}
Soient $C,D\in\ASS$, les équivalences $\overline{C}\to\overline{D}$ correspondent aux $(C^o\otimes_AD)$-modules $E$ tels qu'il existe un $(D^o\otimes_AC)$-module $F$ et des isomorphismes $E\otimes_D F\simeq C$ et $F\otimes_C E\simeq D$, \ie aux $C\z D$-bimodules inversibles.
De tels bimodules sont toujours biprojectifs et de présentation finie.
\end{lemme}
\begin{pr}
Si $E$ est un $C\z D$-bimodule induisant une équivalence $-\otimes_CE:\overline{C}\to\overline{D}$, un quasi-inverse est donné par $N\mapsto \Hom_D(E,N)$. En effet si $F$ désigne un adjoint à droite de $-\otimes_CE$, on déduit de 
$$
\Hom_D(M\otimes_CE,N) \simeq \Hom_C(M,F(N)),
$$
en prenant $M=C$, que $F(N)=\Hom_D(E,N)$.
Comme $-\otimes_CE$ est une équivalence, $F=\Hom_D(E,-)$ l'est aussi et donc commute aux petites colimites filtrantes ainsi qu'aux noyaux, on en déduit que $E$ est projectif de présentation finie comme $D$-module.

Pour l'énoncé analogue sur la structure de $C^o$-module, on remarque que $C^o\otimes_AD\simeq (D^o)^o\otimes_AC^o$ et qu'on a donc une équivalence
$$
\overline{C^o\otimes_AD}\simeq \Hom_{A\z\AB}(\overline{D^o},\overline{C^o}).
$$
Le raisonnement précédent appliqué au même module $E$ assure qu'il est projectif de type sur $C^o$.

Comme $E$ est de présentation finie $\Hom_D(E,N)\simeq N\otimes_D\Hom_D(E,D)$ et $F$ est donc associé au $D\z C^o$-bimodule $\Hom_D(E,D)$.
Le fait que $E$ et $F$ soient des bimodules inversibles est dû aux unités et co-unités de quasi-inversion de $-\otimes_CE$ et $-\otimes_D\Hom_D(E,D)$.
\end{pr}

\begin{lemme}\label{restequivmor}
Soient $u:A\to B\in\kcom$ et $C,D\in A\z\ass$.
On a la commutation (à isomorphisme naturel près) du diagramme suivant :
$$\xymatrix{
\overline{C^o\otimes_AD}\ar[d]^{(-)_u} \ar[r] & \uHom_{A\z\CAT}(\overline{C},\overline{D})\ar[d]_{(-)^u}\\
\overline{C_u^o\otimes_BD_u}\ar[r] & \uHom_{B\z\CAT}(\overline{C_u},\overline{D_u})
}$$
\end{lemme}
\begin{pr}
Si $E\in \overline{C^o\otimes_AD}$ on note $-\otimes_CE$ le foncteur $\overline{C}\to \overline{D}$ associé et $E_u=E\otimes_AB$ le $C_u^o\otimes_AD_u$-module déduit par changement de base ; si $f\in \uHom_{A\z\CAT}(\overline{C},\overline{D})$ on note $f^u$ son image par le changement de base.

Pour un $E$ fixé, on a $(-\otimes_CE)^u:\overline{C}^u\to \overline{D}^u$ et $-\otimes_CE_u:\overline{C_u}\to \overline{D_u}$, l'énoncé est équivalent à démontrer que $(-\otimes_CE)^u\simeq -\otimes_{C_u}E_u$ modulo les équivalences canoniques $\overline{C_u}\simeq \overline{C}^u$ et équivalences canoniques $\overline{D_u}\simeq \overline{D}^u$.

L'équivalence $\overline{C_u}\simeq \overline{C}^u$ associe à $M:C_u\to B\Mod$ le $C$-module $M':C\to C_u\to B\Mod$ où déduit de $C\to C\otimes_AB=C_u$ tiré de l'unité de $B$ ; il est clair que $M'$ a canoniquement une structure de $B$-module. Pour la réciproque, on rappelle que $\B B$ désigne la catégorie à un élément ayant $B$ comme algèbre d'endomorphismes ; $M:\B B\to \overline{C}$ pointe un $C$-module $M'$ dont les valeurs sont en fait des $B$-modules, \ie $M':C\to B\Mod$, qui donne un foncteur  $M'':C_u=C\otimes_AB\to B\Mod$.

Soient $E\in \overline{C^o\otimes_AD}$ et $M:\B B\to \overline{C}$, on a $(-\otimes_CE)^u(M)=\B B\to \overline{C}\to\overline{D}$ ; l'image de l'unique objet de $\B B$ est $M\otimes_CE$, et $(-\otimes_CE)^u(M)$ correspond à l'objet de $\overline{D_u}$ : $y\in D\mapsto \oplus_xM(x)\otimes_AE(x,y)$. En utilisant la structure de $B$-module de $M$ on peut récrire cet objet 
$$
N:y\in D\mapsto \oplus_xM(x)\otimes_BB\otimes_AE(x,y)\simeq \oplus_xM(x)\otimes_BE_u(x,y)\ ;
$$
c'est bien dire que $(-\otimes_CE)^u(M)$ est isomorphe à $M\otimes_{C_u}E_u$.
\end{pr}

\begin{cor}\label{restricabel}
Le foncteur $(-)^u:A\z\CAT\to B\z\CAT$ se restreint en un foncteur $(-)^u:A\z\AB\to B\z\AB$.
\end{cor}
\begin{pr}
Ceci est dû au fait que $\uHom_{A\z\AB}(\overline{C})$ soit l'image essentielle de $\overline{C^o\otimes_AD}\longrightarrow \uHom_{A\z\CAT}(\overline{C})$
\end{pr}

\medskip

On dispose d'un foncteur $I:\uHom_{A\z\CAT}(C,D)\longrightarrow \uHom_{A\z\AB}(\overline{C},\overline{D}) ; f\mapsto f_!$.

\begin{lemme}\label{pleineqbim}
Le foncteur $I$ est pleinement fidèle.
\end{lemme}
\begin{pr}
Soient $f,g:C\to D$ deux foncteurs et $F$ et $G$ les $C\z D$-bimodules associés, il faut prouver que tout morphisme de $F\to G$ provient d'une unique transformation naturelle $f\to g$.

Soit $\alpha:f\to g$ une transformation naturelle, le morphisme $A:F\to G$ associé est donné pour tout $x\in C, y\in D$ par $A_{y,x}: F(y,x) = D(y,f(x)) \to D(y,g(x))=G(y,x) ; u\longmapsto \alpha_x\circ u$.
La commutation à l'action de $D$ est sans condition, celle de $C$ impose que pour tout $u\in D(y,f(x))$ et tout $v\in C(x,x')$ on doit avoir $\alpha_{x'}\circ f(v)\circ u = g(v)\circ\alpha_x\circ u$, ce qui est vrai par hypothèse sur $\alpha$.
Réciproquement,si $A:F\to G$ est un morphisme de bimodules, on définit $A_x:f(x)\to g(x)$ comme $A_{f(x),x}(id_{f(x)})$. La condition précédente prise pour $u=id_{f(x)}$ donne la condition de naturalité sur $\alpha$.
Ces deux correspondances sont clairement inverse l'une de l'autre.
\end{pr}

\bigskip

\begin{defi}
Soient $C$ et $D$ dans $A\z\CAT$. Une {\em équivalence de Morita}\index{equivalence de Morita@{équivalence de Morita}} de $C$ vers $D$ est un foncteur $f:C\to D\in A\z\CAT$ tel que $f_!:\overline{C}\to \overline{D}$ soit une équivalence. 
Une composition d'équivalence de Morita en reste une et on note $Mor(A)$ la sous-catégorie de $A\z\CAT$ formé des équivalences de Morita.
\end{defi}

\begin{lemme}\label{restricmor}
Pour $u:A\to B\in\kcom$, le foncteur $(-)_u$ respecte les équivalences de Morita, \ie envoie $Mor(A)$ dans $Mor(B)$.
\end{lemme}
\begin{pr}
Soit $f:C\to D$ une équivalence de Morita, il agit de vérifier que $f_u:C_u\to D_u$ est encore une équivalence de Morita. On note $E$ le $C\z D$-bimodule associé à $f$, $E_u$ est naturellement muni d'une structure de $C_u\z D_u$-bimodule, en effet, $E_u$ est naturellement un $(C^o\otimes_AD)_u$-module et $(C^o\otimes D)_u\simeq C_u^o\otimes_BD_u$.
Le lemme {\ref{restequivmor}} assure que le bimodule associé à $f_u$ est $E_u$. 
On conclut en remarquant que $E_u$ reste inversible : si $F$ est un inverse pour $E$ alors $F_u$ est un inverse pour $E_u$.
\end{pr}

\begin{lemme}\label{eqmorkar}
Si $W_{A\z\KAR}$ et $W_{A\z\AB}$ désignent respectivement les sous-catégories de $A\z\KAR$ et $A\z\AB$ formées des seules équivalences, le foncteur $Mod^A:A\z\KAR\to A\z\CAT\to A\z\AB$ induit une équivalence $W_{A\z\KAR}\simeq W_{A\z\AB}$.
\end{lemme}
\begin{pr}
On pose les notations suivantes : pour un univers $\UU$, $W_{A\z\KAR_\UU}$ désigne la sous-catégorie de la catégorie des catégories karoubiennes $\UU$-petites formée des équivalences ; pour deux univers $\UU\in\VV$ on a un foncteur pleinement fidèle $W_{A\z\KAR_\UU}\to W_{A\z\KAR_\VV}$ dont on note $W_{A\z\KAR_{\VV,\UU}}$ l'image essentielle. Toujours pour deux univers $\UU\in\VV$, on définit $W_{A\z\AB_{\VV,\UU}}$ comme la catégorie des équivalences catégories abéliennes $\VV$-petites engendrées par des catégories $\UU$-petites. Si on a trois univers $\UU\in\VV\in\WW$ on a un foncteur pleinement fidèle $W_{A\z\AB_{\VV,\UU}}\to W_{A\z\AB_{\WW,\VV}}$ dont on note $W_{A\z\AB_{\WW,\VV,\UU}}$ l'image essentielle.

Si $\UU\in\VV$, on a également un foncteur $c:W_{A\z\AB_{\VV,\UU}}\to W_{A\z\KAR_{\VV,\UU}}$ qui à une catégorie abélienne $C$ associe sa sous-catégorie $C_{ppf}$ des objets projectifs de présentations finies. Pour vérifier que $c$ est bien définit sur les objets, il faut dire que $C_{ppf}$ est une catégorie $\VV$-petite, mais que, comme $C$ est engendrée par une catégorie $\UU$-petite $D$, elle est équivalente à la catégorie $\UU$-petite $\widehat{D}$. Pour la bonne définition de $c$ comme foncteur, elle résulte de ce qu'une équivalence de catégories abéliennes $f:C\to D$ envoie tout objet projectif de présentation fini de $C$ sur un tel de $D$ et donc $C_{ppf}$ dans $D_{ppf}$ et comme $C$ et $D$ étaient équivalentes, il en est de même de $C_{ppf}$ dans $D_{ppf}$ et $c(f):C_{ppf}\to D_{ppf}$ est en fait une équivalence.

Pour trois univers $\UU\in\VV\in\WW$, les considérations précédentes se résument en un diagramme
$$\xymatrix{
W_{A\z\KAR_{\VV,\UU}}\ar[r]^{Mod_\VV^A} & W_{A\z\AB_{\WW,\VV,\UU}}\\
W_{A\z\KAR_\UU}\ar[r]_{Mod_\UU^A}\ar[u]_a^{\textrm{équiv.}} & W_{A\z\AB_{\VV,\UU}}\ar[u]^b_{\textrm{équiv.}}\ar[lu]_c
}$$
commutant à équivalence près. Précisément, pour tout $C\in W_{A\z\KAR_\UU}$, il existe une équivalence naturelle $C\to \overline{C}_{ppf}$ et, pour tout $C\in W_{A\z\AB_{\VV,\UU}}$, il existe une équivalence naturelle $\overline{C_{ppf}}\to C$.

On tire de ce diagramme, un diagramme entre les nerfs de ces catégories dont la commutation à équivalence près, assure qu'il commute à homotopie près. 

$a$ une équivalence, le morphisme déduit $|a|$ sur les nerfs est une équivalence d'homotopie ; on en déduit que $Mod_\UU^A$ est homotopiquement fidèle (injectif sur les $pi_n,\ n\geq0$) et que $c$ est homotopiquement surjectif (surjectif sur les $pi_n,\ n\geq0$).
Similairement, comme $b$ est une équivalence, on déduit que $c$ est homotopiquement fidèle et que 
$Mod_\VV^A$ est homotopiquement surjectif.
$c$ est alors une équivalence d'homotopie, et donc $Mod_\UU^A$ et $Mod_\VV^A$.
\end{pr}

\begin{lemme}\label{cdanscchapeau}
$C\to \widehat{C}$ est une équivalence de Morita.
\end{lemme}
\begin{pr}
On montre que le morphisme naturel $\overline{C}\to \overline{\widehat{C}}$ est une équivalence.
Le morphisme naturel $f:C\to\widehat{C}$ induit deux foncteurs $f_!:\overline{C}\to \overline{\widehat{C}}$ et $f^*:\overline{\widehat{C}}\to\overline{C}$ quasi-inverses l'un de l'autre. En effet, la composition $f^*f_!$ donne l'identité sur les objets de $C$ et donc sur tout $\overline{C}$ et de même la composition $f_!f^*$ donne, grâce à $\widehat{C}\subset\overline{C}$, l'identité sur les objets de $\widehat{C}$ et donc sur tout $\overline{\widehat{C}}$.
\end{pr}

%
%
%
%
%


\section{Structures de modèles}

Le but de cette section est de définir les deux structures de modèles sur les catégorie des catégories linéaires pour lesquelles les équivalences sont les équivalence de catégories et les équivalences de Morita.

\subsection{Adjonctions}

\paragraph{Groupoïdes et catégories}

On note $\catcat$ la catégorie des catégories et $\gpd$ celle des groupoïdes. 
$\catcat$ et $\gpd$ sont des catégories naturellement enrichies sur elles-mêmes et l'une sur l'autre. On note $\uHom_{\catcat}$ et $\uHom_{\gpd}$ les hom dans $\catcat$ et $\uHom_{\catcat}^{int}$ et $\uHom_{\gpd}^{int}$ les hom dans $\gpd$.

L'inclusion $\iota:\gpd\to\catcat$ possède des adjoints à droite et à gauche.
L'adjoint à droite associe à une catégorie $C$ son groupoïde intérieur $C^{int}$ qui est la sous-catégorie formée des seuls isomorphismes de $C$. L'adjoint à gauche associe à $C$ le groupoïde $LC=C[C^{-1}]$ obtenu en inversant toutes les flèches.

On a :
\begin{eqnarray*}
\uHom_{\catcat}^{int}(C,D) &=& \uHom_{\catcat}(C,D)^{int}\ ;\\
\uHom_{\gpd}^{int}(G,H) &=& \uHom_{\gpd}(G,H)^{int}\ ;\\
\uHom_{\gpd}(iG,C) &=& \uHom_{\gpd}(iG,C^{int})\ ;\\
\uHom_{\gpd}(C,iG) &=& \uHom_{\gpd}(LC,G).
\end{eqnarray*}

\paragraph{Nerf et adjoints}

Le foncteur nerf\index{$|\z|$} $|\z|:\catcat\to\sens$ associe à une catégorie l'objet simplicial 
$$
|C| := n\mapsto \Hom_{\catcat}([n],C).
$$
La structure simpliciale et induite par celle cosimpliciale des catégories $[n]$.
Ce foncteur admet un adjoint à gauche\index{$\Pi_1^{cat}$}
$$
\Pi_1^{cat}:\sens\leftrightarrows\catcat
$$
qui à un ensemble simplicial $K$ associe la catégorie ayant pour objets $K_0$ et dont les flèches sont librement engendrées par $K_1$ avec les relations venant de $K_2$ ; la composition est tirée de la flèche $K_2\to K_1$ correspondant à $[1]\to[2] : 0\mapsto 0 ; 1\mapsto 2$.
Par exemple, il est évident que $\Pi_1^{cat}(\Delta^n) = [n]$.
La co-unité de cette adjonction $\Pi_1^{cat}\circ |\z|\to Id_{\catcat}$ est un isomorphisme.

\medskip

Restreint à $\gpd$, $|\z|$ admet comme adjoint à gauche\index{$\Pi_1$} $\Pi_1 = L\circ\Pi_1^{cat}$.
La co-unité de cette adjonction $\Pi_1\circ |\z|\to Id_{\gpd}$ reste un isomorphisme.


\subsection{Structure simpliciale des catégories}
\label{structuresimpliciale}

Les adjonctions précédentes permettent de définir deux structures simpliciale sur $\catcat$, dans la première $\Delta^n$ est représenté dans $\catcat$ par $[n]$ et dans la seconde par $L[n]$. On se limite à cette dernière qui a l'avantage de se restreindre à $\gpd$.

On pose :
$$
C\otimes K := C\times \Pi_1(K).
$$
L'exponentielle est donnée par :
$$
C^K := \uHom_{\catcat}(\Pi_1(K),C),
$$
et les espaces simpliciaux de morphismes entre deux catégories sont
\begin{eqnarray*}
\Hom^\Delta_{\catcat}(C,D) &:=& n\mapsto \Hom_{\catcat}(C\otimes\Delta^n, D)
\end{eqnarray*}
Or, par adjonctions et par $\Pi_1^{cat}(\Delta^n)=[n]$
\begin{eqnarray*}
\Hom_{\catcat}(C\otimes\Delta^n, D) &=& \Hom_{\catcat}(\Pi_1(\Delta^n),\ \uHom_{\catcat}(C,D))\\
 &=& \Hom_{\catcat}([n],\uHom_{\catcat}^{int}(C,D))
\end{eqnarray*}
d'où
\begin{eqnarray*}
\Hom^\Delta_{\catcat}(C,D) &=& |\uHom_{\catcat}^{int}(C,D)|.
\end{eqnarray*}

\subsection{catégories linéaires}
\label{structuresimplicialelineaire}

Soit $A$ un anneau et $A\z\CAT$ la catégorie des catégories linéaires, on a un foncteur d'oubli de la structure linéaire $A\z\CAT\to \catcat$ qui admet un adjoint à gauche : $C\mapsto C\otimes A$ où $C\otimes A$ est la catégorie $A$-linéaire dont les $A$-modules de morphismes sont librement engendrés par les ensembles de morphismes de $C$ : $(C\otimes A)(x,y) := \coprod_{C(x,y)}A$.

Pour $K\in \sens$ on pose $K_A := \Pi_1(K)\otimes A$. $A\z\CAT$ hérite de la structure simpliciale de $\catcat$ et si $C\in A\z\CAT$, $C\otimes K = C\otimes_AK_A$.
On a :
\begin{eqnarray*}
\Hom^\Delta_{A\z \CAT}(C,D) &=& n\mapsto \Hom_{\catcat}(C\times\Delta^n, D) \quad (\in \sens) \\
 &=& n\mapsto \Hom_{A\z \CAT}(C\otimes_A\Delta_A^n, D) \\
&=& |\uHom_{A\z\CAT}^{int}(C,D)|
\end{eqnarray*}
Et en particulier, si on se restreint aux foncteurs qui sont des équivalences :
\begin{eqnarray*}
\Eq^\Delta_{A\z \CAT}(C,D) &=& |\uEq_{A\z\CAT}(C,D)|
\end{eqnarray*}
où $\uEq_{A\z\CAT}(C,D)$ est le groupoïde des équivalences de catégories de $C$ vers $D$ et de leurs isomorphismes naturels.

\begin{lemme}\label{restrictionsplx}
Pour $A\to B\in\kcom$, le foncteur $(-)_u$ vérifie, pour tout $C\in A\z\CAT$, $(C\otimes\Delta^1)_u = C_u\otimes\Delta^1$
\end{lemme}
\begin{pr}
Cela découle du caractère monoïdal de $(-)_B$ (lemme \ref{monoidal}) et du fait que $\Delta^1_B=\Delta^1_A\otimes_AB$.
\end{pr}

\medskip

Si $[n]$ désigne l'ordinal à $n$ flèches, on définit la catégorie $A$-linéaire $[n]_A$ \index{$[n]_A$} comme $[n]\otimes A$. On a $[0]_A = \B A$.

\subsection{Structure de modèles d'équivalences sur $A\z\CAT$}\label{modeleequiv}

Soit $A\z\CAT$ la catégorie des catégories $A$-linéaires, elle est munie d'une structure de catégorie de modèles propre à gauche, engendrée par cofibrations, simpliciale et monoïdale (cf. \cite{bous,hollander,rezk}) pour laquelle les équivalences faibles sont les équivalences $A$-linéaires. L'ensemble des cofibrations génératrices est :
$$
I_A := \left\lbrace\emptyset \to [0]_A,\ [0]_A\coprod [0]_A \to [1]_A,\ [1]'_A \to [1]_A\right\rbrace
$$
où $[1]'_A$ est la catégorie $A$-linéaire à deux objets $x,y$ engendrée par deux flèches $a,b:x\to y$.
Celui des cofibrations triviales génératrices est $J_A:= \{[0]_A \to \Delta^1_A \}$.
Les cofibrations sont les foncteurs induisant une inclusion sur les objets et les fibrations sont les foncteurs surjectifs sur les isomorphismes de leur image essentielle. Tous les objets sont fibrants et cofibrants.

\medskip

Pour la distinguer de la suivante on note $A\z\CATeq$\index{$A\z\CATeq$} la catégorie $A\z\CAT$ munie de cette structure de modèles. On note également $A\z\CATiso$\index{$A\z\CATiso$} la structure de modèles triviale, pour laquelle les équivalences sont les isomorphismes et les sous-catégories de fibrations et cofibrations sont toute la catégorie.
Le fait que ni les fibrations ni les cofibrations de $A\z\CATeq$ ne soient celles de $A\z\CATiso$ lui interdit d'être une localisation de Bousfield de $A\z\CATiso$.

\medskip

Tous les objets étant fibrants et cofibrants, les espaces de morphismes sont donnés par les Hom simpliciaux (cf. {\sffamily \S \ref{structuresimpliciale}}) :
$$
Map_{A\z\CATeq}(C,D) = \Hom^\Delta_{A\z\CAT}(C,D) = |\uHom_{A\z\CAT}(C,D)|.
$$
En particulier si on restreint les objets aux équivalences :
$$
Map_{A\z\CATeq}^{eq}(C,D)= |\uEq_{A\z\CAT}(C,D)|.
$$

\begin{lemme}\label{restrictionquillen}
Pour $u:A\to B$ un morphisme d'anneaux, les foncteurs $(-)_u$ et $(-)^u$ sont de Quillen (à gauche et à droite respectivement) et respectent les objets fibrants et cofibrants.
\end{lemme}
\begin{pr}
Un foncteur entre catégories de modèles est de Quillen à gauche s'il est adjoint à gauche et respecte les cofibrations et les cofibrations triviales. Un foncteur entre catégories de modèles est de Quillen à droite s'il est adjoint à droite et respecte les fibrations et les fibrations triviales. Or, $(-)_u$ et $(-)^u$ sont adjoint à gauche et à droite respectivement et le reste découle du lemme {\ref{compatfibcofib}}.
Le dernier point est trivial puisque tous les objets sont fibrants et cofibrants.
\end{pr}

\subsection{Structure de modèles de Morita sur $A\z\CAT$}\label{modelemor}

On définit une seconde structure de modèle sur $A\z\CAT$, qu'on désignera par $A\z\CATmor$\index{$A\z\CATmor$}.

\bigskip

On rappelle que $Mor(A)$ désigne la sous-catégorie de $A\z\CAT$ formée des équivalences de Morita (cf. {\sf\S \ref{morita}}).
On définit les trois morphismes suivants :
\begin{eqnarray*}
\alpha &:& \emptyset \longrightarrow \widehat{\emptyset}\simeq 0 \\
\beta &:& \B A\coprod \B A \longrightarrow \widehat{\B A\coprod \B A} \\
\gamma &:& \B\frac{A[x]}{(x^2-x)} \longrightarrow \widehat{\B\frac{A[x]}{(x^2-x)}} \\
\end{eqnarray*}

\begin{thm}\label{thmstmodelemorita}
Il existe une structure de modèle sur $A\z\cat$, notée $A\z\CATmor$, pour laquelle les équivalences faibles sont les équivalences de Morita et les objets fibrants les catégories karoubiennes. Cette structure est la localisation de Bousfield à gauche de $A\z\CATeq$ pour la classe des équivalences de Morita.
\end{thm}
\begin{pr}
On définit $A\z\CATmor$ par cofibrations génératrices. On prend pour l'ensemble des cofibrations génératrices celui $I_A$ de $A\z\CATeq$ (cf. {\ref{modeleequiv}}) et pour l'ensemble des cofibrations génératrices triviales l'ensemble 
$$
J'_A := \left\lbrace[0]_A\to \Delta^1_A,\ \alpha,\ \beta,\ \gamma
\right\rbrace
$$
et pour ensemble d'équivalences celui des équivalences de Morita $Mor(A)$.

On vérifie les axiomes du théorème 2.1.19 de \cite{hovey} :
\begin{enumerate}
\item $Mor(A)$ vérifie le trois-pour-deux ;
\item les domaines de $I_A$ sont petits par rapport à $I$-cell = $A\z\CAT$ ;
\item les domaines de $J_A$ sont petits par rapport à $J$-cell ;
\item $J$-cell $\subset Mor(A)\cap I$-cell$=Mor(A)$ ;
\item $I$-inj $= Mor(A)\cap J$-inj.
\end{enumerate}

1. est déduit de ce que les équivalences de Morita soient les images inverses d'isomorphismes par $\overline{(-)}:A\z\CAT\to A\z\CAT$ et que les isomorphisme vérifient le trois-pour-deux.

2. et 3. On montre que tout catégorie est petite. Soit $C$ une catégorie $\UU$-petite et $\sharp C_1$ le cardinal de son ensemble de flèches, alors pour n'importe quel ordinal $\lambda$ $\sharp C_1$-filtré au sens de \cite[def. 2.1.2]{hovey}, $C$ est $\lambda$-petite dans $A\z\CAT$. On utilise pour le montrer un argument similaire à celui de \cite[ex. 2.1.6]{hovey} :
on considère un foncteur $f:C\to X:=\textrm{colim}_{\beta<\lambda}\ X_\beta$ où $X_\beta\in A\z\CAT$.
Les colimites dans $A\z\CAT$ sont telle que l'ensemble $X_1$ des flèches de $X$ est une colimite des ensembles de flèches des $X_\beta$ : $X_1 = \textrm{colim}_{\beta<\lambda}\ (X_\beta)_1$ ; à chaque $x\in C_1$ on associe $\beta(x)$ tel que $f(x)\in X_{\beta(x)}$, par régularité $b:=sup\{\beta(x);x\in C_1\} <\lambda$ ; on en déduit que $f$ se factorise par $X_b$.

4. Comme les éléments de $J$-cell sont des équivalences de Morita, il suffit de montrer que le pushout d'une équivalence de Morita reste une équivalence de Morita.
Soient $A\to B\in Mor(A)$ et $A\to C\in A\z\CAT$, on note $C*_AB$ leur pushout, on veut montrer que $C\to C*_AB$ est une équivalence de Morita. Comme la karoubianisation est un adjoint à gauche on a $\widehat{C*_AB}\simeq \widehat{C}*_{\widehat{A}}\widehat{B}$. Or, parce que $\widehat{A}\simeq \widehat{B}$, on a $\widehat{C}*_{\widehat{A}}\widehat{B}\simeq \widehat{C}$.

5. Les objets de $I$-inj sont les équivalences de catégories surjectives et vérifient donc les conditions de $J$-inj$\cap Mor(A)$.
Pour l'inclusion réciproque, on considère $A\to B\in I_A$, $f:C\to D$ dans $J$-inj$\cap Mor(A)$ et un diagramme
$$\xymatrix{
A\ar[r]\ar[d] & C\ar[r]\ar[d]^f &\widehat{C}\ar[d]^{f_!}\\
B\ar[r]^x\ar@{..>}[ru]^{x'}\ar@{..>}[rru]^{x"} & D \ar[r] &\widehat{D}.
}$$
On cherche à relever $x$ en $x'$. Parce que $f_!$ est une équivalence surjective, on sait que $x"$ existe ; on en déduit $x'$ par pleine fidélité de $C\to \widehat{C}$.

\bigskip

Les objets fibrants de $A\z\CATmor$ sont les objets $C$ tels que pour tout $u:A\to B\in J_A$, tout foncteur $f_u:A\to C$ se prolonge en un foncteur $f'_u:B\to C$. Soit $C$ un objet fibrant, montrons que c'est une catégorie karoubienne.
La condition de prolongement pour $\Delta^0_A\to \Delta^1_A$ est toujours vérifié pour tout catégorie et n'apporte rien. Celle par rapport à $\alpha$ signifie que $C$ possède un objet nul. Celle pour $\beta$ implique que $C$ possède les sommes directes ($f_\beta$ pointe deux objets de $C$ et l'image par $f'_\beta$ de la somme directe des générateurs de $\widehat{\B A\coprod\B A}$ en donne une somme directe dans $C$). Enfin la condition respective à $\gamma$ implique que tout projecteur de $C$ possède un noyau ($f_\gamma$ pointe un projecteur et l'image par $f'_\gamma$ du noyau de $x$ dans $\widehat{\B A[x]/(x^2-x)}$ donne un noyau dans $C$).

\bigskip

Pour montrer que cette structure de modèle est une localisation de Bousfield à gauche, on utilise la définition 3.3.1 de \cite{hirschhorn}. Il s'agit de vérifier que 
\begin{enumerate}
\item les équivalences sont les $Mor(A)$-équivalences ;
\item les cofibrations sont celles de $A\z\CATeq$ ;
\item les fibrations sont les morphismes de relèvement à droite par rapport aux cofibrations de $A\z\CATeq$ qui sont dans $Mor(A)$.
\end{enumerate}
Les deux premiers points sont clairs par définition de la structure $A\z\CATmor$ et le troisième est la définition des fibrations dans une catégorie de modèles fermée.
\end{pr}

\subsection{Comparaison des structures de modèles}\label{compstmod}

On compare les trois structures de modèles existantes sur $A\z\CAT$.

\begin{prop}
On a les foncteurs suivants, surjectifs sur les classes d'équivalences d'objets :
$$\xymatrix{
A\z\CATiso \ar[r]^{\B^A} & A\z\CATeq \ar[r]^{\mathcal{K}^A} & A\z\CATmor
}$$
\end{prop}
\begin{pr}
Ces foncteurs sont en fait l'identité de $A\z\CAT$. La compatibilité avec les équivalences est due aux inclusions $W_{A\z\CATiso}\subset W_{A\z\CATeq}\subset W_{A\z\CATmor}$ et la surjectivité résulte de ce que ces inclusions soient l'identité sur les objets.
\end{pr}

\setcounter{tot}{0}
\newpage
\thispagestyle{empty}
\chapter{Modules des catégories linéaires}\label{chapmodcat}
\thispagestyle{chheadings}

On définit dans ce chapitre les champs auxquels le sujet à conduit à s'intéresser et différents morphismes entre eux. Les preuves de géométricité et le calcul des complexes tangents sont reportées au chapitre suivant.

\medskip

Le diagramme de la section {\ref{morphismes}} résume les objets importants construit dans ce chapitre.

\paragraph{Les différents modules}
Comme dégagé en introduction, tout problème de modules demande des changements de bases et une relation d'identification.
Concernant le premier point, s'il est clair que, pour un anneau $A$, les $A$-points du classifiant des catégories linéaires sont les catégories $A$-linéaires, il apparaît qu'il existe naturellement, pour un morphisme d'anneaux $u:A\to B$ deux foncteurs envoyant les catégories $A$-linéaires dans les catégories $B$-linéaires : $(-)_u$ et $(-)^u$, adjoints respectivement à gauche et à droite du foncteur d'oubli $B\z\CAT\to A\z\CAT$ (cf. {\sf\S\ref{chgtbase}}).
Le problème se pose en fait déjà pour les modules de modules sur un anneau, où il existe aussi deux foncteurs $u^*$ et $u^!$ transportant les $A$-modules dans les $B$-modules, là encore, adjoints à gauche et à droite du foncteur d'oubli. Mais dans ce cadre, seul l'adjoint à gauche est considéré pour construire le champs des modules des modules. Comme il est bien connu, les foncteurs $u^*$ permettent d'associer un préfaisceau sur $Spec(A)$ à tout $A$-module, dont l'interprétation est la décomposition du module en une famille d'espaces vectoriels paramétrée par $Spec(A)$. L'adjoint à droite ne semble pas jouir d'interprétation géométrique.

Au niveau des catégories, l'adjoint à gauche $(-)_u:A\z\CAT\to B\z\CAT$ se compare à celui $u^*$ pour les modules, au sens où, si $C\in A\z\CAT$, les $B$-modules de morphismes de $C_u$ sont les images des $A$-modules par $u^*$. Les foncteurs $(-)_u$ permettent donc de déployer $C$ en un préfaisceau en catégories sur $Spec(A)$ tel qu'étant donné deux objets $x$ et $y$ de $C$, le préfaisceau des morphismes entre ces objets soit celui canoniquement associé au $A$-module des morphismes entre $x$ et $y$ dans $C$. Comme pour les modules, on peut y penser comme à une décomposition de $C$ le long de $Spec(A)$.

Ainsi, poussé par cette intuition de décomposition spectrale, on considérera nos problèmes de modules de catégories en utilisant l'adjoint à gauche comme changement de base. L'adjoint à droite, quand à lui, trouve son utilité dans la construction du champ des catégories abéliennes (cf. infra).

\smallskip

Le choix de ces foncteurs semble également avoir un rapport avec certaines procédure de champification et cela sera étudié dans un prochain travail \cite{anel1}.

\medskip

Concernant les relations d'identification, les catégories linéaires en ont naturellement trois, correspondant à autant de points de vue sur ces objets :
\begin{enumerate}
\item on peut considérer une catégorie comme un anneau à plusieurs objets, l'ensemble des objets devient un invariant et la notion naturelle d'identification est alors l'{\sl isomorphie} ;
\item on peut voir les catégories comme des objets de nature "homotopique" (généralisant les espaces classifiants de groupoïdes), il n'y a plus lieu alors de distinguer deux objets isomorphes et on identifie les catégories via les {\sl équivalences de catégories} ;
\item enfin, on peut les regarder à travers leurs représentations (leurs catégories de modules) et l'équivalence mise en jeu est celle de {\sl Morita}, \ie l'équivalence de leur catégorie de modules.
\end{enumerate}
Il est notable que ces relations sont compatibles entre elles : tout isomorphisme est une équivalence et toute équivalence est une équivalence de Morita. En conséquence, on aura des morphismes surjectifs entre les objets classifiants correspondants.

Chacun des trois points de vue donne un préchamp classifiant qui n'est pas en général un champ ; les raisons en sont multiples :
\begin{itemize}
\item le recollement à isomorphisme près (et a fortiori pour les autres identification) de catégories ayant plusieurs objets peut donner des objets n'ayant plus d'ensemble d'objets globaux ; 
\item plus subtilement, même dans le cas où il n'y a qu'un unique objet, et où le problème précédent ne se pose pas, le recollement de catégories linéaires à équivalence près peut donner encore des objets n'ayant pas d'objet globaux. 
La situation est exactement l'analogue de ce qui se passe lorsqu'on recolle des champs classifiants de groupes : leur recollement donne en général une gerbe non neutre.
\end{itemize}
Dans ces deux cas, les objets recollés n'ayant pas d'ensemble d'objets globaux, ils ne sont pas descriptible par des catégories et ne sont donc pas des points des préchamps classifiants les catégories.

Il faut donc compléter ces préchamps en champs, et, si l'opération est bien définie dans la catégorie des préchamps, il convient néanmoins de décrire ce nouvel objet, qui a priori peut-être très différent du préchamp initial (le morphisme naturel entre les deux, dont on sait qu'il ne va pas être surjectif sur les objets peut a priori ne pas être non plus injectif).

La description des champifiés des classifiants des catégories est mentionnée pour chaque cas mais non démontrée et fera l'objet du travail futur \cite{anel1}. On mentionne simplement que, dans les cas des modules pour les équivalences de catégories, la description attendue du champifié comme classifiant certains champs en catégories linéaires est lié au fait que les changements de bases soient des adjoints à gauche.

\bigskip

En plus des modules de catégories linéaires, on considère également le problème des modules de catégories abéliennes à équivalence près. Comme mentionné plus haut, les changement de base du classifiant des catégories abéliennes sont les foncteurs adjoint à droite de l'oubli, ce qui les rend moins sensible à une intuition géométrique : la fibre d'une catégorie abélienne $C$ en un point $p:A\to B$ de $Spec(A)$ n'est plus la catégorie des morphismes entre les restrictions des objets de $C$ au point, mais la catégorie des $B$-modules dans la catégorie $C$, qu'on peut interpréter -- quoique très imprécisément, mais cela aide à récupérer une certaine intuition géométrique -- comme la catégorie des objets dans $C$ de {\sl support}\footnote{Il faut entendre, ici, {\sl support} en un sens plus naïf que celui que lui donne habituellement la géométrie algébrique.} dans $Spec(B)$.

La raison de ce choix de changement de base tient, d'une part, à ce que, contrairement à l'adjoint à gauche, le changement de base adjoint à droite conserve le caractère abélien (lemme {\ref{restricabel}}), et d'autre, à ce qu'il est le bon changement de base à considérer si on veut un morphisme du classifiant des catégories linéaires vers celui des catégories abélienne qui associé à une catégorie sa catégorie de modules (lemme {\ref{mormod}}).\footnote{On peut noter que, formellement, ce dernier argument est valable aussi pour les modules : le morphisme qui à un $A$-module $M$ associe son dual $M':=\uHom_A(M,A)$, commute aux changements de base si on utilise l'adjoint à droite du côté de $M'$.}

Comme dans le cas linéaire, le problème des modules de catégories abéliennes ne donne pas naturellement un champ et, là encore, on se contente d'indiquer ce que pourrait être ce champ. Il semble que le fait que les changements de bases soient des adjoints à droite permette de décrire ce champ comme classifiant certains cochamps (objets vérifiant la condition duale de celle de descente), cela sera étudié dans le prochain travail \cite{anel1}.

\medskip

La notion de catégorie abélienne à laquelle on se limite est telle qu'elles admettent toujours des petits générateurs, et cela fournit un morphisme surjectif depuis le champ des catégories linéaires à équivalence de Morita près dont on montre que c'est une équivalence (théorème {\ref{thmequivmodulemoritaabel}}).

\paragraph{Les différents champs définis}
Pour chaque problème de modules on définit essentiellement deux champs, l'un noté en majuscules l'autre en minuscules, le second est toujours un sous-champ plein du premier pour lequel on a imposé certaines conditions de finitude sur les objets qu'il classifie. Ce sont ces champs dont il sera montré au chapitre suivant qu'ils sont géométriques. Les champs "majuscules", eux, figurent dans la présentation pour deux raisons : ils servent pour le calcul des champs de lacets des "minuscules" et ils mettent en avant que l'obtention de champs comme réponse aux problèmes de modules est indépendante des conditions de finitude nécessaires au caractère géométrique.

Les définitions sont assez rapides et on ne relève les champs qui sont d'intérêts pour le chapitre {\ref{chapgeom}} que dans le diagramme commutatif de la section {\ref{morphismes}} qui résume les relations entre eux.

\medskip

Les constructions des champs se font à partir de (sous-)préfaisceaux faibles de Quillen et on a essayé de tenir les notations suivantes : pour une lettre $M$ désignant le problème étudié $\widetilde{M}$ représente un (sous-)préfaisceau faible de Quillen ; $M$ son strictifié, qui est un (sous-)préfaisceau de Quillen (cf. {\ref{strictifquillen}}) ; $|M|$ le préfaisceau simplicial nerf des équivalences (cf. {\S\ref{classifquillen}}) et $\underline{M}$ le champifié de $|M|$.

\bigskip

On rappelle qu'on se fixe deux univers $\UU\in \VV$ et un anneau commutatif $\kk$ de $\UU$.
$\kcom$ est la catégorie des algèbres associatives commutatives sur $\kk$ dans $\UU$ et 
$\kaff(\simeq \kcom^o)$ celle des schémas affines dans $\UU$ qu'on munie de la topologie étale : tous les champs construits dans ce chapitre le sont sur ce site.

L'univers $\VV$ ne sert que pour définir certaines catégories de catégories.

\section{Modules à isomorphisme près}\label{modcatiso}

On définit dans cette section le champ classifiant les catégories linéaires à isomorphismes près.

\bigskip

Soit $A\in\kcom$, on considère la catégorie $A\z\CATiso$ des catégories et des foncteurs linéaires munie de la structure de modèle triviale pour laquelle les équivalences sont les isomorphismes.
On définit le préchamp faible de Quillen à gauche suivant :
\begin{eqnarray*}
\tkCATiso : \kaff^o & \longrightarrow & \catcat_\VV\\
	A & \longmapsto & A\z\CATiso \\
	u:A\to B & \longmapsto & (-)_u : A\z\CATiso\to B\z\CATiso
\end{eqnarray*}
Le lemme \ref{fonctofaible} assure que $\tkCATiso$ est un préfaisceau faible. Les structures de modèles $A\z\CATiso$ étant triviales, le caractère de Quillen à gauche, s'établit par le fait que les restrictions soient des adjoints à gauches et qu'elles préservent les isomorphismes (lemme \ref{compatfibcofib}).

On note $\kCATiso$ le strictifié de $\tkCATiso$, $|\kCATiso|$ le préchamp simplicial déduit et $\ukCATiso$ le champ simplicial associé.

\paragraph{Catégories de type fini.}
À $A\in\kcom$ on associe $A\z\catiso$, sous-catégorie pleine de $A\z\CATiso$ formée des catégories dont le graphe sous-jacent est projectif de type fini (cf. {\sf \S \ref{graphes}}). Les foncteurs $(-)_u$ préservant ces conditions (lemme \ref{cdtfinitude}), les $A\z\catiso$ forment un sous-préfaisceau faible de Quillen à gauche $\tkcatiso$ de $\tkCATiso$.

On note $\kcatiso$ le strictifié de $\tkcatiso$ ; $|\kcatiso|$ le préfaisceau simplicial déduit et $\ukcatiso$ le champ simplicial associé.

\begin{prop}
Pour $A\in\kcom$ et $C,D\in A\z\CATiso$ (resp. $C,D\in A\z\CATiso$), le préfaisceau $\Omega_{C,D}|\kCATiso|$ (resp. $\Omega_{C,D}|\kcatiso|$) est un faisceau et classifie les isomorphismes de $C$ vers $D$.
\end{prop}
\begin{pr}
On rappelle que comme $|\kcatiso|$ est plein dans $|\kCATiso|$, si $C$ et $D$ sont deux points de $|\kcatiso|$, on a $\Omega_{C,D}|\kcatiso|\simeq \Omega_{C,D}|\kCATiso|$. On se limite donc au cas où $C,D\in A\z\CATiso$.
Les valeurs de $|\kCATiso|$ étant des nerfs des groupoïdes d'isomorphismes, il est évident que $\Omega_{C,D}|\kCATiso|$ classifie les isomorphismes de $C$ vers $D$. Le fait que ce soit un faisceau est un corollaire de la descente des morphismes d'algèbres \cite[\sf ch. VIII]{sga1}.
\end{pr}

On a le corollaire évident.
\begin{cor}
$\ukCATiso$ et $\ukcatiso$ sont des 1-champs.
\end{cor}

\paragraph{Champification}\label{champifiso}
$|\kCATiso|$ et $|\kcatiso|$ ne sont pas des champs, car, comme dit en introduction, le recollement de catégories à isomorphisme induit un recollement à isomorphisme près des ensembles d'objets qui n'a aucune raison d'être trivial.

Soit $\uEEns\simeq \coprod_n\underline{\B\Sn}$\index{$\uEEns$} le champs sur $\kaff$ classifiant les ensembles $\UU$-petits (champ associé au préchamp constant de valeurs le groupoïde $\ens^{int}$).
L'ensemble d'objets d'une catégorie étant un invariant de sa classe d'isomorphisme, on en déduit un morphisme de champs $ob:\ukCATiso\to \uEEns$.
Soit $X$ un schéma, et un morphisme $C:X\to \ukCATiso$. On sait par des principes généraux que localement sur $X$, $C$ se remonte en un point de $|\kCATiso|$, \ie en une catégorie linéaire et $C$ est obtenu par un recollement à isomorphisme près de ces objets locaux. 
Fixons le nombre $n$ d'objets et considérons la fibre $\ukCATisons$ de $\ukCATiso\to \uEEns$ le long de $n:pt\to \underline{\B\Sn} \to \uEEns$. Le faisceau des objets de $C$ est un faisceau localement constant de degré $n$, il est trivial si et seulement si on a un relèvement 
$$\xymatrix{
&\ukCATisons \ar[d] \ar[r]&pt\ar[d]^n\\
X\ar@{-->}[ru]\ar[r]^C&\ukCATiso \ar[r]^{ob} &\uEEns 
}$$
et l'obstruction est donc donnée par la classe $[ob\circ C]\in H^1(X,\Sn)\simeq \Hom(X,\underline{\B\Sn})$, \ie par le cocycle de recollement des objets.

\bigskip
Il semble clair d'après le raisonnement précédent que le champ $\ukCATiso$ doit être le classifiant des faisceaux en catégories dont le faisceaux des objets est localement constants.
Ce résultat sera démontré dans \cite{anel1}.

\subsection{Modules d'algèbres associatives}\label{chpass}

Soit $A\in\kcom$, on note $A\z\ASS$ la catégorie des $A$-algèbres associatives unitaires, elle s'identifie canoniquement à la sous-catégorie pleine de $A\z\CAT$ formée des catégories ayant un unique objet.
Pour $u:A\to B$, le foncteur $(-)_u$ ne pouvant que diminuer le nombre d'objet d'une catégorie, il envoie $A\z\ASS$ dans $B\z\ASS$. On définit ainsi un sous-champ plein $\ukASS$\index{$\ukASS$} de $\ukCATiso$.

On note $A\z\ass$ la sous-catégorie pleine de $A\z\ASS$ formé des $A$-algèbres dont le module sous-jacent est projectif de type fini. Les foncteurs $(-)_u$ conservant la projectivité des modules, ces catégories définissent un sous-champ plein $\ukass$\index{$\ukass$} de $\ukASS$.

\begin{prop}
$|\kASS|$ est déjà un champ.
\end{prop}
\begin{pr}
C'est un corollaire immédiat de la descente des modules projectifs et des structures d'algèbres \cite[\sf ch. VIII]{sga1}.
\end{pr}

La différence avec la situation de $|\kCATiso|$ et que les recollements d'algèbres n'ont pas de de soucis de recollement des objets.

\medskip

\begin{prop}
$\ukASS$ est un sous-champ ouvert de $\ukCATiso$ et $\ukass$ est un sous-champ ouvert de $\ukcatiso$.
\end{prop}
\begin{pr}
Le nombre d'objets d'une catégorie $C$ est un invariant de sa classe d'isomorphie, conservé par les foncteurs de restriction de $\ukCATiso$. Si $\uEEns$ désigne le champ des ensembles $\UU$-petits sur $\ukaff$, l'association de son nombre d'objet à une catégorie définit un morphisme de champs $o:\ukCATiso\to \uEEns$. Le champ $\uEEns$ est la réunion disjointe des $\underline{\B\Sn}$. En particulier $\{*\}=\underline{\B\mathfrak{S}_1}$ est un sous-champ ouvert de $\uEEns$ ; son image réciproque par $o$ (le sous-champ des catégories ayant un seul objet, \ie le champ $\ukASS$) est donc un sous-champ ouvert de $\ukCATiso$.
La deuxième assertion se déduit de ce que $\ukass \simeq \ukASS\times_{\ukCATiso}\ukcatiso$.
\end{pr}

\subsection{Modules d'algèbres commutatives et schémas affines}\label{chpcom}\label{chpsch}

On définit $A\z\COM$ la sous-catégorie pleine de $A\z\ASS$ formée des algèbres qui sont en plus commutatives. Pour $u:A\to B$, le foncteur $(-)_u$ préserve la nature commutativité des algèbres et les $A\z\COM$ définissent un sous-champ plein $\ukCOM$\index{$\ukCOM$} de $\ukASS$.

On note $A\z\com$ la sous-catégorie pleine de $A\z\COM$ formée des $A$-algèbres dont le module sous-jacent est projectif de type fini. Ces catégories définissent un sous-champs plein $\ukcom$\index{$\ukcom$} de $\ukCOM$.

\begin{prop}
$\ukcom$ est un sous-champ fermé de $\ukass$.
\end{prop}
\begin{pr}
Soit $C:Spec(A)\to \ukass$ une $A$-algèbre, on montre que le lieu sur $Spec(A)$ où $C$ est commutative est un fermé. On considère le module $\uHom_A(C\otimes_AC,C)$ et la section donnée par le commutateur, $C$ est commutatif sur le lieu où cette section est nulle. Ce lieu est fermé car $C$, et donc $\uHom_A(C\otimes_AC,C)$, est projectif de type fini comme $A$-module.
\end{pr}

\paragraph{Schémas} La catégorie $A\z\COM^o$ est équivalente à la catégorie des schémas affines sur $A$ et le champ $\ukCOM^o$ est équivalent au champ $\ukaff$ des schémas affines. Le champ $\ukcom^o$ est, lui, équivalent au champ $\underline{\mathcal{A}f\!f}$ des schémas affines de longueur finie.

\section{Modules à équivalence près}\label{modcateq}

On définit dans cette section le champ associé à la classification des catégories linéaires à équivalence près.
La construction commence avec celle du préchamp classifiant les catégories linéaires où il apparaît un choix pour les foncteurs de restriction menant à deux définitions distinctes. Puis, le classifiant naturel n'étant pas un champ, on discute un peu la champification.

On rappelle qu'on s'est fixé trois univers $\UU\in\VV\in\WW$.

\bigskip

Pour une $A\in\kcom$, on note $A\z\CATeq_\UU$ la catégorie des catégories $A$-linéaires $\UU$-petites et leurs foncteurs munie de sa structure de modèles simpliciale où les équivalences sont les équivalences de catégories (cf. {\S \ref{modeleequiv}}).
Si $u:A\to B$ est un morphisme dans $\kcom$ on rappelle (cf. {\sf\S \ref{operations}}) qu'on a deux foncteurs :
$$
(-)_u\ , (-)^u\ :\ A\z\CATeq_\UU \longrightarrow B\z\CATeq_\UU 
$$
adjoints respectifs à gauche et à droite du foncteur d'oubli : $B\z\CATeq \longrightarrow A\z\CATeq$.

Chacun des foncteurs $(-)_u$ et $(-)^u$ peut servir à définir une structure de préfaisceau faible sur la famille $\{A\z\CATeq_\UU\}_{A\in\kcom}$ et on définit alors 
\begin{eqnarray*}
\tkCATeqg : \kaff^o & \longrightarrow & \catcat_\VV \\
	A & \longmapsto & A\z\CATeq_\UU \\
	u : A\to B & \longmapsto & (-)_u : A\z\CATeq_\UU \to B\z\CATeq_\UU
\end{eqnarray*}
et
\begin{eqnarray*}
\tkCATeqd : \kaff^o & \longrightarrow & \catcat_\WW\\
	A & \longmapsto & A\z\CATeq_\VV \\
	u : A\to B & \longmapsto & (-)^u : A\z\CATeq_\VV \to B\z\CATeq_\VV
\end{eqnarray*}

Le lemmes \ref{fonctofaible}, \ref{restrictionquillen} et \ref{restrictionsplx} établissent que $\tkCATeqg$ et $\tkCATeqd$ sont des préfaisceaux faibles de Quillen simpliciaux, respectivement à gauche et à droite, comme le rappellent les exposants $g$ et $d$.

\subsection{Le champ droit des catégories linéaires}


On note $\kCATeqd$ le strictifié de $\tkCATeqd$ ; par le corollaire \ref{strictifquillen} c'est un préfaisceau de Quillen à droite.
On note $|\kCATeqd|$ le préfaisceau simplicial classifiant les équivalences de $\kCATeqd$ et $\ukCATeqd$ le champ associé.

\medskip
Les univers de références pris pour la définition de $\ukCATeqd$ sont volontairement décalés. Le champ $\ukCATeqd$ ne servira que pour définir le champ $\ukAB$ du {\S \ref{modcatab}}.


\subsection{Le champ gauche des catégories linéaires}


On note $\kCATeqg$ le strictifié de $\tkCATeqg$ ; par le corollaire \ref{strictifquillen}, $\kCATeqg$ est un préfaisceau de Quillen à gauche.
On note $|\kCATeqg|$ le préfaisceau simplicial classifiant les équivalences de $\kCATeqg$ et $\ukCATeqg$ le champ associé.

\subparagraph{Catégories de type fini}
Pour $A\in\kcom$, on note $A\z\cateq$ la sous-catégorie pleine de $A\z\CATeq$ des catégories équivalentes à une catégorie dont le graphe sous-jacent est projectif de type fini (cf. {\sf \S \ref{graphes}}), \ie des catégories ayant un nombre essentiel d'objets fini et des Hom projectifs de type fini sur $A$.
Les foncteurs $(-)_u$ conservent ces conditions (cf. lemme \ref{cdtfinitude}) et ces catégories définissent un sous-préfaisceau faible $\tkcateq$ de $\tkCATeq$.
Les sous-catégories $A\z\cateq$ étant stables par équivalences et tous les objets étant fibrants dans les structures de modèles considérées, $\tkcateq$ est en fait un sous-préfaisceau faible de Quillen {\S \ref{prefQ}} de $\tkCATeq$.

On note $\kcateq$ le strictifié de $\tkcateq$, qui est un sous-préfaisceau de Quillen à gauche.
On note $|\kcateq|$ le préfaisceau simplicial classifiant ses équivalences et $\ukcateq$ le champ associé.

$\kcateq$ étant un sous-préfaisceau de Quillen de $\tkCATeq$, $|\kcateq|$ est un sous-préfaisceau simplicial plein de $|\kCATeqg|$. De même $\ukcateq$ est un sous-champ plein de $\ukCATeq$.

\subparagraph{Catégories à nombre d'objets borné}\label{chpcatcnx}
Pour $A\in\kcom$ et un entier $n$, on note $A\z\CAT^{\leq n}$ (resp. $A\z\cat^{\leq n}$) la sous-catégorie pleine de $A\z\CAT$ (resp. $A\z\cat$) formée des catégories ayant essentiellement moins de $n$ objets. 
Le nombre essentiel d'objets étant conservé par équivalence, et les foncteurs de changement de base $(-)_u$ ne pouvant que diminuer le nombre essentiel d'objets (lemme {\ref{cdtfinitude}}), ces catégories définissent deux  sous-préfaisceau faible de Quillen de $\tkCATeq$, dont le second est un sous-préfaisceau faible de $\tkcateq$. Les champs simpliciaux associés sont notés respectivement $\ukCATeqinfn$\index{$\ukCATeqinfn$} et $\ukcateqinfn$\index{$\ukcateqinfn$}, ce sont des sous-champs pleins de $\ukCATeq$ et $\ukcateq$.

\medskip

Lorsque $n=1$ on abrège $\ukCAT_*$ (resp. $\ukcat_*$\index{$\ukcat_*$}) le champ $\ukCAT^{Eq,\leq1}$ (resp. $\ukcat^{Eq,\leq1}$). On en parle comme du {\em champ des catégories connexes}\index{catégorie connexe} ou des {\em gerbes linéaires}\index{gerbe linéaire}.

\medskip

Si $n\leq m$ on dispose de morphismes pleinement fidèles $\ukCAT^{\leq n}\longrightarrow \ukCAT^{\leq m}$ et $\ukcat^{\leq n}\longrightarrow \ukcat^{\leq m}$.

%


\subsection{Diagonale}\label{diagocateq}

Soit $A\in\kcom$, on considère deux catégories $A$-linéaires $C$ et $D$, de graphes de type fini, vues comme des points de $\kCATeq$. Comme tous les objets sont fibrants et cofibrants, le corollaire {\sf\ref{cheminchampquillen}} donne l'équivalence de champs :
$$
\Omega_{C,D}\ukCATeq  \simeq \underline{\Eq_{\tkCATeq}^\Delta(C,D)}.
$$
Or, les structures simpliciales des $A\z\CAT$ sont telles (cf. {\sf\S\ref{structuresimplicialelineaire}}) que $\Eq_{\tkCATeq}^\Delta(C,D)$ est le préfaisceau simplicial nerf du préfaisceau en groupoïdes :
\begin{eqnarray*}
\uEq(C,D) :  \aff_{/Spec(A)} & \longrightarrow & \gpd \\
	u:A\to B & \longmapsto & \uEq_{A\z\CAT}(C_u,D_u)
\end{eqnarray*}
où $\uEq_{A\z\CAT}(C_u,D_u)$ est le groupoïde des équivalences de catégories $C_u\to D_u$ et de leurs isomorphismes naturels.

\medskip

On en déduit immédiatement la proposition suivante.

\begin{prop}
Soient $A\in\kcom$ et $C\in A\z\CAT$.
\begin{itemize}
\item Le préfaisceau $\hat{\pi}_0(|\kCATeq|)$ classifie les classes d'équivalence de catégories linéaires. 
\item Le préfaisceau $\hat{\pi}_1(|\kCATeq|,C)$ classifie les classes d'isomorphie des auto-équivalences de la catégorie $C$. 
\item Le préfaisceau $\hat{\pi}_2(|\kCATeq|,C)$ classifie les automorphismes du foncteur identité de la catégorie $C$. 
\item Pour $n>2$, les préfaisceaux $\hat{\pi}_n(|\kCATeq|,C)$, et, a fortiori, les faisceaux associés, sont tous triviaux.
\end{itemize}
\end{prop}
%
%
%

$\hat{\pi}_2(|\kCATeq|,C)$ est fait directement un faisceau. Pour $\pi_0(|\kCATeq|]$ et $\pi_1(|\kCATeq|,C)$, conformément à l'interprétation faite au {\sf\S\ref{faischot}}, les faisceaux associés classifient, respectivement, les {\sl formes locales} des classes d'équivalence de catégorie et des classes d'isomorphie d'auto-équivalences de $C$. Il est difficile d'être plus précis, et c'est la raison pour laquelle les objets d'importance sont plus les champs $\ukCATeq$, $\Omega\ukCATeq$ et $\Omega\Omega\ukCATeq$ que leurs invariants d'homotopie.

\subsection{Champification}\label{champifeq}
\label{classifcat}

Les préfaisceaux simpliciaux $|\kCATeqg|$ et $|\kcateq|$ ne sont pas des champs. On a déjà abordé les deux raisons pour lesquels il pouvait ne pas y avoir de d'objets globaux, mais $|\kCATeqg|$ et $|\kcateq|$ ont en plus la particularité que leur préchamps de morphismes ne sont pas des champs : le groupoïde des foncteurs à isomorphisme près n'est pas un champ (cf. infra). Pour les champifier, il va falloir champifier aussi ces champs là.

Ainsi, à strictement parler, le problème des modules des catégories à équivalence près, soit ne se résout pas par un champ, soit se résout par un champ dont ce qu'il classifie n'est pas si clair (puisqu'il contient plus de points et de morphismes que ceux correspondant aux catégories et leurs équivalences).

Toutefois, comme souvent dans ce genre de problèmes, la réponse devance presque la question et on s'attend à ce que le champ $\ukCATeqg$ associé à $\kCATeqg$ soit un champ classifiant certains champs en catégories linéaires, qui localement sont équivalents à un champifié du préchamp issu d'une catégorie linéaire. \'Elément un peu plus subtil, les morphismes de $\ukCATeqg$ doivent être des formes tordues de foncteurs linéaire et, a priori, tout morphisme de champs entre deux points de $\ukCATeqg$ n'est pas de ce type. Cette analyse traduit qu'il doit exister un morphisme du champ $\ukCATeqg$ vers le champs des champs linéaires, dont la remarque sur les morphismes laisse penser qu'il n'est pas plein (mais il doit être fidèle).

\medskip

Nous prévoyons d'étudier ces problèmes dans un travail futur sur la champification \cite{anel1}.

\paragraph{Foncteurs tordus}\index{foncteur tordu}\label{fctrtordu}

Soit $\{U_i\to X\}$ une famille couvrante de $X=Spec(A)\in\kaff$, et soient $C, D\in A\z\CAT$ ; on note $C_i$ et $D_i$ les restrictions de $C$ et $D$ à $U_i$.
On considère une donnée de descente de morphismes de $C$ vers $D$, elle consiste en des foncteurs $f_i:C_i\to D_i$ et des isomorphismes $\phi_{ji}:(f_i)_{ij} \to (f_j)_{ij}$ sur chaque $U_{ij}$ vérifiant la condition de cocycle sur $U_{ijk}$, $\phi_{kj}\phi_{ji}=\phi_{ki}$.

Soit $x\in C$ et $x_i$ ses restrictions aux $C_i$ ; on cherche à définir un objet $f(x)$ de $D$ qui serait son image par le recollement des $f_i$. Le cocycle $\phi_{ij}$ évalué en $x$ donne une condition de recollement pour les objets $f_i(x_i)$ et on peut recoller les $f_i$ seulement si on peut recoller les $f_i(x_i)$ en un objet de $D$ ; un tel objet est un $D$-torseur et n'est pas forcément représentable dans $D$.

En revanche si on sait que le champ associé à $D$ est essentiellement formé des $D$-torseurs, le recollement des $f_i(x_i)$ existera et donc le recollement des $f_i$ en un morphisme de $C$ vers le champ associé à $D$.

\section{Modules des catégories abéliennes}\label{modcatab}

On définit dans cette section le champ classifiant les catégories abéliennes intégrant la théorie des déformations infinitésimales de \cite{lvdb1, lvdb2}.


Pour chaque $A\in\kcom$, on rappelle l'inclusion (non pleine) de catégories $A\z\AB_\UU\to A\z\CAT_\VV$.
Le lemme {\ref{restricabel}} assure que ces sous-catégories sont compatibles avec les foncteurs de restrictions $(-)^u$.
En conséquence on peut définir le sous-préchamp faible non plein de $\tkCATeqd_\VV$
\begin{eqnarray*}
\tkAB : \kaff^o &\longrightarrow & \catcat_\VV\\
A &\longmapsto & A\z\AB\\
u:A\to B &\longmapsto & (-)^u:A\z\AB\to B\z\AB
\end{eqnarray*}

Les valeurs de ce préfaisceau ne sont pas des sous-catégories de modèles, et ce n'est donc pas un sous-préfaisceau faible de Quillen de $\tkCATeqd_\VV$. Néanmoins, $A\z\AB_\UU$ est stable par équivalences dans $A\z\CAT_\VV$ (toute catégorie équivalente à une catégorie abélienne est abélienne).

On note $\kAB$ le strictifié de $\tkAB$, il reste stable par équivalences dans $\kCATeqd$ et, en conséquence, son nerf $|\kAB|$ est plein dans $|\kCATeqd|$, ainsi que le champ simplicial associé $\ukAB$ dans $\ukCATeqd$.

Si $A\z\ab$ désigne la sous-catégorie pleine de $A\z\AB$ formée de l'image essentielle de $A\z\cat$ par $\Mod^A$, le lemme \ref{retricabel} assure que les $A\z\ab$ forment un sous-préfaisceau plein de $\tkAB$ qu'on note $\tkab$.
On note $\kab$ son strictifié ; $|\kab|$ son nerf et $\ukab$ le champ simplicial associé.
$|\kab|$ et $\ukab$ sont pleins dans $|\kAB|$ son nerf et $\ukAB$ respectivement, et aussi dans $|\kCATeqd|$ et $\ukCATeqd$.

%
%
%
%

\bigskip

On rappelle, pour $A\in \kcom$, les foncteurs $\overline{(-)}:A\z\CAT\to A\z\AB$ (cf. {\sf \S \ref{abel}}), que pour les besoins de la section on renomme $Mod_\UU^A$.

\begin{lemme}\label{mormod}
Les foncteurs $Mod_\UU^A$ définissent un morphisme $\MMod_\UU : \tkCATeqg\to \tkAB$\index{$\MMod^d_\UU$} et, en conséquence, un morphisme de préchamps :
$$
\mathcal{M}od^d_\UU : \kCATeqg\longrightarrow \kAB
$$
et de champs 
$$
\mathcal{M}od^d_\UU : \ukCATeqg\longrightarrow \ukAB.
$$
Ces morphismes sont essentiellement surjectifs.
\end{lemme}
\begin{pr}
Il suffit de vérifier que $\MMod^d_\UU$ est définit comme morphisme $\tkCATeqg\to \tkCATeqd$ et il suffit pour cela que les $Mod^d_\UU(A)$ commutent aux restrictions : c'est une application du lemme {\ref{monoidal}}.
La surjectivité est une conséquence de celle de $Mod^A_\UU:A\z\CAT\to A\z\AB$ qui induit un foncteur surjectif entre les préchamps classifiant les équivalences, et donc entre les champs associés.
\end{pr}

\subsection{Diagonale}\label{diagocatab}

\begin{prop}\label{diagcatab}
Soient $C$ et $D$ deux $A$-algèbres associatives projectives de type fini. Le préchamp $\Omega_{\overline{C},\overline{D}}|\kAB|$ est équivalent au champ $\underline{\mathcal{I}nv}(C,D)$ classifiant les $C\z D)$-bimodules inversibles. 
\end{prop}
\begin{pr}
On rappelle que $\mathcal{I}nv(C,D)$ est définit au {\sf\S\ref{modulesinversibles}} où la définition utilise des conditions de finitude sur les catégories mais elle se généralise.

On a
\begin{eqnarray*}
\Omega_{\overline{C},\overline{D}}|\kab| : \aff_A &\longrightarrow & \sens\\
u:A\to B &\longmapsto & \Omega_{\overline{C}^u,\overline{D}^u}|B\z\CAT|\\
A\overset{u}{\to} B\overset{v}{\to} B' &\longmapsto & \omega_v : \Omega_{\overline{C}^u,\overline{D}^u}|B\z\CAT|\to \Omega_{\overline{C}^{vu},\overline{D}^{vu}}|B'\z\CAT|
\end{eqnarray*}
et
\begin{eqnarray*}
\underline{\mathcal{I}nv(C,D)} : \aff_A &\longrightarrow & \sens\\
u:A\to B &\longmapsto & \textrm{Inv}\left(C_u,D_u\right)\\
A\overset{u}{\to} B\overset{v}{\to} B' &\longmapsto & \textrm{Inv}_v : \textrm{Inv}\left(C_u,D_u\right)\to \textrm{Inv}\left(C_{vu},D_{vu}\right).
\end{eqnarray*}
où $\textrm{Inv}\left((C,D)_u\right)$ est le nerf du groupoïde des isomorphismes de $C\z D$-bimodules inversibles.
Il faut donc montrer que les valeurs de ces préchamps sont équivalentes et que les restrictions sont compatibles avec ces équivalences.
Soient $C$ et $D$ deux $A$-algèbres associatives projectives de type fini et $A\overset{u}{\to} B\overset{v}{\to} B'\in\kcom$. On rappelle les isomorphismes canoniques $\overline{C_u}\simeq \overline{C}^u$.
 
Comme $|W_{B\z\ab}|$ est plein dans $|W_{B\z\CAT_\VV}|$, on a le diagramme commutatif :
$$\xymatrix{
\Omega_{\overline{C_u},\overline{D_u}}|W_{B\z\ab}| \ar[d]\ar[r]^-{\sim} & \Omega_{\overline{C_u},\overline{D_u}}|W_{B\z\CAT_\VV}| \ar[d]\\
\Omega_{\overline{C_{vu}},\overline{D_{vu}}}|W_{B\z\ab}| \ar[r]^-{\sim}& \Omega_{\overline{C_{vu}},\overline{D_{vu}}}|W_{B\z\CAT_\VV}|
}$$
où les flèches horizontales sont des équivalences.
Le lemme {\ref{comparclassifmodele}} fournit un diagramme commutatif :
$$\xymatrix{
\Omega_{\overline{C_u},\overline{D_u}}|W_{B\z\CAT_\VV}| \ar[r]^-{\sim}\ar[d] &
\Omega_{\overline{C_u},\overline{D_u}}|\mathcal{G}(B\z\CAT_\VV)| \ar[d]& 
\ar[l]_-{\sim}\Eq^\Delta_{B\z\CAT_\VV}(\overline{C_u},\overline{D_u})\ar[d]\\
\Omega_{\overline{C_{vu}},\overline{D_{vu}}}|W_{B'\z\CAT_\VV}| \ar[r]^-{\sim} &
\Omega_{\overline{C_{vu}},\overline{D_{vu}}}|\mathcal{G}(B'\z\CAT_\VV)| & 
\ar[l]_-{\sim}\Eq^\Delta_{B'\z\CAT_\VV}(\overline{C_{vu}},\overline{D_{vu}})
}$$
où les flèches horizontales sont des équivalences.
La structure simpliciale de $B\z\CAT_\VV$ donne l'isomorphisme :
$$
\Eq^\Delta_{B\z\CAT_\VV}(\overline{C_u},\overline{D_u}) \simeq |\uEq_{B\z\CAT_\VV}(\overline{C_u},\overline{D_u})|
$$
et le lemme {\ref{restequivmor}} donne le diagramme de groupoïdes :
$$\xymatrix{
\uEq_{B\z\CAT_\VV}(\overline{C_u},\overline{D_u}) \ar[d]_{(-)^v}&\ar[l]_-{\sim}  \textrm{Inv}((C_u,D_u)) \ar[d]^{(-)_v}\\
\uEq_{B'\z\CAT_\VV}(\overline{C_{vu}},\overline{D_{vu}}) &\ar[l]_-{\sim}  \textrm{Inv}(C_{vu},D_{vu})) 
}$$
où les flèches horizontales sont des équivalences.

Mis bouts à bouts, ces diagrammes fournissent l'équivalence voulue entre les préchamps.
\end{pr}

On en déduit immédiatement la proposition suivante.

\begin{prop} Soient $A\in\kcom$ et $C\in A\z\CAT$.
\begin{itemize}
\item Le préfaisceau $\hat{\pi}_0(|\kAB|)$ classifie les classes d'équivalence de catégories abéliennes.
\item Le préfaisceau $\hat{\pi}_1(|\kAB|,\overline{C})$ classifie les classes d'isomorphie de $C\z C$-bimodules inversibles.
\item Le préfaisceau $\hat{\pi}_2(|\kAB|,\overline{C})$ classifie les automorphismes de $C$ vu comme $C\z C$-bimodule.
\item Pour $n>2$, les préfaisceaux $\hat{\pi}_n(|\kAB|,\overline{C})$ sont triviaux.
\end{itemize}
\end{prop}

Concernant les faisceaux associés, on a la même remarque qu'au {\sf\S\ref{diagocateq}}.

\subsection{Champification}\label{champifabel}

Comme pour les problèmes de modules précédents, les préchamps $\kAB$ et $\kab$ ne sont pas des champs et il faut les champifier. Mais, différence notable, les préchamps de morphismes sont, ici, des champs (ce qui, en un sens qu'on ne précise pas, correspond au fait que les objets classifiés soient eux-mêmes des champs) et la seule obstruction à ce que $\kAB$ et $\kab$ ne soient pas des champs est leur saturation pour les formes tordues de leurs points.

Tout cela sera étudié en détails dans le futur travail sur la champification \cite{anel1}.

\section{Modules à équivalence de Morita près}\label{modcatmor}




Soit $A\in\kcom$, on rappelle que $A\z\CATmor$ est la catégorie de modèles des catégories linéaires à équivalence de Morita près, c'est la localisation gauche de Bousfield de $A\z\CATeq$ pour laquelle les objets locaux sont les catégories karoubiennes. Un remplacement fibrant dans $A\z\CATmor$ est donné par le foncteur $\widehat{(-)}$ de karoubianisation, c'est un adjoint à gauche (cf. proposition {\ref{univcat}}).

On définit le préchamp faible suivant :
\begin{eqnarray*}
\tkCATmor : \kaff^o & \longrightarrow & \catcat\\
	A & \longmapsto & A\z\CATmor \\
	u:A\to B & \longmapsto & (-)_u : A\z\CATmor\to B\z\CATmor.
\end{eqnarray*}
qui, sauf pour les structures de modèles sur les catégorie points, est identique à $\tkCATeqg$.
Le lemme {\ref{restricmor}} assure que les foncteurs $(-)_u$ sont de Quillen à gauche, et $\tkCATmor$ est donc un préfaisceau faible de Quillen à gauche.

On note $\kCATmor$ son strictifié, $|\kCATmor|$ le préchamp simplicial nerf des équivalences et $\ukCATmor$ le champ simplicial associé.

$\tkcateq$ est un sous-préchamp de $\tkCATmor$ mais il n'est plus de Quillen, car il n'est plu stable par équivalence : la karoubianisation d'une catégorie ayant un nombre fini d'objet a, a priori, un nombre infini d'objets. On définit $\tkcatmor$ comme le plus petit sous-préfaisceau faible de Quillen de $\tkCATmor$ contenant $\tkcateq$. On note $\ukcatmor$ le champ simplicial associé.

\bigskip

On a un morphisme évident $\mathcal{K}:\tkCATeq \longrightarrow\tkCATmor$ qui induit entre les champs associé un morphisme, dit de localisation :
$$
\mathcal{K}:\ukCATeq \longrightarrow \ukCATmor.
$$
De même on récupère un morphisme
$$
\mathcal{K}:\ukcateq \longrightarrow \ukcatmor.
$$

Par construction de la structure de modèles de Morita, le morphisme $\mathcal{K}^A:A\z\CATeq \to A\z\CATmor$ factorise $A\z\CATeq\to A\z\AB$ et, en conséquence, on a une factorisation de morphismes de préfaisceaux de Quillen :
$$\xymatrix{
\kCATeqg \ar[d]\ar[r]^-{\MMod}&\kAB\\
\kCATmor\ar[ru]
}\quad\textrm{ et }\quad
\xymatrix{
\kcateqg \ar[d]\ar[r]^-{\MMod}&\kab\\
\kcatmor\ar[ru]
}$$
Pour le dernier triangle la flèche $\ukcatmor\longrightarrow\ukab$ existe si $C$ n'est que Morita-équivalente à une catégorie de type fini, la catégorie de ses modules est, elle, équivalente à une catégorie engendrée par une catégorie de type fini.


\begin{thm}\label{thmequivmodulemoritaabel}
Les morphismes $\ukCATmor\longrightarrow \ukAB$ et $\ukcatmor\longrightarrow\ukab$ déduit des morphismes précédents sont des équivalences de champs simpliciaux.
\end{thm}
\begin{pr}
On montre que cette équivalence existe déjà au niveau des préchamps simpliciaux.
Pour $Spec(A)\in\kaff$, $|\ukCATmor|(A)=|Mor(A)|$.
Il s'agit de prouver que le morphisme naturel $|Mor(A)|\to |W_{A\z\AB}|$ est une équivalence. On utilise pour cela le diagramme commutatif
$$\xymatrix{
|Mor(A)|\ar[r]^{|Mod^A|} & |W_{A\z\AB}| \\
|W_{A\z\KAR}| \ar[u]^{|\kappa|} \ar[ur]_{|Mod^A\circ \kappa|}&
}$$
$W_{A\z\KAR}$ consiste en les équivalences de $A\z\CATmor$ entre les seuls objets fibrants, qui se réduisent aux équivalences de $A\z\CATeq$ car ces fibrants sont aussi les objets locaux, le lemme {\ref{classifmodele}}) dit alors que $|\kappa|$, qui est l'inclusion naturelle des objets fibrants, est une équivalence.
$|Mod^A\circ \kappa|$ est une équivalence par le lemme {\ref{eqmorkar}}. On conclut en applicant le principe de trois-pour-deux des équivalences dans $\sens$.

Le raisonnement est exactement le même avec les conditions de finitude.
\end{pr}

\begin{rem} 
Les morphismes $\kCATmor\longrightarrow \kAB$ et $\kcatmor\longrightarrow\kab$ ne sont pas des équivalences de préfaisceaux de Quillen, essentiellement parce que les catégories abéliennes possèdent plus de morphismes que les linéaires (au sens où $\MMod$ n'est pas essentiellement surjectif). Néanmoins, si on se restreint aux seules équivalences comme morphismes, les morphismes deviennent des équivalences.
\end{rem}

Pour $A\in\kcom$ on définit $A\z\AB^{proj}$ la sous-catégorie de $A\z\AB$ ayant pour morphismes les seuls foncteurs préservant les objets projectifs de présentation finie. Les restrictions $(-)_u$ conservent ces conditions et ces catégories forment un sous-préchamp $\kAB^{proj}$ (non plein mais essentiel) de $\kAB$.

Par définition de $\kAB^{proj}$, le morphisme $\kCATmor\longrightarrow\kAB$ se factorise par $\kAB^{proj}$.
Si on considère qu'il doit exister une notion de champ en catégories supérieures et qu'on note $\overline{\CAT}^{Mor}$ et $\overline{\AB}^{proj}$ les champs associés aux préchamps $\kCATmor$ et $\kAB^{proj}$, on peut faire la conjecture suivante.

\begin{conj}
$\kCATmor\longrightarrow \kAB^{proj}$ est une équivalence de champs en catégories supérieures.
\end{conj}

\section{Morphismes}\label{morphismes}

Au final, on a le diagramme commutatif suivant entre les champs précédemment définis. 
(On a le même diagramme entre les champs "majuscules", mais il ne servira pas dans la suite.)
$$\xymatrix{
\underline{\mathcal{A}f\!f}_\UU^o\simeq \ukcom_\UU \ar[r]^-c_-{\textrm{{\tiny fermé}}}\ar@/_5pc/[rrrdd]_{QCoh} & \ukass_\UU \ar[d]_{\B^1}^-{\textrm{{\tiny surj.}}} \ar[r]^-\iota_-{\textrm{{\tiny ouvert}}} & \ukcatiso_\UU \ar[d]_\B^-{\textrm{{\tiny surj.}}}  \\
& \ukcat_{*,\UU} \ar[rd]_-{\textrm{{\tiny surj.}}} \ar[r]^-\iota_-{\textrm{{\tiny ouvert}}} & \ukcateqg_\UU \ar[d]_-{\mathcal{K}}^-{\textrm{{\tiny surj.}}} \ar[r]^{\MMod} \ar[rd]^{\MMod} & \ukcateqd_\VV \\
&  & \ukcatmor_\UU  \ar[r]^{\sim} & \ukab_\UU \ar[u]_-{\textrm{{\tiny plein}}} \ar@/^4pc/[llluu]_Z \\
}$$

Le morphisme $QCoh$ est la composition $\MMod\circ\B\circ\iota\circ c$ et le morphisme $Z$ est construit comme suit.

Si $A\in\kcom$ on dispose d'un foncteur $Z^A: A\z\AB \to A\z\COM$ qui à une catégorie $A$-abélienne associé son centre qui est une $A$-algèbre commutative. Il est compatible avec les équivalences naturelles des deux catégories car une équivalence entre deux catégories abéliennes induit un isomorphisme sur leur centre.
Pour qu'il définisse un morphisme de champs, il reste à vérifier sa commutation avec les changements de bases, or, si $u:A\to B\in\kcom$ et $C\in A\z\COM$, on a bien 
$$
Z^A(\overline{C}^u) = Z^A(\overline{C_u})= C_u
$$
(où les égalités sont des isomorphismes canoniques) car le centre d'une catégorie de module sur un anneau commutatif est cet anneau.
On a un isomorphisme canonique de foncteurs $Z^A\circ Mod^A\simeq id$ qui en induit entre les morphismes de champs associés : $Z\circ QCoh\simeq id$.

\begin{prop}\label{abcom}
Si $\ukAB^{com}$ désigne le sous-champ plein de $\ukAB$ image de $QCoh$, le morphisme 
$$
Z: \ukAB^{com}\longrightarrow \ukCOM
$$
est une gerbe ; sa fibre en $C: Spec(A)\to \ukCOM$ est $\mathcal{K}(\underline{Pic}(C),1)=[pt/\underline{Pic}(C)]$ où $\underline{Pic}(C)$ est le champ en groupes\footnote{Pour la notion de champ en groupes, ainsi nommé pour qu'elle soit intuitive, on renvoie à celle de $H_\infty$-champ de \cite[\S 1.4]{chaff}.} des $C$-modules de rang 1 à isomorphisme près ; et où $[pt/\underline{Pic}(C)]$\index{$\mathcal{K}(-,1)$}\index{$[pt/\underline{Pic}(C)]$} désigne le champ quotient du "groupoïde"\footnote{En fait un groupoïde de Segal, notion pour laquelle on renvoie à \cite[\S 1.3.4]{hag2}.} $\underline{Pic}(C)\rightrightarrows pt$.
%
\end{prop}
\begin{pr}
Soit $C: Spec(A)\to \ukCOM$, on note $F$ la fibre de $Z$ en $C$, pour laquelle on suppose qu'il existe un point global $c:Spec(A)\to F$, et on considère la longue suite exacte d'homotopie
\begin{multline*}
0\to \pi_2(F,c) \overset{a}{\longrightarrow} \pi_2(\ukAB^{com},C)\overset{}{\longrightarrow} \pi_2(\ukCOM,C)\overset{}{\longrightarrow} \\
\pi_1(F,c) \overset{}{\longrightarrow} \pi_1(\ukAB^{com},C)\overset{Z_1}{\longrightarrow} \pi_1(\ukCOM,C)\overset{}{\longrightarrow} \\
\pi_0(F) \overset{}{\longrightarrow} \pi_0(\ukAB^{com})\overset{Z_0}{\longrightarrow} \pi_0(\ukCOM).
\end{multline*}

De la relation $Z\circ Qcoh\simeq id$, on déduit que $Z_0$ et $Z_1$ sont surjectifs
Le morphisme $Z_0$ est aussi injectif car deux anneaux commutatifs non-isomorphes ont des catégories de modules non équivalentes (deux anneaux commutatifs Morita-équivalents sont isomorphes).
De plus, $\pi_2(\ukCOM)=0$ car $\ukCOM$ est un 1-champ.

On déduit de tout cela que $F$ est connexe, et que $\Omega_cF$ est équivalent au noyau du morphisme de champs en groupes $\Omega Z:\Omega_C\ukAB^{com}\to \Omega_C\ukCOM$.
$\Omega_C\ukAB^{com}$ est équivalent au champ $\underline{\mathcal{I}nv}(C,C)$ des $C\z C$-bimodules inversibles (cf. proposition {\ref{geominv}}) et $\Omega_C\ukCOM$ est équivalent au schéma des automorphismes de $C$ :
\begin{eqnarray*}
\underline{\mathcal{A}ut}(C) : \kaff_A^o&\longrightarrow & \ens\\
A\to B&\longmapsto & Aut_{B\z\COM}(C\otimes_AB).
\end{eqnarray*}
On décrit le morphisme induit par $Z$ entre ces deux champs.
Soit $A\to B\in \kcom$ et $M$ un $C\otimes_AB\z C\otimes_AB$-bimodule, on note $x$ le générateur canonique de $\overline{C}$ alors $M$, vu comme morphisme, associe à $x$, $M$ vu comme $C\otimes_AB$-module sur sa structure gauche. On en tire un isomorphisme $\alpha : C\otimes_AB= Z(\End(x))\to Z(\End(M))$. De l'autre côté $Z(\End(M))$ est canoniquement identifié à $C\otimes_AB=Z(\End(x))$ en utilisant la structure droite de $M$ ; on note $\beta$ cet isomorphisme.
On a $C\otimes_AB\overset{\alpha}{\to}Z(\End(M))\overset{\beta}{\leftarrow}C\otimes_AB$ et l'isomorphisme associé par $Z$ à $M$ est $\beta^{-1}\circ \alpha$. Cet isomorphisme vaut l'identité si et seulement si $\alpha=\beta$ ce qui revient à ce que les structures droite et gauche de $M$ coïncident.

On définit le sous-champ plein de $\underline{Pic}(C)$ de $\underline{\mathcal{I}nv}(C,C)$ formé des bimodules inversibles dont les structures droite et gauche coïncident, c'est le noyau de $\Omega Z$.
\end{pr}

\bigskip


La fibre de ${\B^1}$ est calculée par la proposition suivante.
\begin{rprop}[{\cite[\S 5.3]{hagdag}}]
Si on a $C\in A\z\ASS$, la fibre de $\B^1$ en $\B C$ est la gerbe pointée $\mathcal{K}(C^*,1)$ définit par
\begin{eqnarray*}
\mathcal{K}(C^*,1) : \kaff^o &\longrightarrow & \sens\\
A\to B &\longmapsto & K((C\otimes_AB)^*,1)
\end{eqnarray*}
où $K((C\otimes_AB)^*,1)$ est un espace classifiant pour le groupe des éléments inversibles de $C\otimes_AB$.
\end{rprop}

\setcounter{tot}{0}
\newpage
\thispagestyle{empty}
\chapter{Géométricité}\label{chapgeom}
\thispagestyle{chheadings}

Ce chapitre reprend, dans ses trois premières sections, les champs du chapitre précédent et prouve leur caractère géométrique. Pour chacun d'entre eux, il est aussi explicité une présentation par un 2-groupoïde (cf. {\sf\S\ref{presentation}}) dont on se sert pour calculer un modèle pour les complexes tangents de ces champs.

La section {\ref{comparaison}} établit que le morphisme naturel $\ukcateq\longrightarrow\ukab$ est étale.

La dernière section, enfin, consiste en un gros diagramme commutatif qui résume tous les résultats.


\section{Modules à isomorphisme près}\label{ch4iso}

On rappelle le champ $\ukcatiso$ sous-champ plein du champ $\ukCATiso$ généré par le préchamp $\kcatiso$ de $\kCATiso$ formé des catégories linéaires dont le nombre d'objets est fini et dont les modules de morphismes sont projectifs de type fini.

On prouve dans cette section que $\ukcatiso$ est un 1-champ 1-géométrique, puis on étudie son tangent.

\subsection{Géométricité}

On rappelle le champ $\underline{Gr}=\coprod_{n\in\mathbb{N}} \underline{Gr}_n$ classifiant les graphes linéaires projectifs de type fini ; les $\underline{Gr}_n\simeq \uvect^{(n)}/\Sn$ classifient les graphes linéaires à $n$ objets (cf. {\sf\S\ref{graphes}}).
L'oubli de la loi d'algèbre fournit un morphisme $\ukcatiso\longrightarrow \underline{Gr}$ et on définit
$\ukcatison$ comme la fibre au-dessus de $\underline{Gr}_n$, il classifie les catégories à $n$ objets.

Avec les notations du {\sf\S\ref{graphes}}, on définit les champs $\ukcatisons$ et $\ukcatisod$ par les deux carrés homotopiquement cartésiens dans $\chaffk$ suivants :
$$\xymatrix{
\ukcatisod \ar@{}[rd]|{\square}\ar[r]\ar[d]\ar@/_3pc/[dd]_{c_d} & \uvectd \ar[d]^{u_d}\ar@/^2pc/[dd]^{v_d} \\
\ukcatisons \ar@{}[rd]|{\square}\ar[r]\ar[d] & \uvect^{(n)} \ar[d]^u \\
\ukcatison \ar[r] & \uvect^{(n)}/\Sn \\
}$$
les morphismes $v_d$, et donc les $c_d$, sont 0-représentables (car $u_d$ est ouvert et $u$ est de fibre $\Sn$).

\begin{lemme}\label{etale}
Ces morphismes sont étales. 
\end{lemme}
\begin{preuve}
$u_d$ est étale car ouvert et $u$ est étale car un c'est un $\Sn$-torseur et $\Sn$ est un groupe étale sur $pt$.
\end{preuve}

\medskip

$\ukcatisons$ est le champ classifiant les catégories linéaires à isomorphisme fixant les objets près, et $\ukcatisod$ en est le sous-champ plein des catégories linéaires de type $d$. La géométricité de $\ukcatison$ est équivalente à celle de tous les $\ukcatisod$.

Sa définition comme produit de $\ubgln$ assure que $\uvectd$ admet $pt$ comme carte ; on l'utilise pour définir le champ $\kcatisod$ des trivialisations de $\ukcatisod$ par le produit fibré homotopique dans $\chaffk$ :
\begin{eqnarray*}
\xymatrix{
\kcatisod \ar@{}[rd]|{\square}\ar[d]_c \ar[r] & pt\ar[d] \\
\ukcatisod \ar[r] & \uvectd \\
}
\end{eqnarray*}

\begin{prop}\label{schcat}
$\kcatisod$ est équivalent au schéma $\kcat^{(d)}$ classifiant les structures de catégorie sur un graphe libre de rang $d$.
\end{prop}

\begin{pr}
On rappelle que $\kcat^{(d)}$ est détaillé à la proposition {\ref{structurecat}}.
On explicite $\kcatisod$ par le produit fibré :
$$\xymatrix{
\kcatisod \ar@{}[rd]|{\square}\ar[d]_c \ar[r] & \left({\uvectd}\right)^{\Delta^1}\ar[d] \\
pt \times\ukcatisod \ar[r] & \uvectd \times \uvectd \\
}$$
Les foncteurs de champification et de strictification commutent aux produits fibrés et on se ramène au calcul du produit fibré au niveau des préfaisceaux faibles.

Si $|G\ell_d|$ désigne le nerf du groupoïde à un objet tiré du groupe $G\ell_d$ et si, pour une algèbre commutative $A$, $|A\z\catiso|_n$ est l'ensemble des $n$-simplexes du nerf des isomorphismes de la catégorie $A\z\cat$, on trouve que les objets de $\kcatisod(A)_n$ sont les couples :
$$
(\alpha:\Delta^1\times\Delta^n\to |G\ell_d|,\ C\in |A\z\catiso|_n)
$$
tels que $\alpha_0:0\times\Delta^n\to |G\ell_d|$ soit le graphe $L$ libre de type $d$ et $\alpha_1:1\times\Delta^n\to |G\ell_d|$ soit l'image de $C$.

En notant à l'identique les objets de $C$ et leurs graphes, on peut représenter $\alpha$ par le diagramme
$$\xymatrix{
&& L\ar[lld]\ar[ld]\ar[rd] & \\
C_0\ar[r] & C_1\ar[r]& \dots\ar[r]& C_n \\
}$$
où toutes les flèches sont des isomorphismes et tous les triangles sont commutatifs.

Si $D$ est une catégorie, on appelle un isomorphisme de graphes $L\to D$ une trivialisation de $D$ et on dit que $D$ est une catégorie trivialisée. Un triangle comme ceux du diagramme ci-dessus définit une notion d'isomorphisme de catégorie trivialisées. Ainsi $\kcatisod(A)$ est simplement le nerf du groupoïde des isomorphismes de catégories trivialisées de $A\z\catiso$.

Chaque trivialisation donne une structure de catégorie sur $L$ qui s'exprime dans la base canonique par des constantes de structures et deux catégories trivialisées sont isomorphes si et seulement si ces structures sont identiques. En d'autres termes, on a une bijection entre les composantes connexes de $\kcatisod(A)$ et les $A$-points de $\kcat^{(d)}$. L'équivalence entre les deux provient de ce que les catégories trivialisées n'ont pas d'automorphismes :
$$\xymatrix{
L\ar[d]_\alpha\ar[rd]^\alpha \\
C\ar[r]_\beta & C \\
}$$
commute si et seulement si $\beta=\alpha\alpha^{-1}=id$.

Le champ $\kcatisod$ est ainsi le champ associé à un préfaisceau faible en groupoïdes, équivalent au
foncteur des points du schéma $\kcat^{(d)}$, il est donc équivalent à ce schéma.
\end{pr}

\begin{cor}\label{catisod}
$\ukcatisod$ est 1-géométrique et $\left(G\ell_d\times \kcat^{(d)} \rightrightarrows \kcat^{(d)}\right)$ en est une présentation par un groupoïde affine.
\end{cor}

\begin{pr}
Par construction $\kcatisod\simeq\kcat^{(d)}\to\ukcatisod$ est un $G\ell_d$-torseur comme tiré en arrière d'un tel et les quotients de groupoïdes étant effectifs dans la catégorie des champs, le groupoïde affine :
$$
G\ell_d\times \kcat^{(d)} \rightrightarrows \kcat^{(d)}
$$
est un modèle local pour $\ukcatisod$.

$\ukcatisod$ est 1-géométrique si le groupoïde est lisse ({\sffamily \S \ref{gpdlisse}}).
Or comme tout pull-back d'un morphisme lisse est lisse : $\kcat^{(d)}\to\ukcatisod$ est lisse et également les flèches 
$$
G\ell_d\times \kcat^{(d)} = \kcat^{(d)}\times_{\ukcatisod} \kcat^{(d)} \rightrightarrows  \kcat^{(d)}.
$$
\end{pr}

\begin{thm}\label{geomiso}
$\ukcatison$ est 0-géométrique localement de présentation finie et une présentation par un groupoïde est donné par :
$$
\Sn\times \coprod_d \left(G\ell_d\times \kcat^{(d)}\right) \rightrightarrows \coprod_d\kcat^{(d)}
$$
(où les réunions sont prises sur toutes les applications $d:E\to\mathbb{N}$ avec $E$ un ensemble fixé de cardinal $n$)
\end{thm}

\begin{pr}
$\coprod_d (G\ell_d\times \kcat^{(d)}) \rightrightarrows \coprod_d\kcat^{(d)}$ définit un modèle local
pour $\ukcatisons$, pour obtenir $\ukcatison$ il suffit d'y rajouter l'action de $\Sn$.
Soit $C^d$ un $A$-point de $\kcat^{(d)}$, un élément $g_d\in G\ell_d$ agit par changement de base dans les $C^d(x,y)$ ; on note $D^d$ le nouveau point de $\kcat^{(d)}$. 
Un élément $\sigma\in\Sn$ agit sur $C^d$ par permutation des objets, ce qui revient, en fixant les objets, à agir sur $d$ ; on note $C^{d\circ \sigma^{-1}}$ l'image de $C$ par $\sigma$. Ces deux actions commutent.
Avec ces notations, l'action du groupoïde est donné sur les $A$-points par : $(\sigma,g^d,C^d) \mapsto D^{d\circ \sigma^{-1}}$.
\end{pr}

\subsection{Tangent}

On renvoie à l'annexe {\ref{hoch}} pour la définition du complexe de Hochschild d'une catégorie.

\begin{thm}\label{tangentiso}
Le complexe tangent au champ $\ukcatiso$ en un point correspondant à une catégorie $A$-linéaire $C$ est donné par un double tronqué du complexe de Hochschild de $C$ :
$$
Der^{\leq1}(C) := HC^1(C)\longrightarrow HZ^2(C)
$$\index{$Der^{\leq1}$}
où $HZ^2(C)$ est en degré 0.
En particulier, le tangent géométrique s'identifie à $HH^2(C)$ et le tangent aux morphismes est $HZ^1(C)$, \ie le module des dérivations de $C$.\footnote{La notation $Der^{\leq1}(C)$ est choisie pour évoquer que ce complexe est naturellement un tronqué-décalé du complexe des dérivations dérivées de $C$.}
\end{thm}
\begin{pr}
La décomposition $\ukcatiso=\coprod_n\ukcatison$ ramène le calcul du tangent à la détermination du tangent à $\ukcatison$. Les applications $c_d:\ukcatisod\to\ukcatison$ étant étales (lemme {\sf\ref{etale}}), le tangent à $\ukcatison$ en un point sera isomorphe à celui de $\ukcatisod$ en un quelconque relevé du point.

Le corollaire {\ref{catisod}} caractérise $\ukcatisod$ comme le champ quotient du groupoïde affine lisse
$$
s,b:G\ell_d\times \kcat^{(d)} \rightrightarrows \kcat^{(d)}.
$$
On utilise la proposition {\ref{tangpd}} pour calculer le complexe tangent en un point $x:X=Spec(A) \to \kcat^{(d)}\to \ukcatisod$.
Comme $G\ell_d\times \kcat^{(d)}\times_{\kcat^{(d)}}X =  G\ell_d\times X$, le complexe tangent est quasi-isomorphe à 
$$
\textrm{T}_{G\ell_d\times X,x}\overset{b^*}{\longrightarrow} \textrm{T}_{\kcat^{(d)},C}.
$$
où $C$ est la catégorie correspondant au point $X\to \kcat^{(d)}$. 
$\textrm{T}_{\kcat^{(d)},C}$ est le $A$-module des déformations à l'ordre 1 de la loi de $C$, \ie $HZ^2(C)$ et $\textrm{T}_{G\ell_d\times X,x}$ est l'algèbre de Lie $g\ell_d=\oplus_{x,y}g\ell_{d(x,y)}\simeq \oplus_{x,y} \Hom_A(C(x,y),C(x,y)) = HC^1(C)$.

Pour définir la différentielle, quelques détails sur l'action de $G\ell_d$ sont nécessaires.
On note $m(\z,\z)$ la multiplication de $C$. Symboliquement, si $\{g_{xy}\}\in G\ell_d(A)$ agit sur les $C(x,y)$ par $f_{xy}\mapsto g_{xy}f_{xy}$, il agit sur :
$$
HC^2(C) = \bigoplus_{(x,y,z)\in E} \Hom_A(C(x,y)\otimes_A C(y,z),C(x,z))
$$
par $M_{xyz}\mapsto g_{xz}M_{xyz}(g_{xy}^{-1},g_{yz}^{-1})$.
L'action dérivée de $\alpha_{xy}\in g\ell_{d(x,y)}$ est donc donnée par
$$
M_{xyz}\mapsto \alpha_{xz} M_{xyz}(\z,\z)-M_{xyz}(\alpha_{xy},\z)-M_{xyz}(\z,\alpha_{yz})),
$$
qu'on reconnaît pour être la différentielle de Hochschild.
\end{pr}


\subsection{Champs des algèbres associatives et commutatives}

On rappelle le champ $\ukass$ classifiant les algèbres associatives à isomorphisme près (cf. {\sf\S\ref{chpass}}), comme c'est un sous-champ plein ouvert de $\ukcatiso$ il est 0-géométrique (corollaire {\ref{ouvfergeom}}) et on retrouve le résultat connu que son tangent en une $A$-algèbre associative $C$ est donné par un tronqué de son complexe de Hochschild.

Dans le cas où $d$ est un type de graphe ayant un seul objet, $d$ est essentiellement un nombre entier ; par souci de clarté on note $\ass^d$ le schéma classifiant les structures de catégorie à un objet, \ie d'algèbre associative unitaire, sur un module libre de rang $d$.

Une présentation de $\ukass$ est donné par le groupoïde
$$
\coprod_{n\in\mathbb{N}} \left(G\ell_n\times \ass^n\right) \rightrightarrows \coprod_{n\in\mathbb{N}}\ass^n.
$$

Le rang $n$ est un invariant de la classe d'isomorphisme d'une algèbre ; on définit $\ukass^n$ comme le sous-champ de $\ukass^n$ dont ce rang est $n$. Il se présente par le groupoïde
$$
G\ell_n\times \ass^n \rightrightarrows \ass^n.
$$
On a la partition :
$$
\ukass = \coprod_n \ukass^n.
$$

\medskip

$\ukaff^o=\ukcom$ (cf. {\sf\S\ref{chpcom}}) étant un sous-champ fermé de $\ukass$, il est également géométrique (corollaire {\ref{ouvfergeom}}).

On déduit de la partition de $\ukass$ la partition
$$
\ukaff = \coprod_n \ukaff^n
$$
où $\ukaff^n\simeq \ukass^n\times_{\ukass}\ukcom$ est le champ des schémas de longueur $n$.

On rappelle de la proposition {\ref{schstcom}} le schéma affine $\kcomm^{(n)}$ classifiant les structures d'algèbres commutatives. Une présentation de $\ukcom^n$ est donnée par le groupoïde
$$
G\ell_n\times \kcomm^{(n)} \rightrightarrows \kcomm^{(n)}.
$$
Le tangent de $\kcomm^{(n)}$ en un $A$-point $C$ est un sous-module du tangent à $\kcat^{(n)}$ en $C$, obtenu en imposant la condition de commutativité. Précisément, comme $\textrm{T}_{\kcat^{(n)},C} = HZ^2(C)$, $\textrm{T}_{\kcomm^{(n)},C}$ est obtenu en ne considérant que les 2-cocycles symétriques en leurs variables dont on note $Z^1_{Harr}(C)$ le module.

Le complexe tangent à $\ukcom$ en $C:Spec(A)\to \ukcom$ est le complexe
$$
Harr^{\leq1}(C) := C^0_{Harr}(C)\to Z^1_{Harr}(C)
$$
où $C^0_{Harr}(C):=HC^1(C) = \End_A(C)$ et où $Z^1_{Harr}(C)$ est en degré 0.
Ce complexe est le tronqué-décalé en degré 1 du complexe de cohomologie de Harrison de $C$ \cite{gerst}.


\section{Modules à équivalence près}\label{ch4eq}

\subsection{Géométricité}

On rappelle que le champ $\ukcateq$ est le champ associé au préchamp classifiant les catégories linéaires ayant un nombre essentiel d'objets fini et des modules de morphismes projectifs de type fini.
On prouve dans ce paragraphe la géométricité de $\ukcateq$. On procède en montrant la géométricité du champ des chemins entre deux points et en explicitant une carte pour $\ukcateq$.

\subsubsection{Diagonale}

Soit $A\in\kcom$, on considère deux catégories $A$-linéaires $C$ et $D$, de graphes finis et libres, vues comme des points de $\ukcateq$. On a déjà vu au {\S\ref{diagocateq}} que $\Omega_{C,D}\ukcateq$ est le champ associé au préfaisceau en groupoïdes :
\begin{eqnarray*}
\uEq(C,D) :  \aff_{/Spec(A)} & \longrightarrow & \gpd \\
	u:A\to B & \longmapsto & \uEq_{A\z\CAT}(C_u,D_u)
\end{eqnarray*}
où $\uEq_{A\z\CAT}(C_u,D_u)$ est le groupoïde des équivalences de catégories $C_u\to D_u$ et de leurs isomorphismes naturels.
$\Omega_{C,D}\ukcateq$ est donc géométrique si le groupoïde est géométrique et lisse.

\begin{prop}
Si $C$ et $D$ sont deux catégories $A$-linéaires, de graphes finis et libres alors le champ $\Omega_{C,D}\ukcateq$ est 1-géométrique et un modèle local est donné par le groupoïde $Eq(C,D)_1\rightrightarrows Eq(C,D)_0$ (cf. prop. {\ref{gpdequiv}}).
\end{prop}
\begin{pr}
\label{repdiageq}
Comme $C$ et $D$ ont été choisies de graphes libres et finis, le groupoïde formé des équivalences et des isomorphismes naturels entre elles est classifié par le groupoïde de la proposition {\ref{gpdequiv}} :
$$
Eq(C,D)_1 \rightrightarrows Eq(C,D)_0.
$$
Comme les deux schémas de ce groupoïde sont affines, la 1-géométricité de $\Omega_{C,D}\ukcateq$ est équivalente à la lissité du morphisme source (cf. {\S\ref{groupoide}}).

Comme les schémas en jeu sont de types finis, on montre la lissité de la projection source par le critère formel.
Pour $B$ une $A$-algèbre et $B'$ une extension infinitésimale de $B$ on considère un carré commutatif
$$
\xymatrix{
Spec(B) \ar[d] \ar[r]^{(f,g,\alpha)} & Eq(C,D)_1 \ar[d]^s \\
Spec(B') \ar[r]_{F} \ar@{-->}[ru]^{(F,\overline{\alpha})} & Eq(C,D)_0
}$$
où $f$ est une équivalence $C_B\to D_B$, $\alpha$ un isomorphisme naturel $f\to g$ (on rappelle que $g$ est entièrement caractérisée par $f$ et $\alpha$) et $F$ une équivalence $C_{B'}\to D_{B'}$. La commutativité assure que $f$ est la restriction de $F$ à $B$. La lissité de $s$ est équivalente à l'existence d'un relèvement $(F,\overline{\alpha})$ de $F$ faisant commuter les deux triangles (ces deux données caractériserons un unique $G$ isomorphe à $F$ par $\overline{\alpha}$).

Comme les restrictions $(-)_B$ conservent les objets, les applications entre les objets de $F$ et $f$ coïncident. 
Pour tout $x\in C_B$, $\alpha_x$ est isomorphisme de $D_B$ de source $f(x)$ et de but un certain $y$ ; on cherche un isomorphisme $\overline{\alpha}_x\in D_{B'}(f(x),y)$ tel que sa restriction donne $\alpha_x$.
Un prolongement $\overline{\alpha}_x$ de $\alpha_x$ existe car $D_{B'}(f(x),y)$ se surjecte sur $D_B(f(x),y)=D_{B'}(f(x),y)\otimes_{B'}B$, il reste à montrer qu'il est encore un isomorphisme. Pour cela on considère la fonction $\delta_x$ de $B'$ définie par le déterminant de la composition par $\overline{\alpha}_x$, elle est inversible ssi $\overline{\alpha}_x$ est un isomorphisme. Or, $\delta_x$ est un prolongement infinitésimal de la fonction $\det(m_{\alpha_x})\in B$ qui est inversible, elle reste donc inversible.
\end{pr}

\subsubsection{Carte}

Les morphismes de catégories de modèles $\B^A:A\z\CATiso\to A\z\CATeq$ se regroupent en un morphisme de préfaisceaux faibles de Quillen et les morphismes $\B^A:A\z\catiso\to A\z\cateq$ en un morphisme de sous-préfaisceaux faibles de Quillen.

Ces morphismes sont surjectifs sur les classes d'équivalences et induisent un morphisme surjectif au niveau des champs associés.
\begin{prop}\label{lissiteB}
$\B: \ukcatiso \longrightarrow \ukcateq$ est un morphisme lisse.
\end{prop}

\begin{pr}
Soit $A$ une $\kk$-algèbre et $A'$ une extension infinitésimale au premier ordre. On considère le carré du type :
$$\xymatrix{
Spec(A) \ar[d] \ar[r]^C & \ukcatiso \ar[d]^{\B} \\
Spec(A') \ar[r]^{D'} \ar@{-->}[ru]^D &\ukcateq
}$$
commutatif à équivalence près. La lissité étant un critère local, il suffit de considérer le cas où $C$ et $D'$ sont des catégories linéaires ; on cherche alors à construire une catégorie $A'$-linéaire $D$ équivalente à $D'$ dont la restriction $A$-linéaire soit isomorphe à $C$.

Les données de commutations du diagramme sont une équivalence de catégorie $f:C\to D'\otimes_{A'}A$. On construit $D$ comme suit : pour chaque paire d'objets $x,y\in C$ on définit $D(x,y)=D(fx,fy)$ où les objets $fx$ et $fy$ sont vus dans $D$ compte tenu que $D'$ et $D'\otimes_{A'}A$ ont les mêmes objets. Comme $f$ est une équivalence, il est clair que $D\otimes_{A'}A$ est isomorphe à $C$.
\end{pr}

\subsubsection{Présentation}

\begin{thm}\label{thmgeomeq}
$\ukcateq$ est un 2-champ 2-géométrique localement de présentation fini. 
Une présentation par un 2-groupoïde affine est donné par (cf. {\S\ref{gpdequiv}}) :
$$
EQ_2 \rightrightarrows EQ_1 \rightrightarrows EQ_0.
$$
\end{thm}
\begin{pr}
$EQ_0$ forme une carte affine de $\ukcatiso$ et donc de $\ukcateq$ par composition avec $\B$, $EQ_1$ est une carte affine de $\Omega_{EQ_0}\ukcateq$ et $EQ_2 \simeq \Omega_{EQ_1}\Omega_{EQ_0}\ukcateq$.
\end{pr}

\begin{rem}\label{chpcatcateq}
 Il est assez naturel qu'il doit exister un champ en 2-catégories enveloppant $\ukcateq$ dont il ne doit être que le champ en groupoïde intérieur. Modulo la définition de tels champs, cet enveloppant devra être encore un champ qui est 2-géométrique par la présentation affine $CL_2 \rightrightarrows CL_1 \rightrightarrows CL_0$.
\end{rem}

\subsection{Tangent}

\begin{thm}\label{tangentequiv}
Soient $X=Spec(A)\in\kaff$, $C:X\to EQ_0$ et $x:X\to EQ_0\to\ukcateq$. Le complexe tangent de $\ukcateq$ en $x$ est donné par le tronqué du complexe de Hochschild :
$$
Hoch^{\leq2}\index{$Hoch^{\leq2}$} := HC^0(C) \to HC^1(C) \to HZ^2(C)
$$\index{$Hoch^{\leq2}$}
(où $HZ^2(C)$ est en degré 0).
\end{thm}
\begin{pr}
On calcule le tangent à $\ukcateq$ en $x$ en utilisant la proposition {\ref{calcultangent}}.
On considère le diagramme de carrés homotopiquement cartésiens
$$\xymatrix{
X\times_{EQ_1}EQ_2 \ar@{}[rd]|{\square}\ar[r]\ar[d] & EQ_2\ar@{}[rd]|{\square}\ar[r]\ar[d]& X\times_{EQ_0}EQ_1 \ar[d]\ar[r] \ar@{}[rd]|{\square}& EQ_1\ar[d] \\
X\ar[r]^{e\circ\sigma} & EQ_1 \ar[r]& \Omega_{X,EQ_0}\ukcateq \ar@{}[rd]|{\square}\ar[r]\ar[d] & \Omega_{EQ_0}\ukcateq \ar[r]\ar[d] \ar@{}[rd]|{\square}&EQ_0\ar[d]\\
& & X \ar@/^1pc/[u]^\sigma \ar[ul]^{e\circ\sigma} \ar[r] & EQ_0\ar[r] & \ukcateq
}$$
où la section $\sigma$ existe par la factorisation de $x$, $e\circ\sigma$ est la composition avec la flèche 'identités' $e:EQ_0\to EQ_1$.

\medskip

On déduit du diagramme par application de la proposition {\ref{calcultangent}} :
$$
\mathbb{T}_{\ukcateq,x} \simeq \textrm{hocolim}\ \left(
\mathbb{T}_{\Omega_{X,EQ_0}\ukcateq,s} \longrightarrow \textrm{T}_{EQ_0,x} \right)
$$
et
$$
\mathbb{T}_{\Omega_{X,EQ_0}\ukcateq,s} \simeq \textrm{hocolim}\ \left(
\textrm{T}_{X\times_{EQ_1}EQ_2,x} \longrightarrow \textrm{T}_{X\times_{EQ_0}EQ_1,x}\right)
$$
d'où le fait que $\mathbb{T}_{\ukcateq,x}$ est quasi-isomorphe au complexe de $A$-modules :
$$
\textrm{T}_{X\times_{EQ_1}EQ_2,x} \longrightarrow \textrm{T}_{X\times_{EQ_0}EQ_1,x}
\longrightarrow \textrm{T}_{EQ_0,x}.
$$

\medskip
Si on note $C$ la catégorie pointée par le morphisme $X\to EQ_0$, on a les interprétations suivantes des objets intervenant dans le complexe précédent :
\begin{itemize}
\item $X\times_{EQ_0}EQ_1$ classifie les équivalences de catégories de source $C$.
\item $X\times_{EQ_1}EQ_0$ classifie les isomorphismes (entres équivalences) de source $id_C$.
\end{itemize}

Le tangent au point $C$ de $EQ_0$ est exactement une déformation au premier ordre de la loi de $C$ (cf. Annexe {\ref{hoch}}) on a donc
$$
\textrm{T}_{EQ_0,x} = HZ^2(C).
$$

Le tangent au point de $X\times_{EQ_0}EQ_1$ correspondant à l'identité de $C$ est le module des équivalences $C[\epsilon]\to D$, où $D$ est une catégorie $A[\epsilon]$-linéaire qui est une déformation au premier ordre de $C$, telles que évaluées en $\epsilon=0$ elles donnent l'identité de $C$. La {\ref{hcun}} caractérise ces objets comme étant les éléments de $HC^1(C)$. D'où 
$$
\textrm{T}_{X\times_{EQ_0}EQ_1,id_C} = HC^1(C).
$$

Le tangent au point de $X\times_{EQ_1}EQ_2$ correspondant à l'identité de $id_C$ est le module des isomorphismes $id_{C[\epsilon]}\to f$ où $f$ est un endofoncteur de $C[\epsilon]$ qui pour $\epsilon=0$ redonne $id_C$. La proposition {\ref{hczero}} caractérise ces objets comme les éléments de $HC^0(C)$. D'où
$$
\textrm{T}_{X\times_{EQ_1}EQ_0,id_{id_C}} = HC^0(C)
$$

La flèche de bord $\textrm{T}_{X\times_{EQ_0}EQ_1,id_C}\to \textrm{T}_{EQ_0,x}$ associe à une équivalence $C[\epsilon]\to D$ la catégorie $D$ et 
la flèche de bord $\textrm{T}_{X\times_{EQ_1}EQ_0,id_{id_C}}\to \textrm{T}_{X\times_{EQ_0}EQ_1,id_C}$ associe à un isomorphisme $\alpha:id_{C[\epsilon]}\to f$ l'équivalence $f$.
La comparaison avec les cobords du complexe de Hochschild (propositions {\ref{hczero}} et {\ref{hcun}}) assure que le complexe tangent de $\ukcateq$ en $x$ 
$$
HC^0(C) \to HC^1(C) \to HZ^2(C)
$$
est exactement un tronqué du complexe de Hochschild.
\end{pr}

\subsection{Stratification par le nombre d'objets}

On rappelle du {\sf\S\ref{chpcatcnx}} les champs $\ukcateqinfn$.

\begin{prop} Les champs $\ukcateqinfn$ sont 2-géométriques et les inclusions $$\ukcateqinfn\longrightarrow \ukcateqinfnplusun$$ sont ouvertes.
En particulier $\iota : \ukcat_*:=\ukcat^{Eq,\leq1}\to \ukcateq$ est une inclusion ouverte.
\end{prop}
\begin{pr}
On utilise la carte du paragraphe précédent et la semi-continuité de la fonction 'nombre essentiel d'objets' (cf. corollaire {\ref{nbessobjets}}).
\end{pr}



\begin{cor} Une déformation infinitésimale d'un point de $\ukcateq$ ne peut pas augmenter le nombre essentiel d'objets, en particulier toute déformation infinitésimale d'une gerbe linéaire reste une gerbe.
\end{cor}


\section{Modules des catégories abéliennes}

On prouve dans cette section la géométricité du champ $\ukab$, puis on montre que le complexe tangent en un point est un tronqué du complexe de Hochschild.

\subsection{Géométricité}

\subsubsection{Carte}

On a les morphismes 
$$
\ukass \overset{\B}{\longrightarrow} \ukcateq \overset{\MMod}{\longrightarrow} \ukab
$$
$\MMod$ est surjectif par définition de $\ukab$ et le composé $\MMod\circ\B$ l'est aussi en vertu du lemme {\ref{surjassab}}.

\begin{prop}\label{carteabel}
Le morphisme $\MMod\circ \B$ est lisse.
\end{prop}
\begin{pr}
Soit $A\in\kcom$ et $A'\to A$ une extension infinitésimale au premier ordre.
Avec les notations de la proposition {\ref{lissequillen}}, il convient de savoir si un objet
$(C,M,N,\alpha,\beta)$ où $C\in A\z\ass$, $M\in A\z\ab$ et $N\in A'\z\ab$ et où $\alpha:\overline{C}\to M$ et $\beta:N\to M$ sont des équivalences, est relié par des équivalences à un objet du type $(D,\overline{D},\overline{D},id,id)$ où $D\in A\z\ass$.
Tout d'abord on remarque que $N$ est du type $\overline{D}$ et que via une équivalence de quintuplets (cf. {\ref{lissequillen}}) on peut toujours supposer que $M=\overline{D_u}$. 
Par cette première équivalence, $C$ peut-être vu par $\alpha$ comme un générateur de $\overline{D_u}$ et le problème se ramène à la construction une déformation de $C$ en un générateur de $\overline{D}$.

L'équivalence $\overline{C}\simeq \overline{D_u}$ est donnée par un $C\z D_u$-bimodule inversible $E$, 
en particulier $\End_{D_u}(E)\simeq C$.
$E$ est projectif de type fini sur $D_u$ (lemme {\ref{eqcateqmor}}) et le lemme {\ref{defmodproj}} nous permet de le déformer en un $D$-module projectif de type fini $F$.
En conséquence, on a $\End_D(F)_u \simeq \End_{D_u}(F_u)=\End_{D_u}(E)\simeq C$ et, en posant $C':=\End_D(F)$, on a $\overline{C'}^u\simeq \overline{C'_u}\simeq \overline{C}$.

On a trouvé une déformation de $\overline{C}$, il reste à voir que $F$ définit une équivalence de Morita $\overline{D}\to \overline{C'}$, il suffit pour cela de vérifier que $F$, vu comme $C'$-module est générateur de $\overline{C'}$. Soit $X\in \overline{C'}$, on montre que $\Hom_{C'}(F,X)=0\Rightarrow X=0$.
Comme $F$ est projectif de type fini $\Hom_{C'}(F,X)_u=\Hom_{C}(F_u,X_u)= 0$ d'où $X_u$ est nul car $F_u=E$ est générateur de $\overline{C}$. $X_u$ est en particulier un $A$-module de type fini et $X$ est en particulier une déformation de $X_u$ en un $A'$-module, il est donc nul par le lemme de Nakayama.
\end{pr}


\medskip

Le fait que $\MMod\circ \B$ soit lisse n'entraîne pas directement la lissité de $\MMod$ car $\B$ n'est pas surjective, néanmoins, la technique de la preuve de la proposition {\ref{carteabel}} s'adapte pour démontrer la lissité de $\MMod$.

\begin{lemme}
Soit $d$ un type de graphe, on note $|d|=\sum_{x,y}d(x,y)$.
Le morphisme de schémas $\cat^{(d)}\to \ass^{|d|}:C\mapsto [C]$ est lisse.
\end{lemme}
\begin{pr}
%
%
Soit $C$ un $A$-point de $\cat^{(d)}$, son image $[C]$ est un anneau muni canoniquement d'une famille complète d'idempotents orthogonaux donnée par les objets de $C$, \ie d'éléments $p_x$, où $x$ est un objet de $C$, telle que $p_xp_y=0$ si $x\not=y$, $p_x^2=p_x$ et $\sum_xp_x=1$.
Les $p_x$ qui permettent de reconstituer $C$ à partir de $[C]$ car $C(x,y)\simeq p_y[C]p_x$.

Le problème de la lissité, dans sa version relèvement des extensions infinitésimales, se reformule alors en le problème de faire suivre ces idempotents le long d'une déformation de $[C]$ en un autre famille complète d'idempotents orthogonaux.

On le démontre par récurrence en indexant les $n$ objets de $C$ par $\{1,\dots n\}$. 
Soit $A\to B\in\kcom$ une extension au premier ordre et soit $D\to C$ une extension au premier ordre de $C$ vu comme $B$-module, dont on note $K$ le noyau. On a déjà vu dans la preuve du lemme {\ref{defmodproj}} qu'on pouvait remonter les projecteurs de $C$ en des projecteurs de $D$, cela permet de remonter $p_1$ en un projecteur $q_1$.
Maintenant on suppose qu'on a remonté $p_1,\dot,p_i$ en des projecteurs orthogonaux $q_1,\dots, q_i$, et on montre qu'on peut remonter $p_{i+1}$ en un élément $q_{i+1}$ tel que $q_jq_{i+1}=0$ pour tout $j\leq i$ et tel que $q_{i+1}^2=q_{i+1}$. On pose $q=\sum_{j\leq i}q_j$ comme ils sont orthogonaux on a, pour tout $j\leq i$ : $q_jq=q_j$, ainsi que $q^2=q$. Soit $q'_{i+1}$ un relevé de $p_{i+1}$, la relation $qq'_{i+1}=0$ est équivalente à $q_jq'_{i+1}=0$ pour tout $j\leq i$. Comme $q'_{i+1}$ ne vérifie pas a priori cette relation, on montre qu'on peut toujours le déformer par un élément de $K$ en un $q''_{i+1}$ qui la vérifie. En effet, l'équation $q(q'_{i+1}+x)=0$ admet $x=qq'_{i+1}$ comme solution et cet $x$ est dans $K$ car son image dans $C$ est $(\sum_{j\leq i}p_i)p_{i+1}=0$ ; on pose donc $q''_{i+1}=qq'_{i+1}$.
On sait par le raisonnement de la preuve du lemme {\ref{defmodproj}} qu'on peut déformer $q''_{i+1}$ en un projecteur $q_{i+1}$ par un élément $x\in K$ et l'hypothèse de récurrence sera montrée si cet élément vérifie encore $qq_{i+1}=0$.
Or la formule obtenue lemme {\ref{defmodproj}} pour $x$ admet un $q''_{i+1}$ en facteur et donc $qx=0$.

On utilise ainsi la récurrence jusqu'à remonter $p_{n-1}$ puis on pose $q_n=1-\sum_{i<n}q_i$.
\end{pr}

\begin{prop}\label{lissitemod}
$\MMod$ est lisse. 
\end{prop}
\begin{pr}
On a un morphisme de schémas affines $[-]:\coprod_d \cat^{(d)}\longrightarrow \coprod_n \mathcal{A}ss^n$ donné par $C\mapsto [C]$. On a un diagramme
$$\xymatrix{
\coprod_n \mathcal{A}ss^n \ar[r]^a & \ukass\ar[r]^{\MMod\circ\B} & \ukab \\
\coprod_d \cat^{(d)}\ar[u]^{[-]} \ar[r]_c&\ukcateq \ar[ru]_{\MMod}
}$$
et le lemme {\ref{surjassab}} assure que $\MMod\circ\B\circ a\circ [-] \simeq \MMod\circ c$.

Comme $\MMod\circ\B\circ a\circ [-]$ est lisse, comme composée de morphismes lisses, et que $c:\coprod_d \cat^{(d)}\longrightarrow \ukcateq$ est lisse, on déduit que $\MMod$ est lisse (par définition de la lissité).
%
%
%
\end{pr}

\subsubsection{Diagonale}

\begin{prop}
Soient $C$ et $D$ deux $A$-algèbres associatives projectives de type fini. Le champ $\Omega_{\overline{C},\overline{D}}\ukab$ est équivalent au champ $\underline{\mathcal{I}nv}(C,D)$ classifiant les $C\z D$-bimodules inversibles. En particulier, c'est un 1-champ 1-géométrique et une présentation est donnée par le groupoïde affine lisse :
$$
\coprod_m Inv(C,D,m)\times G\ell_m \rightrightarrows \coprod_m Inv(C,D,m).
$$
\end{prop}
\begin{pr}
On rappelle que $\mathcal{I}nv(C,D)$ est définit au {\sf\S\ref{modulesinversibles}}.
Le fait que $\Omega_{\overline{C},\overline{D}}\ukab$ soit équivalent à $\mathcal{I}nv(C,D)$ est démontré à la proposition {\ref{diagcatab}}. La géométricité de $\mathcal{I}nv(C,D)$ et sa présentation sont démontrés à la proposition {\ref{geominv}}.
\end{pr}

\subsubsection{Présentation}

\begin{thm}\label{thmgeommod}
$\ukab$ est un 2-champ 1-géométrique localement de présentation finie. 
Une présentation par un 2-groupoïde est donnée par celui des anneaux à équivalence de Morita près :
$$
\coprod_{c,d,m} Inv(c,d,m)\times G\ell_m \rightrightarrows \coprod_{c,d,m}Inv(c,d,m) \rightrightarrows \coprod_n \ass^n.
$$
\end{thm}
\begin{pr}
\end{pr}

\begin{rem} Tout comme $\ukcateq$ (cf. remarque {\ref{chpcatcateq}}), il est assez naturel qu'il doit exister un champ en 2-catégories classifiant les catégories abéliennes, dont une présentation par une 2-catégorie serait :
$$
\coprod_{c,d,m,n}MorBim(c,d,m,n)\rightrightarrows \coprod_{c,d,m} Bim(c,d,m) \rightrightarrows \coprod_n\ass^n.
$$
\end{rem}

\begin{cor}
Une autre présentation de $\ukab$ par un 2-groupoïde est donnée par celui (définit dans la preuve) des catégories à équivalence de Morita près :
$$
\coprod_{c,d,m} Inv'(c,d,m)\times G\ell_m \rightrightarrows \coprod_{c,d,m}Inv'(c,d,m) \rightrightarrows \coprod_d \cat^{(d)}.
$$
où dans les réunions les $c,d,m$ sont des types de graphe.
\end{cor}
\begin{pr}
On définit d'abord le 2-groupoïde en question.
Pour $C$ une catégorie linéaire projective de type fini ayant un nombre fini d'objet, l'association $C\to [C]$ définit un morphisme de schémas $\cat^{(d)}\to \ass^{|d|}$ où $|d|:=\sum_{x,y\in C}d(x,y)$.

Les lemmes {\ref{surjassab}} et {\ref{eqcateqmor}} assurent que les $C\z D$-bimodules inversibles sont les $[C]\z [D]$-bimodules inversibles, en conséquence le classifiant des triplets $(C,D,M)$ où $C$ et $D$ sont des catégories linéaires de graphe libre de types respectifs $c$ et $d$ fixé et où $M$ est un $C\z D$-bimodule inversible de graphe libre de type $m$\footnote{Un bimodule étant en particulier un graphe, son type est celui du graphe sous-jacent.} est donné par 
$$
Inv'(c,d,m) := \left(\cat^{(c)}\times \cat^{(d)}\right)\times_{\left(\ass^{|c|}\times \ass^{|d|}\right)}Inv(|c|,|d|,|m|)
$$
et le schéma classifiant les quintuplets $(C,D,M,N,f)$ où $C$ et $D$ sont des catégories linéaires de graphes libre de type respectifs $c$ et $d$ fixé, où $M$ et $N$ sont des $C\z D$-bimodules inversibles de graphes libres de type $m$ fixé et où $f: M\to N$ est un isomorphisme de bimodules est donné par
$$
Inv'(c,d,m) \times G\ell_m
$$
(car la donnée de $M$ est d'un isomorphisme $f$ caractérise $N$).

Compte tenu de ce que classifient ces schémas il est évident que le 2-graphe
$$
\coprod_{c,d,m} Inv'(c,d,m)\times G\ell_m \rightrightarrows \coprod_{c,d,m}Inv'(c,d,m) \rightrightarrows \coprod_d \cat^{(d)}.
$$
est muni d'une structure de 2-groupoïde.

$EQ_0=\coprod_d \cat^{(d)}$ est une carte pour $\ukab$ d'après la proposition {\ref{lissitemod}} et 
$\Omega_{EQ_0}\ukab$ est le 1-champ classifiant à isomorphisme près les bimodules inversibles entre deux catégories linéaires dont la proposition {\ref{}} assure que $\coprod_{c,d,m} Inv'(c,d,m)\times G\ell_m \rightrightarrows \coprod_{c,d,m}Inv'(c,d,m)$ est une présentation par un groupoïde.
\end{pr}

\begin{nota}\label{notaab}\index{$AB_0$}\index{$AB_1$}\index{$AB_2$} \index{$AB'_0$}\index{$AB'_1$}\index{$AB'_2$}
Pour des besoins de simplicité d'écriture on pose $AB_2 :=\coprod_{c,d,m} Inv(c,d,m)$, $AB_1:=\coprod_{c,d,m}Inv(c,d,m)$ et $AB_0:=\coprod_n \ass^n$ ; ainsi que $AB'_2 :=\coprod_{c,d,m} Inv'(c,d,m)$, $AB'_1:=\coprod_{c,d,m}Inv'(c,d,m)$ et $AB'_0:=\coprod_d \cat^{(d)}$.
\end{nota}

\subsection{Tangent}

Le complexe de Hochschild est définit en annexe {\ref{hoch}}.

\begin{thm}\label{tangentabel}
Soient $X=Spec(A)\in\kaff$, $C:X\to \cat^{(d)}$ et $x:X\to \cat^{(d)}\to \ukab$ le point correspondant de $\ukab$. Le complexe tangent de $\ukab$ en $x$ est donné par un complexe quasi-isomorphe au tronqué $Hoch^{\leq2}$ du complexe de Hochschild de la catégorie $C$ (cf. théorème {\ref{tangentequiv}}).
\end{thm}
%
\begin{pr}
On utilise la même preuve qu'au théorème {\ref{tangentequiv}} avec la présentation $AB'_2\rightrightarrows AB'_1\rightrightarrows AB'_0$ (cf. {\ref{notaab}}) au lieu de $EQ_2\rightrightarrows EQ_1\rightrightarrows EQ_0$.

On en déduit que $\mathbb{T}_{\ukab,x}$ est quasi-isomorphe au complexe de $A$-modules :
$$
\textrm{T}_{X\times_{AB'_1}AB'_2,id_C} \longrightarrow \textrm{T}_{X\times_{AB'_0}AB'_1,C}
\longrightarrow \textrm{T}_{AB'_0,C}.
$$

Le tangent au point $C$ de $EQ_0$ est exactement une déformation au premier ordre de la loi de $C$ (cf. Annexe {\ref{hoch}}) on a donc
$$
\textrm{T}_{EQ_0,x} = HZ^2(C).
$$
Pour les autres on remarque d'abord que 
\begin{itemize}
\item $X\times_{AB'_0}AB'_1$ est le classifiant des bimodules inversibles de source $C$.
\item $X\times_{AB'_1}AB'_2$ est le classifiant des isomorphismes (entres bimodules inversibles) de source $C$ vu comme $C\z C$-bimodule.
\end{itemize}
Avec les notations de {\ref{interpmorhoch}}, on en déduit que le tangent de $X\times_{AB'_0}AB'_1$ au point correspondant à $C$ vu comme $C\z C$-bimodule est le module des $C[\epsilon]\z C[\epsilon]_\alpha$-bimodules inversibles, qui pour $\epsilon=0$ donnent $C$ comme $C\z C$-bimodule. La proposition {\ref{moritatriv}} assure que le morphisme naturel
$$
\textrm{T}_{X\times_{AB'_0}AB'_1,C} \to HZ^2(C)
$$
à pour image le module des 1-cobords de Hochschild, et la proposition {\ref{kertangentbim}} dit que son noyau est 
$$
Z^1_{C\z C}(C).
$$
Toujours avec les notations de {\ref{interpmorhoch}}, le tangent de $X\times_{AB'_1}AB'_2$ au point $id_C$ est le module des isomorphismes $C[\epsilon]\to C[\epsilon]_\mu$, où $\mu$ est une dérivation intérieure, qui pour $\epsilon=0$ donnent $id_C$ ; la proposition {\ref{tangentisobim}} donne :
$$
\textrm{T}_{X\times_{AB'_0}AB'_1,C} = C^0_{C\z C}(C) = \End_A(C).
$$
Le complexe tangent est donc
$$
C^0_{C\z C}(C)\longrightarrow \textrm{T}_{X\times_{AB'_0}AB'_1,C} \longrightarrow HZ^2(C).
$$
Par la proposition {\ref{tangentisobim}} la première flèche est le bord du complexe $C^*_{C\z C}(C)$ et la cohomologie est exactement celle de Hochschild en degré $\leq 2$.
On a même un quasi-isomorphisme de complexes :
$$\xymatrix{
Hoch^{\leq2} \ar[d] &C \ar[r]\ar[d]_{m_0}& \End_A(C) \ar[r]\ar[d]_{m_1}& HZ^2(C) \ar@{=}[d]\\
MHoch^{\leq2}:=&\End_A(C)\ar[r]& \textrm{T}_{X\times_{AB'_0}AB'_1,C} \ar[r]& HZ^2(C)
}$$\index{$MHoch^{\leq2}$}
où $m_1$ associe à une équivalence le bimodule inversible correspondant et où $m_0$ associe à un isomorphisme naturel l'isomorphisme déduit entre les bimodules.
\end{pr}

\medskip

\subsection{Modules de Morita}\label{ch4mor}

On déduit de l'équivalence $\ukcatmor\simeq \ukab$ la proposition suivante

\begin{prop}
Le champ $\ukcatmor$ est un 2-champ 1-géométrique et son complexe tangent en un point $C$ correspondant à une catégorie $A$-linéaire, est quasi-isomorphe à la cohomologie de Hochschild de $C$.
\end{prop}

\section{Comparaisons}\label{comparaison}

\begin{thm}\label{thmetale}
Le morphisme $\MMod : \ukcateq\longrightarrow \ukab$ est étale.
\end{thm}
\begin{pr}
On sait déjà qu'il est lisse (proposition {\ref{lissitemod}}) il suffit de montrer que le morphisme induit entre les tangents est un quasi-isomorphisme.
On a déjà remarqué dans la preuve du théorème {\ref{tangentabel}} que les tangents étaient quasi-isomorphes, il faut montrer que ce quasi-isomorphisme est celui induit par le morphisme $\ukcateq\longrightarrow \ukab$.

Pour le voir, on va montrer que le morphisme $\ukcateq\to ukab$ se relève aux présentations en 2-groupoïdes :
$$\xymatrix{
EQ_2\ar@{}[rd]|2\ar@{=>}[r]\ar[d]_{M_2} &EQ_1\ar@{}[rd]|1 \ar@{=>}[r]\ar[d]_{M_1} &EQ_0\ar@{=}[d]_{M_0}\ar[r]\ar@{}[rd]|0& \ukcateq\ar[d]^{\MMod}\\
AB'_2\ar@{=>}[r] & AB'_1\ar@{=>}[r] &AB'_0 \ar[r]& \ukab
}$$
où $M_0$ est l'identité de $\coprod_d \cat^{(d)}$, $M_1$ associe à une équivalence $f:C\to D$ une structure de $C\z D$-bimodule sur $D$ et $M_2$ associe à un isomorphisme naturel l'isomorphisme déduit entre les bimodules.
(Il est clair que le quasi-isomorphisme de la preuve du théorème {\ref{tangentabel}} provient des différentielles de $M_0$, $M_1$ et $M_2$.)
La commutation des carrés 1 et 2, ainsi que le fait qu'ils forment un morphisme de 2-groupoïdes sont claires compte tenu de ce que classifient les schémas. Celle du carré 0 est déduite de ce que les cartes utilisées de $\ukcateq$ et de $\ukab$ soient les mêmes : $EQ_0=AB'_0=\coprod_d\cat^{(d)}$. 

On déduit de ce dernier fait un morphisme entre les espaces de lacets
$$
M : \Omega_{EQ_0}\ukcateq \longrightarrow \Omega_{AB'_0}\ukab
$$
dont il faut montrer qu'il se relève aux présentations en les morphismes $M_1$ et $M_2$.

Le lemme {\ref{pleineqbim}} assure que $M$ est plein et $P := AB'_1\times_{\Omega_{AB'_0}\ukab} \Omega_{EQ_0}\ukcateq$ est donc un sous-faisceau de $AB'_1$, il classifie les bimodules inversibles entre deux catégories équivalentes et factorise donc $M_1$ en $EQ_1\to P \to AB'_1$ : $M_1$ est bien le morphisme issu de $\MMod$ entre les cartes de $\Omega_{EQ_0}\ukcateq$ et $\Omega_{AB'_0}\ukab$.
On tire de $M_1$ un morphisme de schémas
$$
M' : \Omega_{EQ_1}\Omega_{EQ_0}\ukcateq \longrightarrow \Omega_{AB'_1}\Omega_{AB'_0}\ukab
$$
dont il est facile de voir que c'est $M_2$.
\end{pr}

\begin{prop}\label{zeroetale}
Le morphisme $\ukass\longrightarrow \ukab$ est 0-étale (cf. définition {\ref{defetale}}).
\end{prop}
\begin{pr}
Il est lisse d'après la proposition {\ref{carteabel}} il suffit de montrer que le morphisme induit entre les module de cohomologie tangente est isomorphisme en degré 0.
$\ukass$ admet comme carte $ASS_0:=\coprod_n\ass^n$ où $\ass^n$ classifie les structures d'algèbre associative sur un module libre de rang $n$. $ASS_1:= \left(\coprod_n\ass^n\right)\times {\ukass}\left(\coprod_n\ass^n\right)$ est le schéma classifiant les triplets de deux structures d'algèbre et d'un isomorphisme entre elle. $ASS_1\rightrightarrows ASS_0$ est une présentation de $\ukass$.
On utilise cette fois la carte $AB_2\rightrightarrows AB_1\rightrightarrows AB_0$ de $\ukab$ (cf. {\ref{notaab}}).

Le morphisme $\ukass\to ukab$ se relève aux présentations en 2-groupoïdes :
$$\xymatrix{
&ASS_1\ar@{}[rd]|1 \ar@{=>}[r]\ar[d]_{M_1} &ASS_0\ar@{=}[d]_{M_0}\ar[r]^a\ar@{}[rd]|0& \ukass\ar[d]^{\MMod\circ\B}\\
AB_2\ar@{=>}[r] & AB_1\ar@{=>}[r] &AB_0 \ar[r]^{ab}& \ukab
}$$
En effet, $M_0$ est l'identité de $ASS_0=AB_0$ et la commutativité de 0 vient de ce que $ab=\MMod\circ\B\circ a$ ; cela induit un morphisme 
$$
M : ASS_1=\Omega_{ASS_0}\ukass \longrightarrow \Omega_{AB_0}\ukab
$$
dont il faut montrer qu'il se relève à la carte de $\Omega_{AB_0}\ukab$ donné par $AB_1$. Or $ASS_1$ classifie des isomorphismes entre algèbres de module sous-jacent libre et ceux-ci correspondent bien à des points de $AB_1$ qui classifie les bimodules inversibles entre algèbres de module sous-jacent libre, d'où le relèvement.

On en déduit que le morphisme induit sur les complexes tangents en un point $C:Spec(A)\to \ass^n\to \ukass\to\ukab$ est 
$$\xymatrix{
Der^{\leq1}\ar[d]&&HC^1(C) \ar[r]\ar[d]_{m_1} &HZ^2(C)\ar@{=}[d]_{m_0}\\
MHoch^{\leq2}\ar[d]&\End_A(C)\ar[r]\ar[d]& \textrm{T}_{X\times_{AB'_0}AB'_1,C} \ar[d]\ar[r]& HZ^2(C)\ar@{=}[d]\\
Hoch^{\leq2}& HC^0(C)\ar[r] & HC^1(C)\ar[r] &HZ^2(C)
}$$
où la troisième ligne est le complexe tangent de $\ukcateq$ en $C:Spec(A)\to \ass^n\to \ukass\to\ukcateq$.
La composition $Der^{\leq1}\to MHoch^{\leq2}\to Hoch^{\leq2}$ est l'inclusion évidente qui induit un isomorphisme sur les $H^0$. Comme $MHoch^{\leq2}\to Hoch^{\leq2}$ est un quasi-isomorphisme (théorème {\ref{tangentabel}}), on déduit l'isomorphisme cherché.
\end{pr}

\newpage
\section{Diagramme final}\label{diagrammefinal}

On reprend le diagramme du {\S\ref{morphismes}} en le complétant.
$$
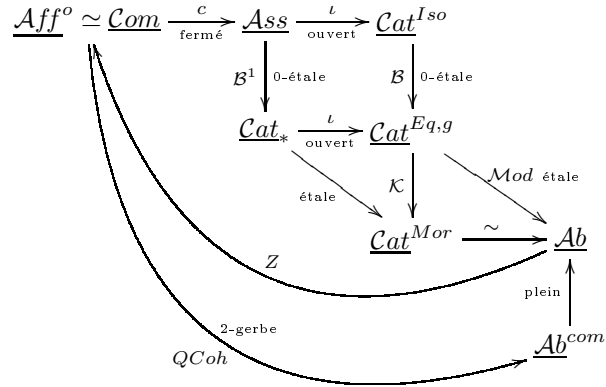
\xymatrix{
\underline{\mathcal{A}f\!f}^o\simeq \ukcom \ar[r]^-c_-{\textrm{{\tiny fermé}}}\ar@/_6pc/[rrrddd]_{QCoh}^{\textrm{{\tiny 2-gerbe}}} & \ukass \ar[d]_{\B^1}^-{\textrm{{\tiny 0-étale}}} \ar[r]^-\iota_-{\textrm{{\tiny ouvert}}} & \ukcatiso \ar[d]_\B^-{\textrm{{\tiny 0-étale}}}\\
& \ukcat_* \ar[rd]_-{\textrm{{\tiny étale}}} \ar[r]^-\iota_-{\textrm{{\tiny ouvert}}} & \ukcateqg \ar[d]_-{\mathcal{K}} \ar[rd]^-{\MMod\ \textrm{{\tiny étale}}}  \\
&  & \ukcatmor  \ar[r]^{\sim} & \ukab\ar@/^5pc/[llluu]_Z\\
&&&\ukab^{com}\ar[u]^-{\textrm{{\tiny plein}}}
}$$

Au niveau des complexes tangents on a schématiquement :
$$
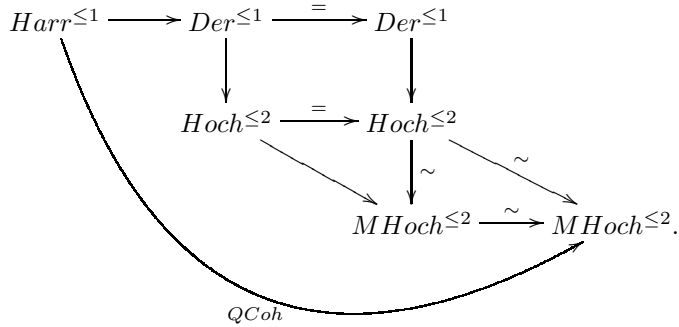
\xymatrix{
Harr^{\leq1} \ar[r] \ar@/_6pc/[rrrdd]_{QCoh} & Der^{\leq1}\ar[d] \ar[r]^-{=} & Der^{\leq1} \ar[d] \\
& Hoch^{\leq2} \ar[rd]\ar[r]^-{=} & Hoch^{\leq2}\ar[d]^-{\sim} \ar[rd]^-{\sim}  \\
&  & MHoch^{\leq2}  \ar[r]^-{\sim} & MHoch^{\leq2}.
}$$

\bigskip

On rappelle que $\ukab^{com}$ est le sous-champ plein de $\ukab$ image de $\ukcom$ par $QCoh$.

\begin{prop}
$\ukAB^{com}$ est un champ géométrique.
\end{prop}
\begin{pr}
La proposition {\ref{abcom}} établit que la fibre de $\ukcom\longrightarrow \ukab^{com}$ en un point $C:Spec(A)\to \com\to \ukcom$ est le champ $\mathcal{K}(\underline{Pic}(C),1)$, comme $\ukcom$ est géométrique, la géométrie $\ukab^{com}$ peut se déduire de celle des $\mathcal{K}(\underline{Pic}(C),1)$ (lemme {\ref{topfibngeom}}).

$\mathcal{K}(\underline{Pic}(C),1)$ admet la présentation par le groupoïde
$$
\underline{Pic}(C)\rightrightarrows pt
$$
et il suffit alors de prouver que $\underline{Pic}(C)$ est géométrique.
$\underline{Pic}(C)$ est le sous-champ plein du 1-champ $\underline{\mathcal{I}nv}(C,C)$ formé des seuls $C\z C$-bimodules dont les structures à droite et à gauche coïncident.
$\underline{\mathcal{I}nv}(C,C)$ admet une présentation par un groupoïde géométrique lisse (cf. proposition {\ref{geominv}}) et il suffit pour avoir une même présentation pour $\underline{Pic}(C)$ de rajouter dans les schémas affines $Bim(C,C,m)$ les équations $\mu(c,x,1)=\mu(1,x,c)$ qui s'écrivent en coordonnées, en utilisant les notations du {\sf\S\ref{modulesinversibles}}, $\mu_{ijk}1^k - 1^i\mu_{ijk}=0$ (où les $1^k$ désignent les coordonnées de l'unité de $C$).
\end{pr}

\appendix

\setcounter{totx}{0}
\newpage

\thispagestyle{empty}
\chapter{Cohomologie de Hochschild}\label{hoch}
\thispagestyle{chheadings}

La première section de l'annexe définit la cohomologie de Hochschild d'une catégorie linéaire et la seconde définit la cohomologie d'un module sur une algèbre associative et compare les deux.

\section{Complexe de Hochschild}

On définit la cohomologie de Hochschild d'une catégorie linéaire comme la cohomologie du complexe de Hochschild, puis on explicite les description traditionnelles des 0-, 1- et 2- cocycles.

\subsection{Définition}

Soit $C\in A\z\CAT$ (cf. {\sf\S\ref{categolin}}). Pour deux objets $x,y\in C$, $C(x,y)$ désigne le $A$-module des morphismes dans $C$ de $x$ vers $y$. Pour trois objets $x,y,z\in C$ on note $m_{xyz}$ la composition
$$
m_{xyz}:C(x,y)\otimes_AC(y,z)\to C(x,z),
$$
c'est une application $A$-bilinéaire.

La {\em cohomologie de Hochschild}\index{cohomologie de Hochschild} de $C$ est définit comme étant la cohomologie du complexe suivant :
$$
HC^0(C)\longrightarrow  HC^1(C)\longrightarrow  HC^2(C)\longrightarrow \dots
$$
où 
\begin{eqnarray*}
&&HC^0(C) = \bigoplus_{x\in C}  C(x,x),\\
&&HC^i(C) = \bigoplus_{x_0,\dots x_i\in C} \Hom_A(C(x_0,x_1)\otimes_AC(x_1,x_2)\otimes_A\dots \otimes_AC(x_{i-1},x_i),C(x_0,x_i))
\end{eqnarray*}
et où la différentielle
$$
d:HC^i(C) \longrightarrow HC^{i+1}(C)
$$
est définit par 
\begin{multline*}
d(f): c_{0,1}\otimes \dots \otimes c_{i,i+1} \longmapsto  m_{0,1,i}(c_{0,1},f(c_{1,2}\otimes\dots\otimes c_{i,i+1})) - f(m_{012}(c_{0,1},c_{1,2})\otimes\dots \otimes c_{i,i+1})) \\
+ f(c_{0,1}\otimes \dots m_{123}(c_{1,2},c_{2,3}) \dots \otimes c_{i,i+1})) + \dots \\
+ (-1)^{i} f(c_{0,1}\otimes \dots \otimes m_{i-1,i,i+1}(c_{i-1,i}),c_{i,i+1}))) 
+ (-1)^{i+1} m_{0,i,i+1}(f(c_{0,1}\otimes \dots \otimes c_{i-1,i})),c_{i,i+1}).
\end{multline*}

\subsection{Cocycles}

Soit $A[\epsilon]$ l'extension au premier ordre de $A$ par $A$, il lui est associé deux morphismes $A\to A[\epsilon]\to A$. On note $C[\epsilon]$ la catégorie $C\otimes_AA[\epsilon]$, dont on rappelle qu'elle a les mêmes objets que $C$. Aux deux morphismes d'anneaux précédents sont associés deux foncteurs $C[\epsilon]\to C$ (évaluation $\epsilon=0$) et $C\to C[\epsilon]$ (prolongement par linéarité).

\medskip

Soit $M$ un $A$-module, une déformation au premier ordre de $M$ est un $A[\epsilon]$-module $\overline{M}$ tel que $\overline{M}\otimes_{A[\epsilon]}A\simeq M$. Dans toute la suite on se limite aux déformations au premier ordre de modules projectifs $M$ qui sont scindée comme $A$-modules. Si en plus on impose une hypothèse de projectivité sur $A[\epsilon]$, le $A$-module sous-jacent à une déformation de $M$ est nécessairement $M\oplus M$.

Une catégorie $A$-linéaire est dite {\em projective de type fini}\index{catégorie projective de type fini}, si tous ses modules de morphismes sont des $A$-modules projectifs de type fini.

Dans la suite, toutes les catégories linéaires considérés seront toujours projectives de type fini.
En particulier les déformations projectives au premier ordre d'une telle catégorie $C$, sont toujours de graphe sous-jacent isomorphe à $C[\epsilon]=C\otimes_AA[\epsilon]$ comme $A$-graphe.

\subsubsection{0-cocycles}

Le module des 0-cocycles de Hochschild consiste en les éléments de $c=\oplus f_t \in \oplus_t C(t,t)$ vérifiant, pour tout $u\in C(x,y)$ :
\begin{eqnarray*}
d(f)(u) &= & m_{xyt}(u,f_t) - m_{txy}(f_t,u) = 0.
\end{eqnarray*}

\begin{defix}
On définit le {\em centre}\index{centre d'une catégorie} de $C$ comme l'algèbre $Z(C)$ formée des automorphismes de $id_{C[\epsilon]}$ qui, pour $\epsilon=0$, donnent $id_{id_C}$,\ie l'algèbre des déformations au premier ordre de $id_{id_C}$. $Id_{id_C}$ étant le neutre de l'algèbre $\End(id_C)$ et $Z(C)$ s'identifie canoniquement à $\End(C)$.
\end{defix}

Un élément $c:id_{C[\epsilon]}\to id_{C[\epsilon]} $ du centre associe à tout objet $x\in C$ un élément $c^0_x+c^1_x\epsilon\in C[\epsilon](x,x)$ tel que pour tout $u\in C[\epsilon](x,y)$ on ait
$$
m_{xxy}(c_x,u)=m_{xyy}(u,c_y)
$$
où $m_{xxy}$ est la loi de $C[\epsilon]$, \ie la $\epsilon$-linéarisation de la loi de $C$.
La condition en $\epsilon=0$ assure que, pour tout $x$, $c^0_x=0$ et l'équation précédente se réduit à 
$$
m_{xxy}(c^1_x,u)=m_{xyy}(u,c^1_y)
$$
Cette condition est exactement celle de 0-cocycle de Hochschild.

On a prouvé la proposition suivante.
\begin{propx}
On a $HH^0(C)=HZ^0(C) = Z(C)$
\end{propx}

\subsubsection{1-cocycles}

Le module des 1-cocycles consiste en les éléments $f\in \oplus_{x,y} End(C(x,y))$ vérifiant, pour tout $c_{x,y}\in C(x,y)$ et pour tout $c_{y,z}\in C(y,z)$ :
\begin{eqnarray*}
d(f)(c_{x,y}\otimes c_{y,z}) & = & m_{xyz}(c_{x,y},f(c_{y,z})) - f(m_{xyz}(c_{x,y}\otimes c_{y,z})) + m_{xyz}(f(c_{x,y}),c_{y,z}) = 0.
\end{eqnarray*}

\begin{defix}
On définit une {\em dérivation}\index{dérivation de catégorie} de $C$ comme un isomorphisme $C[\epsilon]\to C[\epsilon]$ tel que son évaluation en $\epsilon=0$ donne l'identité de $C$, \ie comme une déformation au premier ordre de $id_C$.
$id_C$ étant le neutre de l'algèbre $\End(C)$, le $A$-module des dérivations s'identifie canoniquement à $\End(C)$.
\end{defix}

Un tel foncteur $f$ est nécessairement l'identité sur les objets et se décompose au niveau des modules de morphismes en
$$
f_{xy}=f_{xy}^0\oplus f^1_{xy}\epsilon : C(x,y)[\epsilon]\longrightarrow C(x,y)[\epsilon]
$$
en utilisant la linéarité en $\epsilon$, on se ramène à des applications
$$
f_{xy}=f_{xy}^0\oplus f^1_{xy}\epsilon : C(x,y)\longrightarrow C(x,y)\oplus C(x,y)\epsilon
$$
où $f^0_{xy}$ est l'identité par définition d'une dérivation.
Ainsi, $f$ est entièrement caractérisé par les
$$
f^1_{xy} : C(x,y)\longrightarrow C(x,y).
$$
qui sont des 1-cochaînes de Hochschild. Réciproquement, une 1-cochaîne $f^1$ ne correspond à une dérivation que si $id\oplus f^1\epsilon$ est un foncteur ; ceci se traduit, pour tous $x,y,z\in C$ par les équations :
\begin{eqnarray*}
(id\oplus f^1_{xz})(m_{xyz}(c_{xy},c{yz})) &=& m_{xyz}((id\oplus f^1_{xz})(c_{xy}),(id\oplus f^1_{xz})(c_{yz}))\\
&\iff &\\
f^1_{xz}(m_{xyz}(c_{xy},c{yz}) & = & m_{xyz}(f^1_{xz}(c_{xy}),c_{yz}) + m_{xyz}(c_{xy},f^1_{xz}(c_{yz})).
\end{eqnarray*}
Ces équations sont exactement la condition de 1-cocycle de Hochschild.

On a prouvé la proposition suivante.
\begin{propx}
$HZ^1(C)$ est le $A$-module des dérivations de $C$.
\end{propx}

\subsubsection{2-cocycles}

Le module des 2-cocycles consiste en les éléments $f\in \oplus_{x,y,z} \Hom(C(x,y)\otimes C(y,z),C(x,z))$ vérifiant, pour tous $c_{x,y}\in C(x,y)$, $c_{y,z}\in C(y,z)$ et $c_{z,t}\in C(z,t)$ :
\begin{multline*}
d(f)(c_{x,y}\otimes c_{y,z}\otimes c_{z,t}) = m_{xyt}(c_{x,y},f(c_{y,z},c_{z,t})) - f(m_{xyz}(c_{x,y}\otimes c_{y,z})\otimes c_{z,t}) \\
+ f(c_{x,y}\otimes m_{yzt}(c_{y,z},c_{z,t})) - m_{xzt}(f(c_{x,y}\otimes c_{y,z}),c_{z,t}) = 0.
\end{multline*}

\begin{defix}
Une {\em déformation au premier ordre de la structure de catégorie linéaire de $C$}\index{déformation de catégorie} est une structure de catégorie $A[\epsilon]$-linéaire sur le graphe $C[\epsilon]$ tel que son évaluation en $\epsilon=0$ redonne la structure de $C$.
\end{defix}

Une telle structure se caractérise par des compositions
$$
m_{xyz}^0\oplus m_{xyz}^1\epsilon : C(x,y)[\epsilon]\otimes_{A[\epsilon]} C(y,z)[\epsilon]\longrightarrow C(x,z)[\epsilon]
$$
en utilisant la linéarité en $\epsilon$ on se ramène des applications $A$-linéaires
$$
m_{xyz}^0\oplus m_{xyz}^1\epsilon : C(x,y)\otimes_A C(y,z)\longrightarrow C(x,z)\oplus C(x,z) \epsilon
$$
dont la définition des déformation impose que $m_{xyz}^0=m_{xyz}$, \ie soit la composition de $C$.
Une déformation est donc caractérisé par les applications
$$
m_{xyz}^1 : C(x,y)\otimes_A C(y,z)\longrightarrow C(x,z)
$$
qui sont des 2-cochaîne de Hochschild. Réciproquement, une 2-cochaîne correspond à une déformation si la loi définie par $m_{xyz}\oplus m_{xyz}^1\epsilon$ est associative.

Pour tous $x,y,z,t\in C$, les conditions d'associativité s'écrivent :
\begin{eqnarray*}
(m_{xzt}\oplus m_{xzt}^1\epsilon)((m_{xyz}\oplus m_{xyz}^1\epsilon)(c_{xy},c_{yz}),c_{zt}) &=& (m_{xyt}\oplus m_{xyt}^1\epsilon)(c_{xy},(m_{yzt}\oplus m_{yzt}^1\epsilon)(c_{yz},c_{zt})\\
&\iff &\\
m_{xzt}(m_{xyz}^1(c_{xy},c_{yz}),c_{zt}) + 
m_{xzt}^1(m_{xyz}(c_{xy},c_{yz}),c_{zt}) 
&=& 
m_{xyt}(c_{xy},m_{yzt}^1(c_{yz},c_{zt}) +
m_{xyt}^1(c_{xy},m_{yzt}(c_{yz},c_{zt}).
\end{eqnarray*}
Ces conditions sont exactement celles d'être un 2-cocycle.

On a prouvé la proposition suivante.
\begin{propx}
$HZ^2(C)$ est le $A$-module des déformations au premier ordre de la loi de $C$.
\end{propx}

\subsection{Cobords}

On tire quelques interprétations de l'équation $d^2=0$ dans le complexe de Hochschild.

\begin{defix}
Une {\em dérivation intérieure}\index{dérivation intérieure de catégorie} de $C$ est une dérivation $f$ telle qu'il existe un isomorphisme $\alpha: id_{C[\epsilon]}\to f$, qui induise $id_{id_C}$ pour $\epsilon=0$.
\end{defix}

\begin{propx}\label{hczero}
Les dérivations intérieures sont en bijection avec les cobords de 0-cochaînes de Hochschild ; et les isomorphismes naturels de $\alpha:id_C\to f$, où $f$ est une dérivation de $C$, qui évalués en $\epsilon=0$ donnent $id_{id_C}$ sont en bijection avec les 0-cochaînes. Dans ce cas, $f$ s'identifie au cobord de $\alpha$.
\end{propx}
\begin{pr}
Une déformation de $id_C$ s'écrit toujours sous la forme $id\oplus\beta\epsilon$
Un isomorphisme $\alpha:id_{C[\epsilon]} \to id\oplus \beta\epsilon$ est la donnée pour tout $x\in C$ d'isomorphismes $\alpha_x\in C[\epsilon](x,(id+\beta\epsilon)(x))=C[\epsilon](x,x)$ tels que pour tout $v\in C[\epsilon](x,y)$ 
$$
(id+\beta\epsilon)(v)(\alpha_x\epsilon) = (\alpha_x\epsilon) v
$$
or $(id+d\alpha\epsilon)(v)(\alpha_x\epsilon)= (d\alpha(v)+ v\alpha_x)\epsilon$
et on est ramené à l'égalité $\beta+ v\alpha_x=\alpha_yv$ qui définit $\beta$ comme le cobord de $\alpha$.
\end{pr}

\begin{corx}
$HH^1(C)$ est le $A$-module des dérivations de $C$ modulo les dérivations intérieures.
\end{corx}

\begin{defix}
Une déformation au premier ordre $D$ d'une catégorie $A$-linéaire $C$ est dite {\em triviale} s'il existe une équivalence $f:C[\epsilon]\to D$, où $C[\epsilon]$ désigne la déformation par le 2-cocycle nul (déformation triviale), qui pour $\epsilon=0$ donne l'identité de $C$.
\end{defix}

\begin{propx}\label{hcun}
Les déformations au premier ordre triviales d'une catégorie $C$ sont en bijection avec les cobords de 1-cochaînes de Hochschild ; et les équivalences $f:C[\epsilon]\to D$, où $D$ est une déformation au premier ordre de $C$, qui induit l'identité de $C$ après évaluation en $\epsilon=0$, sont en bijection avec les 1-cochaînes. Dans ce cas, la déformation $D$ s'identifie au cobord de $f$.
\end{propx}
\begin{pr}
Identifiant $D$ à $C[\epsilon]$, une déformation de $C$ s'écrit toujours sous la forme $m\oplus\mu\epsilon$.
La condition sur $f:C[\epsilon]\to D=C[\epsilon]$ en $\epsilon=0$ impose que $f$ soit constante sur les objets, $f$ est donc caractérisé par des isomorphismes $f_{xy}:C[\epsilon](x,y) \to C[\epsilon](x,y)$ qu'on décompose en $f_{xy}=id\oplus f^1_{xy}\epsilon$, qui vérifient la relation de fonctorialité :
$$
(id\oplus f^1_{xz}\epsilon)\circ m_{xyz} = (m_{xyz}\oplus\mu\epsilon)(id\oplus f^1_{xy}\epsilon,id\oplus f^1_{yz}\epsilon).
$$
Or, $(id\oplus f^1_{xz}\epsilon)\circ m_{xyz}= m_{xyz} +  f^1_{xz}\circ m_{xyz} \epsilon$
et $(m_{xyz}\oplus\mu_{xyz}\epsilon)((id\oplus f^1_{xy}\epsilon),(id\oplus f^1_{yz}\epsilon))= m_{xyz}\oplus\mu_{xyz}\epsilon + m_{xyz}(-, f^1_{yz})\epsilon + m_{xyz}( f^1_{xy},-)\epsilon$
d'où
$$
f^1_{xz}\circ m_{xyz} = \mu_{xyz}+ m_{xyz}(-, f^1_{yz})+ m_{xyz}( f^1_{xy},-),
$$
qui définit $\mu$ comme le cobord de $f^1$.
\end{pr}

\begin{corx}
$HH^2(C)$ est le $A$-module des déformations de la structure de $C$ modulo les déformations triviales.
\end{corx}

\section{Cohomologie modulaire}

On définit la cohomologie d'un module, puis on la spécialise au cas des bimodules pour la comparer à la cohomologie de Hochschild.

\bigskip

Soient $A\in\kcom$, $B\in A\z\ASS$ et $M$ un $B$-module.
On définit les modules
$$
C^i_B(M):=\Hom_A(B^{\otimes_A i}\otimes_AM,M)=\Hom_B(B^{\otimes_A i+1}\otimes_AM,M).
$$
Les modules $B^{\otimes_Ai}\otimes_AM$ sont munis d'une structure de complexe {cf. \cite[ch. IX]{CE}}
ceci permet d'obtenir un complexe :
$$
C^0_B(M)\to C^1_B(M)\to \dots
$$
qui calcule la cohomologie de $M$ comme $B$-module, \ie $H^i_B(M):=H^i(C^*_B(M)) = Ext^i_B(M,M)$.

On note $Z_B^i(M)$ le $A$-module des $i$-cocycles de ce complexe et $B^i_B(M)$ le sous-module des $i$-cocycles qui sont des cobords.

\subsection{Cocycles}

\subsubsection{0-cocycles}
$d:\Hom_A(M,M)\to \Hom_A(B\otimes_AM,M)$ est donnée par 
$$
df(b,m) = bf(m)-f(bm)
$$

\begin{propx}\label{morzeroco}
Les 0-cocycles sont exactement les automorphismes de $B[\epsilon]$-module de $M[\epsilon]$ égaux à $id_M$ quand $\epsilon=0$.
\end{propx}
\begin{pr}
Soit $f:M[\epsilon]\to M[\epsilon]$ un tel automorphisme. Par linéarité il se réduit à l'application $B$-linéaire
$$
f=f^0+f^1\epsilon:M\to M\oplus M[\epsilon]
$$
où $f^0=id_M$ par hypothèse. La $b$-linéarité de $f$ impose la condition suivante sur $f^1$, nécessaire et suffisante pour caractériser $f$ en terme de $f^1$ :
$$
f^1(bm) = bf^1(m)
$$
pour tous $b\in B$ et $m\in M$. Cette condition est exactement celle de 0-cocycle.
\end{pr}

\subsubsection{1-cocycles}
$d:\Hom_A(B\otimes_AM,M)\to \Hom_A(B\otimes_AB\otimes_AM,M)$ est donnée par 
$$
df(b',b,m) = b'f(b,m)-f(b'b,m)+f(b',bm)
$$

\begin{propx}\label{morunco}
Les 1-cocycles sont exactement les structures de $B[\epsilon]$-module sur le $B$-module $M[\epsilon]$ qui pour $\epsilon=0$ redonnent la structure de $B$-module de $M$.
\end{propx}
\begin{pr}
Soit $\mu:B[\epsilon]\otimes_{A[\epsilon]} M[\epsilon] \to M[\epsilon]$ une telle structure. Par linéarité elle se réduit à 
$$
\mu=\mu^0+\mu^1\epsilon:B\otimes_A M \to M\oplus M\epsilon
$$
où $\mu^0$ est la structure de $B$-module de $M$.
La condition sur $\mu^1$ pour que $\mu$ soit une structure de module est que, pour tous $b,b'\in B$ et $m\in M$ :
$$
b'\mu^1(b,m)+\mu^1(b',bm) = \mu^1(b'b,m).
$$
Cette condition est exactement celle de 1-cocycle.
\end{pr}

\subsection{Cobords}

Soit $M$ un $B$-module, il existe sur $M[\epsilon]$ une structure canonique de $B[\epsilon]$-module donnée par le prolongement $\epsilon$-linéaire de la structure de $M$. Ceci amène à poser la définition suivante.

\begin{defix}
Une structure de $B[\epsilon]$-bimodule sur $M[\epsilon]$ est dite {\em triviale} si elle est isomorphe à la structure canonique.
\end{defix}

\begin{propx}\label{derintmor}
Les structures triviales de $B[\epsilon]$-bimodule sur $M[\epsilon]$ sont en bijection avec les éléments de 
$B^1_B(M)$. Et les isomorphisme $M[\epsilon]\to M'$ de la structure canonique vers une telle structure induisant $id_M$ en $\epsilon=0$ sont en bijection avec les éléments de $C^0_B(M)$.
\end{propx}
\begin{pr}
Soit $\mu^1:M\to M$ une structure de $B[\epsilon]$-module sur $M[\epsilon]$ et soit un isomorphisme de la structure canonique vers $\mu$. En utilisant la linéarité en $\epsilon$, $f$ se caractérise par 
$$
f:f^0+f^1\epsilon :M\to M\oplus M\epsilon
$$
où $f^0=id_M$ par hypothèse. La $B$-linéarité de $f$ impose sur $f^1 et \mu^1$, que pour tout $b\in B$ et tout $m\in M$, on ait
$$
\mu^1(b,m)+bf^1(m) = f^1(bm)
$$
ce qui est exactement la condition qui décrit $\mu$ comme le cobord de $f^1$.
\end{pr}

\subsection{Catégories linéaires}

Si $C$ est une catégorie $A$-linéaire ayant un nombre fini d'objets et si $M$ un $C\z C$-module, on leur associe le $[C^o\otimes_AC]$-module $[M]$ (cf. {\sf\S\ref{matrice}}).
Dans ce cadre, on abrège $C^*_{C\z C}(M)$ le complexe $C^*_{[C^o\otimes_AC]([M])}$ et on note $Ext_{C\z C}(M,M)$ sa cohomologie.

\subsection{Interprétation de Morita de la cohomologie de Hochschild}\label{interpmorhoch}

\begin{propx}\label{hochmorita}
Soit $C$ une catégorie $A$-linéaire ayant un nombre fini d'objets et donc les modules de morphismes sont projectifs sur $A$, la cohomologie de Hochschild de $C$ se calcule par le complexe 
$C^*_{C\z C}(C)$, \ie $HH^i(C) = Ext^i_{C\z C}(C,C)$.
\end{propx}
\begin{pr}
La cohomologie de Hochschild de $C$ coïncide avec celle de $[C]$, car, définit par $HH^i(C) = Ext^i_{C\z C}(C,C)$, elle ne dépend que de la catégorie $C\Mod\simeq [C]\Mod$.

Si $C$ est un anneau, l'argument pour prouver $HH^*(C) = Ext^*_{C\z C}(C,C)$ consiste à montrer que le complexe de Hochschild provient d'une résolution libre de $C$ comme $C^o\otimes_AC$-module (cf. \cite[ch. IX]{CE}).
\end{pr}

\medskip

Soit $C$ une catégorie $A$-linéaire plate ayant un nombre fini d'objet. On tire de cette proposition les conséquences suivantes.

\begin{defix}
Le {\em centre de Morita de $C$}\index{centre de Morita}, $Z_{mor}(C)$, est l'algèbre des automorphismes de $C[\epsilon]$, vu comme $C[\epsilon]\z C[\epsilon]$-bimodule, qui en $\epsilon=0$ donnent $id_C$ (cette condition est bien stable par composition).
En d'autres termes, $Z_{mor}(C)$ est le module des déformations plates au premier ordre de $id_C$.
\end{defix}

\begin{propx}
$Z_{mor}(C)=\End_{C\z C}(C) = HH^0(C) (=HZ^0(C))$
\end{propx}
\begin{pr}
La première égalité résulte de la proposition {\ref{morzeroco}} et la seconde de la proposition {\ref{hochmorita}}.
\end{pr}

\begin{defix} 
Soit $C$ une catégorie linéaire, une {\em dérivation de Morita de $C$}\index{dérivation de Morita} est un $C[\epsilon]\z C[\epsilon]$-bimodule inversible tel que sa restriction à $\epsilon=0$ soit égale à $C$. En d'autres termes $M$ est une déformation plate au premier ordre de $C$ comme $C\z C$-bimodule.
\end{defix}

\begin{propx}\label{kertangentbim}
Le module des dérivations de Morita de $C$ est isomorphe à $Z^1_{C\z C}(C)$
\end{propx}
\begin{pr}
C'est une application de la proposition {\ref{morunco}}.
\end{pr}

\begin{defix} 
Une dérivation de Morita $M$ de $C$ est dite {\em intérieure}\index{dérivation intérieure de Morita} s'il existe un isomorphisme de $C\z C$-bimodules $f:C\to M$.
\end{defix}

\begin{propx}\label{tangentisobim}
Les dérivations intérieures de Morita de $C$ sont en bijection avec les éléments de $B^1_{C\z C}(C)$ et tout isomorphisme 
$C[\epsilon]\to M$ de $C[\epsilon]$ muni de sa structure canonique de bimodule vers une dérivation intérieure de Morita, qui induit $id_C$ en $\epsilon=0$, est induit par un unique élément de $C^0_{C\z C}(C)$.

Le module des dérivations de Morita de $C$ modulo les dérivations intérieures est $HH^1(C)$.
\end{propx}
\begin{pr}
La première assertion est la proposition {\ref{derintmor}}, la deuxième résulte de ce que $Ext^1_{C\z C}(C,C)=HH^1(C)$.
\end{pr}

\medskip

Pour $\alpha\in HZ^2(C)$, on note $C[\epsilon]_\alpha$ la déformation associée, si $\alpha=0$ on note simplement $C[\epsilon]$.
Une structure de $C[\epsilon]\z C[\epsilon]_\alpha$-bimodule sur $C[\epsilon]$ est donné par une application
$$
\mu : C[\epsilon]_\alpha\otimes_AC[\epsilon]\otimes_AC[\epsilon]\longrightarrow C[\epsilon]
$$
qui, par linéarité en $\epsilon$, se réduit à 
$$
\mu^0+\mu^1\epsilon : C\otimes_AC\otimes_AC\longrightarrow C\oplus C\epsilon
$$
où $\mu^0$ est donné par le produit de $C$ si on veut que la structure déforme celle de $C$.
Dans une telle structure, en notant par concaténation le produit de $C$, la condition sur $\mu^1$ est que, pour tout $a',a,b,c,c'\in C$, on ait :
$$
\alpha(a',a)bcc' = \mu^1((a',abc,c') + a'\mu^1(a,b,c)c' - \mu^1(a'a,b,cc').
$$

Pour $\mu$ une telle structure on note $C[\epsilon]_\mu$ le bimodule associé.

\begin{defix} 
Une déformation $\alpha$ de la structure de $C$ est dite {\em Morita-triviale}\index{déformation Morita-triviale} s'il existe un $C[\epsilon]\z C[\epsilon]_\alpha$-bimodule inversible qui pour $\epsilon=0$ est isomorphe à $C$.
\end{defix}

\begin{propx}\label{moritatriv}
Les déformations Morita-triviales de $C$ sont en bijection avec les cobords des 1-cochaînes de Hochschild, \ie sont en bijection avec les déformations triviales de $C$.
\end{propx}

\medskip

Avant d'attaquer la preuve de la proposition, on établit deux lemmes. Par souci de simplifier la rédaction des preuves on se limite à travailler avec des algèbres associatives plutôt que des catégories, mais les résultats restent vrais en général.

Soit $f:C\to D$ on rappelle (cf. preuve du lemme {\ref{eqcateqmor}}) qu'on associe à $f$ une structure de $C\z D$-bimodule sur $D$ donnée par :
\begin{eqnarray*}
D\otimes_A D\otimes_AC & \longrightarrow & D\\
(a\otimes b\otimes c) &\longmapsto & abf(c)
\end{eqnarray*}
Le cadre d'application ici est celui où $f$ est un isomorphisme et, afin de simplifier les formules, on préférera associer à $f$ le bimodule inverse du précédent donné par :
\begin{eqnarray*}
D\otimes_A D\otimes_AC & \longrightarrow & D\\
(a\otimes b\otimes c) &\longmapsto & f^{-1}(a)bc
\end{eqnarray*}

\begin{lemmex}\label{bimtrivdroite}
Une structure de $C[\epsilon]\z C[\epsilon]_\alpha$-bimodule sur $C[\epsilon]$ qui redonne la structure canonique de $C\z C$-bimodule de $C$ quand $\epsilon=0$ et dont la structure à droite est donné par le produit de $C[\epsilon]$ provient d'un morphisme $C[\epsilon]\to C[\epsilon]_\alpha$ qui pour $\epsilon=0$ est l'identité de $C$.
\end{lemmex}
\begin{pr}
Par hypothèse, pour tous $b,c\in C$, on a 
$$
\pi(1,b,c) = bc \iff \pi^1(1,b,c)=0.
$$
On tire de la condition 
$$
\alpha(a',a)bcc' = \pi^1(a',abc,c') + a'\pi^1(a,b,c)c' - \pi^1(a'a,b,cc')
$$
et de $\alpha(a,1)=\alpha(1,a)=0$ les deux équations :
\begin{eqnarray*}
&\pi^1(1,ab,c') + \pi^1(a,b,1)c' - \pi^1(a,b,c') = 0\\
&\pi^1(a',c,1) + a'\pi^1(1,1,c) - \pi^1(a',1,c) = 0
\end{eqnarray*}
soit en tenant compte de l'hypothèse sur $\pi^1$
\begin{eqnarray*}
&\pi^1(a,b,c) = \pi^1(a,b,1)c \\
&\pi^1(a,b,1) = \pi^1(a,1,b)
\end{eqnarray*}
et finalement, pour tous $a,b,c\in C$ :
$$
\pi^1(a,b,c) = \pi^1(a,1,1)bc.
$$
On pose $\beta(a) := \pi^1(a,1,1)$ la condition sur $\pi^1$ assure que $\alpha=d\beta$ est un 2-cobord de Hochschild et 
le bimodule de départ est celui associé au morphisme $\beta:C[\epsilon]\to C[\epsilon]_{d\beta}$.
\end{pr}

\begin{lemmex}\label{trivstbim}
Une structure de $C[\epsilon]\z C[\epsilon]_\alpha$-bimodule sur $C[\epsilon]$ déformant la structure de $C$
est toujours isomorphe à une structure où la loi droite est donnée par le produit de $C[\epsilon]$.
\end{lemmex}
\begin{pr}
Soit $C[\epsilon]_\mu$ un tel bimodule, on construit un morphisme $f:C[\epsilon]\to C[\epsilon]_\mu$ de $C[\epsilon]$-module droit en posant, pour $x=x^0+x^1\epsilon\in C[\epsilon]$ : 
$$
f(x) = \mu(1,1,x)=x^0 + (x^1+\mu^1(1,1,x^0))\epsilon.
$$
En effet, pour tout $a \in C[\epsilon]$ on a d'un côté 
$$
f(x.a) = f(xa) = \mu(1,1,xa)
$$
et d'autre
$$
f(x).a = \mu(1,f(x),a) = \mu(1,\mu(1,1,x),a) = \mu(1,1,xa).
$$
$f$ est un isomorphisme car son inverse est donné, pour $y=y^0+y^1\epsilon\in C[\epsilon]$, par 
$$
g(y) = y^0 + (y^1-\mu^1(1,1,y^0).
$$
Il suffit, pour conclure, de définir la structure droite sur $C[\epsilon]$ par transport de celle de $\to C[\epsilon]_\mu$.
\end{pr}

\bigskip

\begin{prpr}
Une déformation triviale $C[\epsilon]\to C[\epsilon]_\alpha$ fournit naturellement une structure de bimodule inversible sur $C[\epsilon]$, qui redonne pour $\epsilon=0$ la structure canonique de $C\z C$-bimodule sur $C$ ; la propriété sera donc démontrée si on prouve que toute déformation $\alpha$ Morita-triviale de $C$ est équivalente à la canonique par un bimodule provenant d'une équivalence. 

Soit $\alpha$ une déformation Morita-triviale, par définition il existe $\mu$ une structure de $C[\epsilon]\z C[\epsilon]_\alpha$-bimodule sur $C[\epsilon]$ déformant celle de $C$, dont il suffit, pour que $M$ soit associé à un isomorphisme $C[\epsilon]\to C[\epsilon]_\alpha$, que la structure droite soit donnée par le produit de $C[\epsilon]$ (lemme {\ref{bimtrivdroite}}). A priori, ceci est faux, mais le lemme {\ref{trivstbim}} assure que cela est vrai à isomorphisme de bimodules près. Ce nouveau bimodule donne l'isomorphisme voulu de $C[\epsilon]$ et $C[\epsilon]_\alpha$.
\end{prpr}

\setcounter{totx}{0}
\newpage
\thispagestyle{empty}
\chapter{Espace classifiant de localisation}\label{cmf}
\thispagestyle{chheadings}

\section{Infini-groupoïdes}\index{infini-groupoïde}

On rappelle que les ensembles simpliciaux peuvent être considérés comme des modèles d'infini-groupoïdes.
Sans poser de définition précise, l'intuition sur le sujet considère les infini-catégories comme des catégories possédant des objets, des 1-flèches entre ces objets, des 2-flèches entre les 1-flèches, et plus généralement des $n$-flèches entre les $(n-1)$-flèches, pour tout $n>1$, ces flèches sont munies de compositions pour lesquelles on distingue des flèches neutres ; les infini-groupoïdes sont des infini-catégories dont toutes les flèches possèdent des inverses, éventuellement en un sens affaibli.

Le meilleur outil d'intuition sur les infini-groupoïdes consiste en la théorie homotopique des espaces topologiques où les notions de points, de chemins, d'homotopies de chemins, d'homotopies d'homotopies de chemins, etc. jouent respectivement les rôles des objets, des 1-flèches, des 2-flèches, etc.
On définit donc le {\sl modèle canonique} des infini-groupoïdes comme la catégorie de modèles des espaces topologiques de Kelley ({\cite[def. 2.4.21]{hovey}}) pour laquelle les équivalences sont les équivalences faibles d'homotopie ({\cite[\S 2.4]{hovey}}).

\medskip

La catégorie $\sens$ des ensembles simpliciaux avec la structure de modèle décrite dans {\cite[\S 3.2]{hovey}} est équivalente à celle des espaces topologiques ({\cite[thm. 3.6.7]{hovey}}) et peut donc être utilisée comme catégorie de modèles pour les infini-groupoïdes\footnote{C'est ce qui justifie la définition des champs supérieurs comme champs simpliciaux.}.

\section{Catégories simpliciales}\label{catsplx}

Une catégorie simpliciale\index{catégorie simpliciale} est une catégorie enrichie sur $\sens$. C'est aussi un objet en catégorie $C_1\rightrightarrows C_0$ dans $\sens$ tel que $C_0$ soit un ensemble.
Toute catégorie est canoniquement une catégorie simpliciale dont les espaces de morphismes sont tous discrets ; par opposition aux catégories simpliciales, on les qualifiera de {\em discrètes}\index{catégorie discrète}.
Plus précisément, si on note $\scatcat$ la catégorie des catégories simpliciales et $\catcat$ la catégorie des catégories discrètes, on a une adjonction :
$$
\pi_0 : \scatcat\leftrightarrows\catcat : \iota
$$
où $\pi_0$, est le foncteur qui à une catégorie simpliciale $C$ associe la catégorie discrète ayant les mêmes objets et $\pi_0(C(x,y))$ comme ensembles de morphismes entre deux objets $x,y$.

$\iota$ admet également un adjoint à droite $C\mapsto C_{(0)}$ qui a une catégorie simplicial associe la catégorie discrète ayant les mêmes objets et comme ensemble de morphismes les seules parties de degré 0 des ensembles simpliciaux de morphismes de $C$.


\bigskip

Les catégories simpliciales sont des modèles pour des {\em catégories supérieures}\index{catégorie supérieure} ayant des 1-flèches non forcément inversibles et, pour $n>1$, des $n$-flèches toutes inversibles. En effet, conformément à l'interprétation des ensembles simpliciaux comme modèles d'infini-groupoïdes, les catégories simpliciales peuvent s'interpréter comme des modèles de catégories dont les objets de morphismes sont des infini-groupoïdes.

\medskip

Un morphisme $f:C\to D$ de catégories simpliciales est dit une {\em équivalence} s'il induit une équivalence de catégorie après application du foncteur $\pi_0$ et si, pour tout couple d'objets $x,y\in C$, $f_{xy}:C(x,y)\to D(fx,fy)$ est une équivalence dans $\sens$.
La catégorie homotopique des catégories simpliciales pour leurs équivalences est notée $Ho(\scatcat)$.

\medskip

Une catégorie simpliciale $C$ est dite un {\em groupoïde simplicial}\index{groupoïde simplicial} si $\pi_0C$ est un groupoïde, ce sont des modèles de groupoïdes faibles ; on note $\sgpd$ la sous-catégorie pleine de $\scatcat$ formée des groupoïdes simpliciaux. On a une adjonction :
$$
\iota : \sgpd \leftrightarrows \scatcat : (-)^{int}
$$
où, pour $C\in \scatcat$, $C^{int}$ est la sous-catégorie ayant les mêmes objets et ayant comme ensemble de morphismes de $x$ vers $y$ les images inverses par $C(x,y)\to \pi_0C(x,y)$ des isomorphismes de $\pi_0C(x,y)$. $C^{int}$ est dit le {\em groupoïde intérieur}\index{groupoïde intérieur} de $C$. Dans le cas où $C$ est discrète $C^{int}$ est simplement le groupoïde des isomorphismes de $C$.

\medskip

Le {\em nerf}\index{nerf d'une catégorie simpliciale} $N(C)$\index{$N(C)$} d'une catégorie simpliciale est l'objet simplicial dans $\sens$ : $n\mapsto C_1\times_{C_0} \dots \times_{C_0} C_1\in \sens$ ($n$ facteurs).
La {\em réalisation géométrique}\label{realgeom} $|N(C)|$ (abrégé $|C|$\index{$|C|$}) du nerf d'une catégorie simpliciale $C$ est la colimite homotopique du diagramme $N(C)$ dans $\sens$, elle se réalise par la diagonale de $N(C)$ : $n\mapsto (C_1\times_{C_0} \dots \times_{C_0} C)_n \in\ens$.
Dans le cas où la catégorie $C$ est discrète on a canoniquement $|C|\simeq N(C)$.

On définit $s\sens$ comme la catégorie des ensembles bisimpliciaux. Les constructions précédentes sont fonctorielles :
$$
\scatcat \overset{N}{\longrightarrow} \mathrm{s\sens} \overset{|-|}{\longrightarrow} \sens.
$$

\bigskip

On renvoie à {\cite[ch. IV]{hovey}}, pour les définitions des catégories monoïdales fermées et de modules fermés sur de telles catégories.
La catégorie $\sens$ est monoïdale fermée pour le produit catégoriel. Une catégorie $C$ est dite un {\em $\sens$-module fermé} s'il existe un foncteur 
\begin{eqnarray*}
\otimes : C\times \sens & \longrightarrow& C\\
(x,K) &\longrightarrow & x\otimes K
\end{eqnarray*}
vérifiant certaines conditions d'associativité et d'adjonction (cf. \cite{hovey}).
En particulier un $\sens$-module fermé $C$ permet de définir un enrichissement de $C$ en une catégorie simpliciale dont les espaces de morphismes $\Hom^\Delta_C(x,y)$ sont donnés par 
$$
\Hom^\Delta_C(x,y)_n := \Hom_C(x\otimes \Delta^n,y).
$$
Pour $x,y,z\in C$ le morphisme de composition dans $\sens$ :
$$
\Hom_C^\Delta(x,y)\times \Hom_C^\Delta(y,z)\longrightarrow \Hom_C^\Delta(x,z)
$$
est définit en associant à $(f:x\otimes\Delta^n\to y,g:y\otimes\Delta^n\to z)$ la composition 
$gf:x\otimes\Delta^n\to y\to y\otimes\Delta^n\to z$, où $y\to y\otimes\Delta^n$ provient du morphisme $\Delta^0\to\Delta^n$ pointant $\Delta^n$ en 0.

\section{Localisation et espaces classifiant}

Soit $C$ une catégorie, une {\em sous-catégorie d'équivalences} pour $C$ est simplement une sous-catégorie de $C$, considérée dans une perspective de localisation.
Si $C$ est un catégorie et $W$ une sous-catégorie d'équivalences, la localisation de $C$ par $W$ est notée $C[W^{-1}]$.
$L(C,W)$\index{$L(C,W)$} désigne la localisation simpliciale de $C$ par $W$ (cf. \cite{dk}). On rappelle que c'est une catégorie simpliciale munie d'un morphisme $C\to L(C,W)$ universel dans $Ho(\scatcat)$ (\cite[rem. 2.2.1]{hag1}) pour les morphismes de $C$ vers les catégorie simpliciales $D$ qui envoient $W_C$ dans $D^{int}$.
On a en particulier un isomorphisme :
$$
\pi_0L(C,W)\longrightarrow C[W^{-1}].
$$

\begin{defix}
Pour une catégorie $C$ munie d'une sous-catégorie d'équivalences $W$, on définit un {\em espace classifiant} de $(C,W)$ comme un ensemble simplicial $X$ équivalent à $|L(C,W)^{int}|$.
\end{defix}

Le théorème de délaçage de Segal (\cite{pellissier, segal, tamsamani}) permet de déduire qu'un espace classifiant $X$ vérifie les propriétés suivantes :
\begin{itemize}
\item il existe une bijection de $\pi_0(X)$ vers les classes d'isomorphisme de $C[W^{-1}]$ ;
\item si $x$ est un point de $X$, il existe, pour tout $n\geq1$ des isomorphismes :
$$\pi_n(X,x)\longrightarrow \pi_{n-1}(L(C,W)^{int}(x,x))$$
(si $n\geq 2$, on a $\pi_n(L(C,W)^{int}(x,x)) = \pi_n(L(C,W)(x,x))$).
\end{itemize}

\section{Classifiants de catégories de modèles}\label{quillen}

On renvoie à \cite{hirschhorn, hovey} pour les définitions des catégories de modèles\index{catégorie de modèles}, catégories de modèles simpliciales\index{catégorie de modèles simpliciale} et des adjonctions de Quillen. 

\begin{defix}
On note $\cmg$\index{$\cmg$} (resp. $\cmd$\index{$\cmd$}) la catégorie des catégories de modèles avec les adjoints de Quillen à gauche (resp. à droite) comme morphismes et $\cmsg$\index{$\cmsg$} (resp. $\cmsd$\index{$\cmsd$}) la catégorie des catégories de modèles simpliciales avec comme morphismes les adjoints de Quillen à gauche (resp. à droite) compatible aux structures de $\sens$-modules.
\end{defix}

\subsection{Le cas général}

Si $M$ est une catégorie de modèles, on note $W_M$ la sous-catégorie des équivalences ; 
$M^c$ la sous-catégorie pleine des objets cofibrants ; $M^f$ celle des objets fibrants ; et $M^{cf}=M^c \cap M^f$. 
On note $Ho(M)=M[W_M^{-1}]$ la catégorie homotopique de $M$ et $L(M,W_M)$ la localisée de Dwyer-Kan.
$Ho(M)$ est naturellement enrichie sur $Ho(\sens)$ et $L(M,W_M)$ est un modèle strict pour cet enrichissement.

Pour deux objets $x$ et $y$ de $M$, $Map_M(x,y)$\index{$Map_M(x,y)$} désigne l'ensemble simplicial des morphismes entre $x$ et $y$, calculé à l'aide d'un choix de framing, il est équivalent à l'ensemble simplicial $\Hom_{L(M,W_M)}(x,y)$. $[x,y]$ désigne l'ensemble des morphismes entre $x$ et $y$ dans $Ho(M)$ on a :
$$[x,y] \simeq \pi_0(Map_M(x,y)).$$
$Map^{eq}_M(x,y)$ désigne l'ensemble simplicial des équivalences entre $x$ et $y$, définit comme l'image inverse des isomorphismes de $Ho(M)$ par $Map_M(x,y)\to [x,y]$, il est équivalent à $\Hom_{L(M,W_M)^{int}}(x,y)$.

\begin{defix}\label{sscatmod}
Une {\em sous-catégorie de modèles}\index{sous-catégorie de modèles} $C$ est un couple $(C,M_C)$ tel que :
\begin{itemize}
\item $M_C$ soit une catégorie de modèles et $C$ une sous-catégorie pleine ;
\item $C$ soit saturée par équivalences (tout objet relié par une chaîne d'équivalences à un objet de $C$ est dans $C$).
\end{itemize}
Une {\em sous-catégorie de modèles simpliciale} $C$ est une sous-catégorie de modèles $(M,C)$ d'une catégorie de modèles simpliciale.

\medskip
Soient $C=(C,M_C)$ et $D=(D,M_D)$ deux sous-catégories de modèles, un {\em morphisme de sous-catégories de modèles} $f:C\to D$ est un morphisme $f':M_C\to M_D$ envoyant $C$ dans $D$. On dit que $f$ est un {\em foncteur de Quillen à gauche (resp. à droite)} si $f'$ est un adjoint de Quillen à gauche (resp. à droite). Les morphismes de sous-catégories de modèles simpliciales sont définis de même.
\end{defix}

\begin{defix}
On note $\scmg$\index{$\scmg$} (resp. $\scmd$\index{$\scmd$}) la catégorie des sous-catégories de modèles et des foncteurs de Quillen à gauche (resp. à droite), et $\scmsg$\index{$\scmsd$} (resp. $\scmsg$\index{$\scmsd$}) la catégorie des sous-catégories de modèles simpliciales et des foncteurs de Quillen à gauche (resp. à droite).
\end{defix}

\medskip

Pour une sous-catégorie de modèles $C$, on définit sa {\em catégorie des équivalences} $W_C=C\cap W_M$ et les $C^*=C\cap M^*$ et $W_C^*=W_C\cap M^*$ où $*=c$, $f$ ou $cf$, dite, respectivement catégorie des objets cofibrants, fibrants et fibrants-cofibrants de $C$. Comme $C$ est pleine dans $M$, $W_C$ est pleine dans $W_M$.
On note $Ho(C)$ la sous-catégorie pleine de $Ho(M)$ engendrée par les objets de $C$, $Ho(C) = C[W_C^{-1}]$. On note $L(C,W_C)$ la sous-catégorie simpliciale de $L(M,W_M)$ engendrée par les objets de $C$, c'est aussi la localisation simpliciale de $C$ par rapport à $W_C$.

\medskip
\noindent{\bfseries\sffamily Remarque.} Malgré son nom, le couple $(C,W_C)$ d'une sous-catégorie de modèles ne s'enrichit pas forcément en une structure de modèles (par exemple par ce que $C$ ne possède pas forcément les limites ou colimites).

\medskip

Les équivalences d'une sous-catégorie de modèles sont toujours une classe saturée.

\begin{rpropx}[{\cite[prop. A.0.6]{hag2}}]\label{classifmodele}
Si $C$ une sous-catégorie de modèles, alors $|W_C|$, $|W_C^c|$, $|W_C^f|$ et $|W_C^{cf}|$ sont des espaces classifiant pour $(C,W_C)$.
\end{rpropx}

\subsection{Le cas simplicial}

Ce qui suit est essentiellement l'Appendice A de \cite{hag2}. 
On construit un foncteur $\mathcal{G}:\scms^*\to\scatcat$\label{G}, où *=$g$ ou $d$.

Soit $C$ une sous-catégorie de modèles simpliciale.
On note $\Hom_C^\Delta(x,y)$\label{homsplx} l'ensemble simplicial des morphismes entre deux objets $x$ et $y$ de $C$. 
Si $x$ est cofibrant et $y$ est fibrant, on a une équivalence \cite{hirschhorn, hovey} :
$$Map_M(x,y)\simeq \Hom_C^\Delta(x,y).$$
On note $\Eq_C^\Delta(x,y)$ ou $\G(C)(x,y)$ l'ensemble simplicial des équivalences entre deux objets $x$ et $y$ de $C$, définit par :
$$n\mapsto \Eq_C^\Delta(x,y)_n = \Hom_{W_C}(x\otimes \Delta^n,y)$$

\medskip

Si $x$ est cofibrant et $y$ est fibrant, on a une équivalence \cite{hirschhorn, hovey}
$$\G(C)(x,y)=\Eq_C^\Delta(x,y)\simeq Map^{eq}_M(x,y).$$
Si on se limite à considérer les couples d'objets $(x,y)$ dans $C^{cf}$, les ensembles simpliciaux $\G(C)(x,y)$ définissent une catégorie simpliciale $\G(C)$\label{GC}\index{$\mathcal{G}(C)$} dont la proposition suivante établit son équivalence avec $L(C,W_C)^{int}$. On note $|\G(C)|$ la réalisation géométrique du nerf de $\G(C)$.

\begin{rpropx}[{\cite[prop. A.0.6]{hag2}}]\label{comparclassifmodele}
Pour $C$ une sous-catégorie de modèles, il existe des équivalences $|\G(C)| \simeq |W_C|$ fonctorielles en $C$ (pour les foncteurs respectant les objets fibrants et cofibrants). En particulier, pour $x$ et $y$ dans $C^{cf}$, on a des équivalences, fonctorielles en $(C,x,y)$ :
$$
\Omega_{xy}|W_C| \overset{\sim}{\longrightarrow}\Omega_{xy}|\G(C)|\overset{\sim}{\longleftarrow}\G(C)(x,y).
$$
\end{rpropx}

\subsection{Fonctorialité des espaces de morphismes}\label{fonctomorquillen}

\paragraph{Cas général.} 

On rappelle (cf. {\cite[\S 5.6]{hovey}}) qu'il existe un 2-foncteur faible :
$$
Ho : \cmg \to Ho(\sens)\Mod^g
$$
où $\cmg$ est vue comme 2-catégorie et où $Ho(\sens)\Mod$ est la 2-catégorie des catégories qui sont des modules sur la catégorie monoïdale $Ho(\sens)$, ayant comme morphismes les foncteurs adjoints à gauche et leur transformations naturelles.

On a également un foncteur de 
$$
s : Ho(\sens)\Mod^g \longrightarrow Ho(\sens)\z\catcat^g
$$
vers les catégories enrichies sur $Ho(\sens)$ avec les seuls adjoints à gauche comme foncteurs.
Si $M\in Ho(\sens)\Mod^g$, $sM$ est définie comme ayant les mêmes objets et comme morphismes d'un objet $x$ vers un objet $y$ l'objet $Map_{M}(x,y)\in Ho(\sens)$, définit à isomorphisme unique près comme représentant le foncteur 
\begin{eqnarray*}
Map^{eq}_{C}(x,y) : Ho(\sens) &\longrightarrow & \ens \\
K & \longmapsto & [x\otimes K,y]^{iso}
\end{eqnarray*}
où $[x\otimes K,y]^{iso}$ désigne l'ensemble des éléments de $[x\otimes K,y]$ qui sont des isomorphismes.
	
En conséquence, si $I$ est une petite catégorie et $M:I^o\to \cmg$ est un préfaisceau de catégories de modèles, $M$ fournit, par composition avec $s\circ Ho$, un préfaisceau faible en catégories enrichies sur $Ho(\sens)$, dont le paragraphe suivant en propose un modèle qui est un préfaisceau simplicial.

\paragraph{Cas simplicial.} Soit $f:C\to D\in\scmsg$, $f$ ne préserve pas a priori les objets cofibrants et fibrants, et, en conséquence, $f$ n'induit pas de foncteur $\mathcal{G}(C)\to \mathcal{G}(D)$. On est obligé pour avoir une propriété de fonctorialité pour $\mathcal{G}$ de se restreindre aux foncteurs conservant les objets cofibrants et fibrants. 

\begin{defix}
On note $\scmsg_{cf}$ la sous-catégorie de $\scmsg$ formée de ces foncteurs.
\end{defix}

En effet, si $f:C\to D\in\scmsg_{cf}$, il vérifie, pour tout $x\in C^{cf}$, que $f(x)\in D^{cf}$ et que $f(x\otimes\Delta^n)$ est isomorphe à $f(x)\otimes\Delta^n$, d'où on déduit un morphisme :
$$
\mathcal{G}(f)_{x,y}:\Eq_C^\Delta(x,y)\to \Eq_D^\Delta(f(x),f(y))
$$
compatible avec les compositions, \ie ces morphismes définissent un foncteur $\mathcal{G}(C)\to \mathcal{G}(D)$.

\medskip

On a donc un foncteur $\mathcal{G} : \scmsg_{cf}\longrightarrow \scatcat$ (en fait à valeurs dans $\sgpd$)  qui s'inscrit dans le diagramme commutatif dans les 2-catégories à un 2-isomorphisme naturel près :
$$\xymatrix{
\scmsg_{cf}\ar[r]^-{\mathcal{G}}\ar[d] &\scatcat^g \ar[d]\\
\scmg\ar[r]^-{Ho} & Ho(\sens)\z\catcat
}$$
où on a limité les morphismes de $\scatcat$ à être adjoints à gauche et où le foncteur $\scatcat\to Ho(\sens)\z\catcat$ est simplement le passage aux type d'homotopie des espaces de morphismes.
Seul $Ho$ est un 2-foncteur faible dans ce diagramme.

\section{Résumé}\label{foncteursvaleurschamps}

Au final, dans le cas gauche, on a un diagramme (non strictement commutatif) :
$$
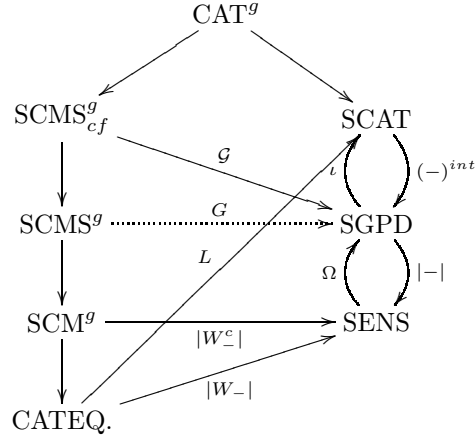
\xymatrix{
&\catcat^g\ar[rd]\ar[ld]\\
\scmsg_{cf}\ar[d]\ar[drr]^{\mathcal{G}} && \scatcat \ar@/^1pc/[d]^{(-)^{int}}\\
\scmsg\ar[d]\ar@{..>}[rr]^{G} && \sgpd\ar@/^1pc/[d]^{|-|}\ar@/^1pc/[u]^{\iota} \\
\scmg\ar[d]\ar[rr]_{|W_-^c|} && \sens\ar@/^1pc/[u]^{\Omega} \\
\catcateq\ar[rru]_{|W_-|}\ar[rruuu]^{L} .
}$$

$\catcat^g$ est la catégorie des catégories avec les seuls foncteurs adjoints à gauche comme morphismes.

$L$ est le foncteur de localisation simpliciale, $\mathcal{G}$ est définit au {\sf \S\ref{G}} foncteur
Le triangle contenant $\mathcal{G}$ et $L$ est homotopiquement commutatif, \ie $\mathcal{G}$ et $L$, vus comme foncteurs $\scmg_f \to \scatcat$, sont des foncteurs homotopes.
L'homotopie de $|L(-)^{int}|$ et de $|W_-|$ est conjecturale.

Les adjoints de Quillen à gauche préservant les équivalences entre objets cofibrants, on dispose aussi d'un foncteurs $|W_-^c|:\cmg\to \sens$, qui est homotope à $|W_-|:\cmg\to \catcateq\to \sens$.

Le morphisme $G$ associe à une catégorie de modèles simpliciale l'intérieur de son enrichissement simplicial, il n'entre dans aucune relation de commutation avec les autres foncteurs ; il faudrait pour cela le "dériver", et c'est en quelque sorte ce que fait le foncteur $\mathcal{G}$, mais au prix de devoir restreindre les morphismes considérés entre catégories de modèles, pour conserver la fonctorialité.

\medskip

On a un diagramme similaire avec les adjoints à droite comme foncteurs.

\section{Produits fibrés de nerfs}\label{profibnerf}

On détaille un cas particulier du théorème de strictification de {\cite[thm. {\sf B.0.7}]{hag2}}.

Soient $M_1$, $M_2$ et $M_3$ trois catégories de modèles engendrées par cofibrations dont on note les sous-catégories d'équivalences entre objets cofibrants respectivement par $W^c_1$, $W^c_2$ et $W^c_3$ ; et soient $a_1:M_1\to M_3 $ et $a_2:M_2\to M_3$ deux morphismes de Quillen à gauche.
On définit une catégorie $M_1\times_{M_3}M_2$ ayant pour objets les quintuplets $(x_1,x_2,x_3,\chi_1,\chi_2)$ où $x_i\in M_i$ et où les $\chi_i:a_i(x_i)\to x_3$ sont des morphismes de $M_3$. Les morphismes de $M_1\times_{M_3}M_2$ entre deux objets $(x_1,x_2,x_3,\chi_1,\chi_2)$ et $(y_1,y_2,y_3,\gamma_1,\gamma_2)$ sont les triplets $(u_1,u_2,u_3)$ où $u_i:x_i\to y_i$ tels que les deux carrés
$$\xymatrix{
a_1(x_1) \ar[r]^{\chi_1} \ar[d]_{a_1(u_1)} & x_3 \ar[d]^{u_2} & a_2(x_2) \ar[d]^{a_3(u_3)}\ar[l]_{\chi_2} \\
a_1(y_1) \ar[r]_{\gamma_1} & y_3 & a_2(y_2) \ar[l]^{\gamma_2}
}$$
commutent.
Comme les structures des $M_i$ sont toutes engendrées par cofibrations, $M_1\times_{M_3}M_2$ hérite également d'une structure de modèles engendrée par cofibration pour laquelle les fibrations et les équivalences sont définies terme à terme \cite{hag2}. 

Pour $i=1,2$, soient $Q_i$ des foncteurs de remplacements cofibrants dans $M_i$, les foncteurs dérivés à gauche des $a_i$ et $b$ sont définis par $\mathbb{L}a_i(x_i) = a_i(Q_i(x_i))$. La flèche naturelle $Q_i(x_i)\to x_i$ donne $\mathbb{L}a_i(x_i) \to a_i(x_i)$.
La catégorie $(M_1\times_{M_3}M_2)_{cart}$ est définie comme la sous-catégorie pleine de $M_1\times_{M_3}M_2$ formé des objets $(x_1,x_2,x_3,\chi_1,\chi_2)$ tels que les composés $\mathbb{L}a_i(x_i) \to a_i(x_i)\to x_3$ soient des équivalences.
On note $W^c_{cart}$ la sous-catégorie des objets cofibrants de $(M_1\times_{M_3}M_2)_{cart}$ et des équivalences entre eux.

\begin{rpropx}[{\cite[\sf thm. B.0.7]{hag2}}]\label{profibnerfprop}
Avec les notations précédentes, on a l'équivalence suivante :
$$
|W^c_1|\times^h_{|W^c_3|}|W^c_2| \simeq |W^c_{cart}|,
$$
fonctorielle en $M_1\to M_3 \leftarrow M_2$.
\end{rpropx}



\printindex

\addcontentsline{toc}{chapter}{Index}

\newpage
\thispagestyle{empty}
~~
\newpage
\thispagestyle{empty}

~~
\vspace{60pt}
\begin{quote}
{\sffamily
\footnotesize
<<~H comme Homotopie : terme mathématique longtemps réservé à la seule opération de déformation continue d'espaces topologiques, il est plus actuellement devenu la dénomination commune d'une repensée de l'algèbre. Tout au long du {\sc xx}\ieme siècle sont apparus de nouvelles notions d'identification d'objets mathématiques différentes de l'isomorphie (équivalence faible, inversibilité homotopique, quasi-isomorphisme, etc.) élargissant celle originale de la topologie, et dont le contexte théorique général est celui des catégories de modèles. En d'autres termes, il est reconnu maintenant que les objets mathématiques ont souvent des processus non triviaux d'identification sous-jacent à leur définition et la relaxe de ces processus conduit à les considérer naturellement comme les objets d'infini-catégories.~>>

\hfill {\em Abécédaire du Pendule}, \verb+http://www.lependule.org/+} 
\end{quote}

~~
\newpage
\thispagestyle{empty}

~~
\newpage
\thispagestyle{empty}
~~
\newpage
\thispagestyle{empty}

\begin{center}
{\large \sffamily \bfseries CHAMPS DE MODULES DES CAT\'EGORIES LIN\'EAIRES ET AB\'ELIENNES\\}

\begin{quote}
{\sffamily \bfseries Résumé} Les {\sl catégories linéaires} ont naturellement plusieurs notions d'identification : l'isomorphie, l'équivalence de catégories et l'équivalence de Morita. On construit les champs classifiant les catégories pour ces trois structures ($\ukcatiso$, $\ukcateq$, $\ukcatmor$) ainsi que le champ classifiant les {\sl catégories abéliennes} ($\ukab$), l'originalité étant que les trois derniers champs sont des {\sl champs supérieurs}.

Le résultat principal de la thèse est que, sous des conditions de finitude des objets classifiés, ces champs sont {\sl géométriques} au sens de C.~Simpson. En particulier, on trouve que les {\sl complexes tangents} de ces champs en une catégorie $C$, \ie les objets classifiant les déformations au premier ordre de $C$, sont donnés par des tronqués du complexe de {\sl cohomologie de Hochschild} de $C$.

En plus, il existe une suite naturelle de morphismes surjectifs de champs :
$$\ukcatiso\overset{a}{\surj} \ukcateq \overset{b}{\surj} \ukcatmor\overset{c}{\surj} \ukab$$
dont on montre que $b$ est étale, et que $c$ est une équivalence.
\end{quote}

\vspace{50pt}
{\large \sffamily \bfseries MODULI STACKS OF LINEAR AND ABELIAN CATEGORIES}
\begin{quote}
{\sffamily \bfseries Abstract} {\sl Linear categories} naturally have several identification relations : isomorphisms, categorical equivalences and Morita equivalences. In this thesis, we construct the classifying stacks for these three relations ($\ukcatiso$, $\ukcateq$, $\ukcatmor$) together with the classifying stack of {\sl abelian categories} ($\ukab$), the originality of the subject being that, apart from the first one, these are {\sl higher stacks}.

The principal result is that, under some finiteness assumptions, these stacks are geometric in the sense of C.~Simpson. In particular, one recover the {\sl Hochschild cohomology} of a category $C$ as the {\sl tangent complex}, \ie the object classifying first order deformations of $C$, of these stacks at the point defined by $C$.

Moreover, there exists a natural sequence of surjective morphisms of stacks :
$$\ukcatiso\overset{a}{\surj} \ukcateq \overset{b}{\surj} \ukcatmor\overset{c}{\surj} \ukab$$
for which we prove that $b$ is étale, and $c$ is an equivalence.
\end{quote}
\vfill
{\sffamily 
Institut de Mathématiques de Toulouse, UMR 5580, UFR MIG\\
Laboratoire \'Emile Picard,\\
Université Paul Sabatier 31062 TOULOUSE Cedex 9
}

\end{center}

\end{document}